\documentclass[11pt]{amsart}

\usepackage[utf8]{inputenc}
\usepackage{lmodern}
\usepackage{amsmath,amssymb,amsthm,latexsym,cite,cancel}
\usepackage[small]{caption}
\usepackage{graphicx,wasysym,overpic,tikz,color}
\usepackage{subfigure}
\usepackage{cite}
\usepackage[colorlinks=true,urlcolor=blue,
citecolor=red,linkcolor=blue,linktocpage,pdfpagelabels,
bookmarksnumbered,bookmarksopen]{hyperref}
\usepackage[english]{babel}
\usepackage{units}
\usepackage[utf8]{inputenc}
\usepackage{enumitem}
\usepackage[left=2.1cm,right=2.1cm,top=2.71cm,bottom=2.71cm]{geometry}
\usepackage[hyperpageref]{backref}
\usepackage{float}
\usepackage[T1]{fontenc}

\usepackage[colorinlistoftodos]{todonotes}
\newtheorem{theorem}{Theorem}[section]

\newtheorem{proposition}[theorem]{Proposition}
\newtheorem{lemma}[theorem]{Lemma}
\newtheorem{remark}[theorem]{Remark}

\newcommand{\abs}[1]{\lvert#1\rvert}
\newcommand{\norm}[1]{\lVert#1\rVert}

\definecolor{revisionpurple}{rgb}{0.5,0,0.5}

\newcommand{\G}{\mathcal{G}}
\newcommand{\K}{\mathcal{K}}

\newcommand{\E}{\mathrm{E}}
\newcommand{\V}{\mathrm{V}}
\newcommand{\R}{\mathbb{R}}

\newcommand{\e}{\mathrm{e}}
\newcommand{\vv}{\mathrm{v}}
\newcommand{\w}{\mathrm{w}}

\definecolor{darkgreen}{rgb}{0,0.45,0}

\newcounter{revisionblock}


\tikzstyle{nodo}=[circle,draw,fill,inner sep=0pt,minimum size=%
1.5mm]
\tikzstyle{bnodo}=[circle,draw,fill,inner sep=0pt,minimum size=%
2mm]

\numberwithin{equation}{section}

\title[The Schr\"odinger equations with nonlinear point defects on noncompact metric graphs]{Infinitely many positive normalized solutions for Schr\"odinger equations with nonlinear point defects on noncompact metric graphs}
\author[Z. He]{Zhentao He}

\address[Z. He]{\newline\indent
	School of Mathematics
	\newline\indent
	East China University of Science and Technology
	\newline\indent
	Shanghai 200237, PR China }
\email{\href{mailto:hezhentao2001@outlook.com}{hezhentao2001@outlook.com}}

\author[C. Ji]{Chao Ji}

\address[C. Ji]{\newline\indent
	School of Mathematics
	\newline\indent
	East China University of Science and Technology
	\newline\indent
	Shanghai 200237, PR China }
\email{\href{mailto:jichao@ecust.edu.cn}{jichao@ecust.edu.cn}}

\author[Y. Tao]{Yifan Tao}

\address[Y. Tao]{\newline\indent
	School of Mathematics
	\newline\indent
	East China University of Science and Technology
	\newline\indent
	Shanghai 200237, PR China }
\email{\href{mailto:taoyf2002@126.com}{taoyf2002@126.com}}

\date{\today}

\begin{document}
\begin{abstract}
In this paper, we study positive normalized solutions of the following Schr\"odinger equation with nonlinear point defects on a metric graph $\mathcal{G}=(\mathrm{V},\mathrm{E})$
\begin{equation*}
\begin{cases}
u''=\lambda u & \text{on every }\mathrm{e}\in\mathrm{E},\\
\int_{\mathcal{G}}|u|^2\,dx=\mu & \\
\displaystyle\sum_{\mathrm{e}\succ\mathrm{v}}u'_\mathrm{e}(\mathrm{v})=-|u(\mathrm{v})|^{p-2}u(\mathrm{v}) & \text{at every }\mathrm{v}\in\mathrm{V}_0,\\
\displaystyle\sum_{\mathrm{e}\succ\mathrm{v}}u'_\mathrm{e}(\mathrm{v})=0 & \text{at every }\mathrm{v}\in\mathrm{V}\setminus\mathrm{V}_0,
\end{cases}
\end{equation*}
where $p>2$, $\varnothing\neq\mathrm{V}_0\subset\mathrm{V}$ is possibly infinite, and $\lambda\in\mathbb{R}$ is unknown. For $p\neq4$, every nonempty finite $W\subset\mathrm{V}_0$ gives $2^{\#W}-1$ positive solutions for every sufficiently small mass $\mu>0$ if $p>4$ and every sufficiently large mass $\mu>0$ if $2<p<4$. If $\#\mathrm{V}_0<+\infty$ and $\mathcal{G}$ is noncompact with the bottom of the spectrum of the Laplacian $\inf\sigma(-d^2/dx^2)=0$, there are exactly $2^{\#\mathrm{V}_0}-1$ positive solutions in these regimes, provided $|\overline{B}_r(o)|=O(r^{4/p})$ as $r\to+\infty$ for some $o\in\mathcal{G}$ when $2<p<4$. If ${\#\{\mathrm{v}\in\mathrm{V}_0:\deg(\mathrm{v})\leq D\}=+\infty}$ for some $D\geq1$, infinitely many positive solutions exist for every sufficiently large mass $\mu>0$ if $2<p<4$ and every mass $\mu>0$ if $p>4$. For $p=4$, under local conditions for $W$, analogous results hold for masses $\mu>0$ slightly below $kq^2/2$, where $1\leq k\leq\#W$ and $\deg(\mathrm{v})=q$ for every $\mathrm{v}\in W$. Finally, if $p>4$, $\#\mathrm{V}_0<+\infty$ and $\mathcal{G}$ has a half-line, every mass $\mu>0$ admits a positive solution, and for $p=4$, we construct a metric tree on which every mass $\mu>0$ admits infinitely many positive solutions.
\end{abstract}
\maketitle

\section{Introduction}
In this paper,  we consider a connected metric graph $\G = (\V, \E)$, where $\E$ is the set of edges and $\V$ is the set of vertices. Each bounded edge $\e$ is identified with a closed and bounded interval $I_\e = [0, \ell_\e]$ with $\ell_\e > 0$, while each unbounded edge $\e$ is identified with a copy of $I_\e = \mathbb{R}^+ = [0, +\infty)$ and is called a half-line.

Following Dovetta, Jeanjean and Serra \cite[Section~1]{Do0}, we denote by $\mathbf{G}$ the class of metric graphs $\G=(\V,\E)$ with the following properties:
\begin{itemize}
\item $\G$ is connected, has no self-loops, and has at most countably many edges;
\item $\deg(\vv)<+\infty$ for every $\vv\in\V$;
\item $\ell_*:=\inf_{\e\in\E}|\,\e\,|>0$, where $|\e|=\ell_\e$ is the length of a bounded edge and $|\e|=+\infty$ for a half-line; if there is no bounded edge of $\G$, we set $\ell_*:=1$.
\end{itemize}

Throughout this paper, we study a metric graph $\G\in \mathbf G$. A connected metric graph is naturally a locally compact metric space, and for any two points  $x,y \in \G$, the distance $\operatorname{dist}(x,y)$ is the length of a shortest path along the edges. Figure~\ref{figure1} shows a graph in $\mathbf G$ with finitely many edges, three half-lines and a pair of parallel edges.

Quantum graphs (metric graphs equipped with differential operators) arise naturally as simplified models in mathematics, physics, chemistry, and engineering when one considers the propagation of waves of various types through a quasi-one-dimensional (e.g., mesoscale or nanoscale) system modeled by a thin neighborhood of a graph. For further details on quantum graphs, see, for example, \cite{Be}.
\begin{figure}[H]
	\begin{tikzpicture}[xscale=0.4, yscale=0.4]
		\draw (-14,3)--(-11,6);
		\draw (0,8)--(-8,6);
		\draw (-11,6)--(-8,6);
		\draw (-22,7)--(-14,3);
		\draw (0,3)--(-8,3);
		\draw (-11,9)--(-11,6);
		\draw (-14,3)--(-14,6);
		\draw (-11,3)--(-11,6);
		\draw (-8,3)--(-8,6);

		\node at (-14,3) [nodo] {};
		\node at (-14,6) [nodo] {};
		\node at (-11,3) [nodo] {};
		\node at (-11,6) [nodo] {};
		\node at (-8,3) [nodo] {};
		\node at (-8,6) [nodo] {};
		\node at (-11,9) [nodo] {};

		\node at (-22,7.5) {$\infty$};
		\node at (0,8.5) {$\infty$};
		\node at (0,3.5) {$\infty$};

		\draw (-11,10.5) circle (1.5);
		\node at (-11,12) [nodo] {};
	\end{tikzpicture}
	\caption{A noncompact metric graph $\G$ with finitely many edges, no self-loops and a nonempty compact core.}
	\label{figure1}
\end{figure}

Extensive research has been devoted to the existence and multiplicity of normalized solutions to the following nonlinear Schr\"{o}dinger equation on a metric graph $\G=(\V,\E)$
\begin{equation}
\label{eqs}
    -u'' - \lambda u= \abs{u}^{p-2}u\quad \text{on every } \e \in \E.
\end{equation}
In the $L^2$-subcritical ($2<p<6$) and $L^2$-critical ($p=6$) cases, the energy functional associated with equation \eqref{eqs} is bounded from below and coercive under the mass constraint (requiring the mass to be below a certain threshold when $p=6$). This structure allows for the use of minimization methods to obtain normalized ground states, i.e., solutions with minimal energy under the mass constraint. The existence of such normalized ground states has been established both on compact metric graphs \cite{CDS, Do} and on noncompact metric graphs \cite{ABD, ACFN, AST, AST2, AST3, Do1, NP, PS}. In addition, the existence of local minimizers was studied by Pierotti, Soave and Verzini in \cite{PSV}. For a noncompact metric graph $\G$ with finitely many edges and a nonempty compact core $\K$, the following modification of equation \eqref{eqs}, in which the nonlinearity acts only on $\K$, was introduced in \cite{Gn,No}
\begin{equation}
\label{eqsloc}
    -u'' - \lambda u= \chi_\K\abs{u}^{p-2}u \quad \text{on every } \e \in \E,
\end{equation}
where $\chi_\K$ is the characteristic function of the compact core $\K$. For the existence and nonexistence of normalized ground states for equation \eqref{eqsloc} in the $L^2$-subcritical case, we refer to \cite{ST, T}; for the $L^2$-critical case, see \cite{DT, DT2}. Moreover, in the $L^2$-subcritical case, Serra and Tentarelli \cite{ST1}, by means of genus theory, proved the existence of multiple bound states with negative energy. 

In the $L^2$-supercritical ($p>6$) case, the energy functional is no longer bounded from below under the mass constraint. Moreover, the scaling technique---commonly used in the analysis in $\mathbb{R}^N$---is no longer applicable, and no corresponding Pohozaev identity is available, which creates further complications in proving the boundedness of Palais--Smale sequences. In addition, the topology of metric graphs further complicates the problem. To address these challenges, several works have focused on the $L^2$-supercritical case. In \cite{Bort0}, Borthwick, Chang, Jeanjean and Soave proved the existence of bounded Palais--Smale sequences carrying Morse index type information for functionals with minimax geometry under the mass constraint. In \cite{Cha}, Chang, Jeanjean, and Soave studied the nonlinear Schr\"{o}dinger equation \eqref{eqs} on compact metric graphs. By combining blow-up analysis, monotonicity methods, and an analysis of Morse index type information for Palais--Smale sequences and bound states, the authors proved the existence of normalized nonconstant bound states for equation \eqref{eqs} when the prescribed mass is sufficiently small. Subsequently, in \cite{Bort}, for arbitrary prescribed mass, Borthwick, Chang, Jeanjean and Soave investigated the existence of normalized solutions for \eqref{eqsloc} on noncompact metric graphs with finitely many edges and a nonempty compact core $\K$. The multiplicity of normalized solutions to equation \eqref{eqsloc} was proved in \cite{Ca} by employing a more general blow-up analysis. In \cite{Do0}, Dovetta, Jeanjean and Serra studied the existence of normalized solutions for the nonlinear Schr\"{o}dinger equation \eqref{eqs} on periodic metric graphs, as well as on noncompact metric graphs with finitely many edges under suitable topological assumptions.

Very recently, Boni, Dovetta, and Serra \cite{BDS} were the first to investigate the existence of normalized ground states for the following Schr\"odinger equation on a metric graph $\G=(\V,\E)$ with nonlinear point defects
\begin{equation} \label{eqpdbds}
    \begin{cases} 
u'' = \lambda u & \text{on every } \e \in \E, \\ 
\int_\G \abs{u}^2\, dx = \mu & \\
\displaystyle \sum_{\e \succ \vv}  u'_\e(\vv) = -\abs{u(\vv)}^{p-2}u(\vv) & \text{at every } \vv \in \V_0, \\ 
\displaystyle \sum_{\e \succ \vv}  u'_\e(\vv) = 0 & \text{at every } \vv \in \V \setminus \V_0,
\end{cases}
\end{equation}
where $\V_0\subset \V$ is nonempty, $\mu>0$ is prescribed, $p \in (2,4)$ (the $L^2$-subcritical case for equation \eqref{eqpdbds}), and $\lambda \in \mathbb{R}$ is an unknown Lagrange multiplier. For metric graphs with finitely many vertices, the authors proved the existence of ground states for every $\mu > 0$ and $p \in (2,4)$.  For metric graphs with infinitely many vertices, focusing on $\mathbb{Z}$-periodic graphs and the two-dimensional square grid, the authors proved a series of results unraveling nontrivial threshold phenomena with respect to both the mass $\mu>0$ and the exponent $p\in(2,4)$. Their study highlights that  the ground state problem depends critically on the degree of periodicity of the metric graph, on the total number of point defects, and on their arrangement on the metric graph. For further results concerning the Schr\"{o}dinger equations with nonlinear point defects on metric graphs or in dimension one, we refer the reader to \cite{BD,BD2,BBDT,BBDT2} and the references therein.

The condition
\[
\sum_{\e \succ \vv} u'_\e(\vv) = 0,\quad \forall \vv \in \V \setminus \V_0
\]
is referred to as the Kirchhoff or natural boundary condition, while
\[
\sum_{\e \succ \vv} u'_\e(\vv) = -\abs{u(\vv)}^{p-2}u(\vv), \quad \forall \vv \in \V_0,
\]
can be interpreted as representing the effect of a deep attractive potential well or a strong attractive defect at each vertex $\vv \in \V_0$. In the literature, such a condition is commonly termed a nonlinear $\delta$-interaction and is regarded as modeling strongly localized, point-like defects or inhomogeneities in the propagation medium (see, e.g., \cite{ABR} for a comprehensive overview of non-Kirchhoff vertex conditions on graphs). Models with concentrated nonlinearities have been proposed, for instance, in semiconductor theory \cite{JPS, Ni}, to describe quantum dynamics in resonant tunneling diodes, as well as the effect of the confinement of charges in small regions.

Formally, solutions to problem \eqref{eqpdbds} correspond to standing waves of the following nonlinear Schr\"odinger equation
\[
i\partial_t \psi + \partial_{xx} \psi + \sum_{\vv \in \V_0} \delta_\vv \abs{\psi}^{p-2} \psi = 0 \quad \text{on } \mathcal{G}
\]
via the ansatz $\psi(t,x) = e^{i\lambda t} u(x)$. The constraint on the $L^2$-norm of $u$ in \eqref{eqpdbds},  which reflects the normalization requirement imposed on the solutions, is a standard condition on the mass of $u$ (sometimes interpreted as the number of particles in the condensate)---a quantity that remains preserved along the evolution in time.

We note that \cite{BDS}  studied the existence of normalized solutions to equation \eqref{eqpdbds} only in the $L^2$-subcritical ($2<p<4$) case. It is therefore natural to investigate the existence of normalized solutions to equation \eqref{eqpdbds} in the $L^2$-supercritical ($p>4$) case and the multiplicity of normalized solutions to equation \eqref{eqpdbds}. To the best of our knowledge, the $L^2$-supercritical problem and the multiplicity of  solutions to \eqref{eqpdbds} have not yet been studied on general noncompact metric graphs. In this paper, we take a first step toward filling this gap.

As in the cases of equations \eqref{eqs} and \eqref{eqsloc}, by \cite[Lemma~2.2]{BDS}, in the $L^2$-subcritical ($2<p<4$) and $L^2$-critical ($p=4$) cases, the energy functional associated with \eqref{eqpdbds} is bounded from below and coercive under the mass constraint (requiring the mass to be below a certain threshold when $p=4$), so minimization methods are available. For $2<p<4$, they give normalized ground states for every $\mu>0$ on graphs with finitely many vertices \cite{BDS}, whereas on some $\mathbb Z$-periodic graphs with infinitely many defects no ground state exists at any mass (see Remark~\ref{ac-no-ground-state-solutions}). The energy functional is unbounded from below under the mass constraint for sufficiently large $\mu>0$ in the $L^2$-critical ($p=4$) case and for every $\mu>0$ in the $L^2$-supercritical ($p>4$) case. This follows from the test functions $w^t$ in the proof of Lemma~\ref{lem3.1} with $\rho=1$ and rules out direct minimization in these regimes. Moreover, in the $L^2$-supercritical minimax argument, the lack of a global mass-preserving scaling and a corresponding Pohozaev identity on general metric graphs makes it difficult to obtain uniform bounds for the normalized solutions of the parameterized problems.

For the multiplicity results, we impose no periodicity on $\G$ or $\V_0$. Even on a periodic metric graph $\G$, translations need not leave $\V_0$ invariant and may change the energy. When $\V_0$ is infinite, the vertex restriction from $H^1(\G)$ to $\ell^p(\V_0)$ need not be compact, so weak $H^1(\G)$ convergence alone need not give convergence of the nonlinear term $\sum_{\vv\in\V_0}\abs{u(\vv)}^p$. Variational methods that require compactness of this restriction or compactness up to translations therefore do not apply directly in our setting for the multiplicity results.

Motivated by \cite{BDS,Bort,Ca,AW}, in this paper, we study the existence and multiplicity of normalized solutions to the following Schr\"odinger equation with nonlinear point defects on a metric graph $\G=(\V,\E)$
\begin{equation}\label{eqpd}
        \begin{cases} 
u'' = \lambda u & \text{on every } \e \in \E, \\ 
\int_\G \abs{u}^2\, dx = \mu & \\
\displaystyle \sum_{\e \succ \vv}  u'_\e(\vv) = -\abs{u(\vv)}^{p-2}u(\vv) & \text{at every } \vv \in \V_0, \\ 
\displaystyle \sum_{\e \succ \vv}  u'_\e(\vv) = 0 & \text{at every } \vv \in \V \setminus \V_0,
\end{cases}
\end{equation}
where $p>2$, $\G\in\mathbf{G}$, $\V_0 \subset \V$ is nonempty, $\mu>0$ is the prescribed mass, $\lambda$ is an unknown Lagrange multiplier, $\e \succ \vv$ means that the edge $\e$ is incident at $\vv$, and the notation $u'_\e(\vv)$ stands for $u'_\e(0)$ or $-u'_\e(\ell_\e)$, according to whether the vertex $\vv$ is identified with $0$ or $\ell_\e$. Throughout this paper, we assume that $\G$ belongs to $\mathbf G$ and has no self-loops. A self-loop may be divided at its midpoint into two parallel edges by inserting a degree-two Kirchhoff vertex outside $\V_0$. This leaves the equation on the metric graph unchanged, and the new lower edge-length bound is at least $\ell_*/2$. Parallel edges and half-lines are allowed.
A solution $(u,\lambda)$ of \eqref{eqpd} is called positive if $u(x)>0$ for every $x\in\G$.

Since every edge $\e$ is identified with the interval $I_\e$, a function $u: \G \to \mathbb{R}$ is actually a family of functions $(u_\e)$, where $u_\e: I_\e \to \mathbb{R}$ is the restriction of $u$ to the edge $\e$. The usual $L^p$ space on the metric graph  $\G$ is defined as follows:
\[
L^p(\mathcal{G}):=\bigoplus_{\e \in \mathrm{E}} L^p(I_\e),
\]
endowed with the norm
\[
\norm{u}_{L^p(\G)}^p := \sum_{\e \in \mathrm{E}}\norm{u_\e}_{L^p(I_\e)}^p, \quad \text { if } p \in[1, \infty) \quad \text { and } \quad \norm{u}_{L^{\infty}(\G)}:=\sup_{\e \in \mathrm{E}}\norm{u_\e}_{L^{\infty}(I_\e)}.
\]
The inner product in $L^2(\mathcal{G})$ is
\begin{equation*}
(u,v)_2= \int_\G uv\,dx.
\end{equation*}
If $u:\G \to\mathbb{R}$ is continuous on $\G$, then we write $u \in C(\G)$. The Sobolev space $H^1(\G)$ consists of all continuous functions $u : \G \to \mathbb{R}$ such that $u_\e \in H^1([0,\ell_\e])$ for every edge $\e\in\E$ (where $\ell_\e$ may equal $+\infty$) and $\sum_{\e\in\E}(\norm{u_\e}_2^2+\norm{u'_\e}_2^2)<+\infty$; the norm in $H^1(\G)$ is defined as
\begin{equation*}
    \norm{u}^2_{H^1(\G)} = \underset{\e \in \E}\sum \left(\norm{u'}^2_{L^2(I_\e)}+\norm{u}^2_{L^2(I_\e)}\right),
\end{equation*}
and the inner product in this space is
\begin{equation*}
(u,v)= \underset{\e \in \E}\sum \int_{I_\e} \left(u'_{\e}v'_{\e}+  u_{\e}v_{\e}\right)\,dx.
\end{equation*}
For convenience, we will abbreviate $\norm{u}_{L^p(\G)}$ as $\norm{u}_p$ and $\norm{u}_{H^1(\G)}$ as $\norm{u}$.

We denote by $E( \cdot,\G): H^1(\G) \to \mathbb{R}$ the energy functional associated with equation \eqref{eqpd}:
\begin{equation}\label{func}
E(u,\G) := \frac{1}{2}\norm{u'}^{2}_2-\frac{1}{p}\sum_{\vv \in \V_0}\abs{u(\vv)}^p.
\end{equation}

For $\G\in\mathbf{G}$, the weighted vertex restriction estimate proved in Lemma~\ref{ac-restriction-extension} implies that the defect series in \eqref{func} is finite for every $u\in H^1(\G)$ and that $E(\cdot,\G)$ is a $C^2$ functional, even when $\V_0$ is infinite. We shall study the existence and multiplicity of critical points of the functional $E( \cdot,\G)$ restricted to the $L^2$-sphere
\[
H^1_\mu(\G):=\{u \in H^1(\G): \norm{u}^2_2=\mu\},
\]
where $\mu > 0$.

We now present our main results, Theorems~\ref{th2} and~\ref{cor-ac-infinite-orbits}--\ref{cor-ac-critical-infinite-orbits}. For all of them the following metric grid $\mathcal G^d$ is a common example ($d\geq2$ for Theorems~\ref{th-ac-critical} and~\ref{cor-ac-critical-infinite-orbits}), with $\V_0$ a finite set satisfying the hypotheses of the statement at hand or an infinite set. 

Let $\mathcal{G}^d = (\V^d, \E^d)$ be the $d$-dimensional metric grid with edge length $1$ given by the subset of $\mathbb{R}^d$ with vertices on $\mathbb{Z}^d$ and edges between every pair of vertices at distance $1$. Each edge $\e \in \E^d$ is identified with a closed and bounded interval $I_\e=[0, 1]$.

We first state the multiplicity result of the solutions to \eqref{eqpd} for $p\ne4$. The constrained Morse index $\hat m_0(u)$ is defined in \eqref{ac-constrained-index-definition} and the unconstrained Morse index $m(u)$ is  defined in \eqref{ac-free-index-definition}.
\begin{theorem}\label{th2}
Let $\G\in\mathbf{G}$, let $p\in(2,4)\cup(4,+\infty)$, and let $\varnothing\ne\V_0\subset\V$. For every nonempty finite set $W\subset\V_0$, there exists $\mu_W>0$, depending only on $p$, $\ell_*$, $\#W$ and $\max_{\vv\in W}\deg(\vv)$, such that, for every prescribed mass satisfying
\[
\begin{cases}
0<\mu<\mu_W,&p>4,\\
\mu>\mu_W,&2<p<4,
\end{cases}
\]
problem \eqref{eqpd} has a family of $2^{\#W}-1$ pairwise distinct positive solutions
\[
\bigl\{(u_{\mu,S},\lambda_{\mu,S}): \varnothing\ne S\subset W\bigr\}
\subset H^1_\mu(\G)\times(0,+\infty).
\]
For each nonempty $S\subset W$, we have
\[
m(u_{\mu,S})=\#S,\qquad
\begin{cases}
E(u_{\mu,S},\G)>0,\quad \hat m_0(u_{\mu,S})=\#S,&p>4,\\
E(u_{\mu,S},\G)<0,\quad \hat m_0(u_{\mu,S})=\#S-1,&2<p<4.
\end{cases}
\]
Moreover, $\lambda_{\mu,S}\to+\infty$ as $\mu\to0^+$ if $p>4$, and as $\mu\to+\infty$ if $2<p<4$.
\end{theorem}
Theorem~\ref{th2} applies without periodicity of $\G$ or $\V_0$, even when the vertex restriction from $H^1(\G)$ to $\ell^p(\V_0)$ is noncompact, whereas the existence results of \cite{BDS} rely on one of these two properties. If $\#\V_0<+\infty$, one may take $W=\V_0$, obtaining at least $2^{\#\V_0}-1$ distinct positive solutions $(u,\lambda)$. If $\V_0$ is infinite, then every fixed nonempty finite set $W\subset\V_0$ gives at least $2^{\#W}-1$ positive solutions for every $0<\mu<\mu_W$ when $p>4$, and for every $\mu>\mu_W$ when $2<p<4$.

The following proposition gives sufficient conditions under which the solutions obtained in Theorem~\ref{th2} are all the positive normalized solutions in the stated mass ranges. 
For $p>4$ and finite $\V_0$, this count is exact for small masses on many noncompact graphs, and for $2<p<4$ it is exact for large masses on noncompact graphs whose balls centered at a fixed point have length $O(r^{4/p})$ as $r\to+\infty$, such as graphs with finitely many edges and the grid $\mathcal G^1$.

For $\G\in\mathbf G$, let
\[
 \lambda(\G):=\inf\bigl\{\norm{\varphi'}_2^2: \varphi\in H^1(\G),\ \norm{\varphi}_2=1\bigr\}
\]
be the bottom of the spectrum of the Laplacian on $\G$ with Kirchhoff conditions at every vertex, that is, $\lambda(\G)=\inf\sigma(-d^2/dx^2)$, where $-d^2/dx^2$ denotes the nonnegative self-adjoint operator in $L^2(\G)$ associated with the quadratic form $\varphi\mapsto\norm{\varphi'}_2^2$ on $H^1(\G)$ and $\sigma(-d^2/dx^2)$ denotes the spectrum of $-d^2/dx^2$ (see \cite[Theorem~1.4.19]{Be} and \cite[Section~1]{DDGST} with $Z=\varnothing$). If $\G$ contains a half-line, then $\lambda(\G)=0$; see \cite[Lemma~A.1]{BDS}. The grids $\mathcal G^d$, $d\geq2$, have no half-line, and also $\lambda(\mathcal G^d)=0$ (see \cite[Lemma~A.1 and Remark~A.2]{BDS}).
For $o\in\G$ and $r\geq0$, write
\[
 \overline B_r(o):=\{x\in\G:\operatorname{dist}(x,o)\leq r\},
\]
and let $\abs{\cdot}$ denote the Lebesgue measure on $\G$.
\begin{proposition}\label{prop-ac-mp-single-vertex}
Let $\G\in\mathbf G$ be noncompact with $\lambda(\G)=0$, and let $\varnothing\ne\V_0\subset\V$ be finite. Assume that $p>4$, or that $2<p<4$ and there are $o\in\G$ and $A>0$ such that
\[
 \abs{\overline B_r(o)}\leq A(1+r^{4/p})\qquad\text{for every }r\geq0.
\]
Then there is $\bar\mu=\bar\mu(p,\G,\V_0)>0$ such that, for every $\mu$ with $0<\mu<\bar\mu$ if $p>4$ and with $\mu>\bar\mu$ if $2<p<4$, the positive solutions $(u,\lambda)\in H^1_\mu(\G)\times\R$ of \eqref{eqpd} are exactly the $2^{\#\V_0}-1$ solutions $(u_{\mu,S},\lambda_{\mu,S})$, $\varnothing\ne S\subset\V_0$, of Theorem~\ref{th2} with $W=\V_0$. For $k=1,\ldots,\#\V_0$, exactly $\binom{\#\V_0}{k}$ of them have constrained Morse index $k$ if $p>4$, and $k-1$ if $2<p<4$.
\end{proposition}
For $2<p<4$, the growth assumption in Proposition~\ref{prop-ac-mp-single-vertex} holds on $\mathcal G^1$ with any $o\in\mathcal G^1$ and $A=2$, since $\abs{\overline B_r(o)}=2r\leq2(1+r^{4/p})$.
It also holds on every noncompact $\G\in\mathbf G$ with finitely many edges: choose $o\in\V$ and set $A:=\sum_{\e\in\E_b}\ell_\e+\#\E_\infty$. Each half-line meets $\overline B_r(o)$ in a segment of length at most $r$, so
\[
 \abs{\overline B_r(o)}
 \leq\sum_{\e\in\E_b}\ell_\e+r\#\E_\infty
 \leq A(1+r^{4/p})
 \qquad\text{for every }r\geq0.
\]
On $\mathcal G^d$ with $d\geq2$, fix $o\in\mathbb Z^d$ and, for $r\geq2d$, put $n:=\lfloor r/d\rfloor$. The ball $\overline B_r(o)$ contains all the $dn(n+1)^{d-1}$ edges joining vertices of $o+\{0,1,\ldots,n\}^d$. Since all these edges have length one and $n\geq r/(2d)$, we have
\[
 \abs{\overline B_r(o)}
 \geq dn(n+1)^{d-1}\geq dn^d
 \geq\frac{r^d}{2^d d^{d-1}}
 \qquad\text{for every }r\geq2d,
\]
and hence, since $d>4/p$,
\[
 \frac{\abs{\overline B_r(o)}}{1+r^{4/p}}
 \geq\frac{r^d}{2^d d^{d-1}(1+r^{4/p})}
 \longrightarrow+\infty
 \qquad\text{as }r\to+\infty.
\]
Thus no $A>0$ gives the required upper bound. 

If $p>4$ and $\G$ contains a half-line, Theorem~\ref{th1} gives, by a minimax argument, a positive solution $(u,\lambda)$ with $\hat m_0(u)\leq1$ for every prescribed mass $\mu>0$. By Proposition~\ref{prop-ac-mp-single-vertex}, for sufficiently small mass $\mu>0$, this solution is $(u_{\mu,S},\lambda_{\mu,S})$ with $\#S=\hat m_0(u)\leq1$, that is, $(u_{\mu,\{\vv\}},\lambda_{\mu,\{\vv\}})$ for some $\vv\in\V_0$. Theorem~\ref{th1} is proved in Appendix~\ref{app-th1}, and Proposition~\ref{prop-ac-mp-single-vertex} is proved in Appendix~\ref{app-exact-count}. 

For $2<p<4$ and $d\in\{2,3,4\}$, Proposition~\ref{prop-grid-examples}\textup{(i)} shows that, for every nonempty finite set of defects, problem \eqref{eqpd} has no positive solution for sufficiently small mass $\mu>0$. Thus the restriction to large masses in Theorem~\ref{th2} cannot be removed for $2<p<4$; for $d=2$, this complements \cite[Theorem~1.3]{BDS}, where ground states on $\mathcal G^2$ are shown not to exist below a critical mass. If, in addition, $\V_0=\{0\}$,  Proposition~\ref{prop-grid-examples}\textup{(i)} also gives at least two positive solutions at every sufficiently large mass $\mu>0$, so the growth assumption in Proposition~\ref{prop-ac-mp-single-vertex} cannot be dropped. For $p>4$, $d\in\{2,3,4\}$ and $\V_0=\{0\}$, problem \eqref{eqpd} has a positive solution for every prescribed mass $\mu>0$, although $\mathcal G^d$ contains no half-line. In contrast, for $p>4$, $d\geq5$ and $\V_0=\{0\}$, problem \eqref{eqpd} on $\mathcal G^d$ has no positive solution for large mass $\mu>0$, so the condition that $\G$ has at least one half-line in Theorem~\ref{th1} cannot be replaced by the condition $\lambda(\G)=0$. These examples are given in Proposition~\ref{prop-grid-examples}, where we use the Fourier transform on $\mathbb Z^d$ to represent the vertex values of solutions and estimate their masses.
\begin{remark}\label{ac-no-ground-state-solutions}
For $2<p<4$, Theorem~\ref{th2} does not require the existence of a ground state. Indeed, \cite[Theorem~1.2(iii)]{BDS} provides $\mathbb Z$-periodic graphs with an infinite defect set $\V_0$ for which the infimum of $E(\cdot,\G)$ on $H^1_\mu(\G)$ is not attained for any $\mu>0$. Fix such a graph $\G$ and such a defect set $\V_0$. For every nonempty finite $W\subset\V_0$, Theorem~\ref{th2} gives $2^{\#W}-1$ distinct positive solutions $(u,\lambda)$ for each $\mu>\mu_W$, and
\[
 \inf_{v\in H^1_\mu(\G)}E(v,\G)<E(u_{\mu,S},\G)<0
 \qquad\text{for }\varnothing\ne S\subset W.
\]
The left inequality is strict because equality would make $u_{\mu,S}$ a ground state. These examples therefore have positive normalized solutions $(u,\lambda)$ of \eqref{eqpd} at every sufficiently large mass $\mu>0$ although they have no ground state at any mass $\mu>0$.
\end{remark}

To prove Theorem~\ref{th2}, we reduce \eqref{eqpd} to the vertex equation $A_tf=\mathcal N(f)$ in \eqref{ac-vertex-equation}, where $t=\sqrt\lambda$ and $f=t^{-1/(p-2)}u|_{\V}$. The normalized Dirichlet-to-Neumann operator $A_t$ converges to the identity in operator norm as $t\to+\infty$. For each nonempty $S\subset W$, the vertex branch of Lemma~\ref{ac-IFT-branch} gives solutions $(u_S(t),t^2)$ of \eqref{eqpd} whose rescaled vertex values converge to $f_S^\infty$, defined by $f_S^\infty(\vv)=\deg(\vv)^{1/(p-2)}\mathbf 1_S(\vv)$ in \eqref{ac-limit-profile}. By Proposition~\ref{prop-ac-branches}, $\mu_S(t):=\norm{u_S(t)}_2^2$ eventually decreases strictly to zero if $p>4$ and increases strictly to $+\infty$ if $2<p<4$. The uniform estimate \eqref{ac-two-sided-mass-bound} then gives the common mass ranges $(0,\mu_W)$ and $(\mu_W,+\infty)$, respectively, in the proof of Theorem~\ref{th-ac-uniform-solutions}, with $\mu_W$ depending only on $p$, $\ell_*$, $\#W$ and $\max_{\vv\in W}\deg(\vv)$. For each $\mu$ in this range, solving $\mu_S(t_S(\mu))=\mu$ on the strictly monotone part of each branch gives the solutions $(u_{\mu,S},\lambda_{\mu,S})=(u_S(t_S(\mu)),t_S(\mu)^2)$, which are pairwise distinct by the separation of the limiting profiles. The inverse parameters $t_S(\mu)$ tend to infinity in the stated mass limits, and Proposition~\ref{prop-ac-branches} gives the energy signs and the Morse indices, completing the proof of Theorem~\ref{th2}.

Normalized multibump solutions and their constrained and unconstrained Morse indices were studied by Ackermann and Weth \cite[Theorems~1.2--1.3 and Section~4]{AW} for $-\Delta u+V(x)u-f(u)=\lambda u$ in $\mathbb R^N$ with a periodic potential. In particular, the semiclassical application in \cite[Theorem~1.5 and Corollary~1.6]{AW} treats the pure power $f(u)=\abs{u}^{p-2}u$ with $p\ne2+4/N$ and $2<p<2^*$, and a positive $C^2$ potential $V$ that is $1$-periodic in each coordinate and has a nondegenerate critical point. Here $2^*:=2N/(N-2)_+$, where we set $2^*=+\infty$ if $N=1,2$. For each prescribed mass $\alpha>0$, there is a sequence $\varepsilon_n\to0$ such that $-\varepsilon_n^2\Delta u+V(x)u-\abs{u}^{p-2}u=\lambda u$ in $\R^N$, with $\norm u_2^2=\alpha$, has infinitely many geometrically distinct positive normalized multibump solutions for every $n$ large enough. Their construction glues well-separated translates of a prescribed fully nondegenerate normalized solution and also gives orbital instability.  In our vertex problem the limiting profiles $f_S^\infty$ have an explicit invertible linearization $L_S^\infty$ (see \eqref{ac-limit-linearization} and \eqref{ac-limit-inverse-coordinate}), and no nondegeneracy assumption on an existing normalized solution is needed.  In the abstract gluing theorem of \cite[Theorem~1.2 and Section~3]{AW}, the total mass $\alpha$ is prescribed, and in the semiclassical application $\varepsilon_n$ is chosen so that the single-bump solution satisfies $\norm{\bar u_{\varepsilon_n}}_2^2=\alpha/n$ \cite[Theorem~1.5 and Corollary~1.6]{AW}. In our construction, we first solve the vertex equation at fixed $t^2$, the Lagrange multiplier in \eqref{eqpd}, and then vary $t$ to attain the prescribed mass by inverting $\mu_S(t)$ on a common mass interval.

The count in Theorem~\ref{th2} refers to functions on the fixed graph. A graph automorphism preserving $\V_0$ may map one solution $u_{\mu,S_1}$ to another $u_{\mu,S_2}$ for distinct $S_1,S_2 \subset W$, so these solutions may represent fewer than $2^{\#W}-1$ equivalence classes. We now define the action of these automorphisms on solutions and the corresponding equivalence relation used in Theorems~\ref{cor-ac-infinite-orbits}--\ref{cor-ac-critical-infinite-orbits}, which give infinitely many pairwise inequivalent solutions for each prescribed mass $\mu>0$ in the stated range, provided that $\V_0$ contains infinitely many vertices of uniformly bounded degree.

Let $\Gamma$ be the group of metric graph automorphisms $\gamma:\G\to\G$ which map vertices to vertices, preserve the relation $\e\succ\vv$, map each edge isometrically onto an edge, and satisfy $\gamma(\V_0)=\V_0$. If an edge $\e$ joins $\vv$ and $\w$, then $\gamma(\e)$ joins $\gamma(\vv)$ and $\gamma(\w)$ and satisfies $\ell_{\gamma(\e)}=\ell_\e$. Its action on functions on $\G$ is
\[
 (\gamma u)(x):=u(\gamma^{-1}x).
\]
Every $\gamma\in\Gamma$ maps solutions of \eqref{eqpd} into solutions of the same \eqref{eqpd}, with the same multiplier and the same mass. This gives the action $\gamma\cdot(u,\lambda):=(\gamma u,\lambda)$ on $H^1_\mu(\G)\times(0,+\infty)$. Two solutions $(u,\lambda)$ and $(v,\nu)$ are called $\Gamma$-equivalent if $(u,\lambda)=(\gamma v,\nu)$ for some $\gamma\in\Gamma$. For group actions on periodic metric graphs and the corresponding action on functions, see \cite[Definition~4.1.1 and Section~4.2.3]{Be}. Compared with \cite[Definition~4.1.1]{Be}, where the group acting is $\mathbb Z^n$ on a periodic graph and continuity, faithfulness, discreteness and co-compactness are required in addition, we require only a group action that preserves the graph structure.

Theorem~\ref{cor-ac-infinite-orbits} follows from Theorem~\ref{th-ac-uniform-solutions} by applying it to two-vertex supports with distinct diameters.
\begin{theorem}\label{cor-ac-infinite-orbits}
Let $\G\in\mathbf G$ and let $p\in(2,4)\cup(4,+\infty)$. Suppose that, for some integer $D\geq1$,
\[
 \#\{\vv\in\V_0:\deg(\vv)\leq D\}=+\infty.
\]
Then there exists $\mu_*>0$, depending only on $p,\ell_*$ and $D$, such that, for every prescribed mass satisfying
\[
\begin{cases}
0<\mu<\mu_*,&p>4,\\
\mu>\mu_*,&2<p<4,
\end{cases}
\]
problem~\eqref{eqpd} has infinitely many pairwise $\Gamma$-inequivalent positive solutions 
\[(u_{\mu,n},\lambda_{\mu,n})\in H^1_\mu(\G)\times(0,+\infty).\]
For every $n\in\mathbb N^+$, we have
\[
 m(u_{\mu,n})=2,\qquad
 \begin{cases}
 E(u_{\mu,n},\G)>0,\quad\hat m_0(u_{\mu,n})=2,&p>4,\\
 E(u_{\mu,n},\G)<0,\quad\hat m_0(u_{\mu,n})=1,&2<p<4.
 \end{cases}
\]
Moreover,
\[
 \inf_{n\geq1}\lambda_{\mu,n}\longrightarrow+\infty
 \quad\text{as }\mu\to0^+\text{ if }p>4,\quad
 \text{and as }\mu\to+\infty\text{ if }2<p<4 .
\]
\end{theorem}
The infinite-defect assumption in Theorem~\ref{cor-ac-infinite-orbits} already implies that $\G$ is noncompact.

For $p>4$, allowing the number of selected vertices to increase extends the infinite multiplicity conclusion to every prescribed mass $\mu>0$.
\begin{theorem}\label{cor-ac-supercritical-all-mass}
Let $\G\in\mathbf G$ and let $p>4$. Suppose that, for some integer $D\geq1$,
\[
 \#\{\vv\in\V_0:\deg(\vv)\leq D\}=+\infty .
\]
For every $\mu>0$, there is an integer $n_\mu\geq1$ such that problem \eqref{eqpd} has a sequence of pairwise $\Gamma$-inequivalent positive solutions $(u_{\mu,n},\lambda_{\mu,n})\in H^1_\mu(\G)\times(0,+\infty)$, indexed by $n\geq n_\mu$, with
\[
 E(u_{\mu,n},\G)>0,\qquad
 m(u_{\mu,n})=\hat m_0(u_{\mu,n})=n,
\]
and
\[
 \lambda_{\mu,n}\longrightarrow+\infty\quad\text{as }n\to+\infty .
\]
Moreover,
\[
 \inf_{n\geq n_\mu}\lambda_{\mu,n}\longrightarrow+\infty
 \qquad\text{as }\mu\to0^+.
\]
\end{theorem}
The proof of Theorem~\ref{cor-ac-supercritical-all-mass} uses a frequency threshold $\widehat T$ independent of the number of selected vertices. For the chosen sets $S$, the masses at that threshold $\mu_S(\widehat{T})$ grow with their cardinality $\#S$, giving the full mass range for $p>4$ (see Remark~\ref{rem-ac-all-mass-limitation}).

We now turn to the $L^2$-critical ($p=4$) case. For $p\ne4$, the leading power $t^{(4-p)/(p-2)}$ in \eqref{ac-mass-asymptotic} determines the eventual monotonicity of the mass curves through \eqref{ac-mass-derivative-asymptotic}. However, at $p=4$, the exponent $(4-p)/(p-2)$ vanishes, and Proposition~\ref{prop-ac-critical-branches} gives $\mu_S(t)\to\mu_S^\infty:=\frac12\sum_{\vv\in S}\deg(\vv)^2$ as $t\to+\infty$. This finite limit depends on $S$ and does not determine the sign of $\mu_S'(t)$. To obtain a common interval below this limiting mass, we select supports with the same limiting mass and show that their mass curves increase to it. 

For $S\subset\V$ and a bounded edge $\e\in\E$, we say that $\e$ is incident at $S$ if at least one of its endpoints belongs to $S$. For vertices $\vv,\w\in\V$, let $\mathcal P(\vv,\w)$ be the set of finite edge paths joining them, including the empty path when $\vv=\w$.  We define the combinatorial distance by
\[
 \operatorname{dist}_{\rm c}(\vv,\w)
 :=\min_{(\e_1,\ldots,\e_m)\in\mathcal P(\vv,\w)}m.
\]
For $W\subset\V$, its set of adjacent vertices outside $W$ is denoted by
\[
 N(W):=\{\w\in\V\setminus W:\w\sim\vv\text{ for some }\vv\in W\},
\]
where $\w\sim\vv$ means that an edge joins the two vertices.
\begin{theorem}\label{th-ac-critical}
Let $\G\in\mathbf G$ be simple and assume that $\G$ is not a star graph whose edges are all half-lines. Let $p=4$ and let $\varnothing\ne\V_0\subset\V$. Suppose that, for some integers $q\geq3$, $k\geq1$ and  some $\ell_+\geq\ell_*$, a  finite set $W\subset\V_0$ satisfies $\#W\geq k$,
\[
 \begin{gathered}
 \deg(\vv)=q \qquad\text{for every }\vv\in W,\\
 \deg(\w)\geq3 \qquad\text{for every }\w\in N(W),\\
 \ell_\e\leq\ell_+ \qquad\text{for every }\text{bounded edge }\e\text{ incident at }W,\\
 \operatorname{dist}_{\rm c}(\vv,\w)\geq3
 \qquad\text{for all distinct }\vv,\w\in W.
 \end{gathered}
\]
There is $\eta_{q,k}\in(0,kq^2/2)$, depending only on $q$, $k$, $\ell_*$ and $\ell_+$, such that, for every prescribed mass satisfying $\frac{kq^2}{2}-\eta_{q,k}<\mu<\frac{kq^2}{2}$, problem \eqref{eqpd} has a family of $\binom{\#W}{k}$ pairwise distinct positive solutions
\[
 \bigl\{(u_{\mu,S},\lambda_{\mu,S}): S\subset W,\ \#S=k\bigr\}
 \subset H^1_\mu(\G)\times(0,+\infty).
\]
For each $S\subset W$ with $\#S=k$, we have
\[
 E(u_{\mu,S},\G)>0,\qquad
 m(u_{\mu,S})=k,\qquad \hat m_0(u_{\mu,S})=k-1.
\]
Moreover, $ \lambda_{\mu,S}\longrightarrow+\infty$ as $\mu\to(\frac{kq^2}{2})^-$.
\end{theorem}
Here simple means that parallel edges are also excluded (the expansions in Proposition~\ref{prop-ac-critical-local-sign} allow parallel edges, and part~\textup{(i)} of Proposition~\ref{prop-ac-critical-sign-conditions} uses simplicity of $\G$ to determine the sign of $\mu_S'$) and $\binom{n}{k}:=n!/[k!(n-k)!]$ for integers $n\geq k\geq0$, with $0!:=1$. The multiplicity asserted in Theorem~\ref{th-ac-critical} counts the $k$-element subsets of $W$:
\[
 \#\{S\subset W:\#S=k\}=\binom{\#W}{k}.
\]
For $p\ne4$, the branches indexed by the nonempty subsets of $W$ in Theorem~\ref{th2} have the same limiting mass by \eqref{ac-mass-asymptotic}, namely zero if $p>4$ and $+\infty$ if $2<p<4$. At $p=4$, Remark~\ref{rem-ac-critical-fixed-multiplier} still gives $2^{\#W}-1$ distinct solutions at a common large multiplier, but their masses may differ. Indeed, when the vertices of $S$ have common degree $q$, Proposition~\ref{prop-ac-critical-branches} gives $\mu_S(t)\to q^2\#S/2$. Theorem~\ref{th-ac-critical} fixes $\#S=k$ such that the selected branches share the limiting mass $kq^2/2$, and inverts their increasing mass curves to obtain at least $\binom{\#W}{k}$ solutions at each mass in a common interval.

On a general metric graph, the mass curve $\mu_S(t)$ in Proposition~\ref{prop-ac-critical-branches} need not be invertible near its limiting mass $\mu_S^\infty$. For every integer $q\geq1$, let $\G_q$ denote the star graph consisting of $q$ half-lines whose origins are identified.
For every $p>2$ and $t>0$, on $\G=\G_q$ with a single defect at the center, $(u(t),t^2)$ is a positive solution of \eqref{eqpd} and satisfies
\begin{equation}\label{ac-star-single-defect-branch}
 \begin{aligned}
 u_\e(t,x)&=(qt)^{1/(p-2)}e^{-tx},\qquad x\geq0,\\
 \norm{u(t)}_2^2&=\frac{q^{p/(p-2)}}{2}\,t^{(4-p)/(p-2)},\\
 E(u(t),\G)&=\frac{p-4}{4p}\,q^{p/(p-2)}t^{p/(p-2)}.
 \end{aligned}
\end{equation}
If $p=4$, then $\norm{u(t)}_2^2=q^2/2$ and $E(u(t),\G)=0$ for every $t>0$. Thus one critical mass carries an entire branch of solutions with different multipliers.

Theorem~\ref{th-ac-critical} treats the case $\mu_S'(t)>0$, whereas Proposition~\ref{prop-ac-critical-sign-conditions} also gives sufficient conditions for $\mu_S'(t)<0$. For instance, if $\G$ is a star graph with at least one bounded edge and a point defect at its center and $S$ consists of that center, then $\mu_S'(t)<0$ and $\mu_S(t)$ decreases to $q^2/2$ from above (see Remark~\ref{rem-ac-critical-star}); if $\G=\mathcal G^d$ with $d \geq 1$ and $S$ is a pair of adjacent vertices, then again $\mu_S'(t)<0$, with $\mu_S(t)$ decreasing to $(2d)^2$ from above (see Remark~\ref{rem-ac-critical-adjacent-supports}). 

Mayteevarunyoo, Malomed and Dong \cite[Section~III]{MMD} explicitly constructed symmetric and asymmetric stationary states on $\R$ with cubic point defects at $-1$ and $1$. For $\V_0=S=\{-1,1\}$ on $\mathcal G^1$, Remark~\ref{rem-ac-critical-adjacent-supports} gives $\mu_S(t)\downarrow4$ as $t\to+\infty$, as also follows from their formula for the symmetric branch after the change of variables $u=\sqrt2\,\phi$. In our proof, this direction of approach follows from the sign of $\Xi_S$ in Proposition~\ref{prop-ac-critical-sign-conditions}.

To the best of our knowledge, Theorem~\ref{th-ac-critical} is the first result to establish the existence of a positive normalized solution $(u,\lambda)$ for every prescribed mass in a nonempty open interval for  problem \eqref{eqpd} at the $L^2$-critical exponent $p=4$ on metric grids $\mathcal G^d$, $d\geq2$. For the standard edgewise NLS with Kirchhoff conditions on $\mathcal G^2$, the positive action ground states at the $L^2$-critical exponent $p=6$ have masses tending to infinity as the frequency tends to zero \cite[Section~2 and Proposition~5.4]{DoGrid}.

In the proof of Theorem~\ref{th-ac-critical}, at $p=4$, the leading term of $\mu_S(t)$ is independent of $t$, and the first-order terms vanish for the separated sets $S$. Then, the second-order expansions of $\mu_S(t)$ and $\mu_S'(t)$ in Proposition~\ref{prop-ac-critical-local-sign}, together with the uniform remainder estimates, show that the mass curves increase to $kq^2/2$ above a common frequency threshold, so they can be inverted on a common mass interval below $kq^2/2$.

The mass interval in Theorem~\ref{th-ac-critical} is independent of the positions of the selected vertices. This uniformity also gives infinitely many pairwise $\Gamma$-inequivalent solutions in the following theorem.
\begin{theorem}\label{cor-ac-critical-infinite-orbits}
Let $\G\in\mathbf G$ be simple. Let $p=4$ and let $\varnothing\ne\V_0\subset\V$. Suppose that, for some integer $q\geq3$ and some $\ell_+\geq\ell_*$, there is an infinite set $W_\infty\subset\V_0$ such that
\[
 \begin{gathered}
 \deg(\vv)=q \qquad\text{for every }\vv\in W_\infty,\\
 \deg(\w)\geq3 \qquad\text{for every }\vv\in W_\infty
                      \text{ and }\w\in N(\{\vv\}),\\
 \ell_\e\leq\ell_+ \qquad\text{for every }\text{bounded edge }\e\text{ incident at }W_\infty.
 \end{gathered}
\]
For every integer $k\geq2$ and every prescribed mass satisfying \( \frac{kq^2}{2}-\eta_{q,k}<\mu<\frac{kq^2}{2} \), where $\eta_{q,k}$ is as in Theorem~\ref{th-ac-critical}, problem \eqref{eqpd} has infinitely many pairwise $\Gamma$-inequivalent positive solutions
\[(u_{\mu,n},\lambda_{\mu,n})\in H^1_\mu(\G)\times(0,+\infty).\]
For every $n \in \mathbb{N}^+$, we have
\[
 E(u_{\mu,n},\G)>0,\qquad
 m(u_{\mu,n})=k,\qquad \hat m_0(u_{\mu,n})=k-1.
\]
Moreover,
\[
 \inf_{n\geq1}\lambda_{\mu,n}\longrightarrow+\infty
 \qquad\text{as }\mu\to(kq^2/2)^-.
\]
\end{theorem}
We conclude this discussion with a variational multiplicity result for \eqref{eqpd} on compact metric graphs.
\begin{remark}\label{th3}
On a compact metric graph $\G\in\mathbf G$ with $\varnothing\ne\V_0\subset\V$, compactness of $H^1(\G)\hookrightarrow L^2(\G)$ gives additional normalized solutions $(u,\lambda)$ without a positivity requirement on $u$. For $p>4$, as in \cite[Theorem~1.4]{DAHJ}, for every integer $m\geq1$ there is $\mu_m=\mu_m(p,\G,\V_0)>0$ such that every $0<\mu<\mu_m$ admits at least $m$ distinct pairs of normalized solutions, with multipliers below $1$.

For $2<p<4$, every prescribed mass $\mu>0$ admits infinitely many distinct pairs of normalized solutions on a compact graph, with energy tending to $+\infty$. For $p=4$, there is $\mu^*=\mu^*(\G,\V_0)>0$ such that the same conclusion holds for every $0<\mu<\mu^*$. See \cite[Theorem~1.4]{DAHJ} for the corresponding results for edgewise magnetic NLS.

For $p\ne4$, the variational constructions thus yield more solution pairs than the $2^{\#\V_0}-1$ pairs obtained from Theorem~\ref{th2} with $W=\V_0$ on the common mass ranges, taking $m>2^{\#\V_0}-1$ if $p>4$.
\end{remark}

The rest of the paper is organized as follows. Section~\ref{sect5} develops the weighted vertex framework and constructs branches of positive solutions $(u_S(t),t^2)$. Sections~\ref{sect-branches} and \ref{sect-critical} invert the mass curves and prove normalized multiplicity for $p\ne4$ and $p=4$, respectively. Appendix~\ref{app-th1} proves Theorem~\ref{th1}, which gives a positive solution of \eqref{eqpd} for every prescribed mass $\mu>0$ by a minimax argument. Appendix~\ref{app-exact-count} proves Proposition~\ref{prop-ac-mp-single-vertex} and studies the metric grids $\mathcal G^d$, $d\geq2$, mainly with one defect vertex (see Proposition~\ref{prop-grid-examples}). Appendix~\ref{app-ac-critical-all-mass-tree} constructs a metric tree on which, for $p=4$, every prescribed mass carries infinitely many positive solutions (see Proposition~\ref{prop-ac-critical-all-mass-tree}).

\section{The construction of the branches}
\label{sect5}
In this section, we develop the framework used in the subsequent proofs: weighted restrictions and extensions for $H^1(\G)$, the associated Dirichlet-to-Neumann operator, and the construction of branches of positive solutions $(u_S(t),t^2)$ near the limiting profiles.

Throughout Section~\ref{sect5}, we assume that $\G=(\V,\E)\in\mathbf G$, $p>2$ and $\varnothing\ne\V_0\subset\V$.  Write $\E_b$ for the set of bounded edges and $\E_\infty:=\E\setminus\E_b$ for the set of half-lines; recall that $\ell_*$ denotes the uniform lower bound for the lengths of the bounded edges, as introduced in the definition of $\mathbf{G}$.  Given a Lagrange multiplier $\lambda>0$, we write $\lambda=t^2$ with $t>0$; the edge equation in \eqref{eqpd} then reads $-u''+t^2u=0$.

\subsection{Weighted restrictions and harmonic extension}

Let $C(\V)$ be the space of real-valued functions on $\V$.  Define
\[
 X:=\ell^2(\V,\deg)
   =\Bigl\{f\in C(\V):\sum_{\vv\in\V}\deg(\vv)\abs{f(\vv)}^2<\infty\Bigr\},
 \qquad
 \langle f,g\rangle_X:=\sum_{\vv\in\V}\deg(\vv) f(\vv)g(\vv).
\]
The map $f\mapsto\bigl(\sqrt{\deg(\vv)}\,f(\vv)\bigr)_{\vv\in\V}$ is an isometric isomorphism of $X$ onto $\ell^2(\V)$, so $X$ is a Hilbert space.  Since $\deg(\vv)\geq1$,
\begin{equation}\label{ac-norm-ordering}
 \norm{f}_{\ell^\infty(\V)}\leq\norm{f}_{\ell^2(\V)}\leq\norm{f}_X .
\end{equation}
For normed spaces $X_1$ and $X_2$, $\mathcal L(X_1,X_2)$ denotes the space of bounded linear maps from $X_1$ to $X_2$, and $\mathcal L(X_1):=\mathcal L(X_1,X_1)$.

For $t>0$ define the bounded symmetric bilinear form
\[
 \mathfrak b_t:H^1(\G)\times H^1(\G)\to\R,\qquad
 \mathfrak b_t(h,k):=(h',k')_2+t^2(h,k)_2,
\]
where $(h,k)_2:=\int_\G hk\,dx$.  With $\norm{h}_{H^1(\G)}^2=\norm{h}_2^2+\norm{h'}_2^2$ we have
\begin{equation}\label{ac-form-coercivity}
 \mathfrak b_t(h,h)\geq\min\{1,t^2\}\norm{h}_{H^1(\G)}^2 .
\end{equation}
The subspace of functions with zero vertex restriction is
\[
 H^1_0(\G;\V):=\{h\in H^1(\G):h(\vv)=0\text{ for every }\vv\in\V\}.
\]

We first introduce the explicit form of the harmonic extension of prescribed vertex values, which will be used repeatedly.  For a bounded edge $\e\in\E_b$ with endpoints $\vv,\w\in\V$, identify $\e$ with the interval $I_\e=[0,\ell_\e]$ in such a way that $\vv$ is identified with $0$ and $\w$ with $\ell_\e$, and set
\begin{equation}\label{ac-edge-extension}
 (U_tf)_\e(x):=
 f(\vv)\frac{\sinh\bigl(t(\ell_\e-x)\bigr)}{\sinh(t\ell_\e)}
 +f(\w)\frac{\sinh(tx)}{\sinh(t\ell_\e)},
 \qquad x\in[0,\ell_\e].
\end{equation}
Note that the definition of $(U_tf)_\e$ does not depend on the orientation chosen for $\e$.
Likewise, identify a half-line $\e\in\E_\infty$ attached at $\vv\in\V$ with $[0,\infty)$ in such a way that $\vv$ is identified with $0$, and set
\begin{equation}\label{ac-halfline-extension}
 (U_tf)_\e(x):=f(\vv)e^{-tx},
 \qquad x\in[0,\infty).
\end{equation}
Each term in \eqref{ac-edge-extension} solves $-y''+t^2y=0$, and their sum takes the prescribed boundary values $f(\vv)$, $f(\w)$.  Likewise, \eqref{ac-halfline-extension} is a solution of $y''=t^2y$ on $(0,\infty)$ with $(U_tf)_\e(0)=f(\vv)$ and $(U_tf)_\e\in H^1(0,\infty)$. By \cite[Remark~4 in Section~3.3]{AD}, a function $h\in H^1(\G)$ with restriction $f$ and $h''=t^2h$ on every edge is unique, and it is given by \eqref{ac-edge-extension} and \eqref{ac-halfline-extension}.  Moreover, for $\e=\vv\w\in\E_b$, where $\vv$ is identified with $0$ and $\w$ with $\ell_\e$, we have
\begin{equation}\label{ac-edge-derivatives}
 \begin{aligned}
 (U_tf)'_\e(x)
   &=f(\vv)\frac{-t\cosh\bigl(t(\ell_\e-x)\bigr)}{\sinh(t\ell_\e)}
     +f(\w)\frac{t\cosh(tx)}{\sinh(t\ell_\e)},\qquad x\in[0,\ell_\e],\\
 (U_tf)'_\e(\vv)
   &=t\bigl[-\coth(t\ell_\e)f(\vv)+\operatorname{csch}(t\ell_\e)f(\w)\bigr],\\
 (U_tf)'_\e(\w)
   &=t\bigl[\operatorname{csch}(t\ell_\e)f(\vv)-\coth(t\ell_\e)f(\w)\bigr].
 \end{aligned}
\end{equation}

\begin{lemma}\label{ac-restriction-extension}\label{ac-trace-extension}
The restriction map $\mathcal R:H^1(\G)\to X$, $\mathcal R h=h|_\V$, is bounded and onto.  For every $t>0$, \eqref{ac-edge-extension} and \eqref{ac-halfline-extension} define a bounded linear operator $U_t:X\to H^1(\G)$, and $U_tf$ is the unique function satisfying $\mathcal R U_tf=f$ and $(U_tf)''=t^2U_tf$ on every edge.  In particular, $\mathcal R U_t=I_X$.  Moreover,
\begin{equation}\label{ac-extension-orthogonality}
 \mathfrak b_t(U_tf,z)=0\qquad\text{for all }f\in X,\ z\in H^1_0(\G;\V),
\end{equation}
and for every $t_0>0$ there is $C=C(t_0,\ell_*)>0$ such that, for $t\geq t_0$,
\begin{equation}\label{ac-extension-estimates}
 \norm{U_tf}_2^2\leq\frac Ct\norm{f}_X^2,
 \qquad
 \norm{(U_tf)'}_2^2\leq Ct\norm{f}_X^2 .
\end{equation}
The map $t\mapsto U_t$ belongs to $C^1((0,\infty),\mathcal L(X,H^1(\G)))$, and $\partial_tU_tf\in H^1_0(\G;\V)$ for every $f\in X$.
\end{lemma}

\begin{proof}
Fix $0<\delta<\min\{1,\ell_*/2\}$ and $h\in H^1(0,\delta)$.  By the fundamental theorem of calculus and the Cauchy--Schwarz inequality,
\[
 \abs{h(0)}^2
 =\Bigl|h(x)-\int_0^x h'(s)\,ds\Bigr|^2
 \leq 2\abs{h(x)}^2+2x\int_0^x\abs{h'(s)}^2\,ds
 \leq 2\abs{h(x)}^2+2\delta\int_0^\delta\abs{h'}^2\,ds
\]
for a.e. $x\in(0,\delta)$, hence
\[
 \abs{h(0)}^2\leq\frac2\delta\int_0^\delta\abs{h}^2\,dx
   +2\delta\int_0^\delta\abs{h'}^2\,dx .
\]
Summing over all edge endpoints gives
\begin{equation}\label{ac-weighted-restriction-bound}
 \sum_{\vv\in\V}\deg(\vv)\abs{h(\vv)}^2
   \leq \frac2\delta\norm{h}_2^2+2\delta\norm{h'}_2^2 .
\end{equation}
This proves that $\mathcal R$ is bounded on $H^1(\G)$; the surjectivity of $\mathcal R$ follows from the explicit construction of the extension $U_tf$ above.
Integrating by parts with $u:=(U_tf)_\e$ and using \eqref{ac-edge-derivatives},
\begin{equation}\label{ac-edge-energy}
 \begin{aligned}
 \int_\e\bigl(\abs{(U_tf)'}^2+t^2\abs{U_tf}^2\bigr)\,dx
 &=\bigl[u u'\bigr]_0^{\ell_\e}
   =f(\w)u'(\ell_\e)-f(\vv)u'(0)\\
 &=-f(\vv)(U_tf)'_\e(\vv)-f(\w)(U_tf)'_\e(\w)\\
 &= t\Bigl[\coth(t\ell_\e)\bigl(f(\vv)^2+f(\w)^2\bigr)
      -2\operatorname{csch}(t\ell_\e)f(\vv)f(\w)\Bigr],
 \end{aligned}
\end{equation}
while on a half-line $\e\in\E_\infty$, we have
\begin{equation}\label{ac-halfline-energy}
 \int_\e\bigl(\abs{(U_tf)'}^2+t^2\abs{U_tf}^2\bigr)\,dx
 =2t^2f(\vv)^2\int_0^\infty e^{-2tx}\,dx
 =tf(\vv)^2 .
\end{equation}
For a bounded edge $\e=\vv\w\in\E_b$ with $x=f(\vv)$, $y=f(\w)$ and $s=t\ell_\e$, the last line of \eqref{ac-edge-energy} factors as
\[
 \coth s\,(x^2+y^2)-2\operatorname{csch}s\,xy
 =\frac{\coth s-\operatorname{csch}s}{2}(x+y)^2
  +\frac{\coth s+\operatorname{csch}s}{2}(x-y)^2 .
\]
Since $\coth s-\operatorname{csch}s=\tanh(s/2)$ and $\coth s+\operatorname{csch}s=\coth(s/2)$, and $(x+y)^2+(x-y)^2=2(x^2+y^2)$, this gives
\begin{equation}\label{ac-edge-matrix-bounds}
 t\tanh(t\ell_\e/2)\bigl(f(\vv)^2+f(\w)^2\bigr)
 \leq \int_\e\bigl(\abs{(U_tf)'}^2+t^2\abs{U_tf}^2\bigr)\,dx
 \leq t\coth(t\ell_\e/2)\bigl(f(\vv)^2+f(\w)^2\bigr).
\end{equation}
Fix $t_0> 0$. Summing \eqref{ac-edge-matrix-bounds} over the bounded edges, using \eqref{ac-halfline-energy} on the half-lines, and recalling $\ell_\e\geq\ell_*$ gives, for $t\geq t_0$,
\[
 t\tanh(t_0\ell_*/2)\norm{f}_X^2
 \leq\mathfrak b_t(U_tf,U_tf)
 \leq t\coth(t_0\ell_*/2)\norm{f}_X^2 .
\]
Since $\norm{(U_tf)'}_2^2\leq\mathfrak b_t(U_tf,U_tf)$ and $t^2\norm{U_tf}_2^2\leq\mathfrak b_t(U_tf,U_tf)$, the two assertions in \eqref{ac-extension-estimates} follow with $C:=\coth(t_0\ell_*/2)$.  In particular, $U_t:X\to H^1(\G)$ is a bounded linear operator for every $t>0$.  Moreover, since $(U_tf)''=t^2U_tf$, we know that $\mathfrak b_t(U_tf,z)=0$ for any $z\in H^1_0(\G;\V)$, which is \eqref{ac-extension-orthogonality}.

Finally, we prove that $t\mapsto U_t$ belongs to $C^1((0,\infty),\mathcal L(X,H^1(\G)))$.
Let $I\subset(0,\infty)$ be a compact interval and $s,t\in I$, and let $C=C(\min I,\ell_*)$ be the constant of \eqref{ac-extension-estimates}.  For fixed $f\in X$ set $W_{s,t}:=(U_s-U_t)f$.  Since $\mathcal R W_{s,t}=\mathcal R U_sf- \mathcal R U_tf=0$, we have $W_{s,t}\in H^1_0(\G;\V)$, and thus, by \eqref{ac-extension-orthogonality},
\begin{equation}\label{ac-extension-difference}
 \mathfrak b_t(W_{s,t},W_{s,t})=\mathfrak b_t(U_sf,W_{s,t})
 =(t^2-s^2)(U_sf,W_{s,t})_2.
\end{equation}
Since $\min I>0$, $\mathfrak b_t(h,h)\geq\min\{1,(\min I)^2\}\norm{h}_{H^1(\G)}^2$ for all $t\in I$, and \eqref{ac-extension-difference} together with \eqref{ac-extension-estimates} at $t_0=\min I$  gives 
\begin{equation*}
\norm{W_{s,t}}_{H^1(\G)}^2\leq \frac{\sqrt{C}}{\min\{1,t^2\}\sqrt{s}}\abs{t^2-s^2}\norm{f}_X\norm{W_{s,t}}_{H^1(\G)},
\end{equation*}
which implies
\begin{equation}\label{ac-extension-lipschitz}
 \norm{W_{s,t}}_{H^1(\G)}\leq \frac{\sqrt{C}(t+s)}{\min\{1,t^2\}\sqrt{s}}\abs{s-t}\norm{f}_X .
\end{equation}
For fixed $t>0$ and $f \in X$, the Lax--Milgram theorem provides a unique $Z_tf\in H^1_0(\G;\V)$ such that
\begin{equation}\label{ac-extension-derivative}
 \mathfrak b_t(Z_tf,z)=-2t(U_tf,z)_2\qquad \text{for any }z\in H^1_0(\G;\V).
\end{equation}
Taking $z=Z_tf$ in \eqref{ac-extension-derivative} and using \eqref{ac-form-coercivity} and \eqref{ac-extension-estimates} gives $\min\{1,t^2\}\norm{Z_tf}_{H^1(\G)}^2\leq2t\norm{U_tf}_2\norm{Z_tf}_{H^1(\G)}$, and hence
\begin{equation}\label{ac-extension-Z-bound}
 \norm{Z_tf}_{H^1(\G)}\leq\frac{2\sqrt{tC}}{\min\{1,t^2\}}\norm{f}_X .
\end{equation}
For $s\ne t$ in $I$,
\[
 \mathfrak b_t\Bigl(\frac{W_{s,t}}{s-t}-Z_tf,z\Bigr)
 =-(s+t)(U_sf,z)_2+2t(U_tf,z)_2
 =-(s-t)(U_sf,z)_2-2t(W_{s,t},z)_2 ,
\]
so \eqref{ac-form-coercivity} and \eqref{ac-extension-lipschitz} yield
\[
 \Bigl\|\frac{(U_s-U_t)f}{s-t}-Z_tf\Bigr\|_{H^1(\G)}
 \leq\frac{\sqrt{C}\abs{s-t}}{\min\{1,t^2\}\sqrt{s}}\Bigl(1+\frac{2t(t+s)}{\min\{1,t^2\}}\Bigr)\norm{f}_X .
\]
Similarly, since $\mathfrak b_t-\mathfrak b_s=(t^2-s^2)(\cdot,\cdot)_2$,
\begin{equation*}
 \begin{aligned}
 \mathfrak b_t\bigl((Z_s-Z_t)f,z\bigr)
 &=\mathfrak b_t(Z_sf,z)-\mathfrak b_t(Z_tf,z)\\
 &=\mathfrak b_s(Z_sf,z)+(t^2-s^2)(Z_sf,z)_2-\mathfrak b_t(Z_tf,z)\\
 &=-2s(U_sf,z)_2+2t(U_tf,z)_2+(t^2-s^2)(Z_sf,z)_2\\
 &=2(t-s)(U_sf,z)_2+2t\bigl((U_t-U_s)f,z\bigr)_2+(t^2-s^2)(Z_sf,z)_2\\
 &=2(t-s)(U_sf,z)_2-2t(W_{s,t},z)_2+(t^2-s^2)(Z_sf,z)_2,
 \end{aligned}
\end{equation*}
which, together with \eqref{ac-extension-derivative} and \eqref{ac-extension-Z-bound}, leads to 
\begin{equation*}
\norm{(Z_s-Z_t)f}_{H^1(\G)}\leq\frac{\sqrt{C}\abs{s-t}}{\min\{1,t^2\}}\Bigl(\frac{2}{\sqrt{s}}+\frac{2t(t+s)}{\min\{1,t^2\}\sqrt{s}}+\frac{2(t+s)\sqrt{s}}{\min\{1,s^2\}}\Bigr)\norm{f}_X.
\end{equation*}
Hence $\partial_tU_t=Z_t$ in $\mathcal L(X,H^1(\G))$ and $t\mapsto U_t$ is $C^1$.  Since $Z_tf\in H^1_0(\G;\V)$ by construction, this also shows that $\partial_tU_tf\in H^1_0(\G;\V)$ for every $f\in X$.
\end{proof}

By Lemma~\ref{ac-restriction-extension}, every $h\in H^1(\G)$ splits uniquely as
\begin{equation}\label{ac-restriction-decomposition}
 h=z+U_tg,\qquad g=\mathcal R h\in X,\qquad
 z=h-U_t\mathcal R h\in H^1_0(\G;\V),
\end{equation}
with $\mathfrak b_t(z,U_tg)=0$.  Indeed, $\mathcal R U_t=I_X$ produces the prescribed restrictions and \eqref{ac-extension-orthogonality} gives the orthogonality.

Integrating by parts with $(U_tf)''=t^2U_tf$ and substituting the inward derivatives \eqref{ac-edge-derivatives} gives, for $f,g\in X$, $\mathfrak b_t(U_tf,U_tg)=t\,a_t(f,g)$ with
\begin{equation}\label{ac-DtN-form}
 \begin{aligned}
 a_t(f,g)
 :={}&\sum_{\e\in\E_b}
   \coth(t\ell_\e)\bigl[f(\vv_\e)g(\vv_\e)+f(\w_\e)g(\w_\e)\bigr]\\
 &-\sum_{\e\in\E_b}\operatorname{csch}(t\ell_\e)
   \bigl[f(\vv_\e)g(\w_\e)+f(\w_\e)g(\vv_\e)\bigr]
 +\sum_{\e\in\E_\infty}f(\vv_\e)g(\vv_\e),
 \end{aligned}
\end{equation}
where $\vv_\e,\w_\e$ are the endpoints of $\e$ for $\e\in\E_b$ and $\vv_\e$ is the endpoint of $\e$ for $\e\in\E_\infty$.  Moreover, \eqref{ac-edge-matrix-bounds} gives
\begin{equation}\label{ac-DtN-bound}
 \abs{a_t(f,g)}\leq\coth(t\ell_*/2)\norm{f}_X\norm{g}_X .
\end{equation}
By the Riesz representation theorem there is a unique bounded linear operator
\[
 A_t:X\longrightarrow X,\qquad
 \langle A_tf,g\rangle_X=a_t(f,g)\qquad \text{for any }f,g\in X.
\]
Its vertex values are determined in \eqref{ac-DtN-coordinate}, and $a_t(f,g)=a_t(g,f)$ shows that $A_t$ is self-adjoint on $X$.  Moreover, we have
\begin{equation}\label{ac-DtN-identity}
 \mathfrak b_t(U_tf,U_tg)=t\,a_t(f,g)=t\langle A_tf,g\rangle_X \qquad \text{for any }f,g\in X.
\end{equation}

For a bounded edge $\e\in\E_b$ incident at $\vv$, let $\w_{\e,\vv}$ be the other endpoint of $\e$.  For any $f \in X$, testing $ \langle A_tf,g\rangle_X=a_t(f,g)$ with the coordinate function $g=\mathbf 1_{\{\vv\}}$ yields
\begin{equation}\label{ac-DtN-coordinate}
 \begin{aligned}
 \deg(\vv)(A_tf)(\vv)
 &=\langle A_tf,\mathbf 1_{\{\vv\}}\rangle_X\\
 &=a_t(f,\mathbf 1_{\{\vv\}})\\
 &=\sum_{\substack{\e\in\E_b\\\e\succ\vv}}
   \bigl[\coth(t\ell_\e)f(\vv)
    -\operatorname{csch}(t\ell_\e)f(\w_{\e,\vv})\bigr]
   +\sum_{\substack{\e\in\E_\infty\\\e\succ\vv}}f(\vv).
 \end{aligned}
\end{equation}
Substituting \eqref{ac-edge-derivatives} into \eqref{ac-DtN-coordinate} gives, for every vertex $\vv\in\V$,
\begin{equation}\label{ac-DtN-flux}
 (A_tf)(\vv)=-\frac{1}{t\deg(\vv)}\sum_{\e\succ\vv}(U_tf)'_\e(\vv).
\end{equation}
The operator $A_t$ is the normalized Dirichlet-to-Neumann operator of the problem. The factor $1/t$ in \eqref{ac-DtN-flux} removes the factor $t$ in the endpoint derivatives \eqref{ac-edge-derivatives}.

We next prove that $A_t$ is coercive and that $A_t\to I$ in $\mathcal L(X)$ as $t\to\infty$. We also express $\norm{U_tf}_2^2$ in terms of $A_t$ and $\partial_tA_t$, which will be used to compute the mass of the solution branches.

\begin{lemma}\label{ac-DtN-estimates}
The map $t\mapsto A_t$ belongs to $C^2((0,\infty),\mathcal L(X))$.  For $t>0$,
\begin{equation}\label{ac-DtN-coercivity}
 \langle A_tf,f\rangle_X\geq\tanh(t\ell_*/2)\norm{f}_X^2 .
\end{equation}
There are constants $C_j>0$, depending only on $\ell_*$, such that for all $t\geq1$,
\begin{equation}\label{ac-DtN-exponential}
 \norm{A_t-I_X}_{\mathcal L(X)}\leq C_0e^{-t\ell_*/2},
 \qquad
 \norm{\partial_t^jA_t}_{\mathcal L(X)}\leq C_je^{-t\ell_*/2},
 \qquad j=1,2.
\end{equation}
Finally, for every $f\in X$,
\begin{equation}\label{ac-mass-fixed-restriction}
 \norm{U_tf}_2^2=\frac{1}{2t}\bigl\langle (A_t+t\partial_tA_t)f,f\bigr\rangle_X .
\end{equation}
\end{lemma}

\begin{proof}
For $f\in X$, since $\frac1t\mathfrak b_t(U_tf,U_tf)=\langle A_tf,f\rangle_X$ by \eqref{ac-DtN-identity}, the bounds \eqref{ac-edge-matrix-bounds} on the bounded edges and \eqref{ac-halfline-energy} on the half-lines, together with $\ell_\e\geq\ell_*$, the monotonicity of $s\mapsto\tanh(s/2)$ on $(0,\infty)$, and $\tanh(t\ell_*/2)\leq1$, give
\[
 \begin{aligned}
 \frac1t\mathfrak b_t(U_tf,U_tf)
 &=\sum_{\e\in\E_b}\frac1t\int_\e\bigl(\abs{(U_tf)'}^2+t^2\abs{U_tf}^2\bigr)
  +\sum_{\e\in\E_\infty}\frac1t\int_\e\bigl(\abs{(U_tf)'}^2+t^2\abs{U_tf}^2\bigr)\\
 &\geq\sum_{\e\in\E_b}\tanh(t\ell_\e/2)\bigl(f(\vv_\e)^2+f(\w_\e)^2\bigr)
  +\sum_{\e\in\E_\infty}f(\vv_\e)^2\\
 &\geq\tanh(t\ell_*/2)\Bigl[\sum_{\e\in\E_b}\bigl(f(\vv_\e)^2+f(\w_\e)^2\bigr)+\sum_{\e\in\E_\infty}f(\vv_\e)^2\Bigr]
 =\tanh(t\ell_*/2)\norm{f}_X^2,
 \end{aligned}
\]
which is \eqref{ac-DtN-coercivity}.

We next verify that $t\mapsto A_t$ is differentiable in the norm $\norm{\cdot}_{\mathcal L(X)}$. Since $\coth^2(t\ell)-\operatorname{csch}^2(t\ell)=1$, repeated differentiation gives
\begin{equation}\label{ac-DtN-derivatives}
 \begin{aligned}
 \partial_t\coth(t\ell)&=-\ell\operatorname{csch}^2(t\ell),\\
 \partial_t\operatorname{csch}(t\ell)&=-\ell\coth(t\ell)\operatorname{csch}(t\ell),\\
 \partial_t^2\coth(t\ell)&=2\ell^2\coth(t\ell)\operatorname{csch}^2(t\ell),\\
 \partial_t^2\operatorname{csch}(t\ell)&=\ell^2\operatorname{csch}(t\ell)\bigl(1+2\operatorname{csch}^2(t\ell)\bigr),\\
 \partial_t^3\coth(t\ell)&=-2\ell^3\operatorname{csch}^2(t\ell)\bigl(2+3\operatorname{csch}^2(t\ell)\bigr),\\
 \partial_t^3\operatorname{csch}(t\ell)&=-\ell^3\coth(t\ell)\operatorname{csch}(t\ell)\bigl(1+6\operatorname{csch}^2(t\ell)\bigr).
 \end{aligned}
\end{equation}
Let $I\subset(0,\infty)$ be a compact interval and $\sigma:=\min I$. For $t\in I$ and $\ell\geq\ell_*$, $t\ell\geq\sigma\ell_*>0$, so
\begin{equation}\label{ac-coth-csch-bounds}
 \coth(t\ell)-1\leq\frac{2e^{-2\sigma\ell}}{1-e^{-2\sigma\ell_*}},
 \qquad
 \operatorname{csch}(t\ell)\leq\frac{2e^{-\sigma\ell}}{1-e^{-2\sigma\ell_*}},
\end{equation}
and therefore, with $j$ denoting the number of derivatives, the eight quantities $\coth(t\ell)-1$, $\operatorname{csch}(t\ell)$, $\partial_t\coth(t\ell)$, $\partial_t\operatorname{csch}(t\ell)$, $\partial_t^2\coth(t\ell)$, $\partial_t^2\operatorname{csch}(t\ell)$, $\partial_t^3\coth(t\ell)$, $\partial_t^3\operatorname{csch}(t\ell)$ are bounded on $I\times[\ell_*,\infty)$ by $C(\sigma,j)\ell^je^{-\sigma\ell}$, with $C(\sigma,j)$ depending only on $\sigma$, $\ell_*$ and $j$.  For $j=1,2,3$, the derivative of $a_t$ is given by
\[
 \partial_t^j a_t(f,g)=\sum_{\E_b}\Bigl[\partial_t^j\coth(t\ell_\e)\bigl(f(\vv_\e)g(\vv_\e)+f(\w_\e)g(\w_\e)\bigr)-\partial_t^j\operatorname{csch}(t\ell_\e)\bigl(f(\vv_\e)g(\w_\e)+f(\w_\e)g(\vv_\e)\bigr)\Bigr].
\]
Here the series is differentiated term by term. This is justified because, for fixed $f,g\in X$, the second inequality below bounds the $\e$-th term of the differentiated series, uniformly for $t\in I$, by $C(\sigma,j)\ell_\e^je^{-\sigma\ell_\e}(\abs{f(\vv_\e)}+\abs{f(\w_\e)})(\abs{g(\vv_\e)}+\abs{g(\w_\e)})$, and these bounds are summable by the last step; hence the differentiated series converge uniformly on $I$. Then, for $t\in I$ and $f,g\in X$, \eqref{ac-coth-csch-bounds} and \eqref{ac-DtN-derivatives} give
\[
 \begin{aligned}
 \abs{\partial_t^ja_t(f,g)}
 &\leq\sum_{\E_b}\Bigl[\abs{\partial_t^j\coth(t\ell_\e)}\bigl(\abs{f(\vv_\e)}\abs{g(\vv_\e)}+\abs{f(\w_\e)}\abs{g(\w_\e)}\bigr)\\
 &\qquad+\abs{\partial_t^j\operatorname{csch}(t\ell_\e)}\bigl(\abs{f(\vv_\e)}\abs{g(\w_\e)}+\abs{f(\w_\e)}\abs{g(\vv_\e)}\bigr)\Bigr]\\
 &\leq C(\sigma,j)\sum_{\E_b}\ell_\e^je^{-\sigma\ell_\e}\bigl(\abs{f(\vv_\e)}+\abs{f(\w_\e)}\bigr)\bigl(\abs{g(\vv_\e)}+\abs{g(\w_\e)}\bigr)\\
 &\leq C(\sigma,j)\Bigl(\sup_{\ell\geq\ell_*}\ell^je^{-\sigma\ell}\Bigr)\sum_{\E_b}\bigl(\abs{f(\vv_\e)}+\abs{f(\w_\e)}\bigr)\bigl(\abs{g(\vv_\e)}+\abs{g(\w_\e)}\bigr)\\
 &\leq 2C(\sigma,j)\Bigl(\sup_{\ell\geq\ell_*}\ell^je^{-\sigma\ell}\Bigr)\norm f_X\norm g_X,
 \end{aligned}
\]
where the last step uses $(\abs{f(\vv_\e)}+\abs{f(\w_\e)})^2\leq2(\abs{f(\vv_\e)}^2+\abs{f(\w_\e)}^2)$ and $\sum_{\E_b}(\abs{f(\vv_\e)}^2+\abs{f(\w_\e)}^2)\leq\norm f_X^2$.  Writing $\tilde C(\sigma,j):=2C(\sigma,j)\bigl(\sup_{\ell\geq\ell_*}\ell^je^{-\sigma\ell}\bigr)$, we have
\begin{equation}\label{ac-DtN-derivative-bound}
 \abs{\partial_t^ja_t(f,g)}\leq\tilde C(\sigma,j)\norm f_X\norm g_X, \qquad t\in I,\ j=1,2,3,
\end{equation}
so that $\partial_t^ja_t$ is a bounded bilinear form on $X\times X$.
By the Riesz representation theorem there are therefore operators $\partial_tA_t,\partial_t^2A_t\in\mathcal L(X)$ such that $\partial_t^j a_t(f,g)=\langle\partial_t^jA_tf,g\rangle_X$ for all $f,g\in X$.

For $t\in I$ and $h\in\R$ with $t+h\in I$, Taylor's formula with integral remainder, applied to $s\mapsto a_s(f,g)$ and to $s\mapsto\partial_ta_s(f,g)$ on the interval with endpoints $t$ and $t+h$, and the fundamental theorem of calculus applied to $s\mapsto\partial_t^2a_s(f,g)$, give
\[
 \begin{aligned}
 a_{t+h}(f,g)&=a_t(f,g)+h\,\partial_ta_t(f,g)+\int_t^{t+h}(t+h-s)\,\partial_t^2a_s(f,g)\,ds,\\
 \partial_ta_{t+h}(f,g)&=\partial_ta_t(f,g)+h\,\partial_t^2a_t(f,g)+\int_t^{t+h}(t+h-s)\,\partial_t^3a_s(f,g)\,ds,\\
 \partial_t^2a_{t+h}(f,g)&=\partial_t^2a_t(f,g)+\int_t^{t+h}\partial_t^3a_s(f,g)\,ds.
 \end{aligned}
\]
By \eqref{ac-DtN-derivative-bound}, together with $\Bigl|\int_t^{t+h}(t+h-s)\,ds\Bigr|=\abs h^2/2$ and $\Bigl|\int_t^{t+h}ds\Bigr|=\abs h$, this gives
\[
 \begin{aligned}
 \abs{a_{t+h}(f,g)-a_t(f,g)-h\,\partial_ta_t(f,g)}&\leq\tfrac12\tilde C(\sigma,2)\abs h^2\norm f_X\norm g_X,\\
 \abs{\partial_ta_{t+h}(f,g)-\partial_ta_t(f,g)-h\,\partial_t^2a_t(f,g)}&\leq\tfrac12\tilde C(\sigma,3)\abs h^2\norm f_X\norm g_X,\\
 \abs{\partial_t^2a_{t+h}(f,g)-\partial_t^2a_t(f,g)}&\leq\tilde C(\sigma,3)\abs h\norm f_X\norm g_X.
 \end{aligned}
\]
Taking the supremum over $\norm f_X=\norm g_X=1$ in these three estimates gives, by the definition of the operator norm,
\[
 \begin{aligned}
 \norm{A_{t+h}-A_t-h\partial_tA_t}_{\mathcal L(X)}&\leq\tfrac12\tilde C(\sigma,2)\abs h^2,\\
 \norm{\partial_tA_{t+h}-\partial_tA_t-h\partial_t^2A_t}_{\mathcal L(X)}&\leq\tfrac12\tilde C(\sigma,3)\abs h^2,\\
 \norm{\partial_t^2A_{t+h}-\partial_t^2A_t}_{\mathcal L(X)}&\leq\tilde C(\sigma,3)\abs h.
 \end{aligned}
\]
In particular, $t\mapsto A_t$ belongs to $C^2((0,\infty),\mathcal L(X))$.

For $t\geq1$ and $\ell\geq\ell_*$ we have $t\ell\geq\ell_*$, hence
\[
 \coth(t\ell)-1\leq\frac{2e^{-2t\ell}}{1-e^{-2\ell_*}},
 \qquad
 \operatorname{csch}(t\ell)\leq\frac{2e^{-t\ell}}{1-e^{-2\ell_*}},
 \qquad
 \coth(t\ell)\leq\coth(\ell_*).
\]
Substituting \eqref{ac-DtN-derivatives} and using $e^{-2t\ell}\leq e^{-t\ell}$, we see that, for $j=0,1,2,3$ and with $\partial_t^0$ acting as the identity,
\[
\max\bigl\{\abs{\partial_t^j\bigl(\coth(t\ell)-1\bigr)},\abs{\partial_t^j\operatorname{csch}(t\ell)}\bigr\}\leq C(1,j)\ell^je^{-t\ell},
\]
where $C(1,j)$ is the constant $C(\sigma,j)$ of the preceding display with $\sigma=1$, so that it depends only on $\ell_*$ and $j$.  Repeating the argument of \eqref{ac-DtN-derivative-bound} with $\coth(t\ell_\e)-1$ in place of $\coth(t\ell_\e)$ yields
\[
 \begin{aligned}
 \norm{\partial_t^j(A_t-I_X)}_{\mathcal L(X)}
 &\leq 2C(1,j)\sup_{\ell\geq\ell_*}\ell^je^{-t\ell}
 \leq 2C(1,j)e^{-t\ell_*/2}\sup_{\ell\geq\ell_*}\ell^je^{-\ell/2}
 \leq C_je^{-t\ell_*/2},
 \end{aligned}
\]
where $C_j:=2C(1,j)\sup_{\ell\geq\ell_*}\ell^je^{-\ell/2}$ depends only on $\ell_*$ and $j$.  This proves \eqref{ac-DtN-exponential}.

For $f\in X$, differentiating \eqref{ac-DtN-identity} with $g=f$ in $t$ gives
\begin{equation}\label{ac-mass-derivative}
 \frac{d}{dt}\mathfrak b_t(U_tf,U_tf)=\frac{d}{dt}\bigl(t\langle A_tf,f\rangle_X\bigr)=\langle(A_t+t\partial_tA_t)f,f\rangle_X .
\end{equation}
Since $\partial_tU_tf\in H^1_0(\G;\V)$ by Lemma~\ref{ac-restriction-extension}, \eqref{ac-extension-orthogonality} gives $\mathfrak b_t(U_tf,\partial_tU_tf)=0$; hence, by \eqref{ac-mass-derivative},
\[
 \langle(A_t+t\partial_tA_t)f,f\rangle_X
 =2t\norm{U_tf}_2^2+2\mathfrak b_t(U_tf,\partial_tU_tf)
 =2t\norm{U_tf}_2^2,
\]
which proves \eqref{ac-mass-fixed-restriction}.
\end{proof}

\subsection{Branches near the limiting profiles}
For $p>2$, define $\mathcal N:X\to X$ by
\[
 (\mathcal N(f))(\vv):=\frac{\mathbf1_{\V_0}(\vv)}{\deg(\vv)}
   \abs{f(\vv)}^{p-2}f(\vv),\qquad\vv\in\V.
\]
Since $\deg(\vv)\geq1$ and $X\hookrightarrow\ell^\infty(\V)$ by \eqref{ac-norm-ordering},
\[
 \begin{aligned}
 \norm{\mathcal N(f)}_X^2
 &=\sum_{\vv\in\V_0}\frac{1}{\deg(\vv)}\abs{f(\vv)}^{2(p-1)}\leq\sum_{\vv\in\V_0}\abs{f(\vv)}^{2(p-1)}
 \leq\norm f_\infty^{2(p-2)}\sum_{\vv\in\V_0}\abs{f(\vv)}^2
 \leq\norm f_\infty^{2(p-2)}\norm f_X^2.
 \end{aligned}
\]
For $f\in X$ define $\mathcal N'(f)\in\mathcal L(X)$ by
\begin{equation}\label{ac-N-derivative}
 \bigl(\mathcal N'(f)g\bigr)(\vv)
 =(p-1)\frac{\mathbf1_{\V_0}(\vv)}{\deg(\vv)}\abs{f(\vv)}^{p-2}g(\vv),
 \qquad \forall g\in X,
\end{equation}
so that
\[
 \norm{\mathcal N'(f)g}_X\leq(p-1)\norm f_\infty^{p-2}\norm g_X,
 \qquad \forall g\in X.
\]
We show first that $\mathcal N$ has a Fr\'echet derivative at every $f\in X$, namely $\mathcal N'(f)$. For $R>0$ and $s_1,s_2\in[-R,R]$, the inequality $\bigl|\abs{a}^{\alpha}-\abs{b}^{\alpha}\bigr|\leq\abs{a-b}^{\alpha}$ with $\alpha=p-2$ if $2<p\leq3$, and the mean value theorem if $p>3$, give
\begin{equation}\label{ac-N-mvt}
 \bigl|\abs{s_1}^{p-2}-\abs{s_2}^{p-2}\bigr|
 \leq C_R|s_1-s_2|^{\min\{1,p-2\}},
\end{equation}
where $C_R:=(p-1)\max\{1,(p-2)R^{p-3}\}$ depends only on $p$ and $R$; the factor $p-1$ makes \eqref{ac-N-remainder} and \eqref{ac-N-lipschitz} below hold with the same constant. For $f,h\in X$ with $\norm f_X+\norm h_X\leq R$, the coordinatewise fundamental theorem of calculus gives
\[
 \mathcal N(f+h)-\mathcal N(f)-\mathcal N'(f)h=\int_0^1\bigl[\mathcal N'(f+sh)-\mathcal N'(f)\bigr]h\,ds,
\]
and hence, by \eqref{ac-N-mvt},
\begin{equation}\label{ac-N-remainder}
 \norm{\mathcal N(f+h)-\mathcal N(f)-\mathcal N'(f)h}_X
 \leq C_R\norm h_X^{1+\min\{1,p-2\}}=o(\norm h_X),
\end{equation}
since $p>2$. This proves that $\mathcal N$ has a Fr\'echet derivative at every $f\in X$.
It remains to prove that $f\mapsto\mathcal N'(f)$ is continuous. If $f,\widetilde f\in X$ with $\norm f_X,\norm{\widetilde f}_X\leq R$, then \eqref{ac-N-mvt} gives
\begin{equation}\label{ac-N-lipschitz}
 \norm{\mathcal N'(f)-\mathcal N'(\widetilde f)}_{\mathcal L(X)}
 \leq C_R\norm{f-\widetilde f}_X^{\min\{1,p-2\}}.
\end{equation}
Hence $\mathcal N\in C^1(X,X)$.

For $f\in X$ and $t>0$, set $u:=t^{1/(p-2)}U_tf\in H^1(\G)$.  For any $h\in H^1(\G)$ written as in \eqref{ac-restriction-decomposition}, with vertex restriction $g=\mathcal R h$ and interior part $z=h-U_tg\in H^1_0(\G;\V)$, the orthogonality \eqref{ac-extension-orthogonality} and the identity \eqref{ac-DtN-identity} yield
\begin{equation}\label{ac-weak-vertex-form}
 \begin{aligned}
 E'(u,\G)[h]+t^2(u,h)_2
 &=t^{1/(p-2)}\bigl[\mathfrak b_t(U_tf,U_tg)+\mathfrak b_t(U_tf,z)\bigr]
 -t^{(p-1)/(p-2)}\sum_{\vv\in\V_0}\abs{f(\vv)}^{p-2}f(\vv)g(\vv)\\
 &=t^{1/(p-2)}\mathfrak b_t(U_tf,U_tg)
 -t^{(p-1)/(p-2)}\sum_{\vv\in\V_0}\abs{f(\vv)}^{p-2}f(\vv)g(\vv)\\
 &=t^{(p-1)/(p-2)}\bigl\langle A_tf-\mathcal N(f),g\bigr\rangle_X,
 \end{aligned}
\end{equation}
where the first equality uses $h=z+U_tg$, the second $\mathfrak b_t(U_tf,z)=0$ from \eqref{ac-extension-orthogonality}, and the third $\mathfrak b_t(U_tf,U_tg)=t\langle A_tf,g\rangle_X$ from \eqref{ac-DtN-identity} together with the definition of $\mathcal N$.  Since $\mathcal R$ is onto, \eqref{ac-weak-vertex-form} shows that, for $u\in H^1(\G)$, the weak equation for \eqref{eqpd} is equivalent to
\begin{equation}\label{ac-vertex-equation}
 A_tf=\mathcal N(f)\qquad\text{in }X .
\end{equation}
More explicitly, every solution $f$ of \eqref{ac-vertex-equation} yields a solution $u=t^{1/(p-2)}U_tf$ of \eqref{eqpd} with Lagrange multiplier $t^2$.  Conversely, if $u\in H^1(\G)$ is a solution of \eqref{eqpd} with multiplier $t^2$ and $f:=t^{-1/(p-2)}\mathcal R u$, then $u=t^{1/(p-2)}U_tf$ by the uniqueness in Lemma~\ref{ac-restriction-extension}.  Imposing the prescribed mass $\mu$ is therefore equivalent to the scalar condition
\begin{equation}\label{ac-mass-scalar}
 t^{2/(p-2)}\norm{U_tf}_2^2=\mu ,
\end{equation}
so that, for each $t>0$, solutions of \eqref{eqpd} with multiplier $t^2$ correspond exactly to solutions of \eqref{ac-vertex-equation} satisfying \eqref{ac-mass-scalar}.

Since $A_t\to I_X$, we seek solutions of \eqref{ac-vertex-equation} near solutions of $f=\mathcal N(f)$ supported on a finite subset of $\V_0$. The next lemma constructs a branch for each nonempty finite subset of $\V_0$ and proves that the corresponding solutions are positive. For $S\subset\V_0$ finite, set
\begin{equation}\label{ac-limit-profile}
 f_S^\infty(\vv):=\deg(\vv)^{1/(p-2)}\mathbf 1_S(\vv),\qquad\vv\in\V .
\end{equation}

\begin{lemma}\label{ac-IFT-branch}
For every finite nonempty $S\subset\V_0$ there are constants $T_S>1$, $0<\delta_S\leq\frac12$ and $C>0$, depending only on $p$, $\ell_*$, $\#S$ and $\max_{\vv\in S}\deg(\vv)$, and a map $f_S\in C^1((T_S,\infty),X)$ such that, for every $t>T_S$, $f_S(t)$ solves \eqref{ac-vertex-equation}, $U_tf_S(t)>0$ on $\G$, and
\begin{equation}\label{ac-branch-estimates}
 \norm{f_S(t)-f_S^\infty}_X+\norm{f_S'(t)}_X\leq Ce^{-t\ell_*/2}.
\end{equation}
For every $t>T_S$, the only solution $f\in X$ of \eqref{ac-vertex-equation} with $\norm{f-f_S^\infty}_X<\delta_S$ is $f=f_S(t)$.
\end{lemma}

\begin{proof}
By \eqref{ac-DtN-exponential}, $A_t\to I_X$ in operator norm as $t\to\infty$, so the limiting form of \eqref{ac-vertex-equation} is $f=\mathcal N(f)$.  Since $S$ is finite, $f_S^\infty$ defined as in \eqref{ac-limit-profile} belongs to $X$ and satisfies
\[
 \norm{f_S^\infty}_X^2
   =\sum_{\vv\in S}\deg(\vv)^{p/(p-2)}<\infty,
 \qquad \mathcal N(f_S^\infty)=f_S^\infty .
\]

Substituting $f_S^\infty$ into \eqref{ac-N-derivative} gives
$(\mathcal N'(f_S^\infty)g)(\vv)=(p-1)\mathbf 1_S(\vv)g(\vv)$.  We define the limiting linearization $L_S^\infty:X\to X$ by
\begin{equation}\label{ac-limit-linearization}
 L_S^\infty:=I_X-\mathcal N'(f_S^\infty),\qquad
 (L_S^\infty g)(\vv)=
 \begin{cases}
  (2-p)g(\vv),&\vv\in S,\\
  g(\vv),&\vv\notin S,
 \end{cases}
\end{equation}
whose inverse is given on all of $X$ by
\begin{equation}\label{ac-limit-inverse-coordinate}
 \bigl((L_S^\infty)^{-1}g\bigr)(\vv)=
 \begin{cases}
  (2-p)^{-1}g(\vv),&\vv\in S,\\
  g(\vv),&\vv\notin S.
 \end{cases}
\end{equation}
Thus
\begin{equation}\label{ac-limit-inverse-bound}
 \norm{(L_S^\infty)^{-1}}_{\mathcal L(X)}
 \leq\max\{1,(p-2)^{-1}\}<+\infty.
\end{equation}
We now introduce the subspace of vertex values carried by $S$ and its orthogonal projection,
\begin{equation}\label{ac-coordinate-subspace}
 X_S:=\{g\in X:g(\vv)=0\text{ for }\vv\notin S\},\qquad
 P_S:X\to X_S,\qquad P_Sg:=\mathbf 1_Sg .
\end{equation}
In this notation
\[
 L_S^\infty=(2-p)P_S+(I_X-P_S).
\]
Thus $L_S^\infty$ is the negative multiple $(2-p)$ of the identity on the $\#S$-dimensional space $X_S$, and the identity on $X_S^\perp$.  In particular,
\begin{equation}\label{ac-limit-form}
 \langle L_S^\infty g,g\rangle_X
 =(2-p)\norm{P_Sg}_X^2+\norm{(I_X-P_S)g}_X^2
 \qquad\forall g\in X.
\end{equation}
Set $\widetilde A_\tau:=A_{1/\tau}$ for $\tau>0$ and $\widetilde A_\tau:=I_X$ for $\tau\leq0$. For $\tau>0$,
\[
 \widetilde A_\tau'=-\tau^{-2}\partial_tA_{1/\tau}.
\]
By \eqref{ac-DtN-exponential}, as $\tau\to0^+$,
\[
 \frac{\norm{\widetilde A_\tau-I_X}_{\mathcal L(X)}}{\tau}
 +\norm{\widetilde A_\tau'}_{\mathcal L(X)}
 \leq C_0\frac{e^{-\ell_*/(2\tau)}}{\tau}
 +C_1\frac{e^{-\ell_*/(2\tau)}}{\tau^2}\longrightarrow0.
\]
Hence $\tau \mapsto \widetilde A_\tau$ is $C^1$ near zero. Define the $C^1$ map
\[
 \Phi:(-1,1)\times X\longrightarrow X,\qquad
 \Phi(\tau,f):=\widetilde A_\tau f-\mathcal N(f).
\]

We now show that the solutions of $\Phi(\tau,f)=0$ near $f_S^\infty$ form a locally unique $C^1$ curve.
Set $R:=(\#S)^{1/2}\bigl(\max_{\vv\in S}\deg(\vv)\bigr)^{p/(2(p-2))}+1$ and let $C_R>0$ be the constant in \eqref{ac-N-mvt}.
Since
\begin{equation}\label{ac-profile-norm-bound}
 \norm{f_S^\infty}_X^2
 =\sum_{\vv\in S}\deg(\vv)^{p/(p-2)}
 \leq\#S\left(\max_{\vv\in S}\deg(\vv)\right)^{p/(p-2)},
\end{equation}
we have $\norm{f_S^\infty}_X+1\leq R$, and $C_R$ depends only on $p$, $\#S$ and $\max_{\vv\in S}\deg(\vv)$.
Choose
\[
 \delta_S:=\min\left\{\frac12,
 \left(\frac{\min\{1,p-2\}}{4C_R}\right)^{1/\min\{1,p-2\}}\right\}>0,
\]
and choose $\varepsilon_S\in(0,1)$ so that, for every $0<\tau<\varepsilon_S$,
\[
 C_0e^{-\ell_*/(2\tau)}
 <\frac{\min\{1,p-2\}}4
   \min\left\{1,\frac{\delta_S}{\norm{f_S^\infty}_X}\right\},
\]
where $C_0>0$ is the constant in \eqref{ac-DtN-exponential}, which depends only on $\ell_*$. Moreover, $\delta_S$ and $\varepsilon_S$ may be chosen to depend only on $p$, $\ell_*$, $\#S$ and $\max_{\vv\in S}\deg(\vv)$.

 By \eqref{ac-DtN-exponential} and the definition of $\widetilde A_\tau$, for every $\abs\tau<\varepsilon_S$,
\begin{equation}\label{ac-tau-smallness}
 \norm{\widetilde A_\tau-I_X}_{\mathcal L(X)}<\frac14\min\{1,p-2\},
 \qquad
 \norm{\widetilde A_\tau-I_X}_{\mathcal L(X)}\norm{f_S^\infty}_X
 <\frac{\delta_S}{4}\min\{1,p-2\}.
\end{equation}
For each $\abs\tau<\varepsilon_S$, define the map $\mathcal T_{\tau,S}:\{h\in X:\norm h_X\leq\delta_S\}\to X$ by
\[
 \mathcal T_{\tau,S}(h):=-(L_S^\infty)^{-1}\Bigl[(\widetilde A_\tau-I_X)(f_S^\infty+h)-\bigl(\mathcal N(f_S^\infty+h)-\mathcal N(f_S^\infty)
           -\mathcal N'(f_S^\infty)h\bigr)\Bigr].
\]
The fixed-point equation $h=\mathcal T_{\tau,S}(h)$ is equivalent to $\Phi(\tau,f_S^\infty+h)=0$. Indeed, writing $\Theta(h):=\mathcal N(f_S^\infty+h)-\mathcal N(f_S^\infty)-\mathcal N'(f_S^\infty)h$ and using $f_S^\infty=\mathcal N(f_S^\infty)$, $L_S^\infty=I_X-\mathcal N'(f_S^\infty)$ and $\widetilde A_\tau=I_X+(\widetilde A_\tau-I_X)$ gives
\[
\begin{aligned}
   \Phi(\tau,f_S^\infty+h)
 &=\widetilde A_\tau f_S^\infty-\mathcal N(f_S^\infty)
  +\widetilde A_\tau h-\bigl(\mathcal N(f_S^\infty+h)-\mathcal N(f_S^\infty)\bigr)\\
 &=L_S^\infty h+(\widetilde A_\tau-I_X)(f_S^\infty+h)-\Theta(h),
\end{aligned}
\]
and thus, since $L_S^\infty$ is invertible by \eqref{ac-limit-inverse-bound}, $\Phi(\tau,f_S^\infty+h)=0$ if and only if $h=\mathcal T_{\tau,S}(h)$. For $\norm{h_1}_X,\norm{h_2}_X\leq\delta_S$, \eqref{ac-N-lipschitz} and the fundamental theorem of calculus give
\[
 \begin{aligned}
   \Theta(h_1)-\Theta(h_2)
 &=\int_0^1\bigl[\mathcal N'\bigl(f_S^\infty+h_2+s(h_1-h_2)\bigr)-\mathcal N'(f_S^\infty)\bigr]ds\;\,(h_1-h_2),\\
   \norm{\Theta(h_1)-\Theta(h_2)}_X
 &\leq\int_0^1\norm{\mathcal N'\bigl(f_S^\infty+h_2+s(h_1-h_2)\bigr)-\mathcal N'(f_S^\infty)}_{\mathcal L(X)}ds\;\norm{h_1-h_2}_X\\
 &\leq C_R\Bigl(\sup_{0\leq s\leq1}\norm{h_2+s(h_1-h_2)}_X\Bigr)^{\min\{1,p-2\}}\norm{h_1-h_2}_X\\
 &\leq C_R\delta_S^{\min\{1,p-2\}}\norm{h_1-h_2}_X,
 \end{aligned}
\]
and hence, by \eqref{ac-limit-inverse-bound} and \eqref{ac-tau-smallness},
\[
 \begin{aligned}
 \norm{\mathcal T_{\tau,S}(h_1)-\mathcal T_{\tau,S}(h_2)}_X
 &=\norm{(L_S^\infty)^{-1}\bigl[(\widetilde A_\tau-I_X)(h_1-h_2)-\bigl(\Theta(h_1)-\Theta(h_2)\bigr)\bigr]}_X\\
 &\leq\max\{1,(p-2)^{-1}\}\Bigl(\norm{\widetilde A_\tau-I_X}_{\mathcal L(X)}\norm{h_1-h_2}_X+\norm{\Theta(h_1)-\Theta(h_2)}_X\Bigr)\\
 &\leq\frac{\norm{\widetilde A_\tau-I_X}_{\mathcal L(X)}
       +C_R\delta_S^{\min\{1,p-2\}}}{\min\{1,p-2\}}
       \norm{h_1-h_2}_X\\
 &\leq\frac12\norm{h_1-h_2}_X,
  \end{aligned}
 \]
 and
\[
 \norm{\mathcal T_{\tau,S}(0)}_X
=\norm{(L_S^\infty)^{-1}(\widetilde A_\tau-I_X)f_S^\infty}_X
 \leq\frac{\norm{\widetilde A_\tau-I_X}_{\mathcal L(X)}\norm{f_S^\infty}_X}{\min\{1,p-2\}}
 \leq\frac{\delta_S}{4}.
\]
Thus $\mathcal T_{\tau,S}$ maps $\{h\in X:\norm h_X\leq\delta_S\}$ into $\{h\in X:\norm h_X\leq\frac34\delta_S\}$ and is a $\frac12$-contraction there, so it has exactly one fixed point $h_{\tau,S}:=\lim_{n\to\infty}\mathcal T_{\tau,S}^n(0)$ in $\{h\in X:\norm h_X\leq\delta_S\}$. Moreover, since \(
\|h_{\tau,S}\|_X
=\|\mathcal T_{\tau,S}(h_{\tau,S})\|_X
\leq \tfrac12\|h_{\tau,S}\|_X+\|\mathcal T_{\tau,S}(0)\|_X
\), we have
\begin{equation}\label{ac-uniform-fixedpoint-error}
 \norm{h_{\tau,S}}_X\leq2\norm{\mathcal T_{\tau,S}(0)}_X
 \leq\min\left\{\frac{\delta_S}2,
 \frac{2\norm{f_S^\infty}_X}{\min\{1,p-2\}}
 \norm{\widetilde A_\tau-I_X}_{\mathcal L(X)}\right\}.
\end{equation}
Let
\[
 K:=(L_S^\infty)^{-1}D_f\Phi(\tau,f_S^\infty+h_{\tau,S})-I_X.
\]
Then we have $D_f\Phi(\tau,f_S^\infty+h_{\tau,S})=L_S^\infty(I_X+K)$.
For every $g\in X$ with $\norm g_X\leq\delta_S$, since $\widetilde A_\tau-\mathcal N'(f_S^\infty+g)-L_S^\infty=(\widetilde A_\tau-I_X)-\bigl(\mathcal N'(f_S^\infty+g)-\mathcal N'(f_S^\infty)\bigr)$, \eqref{ac-limit-inverse-bound}, \eqref{ac-N-lipschitz}, \eqref{ac-tau-smallness} and the choice of $\delta_S$ give
\[
 \begin{aligned}
   &\norm{(L_S^\infty)^{-1}
          [\widetilde A_\tau-\mathcal N'(f_S^\infty+g)-L_S^\infty]}_{\mathcal L(X)}\\
   &\qquad\leq\norm{(L_S^\infty)^{-1}}_{\mathcal L(X)}
        \Bigl[\norm{\widetilde A_\tau-I_X}_{\mathcal L(X)}
        +\norm{\mathcal N'(f_S^\infty+g)-\mathcal N'(f_S^\infty)}_{\mathcal L(X)}\Bigr]\\
   &\qquad\leq\frac{\norm{\widetilde A_\tau-I_X}_{\mathcal L(X)}
        +C_R\norm g_X^{\min\{1,p-2\}}}{\min\{1,p-2\}}\leq\frac{\frac14\min\{1,p-2\}+\frac14\min\{1,p-2\}}{\min\{1,p-2\}}=\frac12.
 \end{aligned}
\]
Thus, since $\norm{h_{\tau,S}}_X\leq\frac12\delta_S$ by \eqref{ac-uniform-fixedpoint-error}, we obtain $\norm K_{\mathcal L(X)}\leq\frac12$. Since $I_X+K=I_X-(-K)$ is invertible by \cite[Exercise~6.14]{Ha} and $L_S^\infty$ is invertible by \eqref{ac-limit-inverse-bound}, we know that $D_f\Phi(\tau,f_S^\infty+h_{\tau,S})$ is invertible. 

For each $\tau$ with $\abs\tau<\varepsilon_S$, applying the implicit function theorem at $h_{\tau,S}$ gives a $C^1$ map on a neighborhood of $\tau$. Since $\norm{h_{\tau,S}}_X\leq\delta_S/2$, each local solution of $\Phi(\tau,f_S^\infty+h)=0$ in the ball $\{h:\norm h_X\leq\delta_S\}$ is a fixed point of $\mathcal T_{\tau,S}$ in that ball, where the fixed point is unique. Hence the local maps obtained for different $\tau$ with $\abs\tau<\varepsilon_S$ agree and combine into a single $C^1$ map $\psi_S:(-\varepsilon_S,\varepsilon_S)\to X$ such that
\[
 \Phi(\tau,\psi_S(\tau))=0, \qquad\abs\tau<\varepsilon_S,
\]
and we conclude the uniqueness of the solution of $\Phi(\tau,f)=0$ in the ball $\norm{f-f_S^\infty}_X<\delta_S$:
\begin{equation}\label{ac-IFT-local-uniqueness}
 f=\psi_S(\tau)\quad\text{if}\quad \Phi(\tau,f)=0,\ \ \abs\tau<\varepsilon_S,\ \ \norm{f-f_S^\infty}_X<\delta_S .
\end{equation}
Since $\psi_S(\tau)=f_S^\infty+h_{\tau,S}$ and $\norm{h_{\tau,S}}_X\leq\delta_S/2$ by \eqref{ac-uniform-fixedpoint-error},
\begin{equation}\label{ac-psi-norm-bound}
 \norm{\psi_S(\tau)}_X\leq\norm{f_S^\infty}_X+\frac{\delta_S}{2},
 \qquad\abs\tau<\varepsilon_S.
\end{equation}
Setting $\widetilde T_S:=2\varepsilon_S^{-1}$ and $f_S(t):=\psi_S(1/t)$ for $t>\widetilde T_S$ gives $A_tf_S(t)=\mathcal N(f_S(t))$ and $f_S(t)\to f_S^\infty$ in $X$ as $t\to+\infty$, since $\psi_S$ is continuous at $\tau=0$ and $\psi_S(0)=f_S^\infty$ by \eqref{ac-IFT-local-uniqueness}.
By \eqref{ac-psi-norm-bound}, \eqref{ac-profile-norm-bound} and $\delta_S\leq1/2$,
\[
 \sup_{s>\widetilde T_S}\norm{f_S(s)}_X
 \leq\norm{f_S^\infty}_X+\frac{\delta_S}{2}
 \leq(\#S)^{1/2}\left(\max_{\vv\in S}\deg(\vv)\right)^{p/(2(p-2))}+\frac14.
\]
Set
\[
 \begin{aligned}
 C:= \frac{4\max\{C_0,C_1\}\max\{1,(p-2)^{-1}\}}{\min\{1,p-2\}}\times\left[(\#S)^{1/2}\left(\max_{\vv\in S}\deg(\vv)\right)^{p/(2(p-2))}+\frac14\right]>0.
 \end{aligned}
\]
Here $C_0$ and $C_1$ are the constants in \eqref{ac-DtN-exponential}, depending only on $\ell_*$. Thus $C$ depends only on $p$, $\ell_*$, $\#S$ and $\max_{\vv\in S}\deg(\vv)$.
The remainder estimate \eqref{ac-N-remainder} for $\mathcal N$, with $R$ chosen above, gives, for every $t>\widetilde T_S$,
\[
 \norm{\mathcal N(f_S(t))-\mathcal N(f_S^\infty)
   -\mathcal N'(f_S^\infty)(f_S(t)-f_S^\infty)}_X
 \leq C_R\norm{f_S(t)-f_S^\infty}_X^{1+\min\{1,p-2\}}.
\]
Writing $A_t=I_X+(A_t-I_X)$ in $A_tf_S(t)=\mathcal N(f_S(t))$ gives
\[
 f_S(t)=\mathcal N(f_S(t))-(A_t-I_X)f_S(t).
\]
  Since $L_S^\infty=I_X-\mathcal N'(f_S^\infty)$, it follows that
\[
 \begin{aligned}
 L_S^\infty(f_S(t)-f_S^\infty)
 ={}&f_S(t)-f_S^\infty-\mathcal N'(f_S^\infty)(f_S(t)-f_S^\infty)\\
 ={}&\mathcal N(f_S(t))-\mathcal N(f_S^\infty)
   -\mathcal N'(f_S^\infty)(f_S(t)-f_S^\infty)
   -(A_t-I_X)f_S(t).
 \end{aligned}
\]
Using \eqref{ac-DtN-exponential}, the boundedness of $f_S(t)$ and \eqref{ac-N-remainder}, for every $t>\widetilde T_S$ we obtain
\[
 \begin{aligned}
 \norm{f_S(t)-f_S^\infty}_X
 &\leq\norm{(L_S^\infty)^{-1}}_{\mathcal L(X)}
   \Bigl(\norm{(A_t-I_X)f_S(t)}_X
   +\bigl\|\mathcal N(f_S(t))-\mathcal N(f_S^\infty)-\mathcal N'(f_S^\infty)(f_S(t)-f_S^\infty)\bigr\|_X\Bigr)\\
 &\leq\max\{1,(p-2)^{-1}\}
   \Bigl(C_0e^{-t\ell_*/2}\sup_{s>\widetilde T_S}\norm{f_S(s)}_X
   +C_R\norm{f_S(t)-f_S^\infty}_X^{1+\min\{1,p-2\}}\Bigr)\\
 &\leq\frac{C}{4}e^{-t\ell_*/2}+\tfrac14\norm{f_S(t)-f_S^\infty}_X,
 \end{aligned}
\]
where the last step uses the bound $\norm{f_S(t)-f_S^\infty}_X=\norm{h_{1/t,S}}_X\leq\delta_S/2$ from \eqref{ac-uniform-fixedpoint-error} and $C_R\delta_S^{\min\{1,p-2\}}\leq\tfrac14\min\{1,p-2\}$ from the choice of $\delta_S$; hence, for every $t>\widetilde T_S$,
\begin{equation}\label{ac-branch-limit-estimate}
 \norm{f_S(t)-f_S^\infty}_X\leq\frac{C}{2}e^{-t\ell_*/2}.
\end{equation}
Differentiating \eqref{ac-vertex-equation} with respect to $t$ gives
\begin{equation}\label{ac-branch-derivative}
 \bigl[A_t-\mathcal N'(f_S(t))\bigr]f_S'(t)=-\partial_tA_tf_S(t).
\end{equation}
Since $L_S^\infty$ is invertible, with $\norm{(L_S^\infty)^{-1}}_{\mathcal L(X)}\le\max\{1,(p-2)^{-1}\}$ by \eqref{ac-limit-inverse-bound}, we have
\[
 \norm{L_S^\infty g}_X\geq\min\{1,p-2\}\norm g_X\qquad\forall g\in X.
\]
By \eqref{ac-DtN-exponential}, and by \eqref{ac-N-lipschitz} combined with \eqref{ac-branch-limit-estimate}, we have, for every $t$ such that
\[t>\max\Bigl\{\widetilde T_S,\ \frac{2}{\min\{1,p-2\}\ell_*}\log\frac{2\bigl(C_0+C_R(C/2)^{\min\{1,p-2\}}\bigr)}{\min\{1,p-2\}}\Bigr\},
\]
\[
 \begin{aligned}
 \norm{A_t-\mathcal N'(f_S(t))-L_S^\infty}_{\mathcal L(X)}
 &=\norm{(A_t-I_X)-\bigl(\mathcal N'(f_S(t))-\mathcal N'(f_S^\infty)\bigr)}_{\mathcal L(X)}\\
 &\leq\norm{A_t-I_X}_{\mathcal L(X)}+\norm{\mathcal N'(f_S(t))-\mathcal N'(f_S^\infty)}_{\mathcal L(X)}\\
 &\leq C_0e^{-t\ell_*/2}+C_R\norm{f_S(t)-f_S^\infty}_X^{\min\{1,p-2\}}\\
 &\leq C_0e^{-t\ell_*/2}+C_R\bigl(C/2\bigr)^{\min\{1,p-2\}}e^{-\min\{1,p-2\}t\ell_*/2}\\
 &\leq\bigl(C_0+C_R(C/2)^{\min\{1,p-2\}}\bigr)e^{-\min\{1,p-2\}t\ell_*/2}\leq\tfrac12\min\{1,p-2\}.
 \end{aligned}
\]
The triangle inequality therefore yields
\[
 \norm{[A_t-\mathcal N'(f_S(t))]g}_X
 \geq\tfrac12\min\{1,p-2\}\norm g_X
 \qquad\forall g\in X.
\]
Combining this with \eqref{ac-branch-derivative} and \eqref{ac-DtN-exponential}, we obtain
\[
 \tfrac12\min\{1,p-2\}\norm{f_S'(t)}_X \leq\norm{\partial_tA_tf_S(t)}_X\leq\norm{\partial_tA_t}_{\mathcal L(X)}\norm{f_S(t)}_X
 \leq C_1e^{-t\ell_*/2}\sup_{s>\widetilde T_S}\norm{f_S(s)}_X,
  \]
which proves $\norm{f_S'(t)}_X\leq\frac{C}{2}e^{-t\ell_*/2}$.  Summing this with \eqref{ac-branch-limit-estimate} gives \eqref{ac-branch-estimates}.

Finally, we prove that $U_tf_S(t)>0$ on $\G$ for every $t>T_S$, where the threshold $T_S\geq\widetilde T_S$ will be chosen below.  Write $f^+:=\max\{f,0\}$ and $f^-:=\min\{f,0\}$, so that $f=f^++f^-$ and $f^+(\vv)f^-(\vv)=0$ for every $\vv\in\V$.  Since $f_S^\infty\geq0$, by \eqref{ac-norm-ordering} and \eqref{ac-branch-limit-estimate},
\begin{equation}\label{ac-negative-part-decay}
 \norm{f_S(t)^-}_\infty
 \leq\norm{f_S(t)-f_S^\infty}_\infty
 \leq\norm{f_S(t)-f_S^\infty}_X\leq\frac{C}{2}e^{-t\ell_*/2}.
\end{equation}
Since $f_S(t)^+(\vv)f_S(t)^-(\vv)=0$ for every $\vv\in\V$,
\[
 \sum_{\e\in\E_b}\coth(t\ell_\e)
   \bigl[f_S(t)^+(\vv_\e)f_S(t)^-(\vv_\e)+f_S(t)^+(\w_\e)f_S(t)^-(\w_\e)\bigr]=0,
 \qquad
 \sum_{\e\in\E_\infty}f_S(t)^+(\vv_\e)f_S(t)^-(\vv_\e)=0.
\]
 Recall that \eqref{ac-DtN-form} reads, in the endpoint notation,
\begin{equation*}
 \begin{aligned}
 \bigl\langle A_tf_S(t)^+,f_S(t)^-\bigr\rangle_X
 =&a_t\bigl(f_S(t)^+,f_S(t)^-\bigr)\\
 =&\sum_{\e\in\E_b}\coth(t\ell_\e)\bigl[f_S(t)^+(\vv_\e)f_S(t)^-(\vv_\e)+f_S(t)^+(\w_\e)f_S(t)^-(\w_\e)\bigr]\\
 &-\sum_{\e\in\E_b}\operatorname{csch}(t\ell_\e)\bigl[f_S(t)^+(\vv_\e)f_S(t)^-(\w_\e)+f_S(t)^+(\w_\e)f_S(t)^-(\vv_\e)\bigr]\\
 &+\sum_{\e\in\E_\infty}f_S(t)^+(\vv_\e)f_S(t)^-(\vv_\e)\\
 =&-\sum_{\e\in\E_b}\operatorname{csch}(t\ell_\e)\bigl[f_S(t)^+(\vv_\e)f_S(t)^-(\w_\e)+f_S(t)^+(\w_\e)f_S(t)^-(\vv_\e)\bigr]\geq0.
 \end{aligned}
\end{equation*}
Since $A_tf_S(t)=\mathcal N(f_S(t))$, we obtain
\[
 \sum_{\vv\in\V_0}\abs{f_S(t)^-(\vv)}^p
 =\bigl\langle A_tf_S(t)^+,f_S(t)^-\bigr\rangle_X+\bigl\langle A_tf_S(t)^-,f_S(t)^-\bigr\rangle_X .
\]
Moreover, by \eqref{ac-DtN-coercivity}, it follows that
\[
 \begin{aligned}
 \tanh(t\ell_*/2)\norm{f_S(t)^-}_X^2
 \leq\bigl\langle A_tf_S(t)^-,f_S(t)^-\bigr\rangle_X
 &=\sum_{\vv\in\V_0}\abs{f_S(t)^-(\vv)}^p-\bigl\langle A_tf_S(t)^+,f_S(t)^-\bigr\rangle_X\\
 &\leq\sum_{\vv\in\V_0}\abs{f_S(t)^-(\vv)}^p
 \leq\norm{f_S(t)^-}_\infty^{p-2}\norm{f_S(t)^-}_X^2 .
 \end{aligned}
\]
By \eqref{ac-branch-estimates} and \eqref{ac-negative-part-decay}, we can choose $T_S\geq\widetilde T_S$ so large that, for every $t>T_S$, $\norm{f_S(t)-f_S^\infty}_X\leq Ce^{-t\ell_*/2}$, $\norm{f_S'(t)}_X\leq Ce^{-t\ell_*/2}$, $\norm{f_S(t)^-}_\infty^{p-2}<\tanh(t\ell_*/2)$ and
\[
 Ce^{-t\ell_*/2}\leq\frac12
 \leq\frac12\left(\sum_{\vv\in S}\deg(\vv)^{p/(p-2)}\right)^{1/2}
 =\frac12\norm{f_S^\infty}_X.
\]
Then the inequality $\norm{f_S(t)^-}_\infty^{p-2}<\tanh(t\ell_*/2)$ forces $f_S(t)^-=0$ and hence $f_S(t)\geq0$ for every $t>T_S$. Moreover,  $\norm{f_S(t)}_X\geq\norm{f_S^\infty}_X-Ce^{-t\ell_*/2}\geq\tfrac12\norm{f_S^\infty}_X>0$, which implies that $f_S(t) \not \equiv 0$. Since $t>0$ and $\ell_\e>0$, the coefficients $\sinh(t(\ell_\e-x))/\sinh(t\ell_\e)$ and $\sinh(tx)/\sinh(t\ell_\e)$ in \eqref{ac-edge-extension} are nonnegative for $x\in[0,\ell_\e]$, and the coefficient $e^{-tx}$ in \eqref{ac-halfline-extension} is nonnegative for $x\geq0$. If $f_S(t)(\vv)=0$ at some vertex $\vv$, then $\mathcal N(f_S(t))(\vv)=0$, so \eqref{ac-vertex-equation} and \eqref{ac-DtN-coordinate} give $0=\deg(\vv)(A_tf_S(t))(\vv)=-\sum_{\e\in\E_b,\,\e\succ\vv}\operatorname{csch}(t\ell_\e)f_S(t)(\w_{\e,\vv})$. Since $f_S(t)\geq0$, $f_S(t)$ vanishes at every neighbor of $\vv$, and, since $\G$ is connected, $f_S(t)\equiv0$, which is impossible. Hence $f_S(t)>0$ at every vertex. Since the two coefficients in \eqref{ac-edge-extension} do not vanish simultaneously and $e^{-tx}>0$, we conclude that $U_tf_S(t)>0$ on $\G$ for every $t>T_S$.
Since $0<1/t<\varepsilon_S$ for every $t>T_S\geq\widetilde T_S=2\varepsilon_S^{-1}$, \eqref{ac-IFT-local-uniqueness} with $\tau=1/t$ proves the uniqueness assertion of the lemma.
\end{proof}

For the branches in Lemma~\ref{ac-IFT-branch}, we shall use the following functions and their masses:
\begin{equation}\label{ac-physical-branch}
 u_S(t):=t^{1/(p-2)}U_tf_S(t),\qquad
 \mu_S(t):=\norm{u_S(t)}_2^2 .
\end{equation}
By Lemma~\ref{ac-IFT-branch}, the rescaled vertex restriction $t^{-1/(p-2)}\mathcal Ru_S(t)=f_S(t)$ converges in $X$ to $f_S^\infty$, which is supported on $S$.

\section{Proofs of Theorems~\ref{th2}, \ref{cor-ac-infinite-orbits} and \ref{cor-ac-supercritical-all-mass}}
\label{sect-branches}

The main result of this section is Theorem~\ref{th-ac-uniform-solutions}, a uniform version of Theorem~\ref{th2}. We prove it by inverting the mass curves with a common threshold: Proposition~\ref{prop-ac-branches} computes these curves asymptotically as $t\to+\infty$, and inverting them turns the parameter $t$ back into the prescribed mass $\mu$. Theorem~\ref{th2} and Theorem~\ref{cor-ac-infinite-orbits} then follow. At the end of the section, we prove Theorem~\ref{cor-ac-supercritical-all-mass} by choosing a frequency threshold independent of the number of selected vertices.
Throughout this section, $\G=(\V,\E)\in\mathbf{G}$, $p>2$ and $\varnothing\ne\V_0\subset\V$, and we keep the notation of Section~\ref{sect5}.

By Lemmas~\ref{ac-restriction-extension} and~\ref{ac-IFT-branch}, $u_S\in C^1((T_S,\infty),H^1(\G))$, with
\begin{equation}\label{ac-physical-branch-derivative}
 u_S'(t)
 =\frac{t^{1/(p-2)-1}}{p-2}U_tf_S(t)
   +t^{1/(p-2)}\bigl(\partial_tU_tf_S(t)+U_tf_S'(t)\bigr)\in H^1(\G),
\end{equation}
and
\[
 A_t\bigl(f_S(t)\bigr)
 =\mathcal N\bigl(f_S(t)\bigr) \qquad \text{in }X.
\]
Then, by \eqref{ac-weak-vertex-form}, $u_S(t)$ satisfies the weak equation of \eqref{eqpd} with multiplier $t^2$:
\begin{equation}\label{ac-branch-weak-identity}
 E'(u_S(t),\G)[h]+t^2(u_S(t),h)_2
 =0 \qquad\text{for every }h\in H^1(\G).
\end{equation}
Thus each fixed $S$ gives a branch $\{(u_S(t),t^2)\}_{t>T_S}$ of positive solutions of \eqref{eqpd} with mass $\mu_S(t)$.

For any weak solution $(u,\lambda)\in H^1_\mu(\G)\times\R$ to \eqref{eqpd}, define $Q_u:H^1(\G)\times H^1(\G)\to\R$ as the Hessian of $E(\cdot,\G)+\frac{\lambda}{2}\norm{\cdot}_2^2$:
\begin{equation}\label{ac-Q-definition}
 \begin{aligned}
Q_u[h,k]:=(h',k')_2+\lambda(h,k)_2-(p-1)\sum_{\vv\in\V_0}\abs{u(\vv)}^{p-2}h(\vv)k(\vv).
 \end{aligned}
\end{equation}
Along the branches in \eqref{ac-physical-branch}, we have $\lambda=t^2$.
Moreover, denote the \emph{unconstrained Morse index} by
\begin{equation}\label{ac-free-index-definition}
 m(u):=\sup\bigl\{\dim L: L \text{ is a subspace of } H^1(\G) \text{ such that } Q_u[\varphi,\varphi]<0 \text{ for every } \varphi\in L\setminus\{0\}\bigr\},
\end{equation}
and recall the definition of the \emph{constrained Morse index} \cite[Definition~1.4]{Bort0}:
\begin{equation}\label{ac-constrained-index-definition}
 \hat m_0(u)=\sup\bigl\{\dim L: L \text{ is a subspace of } T_uH^1_\mu(\G) \text{ such that } Q_u[\varphi,\varphi]<0 \text{ for every } \varphi\in L\setminus\{0\}\bigr\}.
\end{equation}
Here $T_uH^1_\mu(\G):=\{v\in H^1(\G): (u,v)_2=0\}$ is the tangent space of $H^1_\mu(\G)$ at $u$. Directly from the definitions, we have $\hat m_0(u)\leq m(u)$.

 From now on we assume in addition that $p\ne4$. In Proposition~\ref{prop-ac-branches}, we determine the mass, the energy and the Morse indices along $u_S(t)$. The critical case $p=4$ is treated in Section~\ref{sect-critical}.
\begin{proposition}\label{prop-ac-branches}
Assume $p>2$ and $p\ne4$. For every finite nonempty $S\subset\V_0$, there is $C_S>0$, depending only on $p$, $\ell_*$, $\#S$ and $\max_{\vv\in S}\deg(\vv)$, such that for every $t>T_S$ 
\begin{align}
 \mu_S(t)&=\frac12
   \Bigl(\sum_{\vv\in S}\deg(\vv)^{p/(p-2)}\Bigr)
   t^{(4-p)/(p-2)}\bigl(1+\varepsilon_{S,1}(t)\bigr),
                                      \label{ac-mass-asymptotic}\\
 \mu_S'(t)&=-\frac{p-4}{2(p-2)}
   \Bigl(\sum_{\vv\in S}\deg(\vv)^{p/(p-2)}\Bigr)
   t^{(4-p)/(p-2)-1}\bigl(1+\varepsilon_{S,2}(t)\bigr),
                               \label{ac-mass-derivative-asymptotic}\\
 E(u_S(t),\G)&=\frac{p-4}{4p}
   \Bigl(\sum_{\vv\in S}\deg(\vv)^{p/(p-2)}\Bigr)
   t^{p/(p-2)}\bigl(1+\varepsilon_{S,3}(t)\bigr),
                                      \label{ac-energy-asymptotic}
\end{align}
where $\abs{\varepsilon_{S,1}(t)}\leq C_St e^{-t\ell_*/2}$, $\abs{\varepsilon_{S,2}(t)}\leq C_St^2e^{-t\ell_*/2}$ and $\abs{\varepsilon_{S,3}(t)}\leq C_St e^{-t\ell_*/2}$ and $T_S$ is the threshold of Lemma~\ref{ac-IFT-branch}.
There is $T_S^*\geq T_S$, depending only on $p$, $\ell_*$, $\#S$ and $\max_{\vv\in S}\deg(\vv)$, such that, for every $t\geq T_S^*$,
\begin{equation}\label{ac-index-equality}
 m(u_S(t))=\#S,\qquad
 \hat m_0(u_S(t))=
 \begin{cases}
 \#S,&p>4,\\
 \#S-1,&2<p<4.
 \end{cases}
\end{equation}
and $\mu_S'(t)>0$ if $2<p<4$ and $\mu_S'(t)<0$ if $p>4$. Moreover, for every $t\geq T_S^*$, the energy $E(u_S(t),\G)$ is positive if $p>4$ and negative if $2<p<4$. Furthermore, $\mu_S:[T_S^*,\infty)\to(0,\mu_S^*]$ is a bijection if $p>4$, and $\mu_S:[T_S^*,\infty)\to[\mu_S^*,+\infty)$ is a bijection if $2<p<4$, where $\mu_S^*:=\mu_S(T_S^*)$. Consequently, with $t_S$ the inverse of $\mu_S$, the branch $\{u_S(t)\}_{t>T_S}$ gives a positive solution $(u_S(t_S(\mu)),t_S(\mu)^2)$ of mass $\mu$ for every $0<\mu<\mu_S^*$ if $p>4$, and for every $\mu>\mu_S^*$ if $2<p<4$.
\end{proposition}

\begin{proof}
Throughout the proof $C_S>0$ denotes a constant, depending only on $p$, $\ell_*$, $\#S$ and $\max_{\vv\in S}\deg(\vv)$, whose value may change from line to line; it may be taken larger than the constants of Lemmas~\ref{ac-DtN-estimates} and~\ref{ac-IFT-branch} and than the bound $\norm{f_S^\infty}_X^2=\sum_{\vv\in S}\deg(\vv)^{p/(p-2)}\leq\#S(\max_{\vv\in S}\deg(\vv))^{p/(p-2)}$. Throughout the proof we also take $t>T_S$, so that the estimates of Lemma~\ref{ac-IFT-branch} hold.

We first prove \eqref{ac-mass-asymptotic}--\eqref{ac-energy-asymptotic}.
By \eqref{ac-mass-fixed-restriction},
\begin{equation}\label{ac-mass-exact}
 \mu_S(t)
 =t^{2/(p-2)}\norm{U_tf_S(t)}_2^2=\frac12t^{(4-p)/(p-2)}
   \bigl\langle (A_t+t\partial_tA_t)f_S(t),f_S(t)\bigr\rangle_X .
\end{equation}
By \eqref{ac-DtN-exponential}, we see that $\norm{A_t+t\partial_tA_t-I_X}_{\mathcal L(X)}\leq(C_0+C_1t)e^{-t\ell_*/2}$, so $\bigl|\bigl\langle(A_t+t\partial_tA_t-I_X)f_S(t),f_S(t)\bigr\rangle_X\bigr|\leq C_St e^{-t\ell_*/2}$. Moreover, by \eqref{ac-branch-estimates}, we obtain $\bigl|\norm{f_S(t)}_X^2-\norm{f_S^\infty}_X^2\bigr|\leq C_Se^{-t\ell_*/2}$. Hence, setting
\[
 \eta_S(t):=\bigl\langle (A_t+t\partial_tA_t)f_S(t),f_S(t)\bigr\rangle_X-\norm{f_S^\infty}_X^2,
 \qquad \abs{\eta_S(t)}\leq C_St e^{-t\ell_*/2},
\]
we have $\bigl\langle (A_t+t\partial_tA_t)f_S(t),f_S(t)\bigr\rangle_X=\norm{f_S^\infty}_X^2+\eta_S(t)$, and this proves \eqref{ac-mass-asymptotic} with $\varepsilon_{S,1}(t):=\eta_S(t)/\norm{f_S^\infty}_X^2$. Moreover, it follows from \eqref{ac-DtN-exponential}, \eqref{ac-branch-estimates} and $\bigl|\norm{f_S(t)}_X^2-\norm{f_S^\infty}_X^2\bigr|\leq C_Se^{-t\ell_*/2}$ that
\begin{equation}\label{ac-mass-quadratic-derivative}
 \begin{aligned}
 \frac{d}{dt}\bigl\langle (A_t+t\partial_tA_t)f_S(t),f_S(t)\bigr\rangle_X
 &=\bigl\langle (2\partial_tA_t+t\partial_t^2A_t)f_S(t),f_S(t)\bigr\rangle_X+2\bigl\langle (A_t+t\partial_tA_t)f_S(t),f_S'(t)\bigr\rangle_X\\
 &=\eta_S'(t),\qquad \abs{\eta_S'(t)}\leq C_St e^{-t\ell_*/2}.
 \end{aligned}
\end{equation}
Thus, differentiating \eqref{ac-mass-exact} gives
\[
 \begin{aligned}
 \mu_S'(t)
 ={}&-\frac{p-4}{2(p-2)}t^{(4-p)/(p-2)-1}
   \bigl\langle (A_t+t\partial_tA_t)f_S(t),f_S(t)\bigr\rangle_X\\
 &+\frac12t^{(4-p)/(p-2)}
   \frac{d}{dt}\bigl\langle (A_t+t\partial_tA_t)f_S(t),f_S(t)\bigr\rangle_X\\
 ={}&-\frac{p-4}{2(p-2)}
   t^{(4-p)/(p-2)-1}\norm{f_S^\infty}_X^2\bigl(1+\varepsilon_{S,2}(t)\bigr),
 \end{aligned}
\]
with
\[
 \varepsilon_{S,2}(t):=\frac{\eta_S(t)-\frac{p-2}{p-4}\,t\,\eta_S'(t)}{\norm{f_S^\infty}_X^2},\qquad \abs{\varepsilon_{S,2}(t)}\leq C_St^2e^{-t\ell_*/2},
\]
which proves \eqref{ac-mass-derivative-asymptotic}. 

Testing the weak equation of \eqref{eqpd} with $u_S(t)=t^{1/(p-2)}U_tf_S(t)$ gives the identity
\begin{equation}\label{ac-solution-energy-identity}
 \norm{u_S(t)'}_2^2+t^2\mu_S(t)
 =\sum_{\vv\in\V_0}\abs{u_S(t)(\vv)}^p
 =t^{p/(p-2)}\sum_{\vv\in\V_0}\abs{f_S(t)(\vv)}^p .
\end{equation}
For $g\in X$, \eqref{ac-norm-ordering} gives $\norm g_{\ell^p(\V)}^p\leq\norm g_\infty^{p-2}\norm g_{\ell^2(\V)}^2\leq\norm g_X^p$. Since $\norm{f_S(t)-f_S^\infty}_X\leq Ce^{-t\ell_*/2}$ by \eqref{ac-branch-estimates} and $\norm{f_S(t)}_{\ell^p(\V)}\leq\norm{f_S(t)}_X\leq\norm{f_S^\infty}_X+Ce^{-t\ell_*/2}$ for $t>T_S$, while $\sum_{\vv\in S}\deg(\vv)^{p/(p-2)}=\norm{f_S^\infty}_{\ell^p(\V)}^p$, the mean value 
inequality 
$\abs{a^p-b^p}\leq p\max\{a,b\}^{p-1}\abs{a-b}$, $a,b\geq0$, and H\"older's inequality give
\[
 \Bigl|\sum_{\vv\in\V_0}\abs{f_S(t)(\vv)}^p-\sum_{\vv\in S}\deg(\vv)^{p/(p-2)}\Bigr|\leq p\sum_{\vv\in\V_0}\max\bigl\{\abs{f_S(t)(\vv)},\abs{f_S^\infty(\vv)}\bigr\}^{p-1}\abs{f_S(t)(\vv)-f_S^\infty(\vv)}\leq C_Se^{-t\ell_*/2}.
\]
Using \eqref{ac-solution-energy-identity} and \eqref{ac-mass-asymptotic}, we obtain
\[
 \begin{aligned}
 E(u_S(t),\G)
 &=\Bigl(\frac12-\frac1p\Bigr)
    \sum_{\vv\in\V_0}\abs{u_S(t)(\vv)}^p-\frac{t^2}{2}\mu_S(t)\\
 &=\frac{p-4}{4p}
   \Bigl(\sum_{\vv\in S}\deg(\vv)^{p/(p-2)}\Bigr)t^{p/(p-2)}\bigl(1+\varepsilon_{S,3}(t)\bigr),\qquad \abs{\varepsilon_{S,3}(t)}\leq C_St e^{-t\ell_*/2},
 \end{aligned}
\]
which is \eqref{ac-energy-asymptotic}.

We now compute the Morse indices of the branch.   Since $\mathcal N'(f)$ is diagonal in the weighted space $X$, the operator $\mathcal N'(f)$ has operator norm
\[
 \norm{\mathcal N'(f)}_{\mathcal L(X)}=(p-1)\sup_{\vv\in\V_0}\frac{\abs{f(\vv)}^{p-2}}{\deg(\vv)}\leq(p-1)\norm f_X^{p-2},
\]
and thus, by \eqref{ac-branch-estimates},
\[
 \norm{\mathcal N'(f_S(t))}_{\mathcal L(X)}\leq(p-1)\bigl(\norm{f_S^\infty}_X+Ce^{-t\ell_*/2}\bigr)^{p-2}\qquad\text{for every }t>T_S .
\]
Define the bounded self-adjoint operator
\[
 L_S(t):=A_t-\mathcal N'(f_S(t)):X\to X .
\]
 Decompose $h\in H^1(\G)$ as in \eqref{ac-restriction-decomposition}, with restriction $g=\mathcal R h$ and interior part $z=h-U_tg$.  Using \eqref{ac-extension-orthogonality}, \eqref{ac-DtN-identity} and $\abs{u_S(t)(\vv)}^{p-2}=t\abs{f_S(t)(\vv)}^{p-2}$, the form \eqref{ac-Q-definition} becomes
\begin{equation}\label{ac-Hessian-splitting}
 \begin{aligned}
 Q_{u_S(t)}[h,h]
 &=\mathfrak b_t(z,z)+\mathfrak b_t(U_tg,U_tg)
   -t(p-1)\sum_{\vv\in\V_0}\abs{f_S(t)(\vv)}^{p-2}g(\vv)^2\\
 &=\mathfrak b_t(z,z)+t\bigl\langle L_S(t)g,g\bigr\rangle_X .
 \end{aligned}
\end{equation}
In particular, on $H^1_0(\G;\V)$ the form $Q_{u_S(t)}$ coincides with $\mathfrak b_t$ and is strictly positive; hence a negative direction must have a nonzero vertex restriction.  Recalling the subspace $X_S$ and the projection $P_S$ from \eqref{ac-coordinate-subspace}, we see that $\dim X_S=\#S$.  Since $\norm{A_t-I_X}_{\mathcal L(X)}\leq C_0e^{-t\ell_*/2}$ by \eqref{ac-DtN-exponential} and $\norm{\mathcal N'(f_S(t))-\mathcal N'(f_S^\infty)}_{\mathcal L(X)}\leq C_Se^{-\min\{1,p-2\}t\ell_*/2}$ for $t>T_S$ by \eqref{ac-N-lipschitz} and \eqref{ac-branch-estimates}, we have
\begin{equation}\label{ac-LS-uniform}
 \norm{L_S(t)-L_S^\infty}_{\mathcal L(X)}\leq C_Se^{-\min\{1,p-2\}t\ell_*/2}.
\end{equation}
By \eqref{ac-LS-uniform}, take $\widehat T_S\geq T_S$ so large that $e^{-\min\{1,p-2\}\widehat T_S\ell_*/2}\leq\min\{1,p-2\}/(2C_S)$. This $\widehat T_S$ depends only on $p$, $\ell_*$, $\#S$ and $\max_{\vv\in S}\deg(\vv)$, and $\norm{L_S(t)-L_S^\infty}_{\mathcal L(X)}\leq\tfrac12\min\{1,p-2\}$ for every $t>\widehat T_S$. Then \eqref{ac-limit-form} gives, for $g\in X_S$,
\[
 \begin{aligned}
 \bigl\langle L_S(t)g,g\bigr\rangle_X
 &\leq-(p-2)\norm g_X^2+\tfrac12\min\{1,p-2\}\norm g_X^2\\
 &\leq-\frac{p-2}{2}\norm g_X^2,
 \end{aligned}
\]
and, for $g\in X_S^\perp$,
\begin{equation}\label{ac-positive-part-lower-bound}
 \bigl\langle L_S(t)g,g\bigr\rangle_X
 \geq\bigl(1-\tfrac12\min\{1,p-2\}\bigr)\norm g_X^2
 \geq\tfrac12\norm g_X^2.
\end{equation}
Since $\mathcal RU_t=I_X$, the operator $U_t$ is injective, hence $\dim U_tX_S=\#S$. Thus, \eqref{ac-Hessian-splitting} shows that $Q_{u_S(t)}$ is negative definite on $U_tX_S$.  By \eqref{ac-free-index-definition}, $m(u_S(t))\geq\#S$.

Conversely, let $Y\subset H^1(\G)$ be a subspace on which $Q_{u_S(t)}$ is negative definite.  If $h\in Y$ is nonzero and $P_S\mathcal Rh=0$, then $g=\mathcal Rh\in X_S^\perp$, and \eqref{ac-Hessian-splitting} and \eqref{ac-positive-part-lower-bound} give
\[
 Q_{u_S(t)}[h,h]\geq\mathfrak b_t(z,z)+\frac t2\norm{g}_X^2>0,
\]
a contradiction.  Hence no nonzero $h\in Y$ satisfies $P_S\mathcal Rh=0$, so the linear map $(P_S\mathcal R)|_Y:Y\to X_S$ is injective and $\dim Y\leq\dim X_S=\#S$.  Taking the supremum over $Y$ proves $m(u_S(t))\leq\#S$, and thus, we conclude that  $m(u_S(t))=\#S$.

For fixed $h\in H^1(\G)$,  since $u_S$ is $C^1$ in $H^1(\G)$ and $E$ is $C^2$ for $p>2$, differentiating $E'(u_S(t),\G)[h]+t^2(u_S(t),h)_2=0$ with respect to $t$ gives
\[
E''(u_S(t),\G)[u_S'(t),h]+t^2(u_S'(t),h)_2+2t(u_S(t),h)_2=0,
\] 
hence,
\begin{equation}\label{ac-Q-branch-derivative}
 Q_{u_S(t)}[u_S'(t),h]=-2t(u_S(t),h)_2
 \qquad\forall h\in H^1(\G).
\end{equation}
Then, $2(u_S(t),u_S'(t))_2=\mu_S'(t)$ leads to
\begin{equation}\label{ac-positive-branch-direction}
 Q_{u_S(t)}[u_S'(t),u_S'(t)]=-t\mu_S'(t).
\end{equation}
For $t>\widehat T_S$ with $C_St^2e^{-t\ell_*/2}\leq\frac12$, we have $\abs{\varepsilon_{S,2}(t)}\leq\frac12$, hence $\mu_S'(t)\ne0$ by \eqref{ac-mass-derivative-asymptotic}. Moreover, the sign of $\mu_S'(t)$ is that of $4-p$. Define
\begin{equation}\label{ac-exact-tangent-correction}
 \Pi_t:H^1(\G)\longrightarrow T_{u_S(t)}H^1_{\mu_S(t)}(\G), \qquad
 \Pi_th:=h-\frac{2(u_S(t),h)_2}{\mu_S'(t)}u_S'(t).
\end{equation}
The identity $2(u_S(t),u_S'(t))_2=\mu_S'(t)$ gives $(u_S(t),\Pi_th)_2=0$. Moreover, $\Pi_th=h$ for every $h\in T_{u_S(t)}H^1_{\mu_S(t)}(\G)$, and $\ker\Pi_t=\operatorname{span}\{u_S'(t)\}$. Since $\mu_S'(t)\ne0$ and $2(u_S(t),u_S'(t))_2=\mu_S'(t)$, we have $u_S'(t)\notin T_{u_S(t)}H^1_{\mu_S(t)}(\G)$, and hence
\[
 H^1(\G)=T_{u_S(t)}H^1_{\mu_S(t)}(\G)
 \oplus\operatorname{span}\{u_S'(t)\}.
\]
Equivalently, by \eqref{ac-exact-tangent-correction}, every $h\in H^1(\G)$ splits as
\[
 h=\Pi_th+\frac{2(u_S(t),h)_2}{\mu_S'(t)}u_S'(t).
\]
We emphasize that here $\oplus$ denotes the algebraic direct sum: the decomposition is not necessarily orthogonal, and $\operatorname{span}\{u_S'(t)\}$ is not necessarily orthogonal to $T_{u_S(t)}H^1_{\mu_S(t)}(\G)$, neither for the $L^2$ nor for the $H^1$ inner product, so $\Pi_t$ is not necessarily an orthogonal projection. However, the two summands are orthogonal for the form $Q_{u_S(t)}$: by \eqref{ac-Q-branch-derivative},
\[
 Q_{u_S(t)}[\Pi_th,u_S'(t)]=-2t(u_S(t),\Pi_th)_2=0.
\]
Consequently,
\begin{equation}\label{ac-correction-formula}
 Q_{u_S(t)}[\Pi_th,\Pi_th]
 =Q_{u_S(t)}[h,h]+\frac{4t(u_S(t),h)_2^2}{\mu_S'(t)}.
\end{equation}

If $p>4$, then $\mu_S'(t)<0$, hence \eqref{ac-correction-formula} implies
\[
 Q_{u_S(t)}[\Pi_th,\Pi_th]\leq Q_{u_S(t)}[h,h]<0
 \qquad\text{for }h\in U_tX_S\setminus\{0\}.
\]
Thus, $\Pi_t$ is injective on $U_tX_S$, and $Q_{u_S(t)}$ is negative definite on its image $\Pi_t(U_tX_S)$, a subspace of the tangent space $T_{u_S(t)}H^1_{\mu_S(t)}(\G)$ of dimension $\#S$. Since $\hat m_0(u_S(t))\leq m(u_S(t))=\#S$, it follows that $\hat m_0(u_S(t))=\#S$.

If $2<p<4$, then $\mu_S'(t)>0$, and \eqref{ac-positive-branch-direction} shows that $Q_{u_S(t)}$ is negative definite on $\operatorname{span}\{u_S'(t)\}$. On the one hand, since $Q_{u_S(t)}$ is negative definite on $U_tX_S$ and $T_{u_S(t)}H^1_{\mu_S(t)}(\G)$ has codimension one in $H^1(\G)$,
\[
 \hat m_0(u_S(t))
 \geq\dim\bigl(U_tX_S\cap T_{u_S(t)}H^1_{\mu_S(t)}(\G)\bigr)
 \geq\#S-1.
\]
On the other hand, for every subspace $Y$ of the tangent space $T_{u_S(t)}H^1_{\mu_S(t)}(\G)$ on which $Q_{u_S(t)}$ is negative definite, the identities \eqref{ac-Q-branch-derivative} and \eqref{ac-positive-branch-direction}, together with $(u_S(t),h)_2=0$ for $h\in Y$, give
\[
 \begin{aligned}
 Q_{u_S(t)}[h+cu_S'(t),h+cu_S'(t)]
 &=Q_{u_S(t)}[h,h]-c^2t\mu_S'(t)<0
 \end{aligned}
\]
whenever $(h,c)\in Y\times\R\setminus\{(0,0)\}$. Hence $Y\oplus\operatorname{span}\{u_S'(t)\}$ is a subspace on which $Q_{u_S(t)}$ is negative definite, of dimension $\dim Y+1$, and
\[
 \dim Y+1\leq m(u_S(t))=\#S.
\]
Taking the supremum over $Y$ yields $\hat m_0(u_S(t))\leq\#S-1$. This proves \eqref{ac-index-equality} in both cases.

Since $\widehat T_S\geq T_S$, the positivity $U_tf_S(t)>0$ (see Lemma~\ref{ac-IFT-branch}) and the error bounds in \eqref{ac-mass-asymptotic}--\eqref{ac-energy-asymptotic} hold for every $t>\widehat T_S$; the equality \eqref{ac-index-equality} and the sign of $\mu_S'$ hold for every $t>\widehat T_S$ with $C_St^2e^{-t\ell_*/2}\leq\frac12$. Hence there is $T_S^*\geq\widehat T_S$, depending only on $p$, $\ell_*$, $\#S$ and $\max_{\vv\in S}\deg(\vv)$, such that all of them hold for every $t\geq T_S^*$. Put $\mu_S^*:=\mu_S(T_S^*)$. Then $\mu_S$ is continuous  and strictly monotone on $[T_S^*,\infty)$, and \eqref{ac-mass-asymptotic} gives $\mu_S(t)\to0$ if $p>4$ and $\mu_S(t)\to+\infty$ if $2<p<4$ as $t\to+\infty$, and thus, the continuous map
\begin{equation}\label{ac-mass-bijection}
 \begin{cases}
 \mu_S:[T_S^*,\infty)\longrightarrow(0,\mu_S^*],&p>4,\\
 \mu_S:[T_S^*,\infty)\longrightarrow[\mu_S^*,\infty),&2<p<4
 \end{cases}
\end{equation}
is a bijection in either case. This establishes the two bijections stated in Proposition~\ref{prop-ac-branches}, and hence the existence of the corresponding normalized solutions $u_S(t_S(\mu))$.
\end{proof}

We now provide the main result of this section. The following theorem constructs families of positive solutions $(u_{\mu,S},\lambda_{\mu,S})$ of \eqref{eqpd} at the same prescribed mass, indexed by the nonempty finite subsets $S\subset\V_0$. For prescribed bounds $k$ and $D$ on the number and degrees of the selected vertices, it gives a common mass range, determines the Morse indices, and provides uniform estimates for the multipliers and vertex profiles. It is the main step in the proofs of Theorem~\ref{th2} and Theorem~\ref{cor-ac-infinite-orbits}.

\begin{theorem}\label{th-ac-uniform-solutions}
Let $\G=(\V,\E)\in\mathbf G$, let $p\in(2,4)\cup(4,+\infty)$ and let $\varnothing\ne\V_0\subset\V$. For every integer $k\geq1$ and every integer $D\geq1$, there are constants $T_{k,D}>1$ and $\mu_{k,D}>0$, depending only on $p$, $\ell_*$, $k$ and $D$, such that, for every nonempty finite $W\subset\V_0$ satisfying $\#W\leq k$ and $\max_{\vv\in W}\deg(\vv)\leq D$, and every prescribed mass $\mu>0$ in the range
\[
 \begin{cases}
  0<\mu<\mu_{k,D},&p>4,\\
  \mu>\mu_{k,D},&2<p<4,
 \end{cases}
\]
problem \eqref{eqpd} has a family of $2^{\#W}-1$ pairwise distinct positive solutions
\[
 \bigl\{(u_{\mu,S},\lambda_{\mu,S}): \varnothing\ne S\subset W\bigr\}
 \subset H^1_\mu(\G)\times(0,+\infty).
\]
For every nonempty $S\subset W$, there is a unique $t_S(\mu)>T_{k,D}$ satisfying $\mu_S(t_S(\mu))=\mu$. More explicitly, with $u_{\mu,S}=u_S(t_S(\mu))$, $\lambda_{\mu,S}=t_S(\mu)^2$ and $\mu=\mu_S(t_S(\mu))$, these solutions satisfy
\[
 t_S(\mu)>T_{k,D},\qquad
 \bigl\|\lambda_{\mu,S}^{-1/(2(p-2))}\mathcal Ru_{\mu,S}-f_S^\infty\bigr\|_X<\frac13,
\]
as well as
\[
 m(u_{\mu,S})=\#S,\qquad
 \begin{cases}
  E(u_{\mu,S},\G)>0,\quad\hat m_0(u_{\mu,S})=\#S,&p>4,\\
  E(u_{\mu,S},\G)<0,\quad\hat m_0(u_{\mu,S})=\#S-1,&2<p<4.
 \end{cases}
\]
Moreover, there is a constant $C=C(p,\ell_*,k,D)>0$ such that, for every nonempty finite $S\subset\V_0$ satisfying $\#S\leq k$ and $\max_{\vv\in S}\deg(\vv)\leq D$, and every $\mu$ in the corresponding mass range above,
\[
 \lambda_{\mu,S}\geq
 \begin{cases}
  \displaystyle\left(\frac{1}{4\mu}\right)^{2(p-2)/(p-4)},&p>4,\\[2mm]
  \displaystyle\left(\frac{\mu}{2kD^{p/(p-2)}}\right)^{2(p-2)/(4-p)},&2<p<4,
 \end{cases}
\]
and
\[
 \bigl\|\lambda_{\mu,S}^{-1/(2(p-2))}\mathcal Ru_{\mu,S}-f_S^\infty\bigr\|_X
 \leq
 \begin{cases}
  \displaystyle C\exp\!\left[-\frac{\ell_*}{2}\left(\frac{1}{4\mu}\right)^{(p-2)/(p-4)}\right],&p>4,\\[2mm]
  \displaystyle C\exp\!\left[-\frac{\ell_*}{2}\left(\frac{\mu}{2kD^{p/(p-2)}}\right)^{(p-2)/(4-p)}\right],&2<p<4.
 \end{cases}
\]

\end{theorem}
\begin{proof}
Fix integers $k\geq1$ and $D\geq1$. For every nonempty finite $S\subset\V_0$ with $\#S\leq k$ and $\max_{\vv\in S}\deg(\vv)\leq D$,
\begin{equation}\label{ac-fS-bound}
 1\leq\norm{f_S^\infty}_X^2
 =\sum_{\vv\in S}\deg(\vv)^{p/(p-2)}
 \leq kD^{p/(p-2)}.
\end{equation}
The constants and thresholds in Lemma~\ref{ac-IFT-branch} and Proposition~\ref{prop-ac-branches} can be bounded using only $p$, $\ell_*$, $\#S$ and $\max_{\vv\in S}\deg(\vv)$. Moreover, the pairs $(\#S,\max_{\vv\in S}\deg(\vv))$ range over the finitely many integer pairs $(m,d)$ with $1\leq m\leq k$ and $1\leq d\leq D$. Hence there are constants $C=C(p,\ell_*,k,D)>0$ and $T_0=T_0(p,\ell_*,k,D)>1$ such that $T_0 \geq T_S^*$ and, for every such $S$ and every $t\geq T_0$,
\begin{equation}\label{ac-S-uniform-estimates}
 \norm{f_S(t)-f_S^\infty}_X\leq Ce^{-t\ell_*/2},\qquad
 \abs{\varepsilon_{S,1}(t)}\leq Ct e^{-t\ell_*/2}.
\end{equation}
Choose $T_{k,D}>\max\{T_0,1,2/\ell_*\}$ so large that $CT_{k,D}e^{-T_{k,D}\ell_*/2}<1/3$. Since the functions $t\mapsto e^{-t\ell_*/2}$ and $t\mapsto te^{-t\ell_*/2}$ are both decreasing for $t>2/\ell_*$, and $t\geq T_{k,D}>2/\ell_*$, we know that, for every $t\geq T_{k,D}$,
\begin{equation}\label{ac-S-epsilon-estimates}
 \norm{f_S(t)-f_S^\infty}_X<\frac13,\qquad
 \abs{\varepsilon_{S,1}(t)}<\frac12.
\end{equation}
Substituting \eqref{ac-S-epsilon-estimates} and \eqref{ac-fS-bound} into \eqref{ac-mass-asymptotic} gives, for every $t\geq T_{k,D}$,
\begin{equation}\label{ac-two-sided-mass-bound}
 \frac14t^{(4-p)/(p-2)}
 \leq\mu_S(t)
 \leq2kD^{p/(p-2)}t^{(4-p)/(p-2)}.
\end{equation}
Define
\begin{equation}\label{ac-common-mass-threshold}
 \mu_{k,D}:=
 \begin{cases}
  \displaystyle\frac14T_{k,D}^{(4-p)/(p-2)},&p>4,\\[2mm]
  \displaystyle2kD^{p/(p-2)}T_{k,D}^{(4-p)/(p-2)},&2<p<4.
 \end{cases}
\end{equation}
Both $T_{k,D}$ and $\mu_{k,D}$ depend only on $p$, $\ell_*$, $k$ and $D$, and $0<\mu_{k,D}<+\infty$.

Let $W\subset\V_0$ be nonempty and finite with $\#W\leq k$ and $\max_{\vv\in W}\deg(\vv)\leq D$, and let $\mu>0$ satisfy $0<\mu<\mu_{k,D}$ if $p>4$ and $\mu>\mu_{k,D}$ if $2<p<4$. For every nonempty $S\subset W$, Proposition~\ref{prop-ac-branches} gives strict monotonicity of $\mu_S$ on $[T_{k,D},+\infty)$. If $p>4$, then $\mu_S(T_{k,D})\geq\mu_{k,D}>\mu$ and $\mu_S(t)\to0$; if $2<p<4$, then $\mu_S(T_{k,D})\leq\mu_{k,D}<\mu$ and $\mu_S(t)\to+\infty$. Consequently, there is a unique $t_S(\mu)>T_{k,D}$ such that $\mu_S(t_S(\mu))=\mu$. Set
\[
 u_{\mu,S}:=u_S(t_S(\mu)),\qquad
 \lambda_{\mu,S}:=t_S(\mu)^2.
\]
By \eqref{ac-branch-weak-identity}, $(u_{\mu,S},\lambda_{\mu,S})\in H^1_\mu(\G) \times (0,+\infty)$ is a positive solution of \eqref{eqpd} with $\norm{u_{\mu,S}}_2^2=\mu$. Since $t_S(\mu)>T_{k,D}>T_S^*$, its energy sign and both Morse indices are those given in Proposition~\ref{prop-ac-branches}.

The two-sided mass bound \eqref{ac-two-sided-mass-bound}, evaluated at $t=t_S(\mu)$ with $\mu_S(t_S(\mu))=\mu$, gives
\begin{equation}
\label{ac-tS-lower-bound}
 t_S(\mu)\geq
 \begin{cases}
  \displaystyle\left(\frac{1}{4\mu}\right)^{(p-2)/(p-4)},&p>4,\\[2mm]
  \displaystyle\left(\frac{\mu}{2kD^{p/(p-2)}}\right)^{(p-2)/(4-p)},&2<p<4.
 \end{cases}
\end{equation}
Since $\lambda_{\mu,S}=t_S(\mu)^2$, we obtain
\[
 \lambda_{\mu,S}\geq
 \begin{cases}
  \displaystyle\left(\frac{1}{4\mu}\right)^{2(p-2)/(p-4)},&p>4,\\[2mm]
  \displaystyle\left(\frac{\mu}{2kD^{p/(p-2)}}\right)^{2(p-2)/(4-p)},&2<p<4.
 \end{cases}
\]
By \eqref{ac-physical-branch} and $\mathcal RU_t=I_X$, we have $\lambda_{\mu,S}^{-1/(2(p-2))}\mathcal Ru_{\mu,S}=f_S(t_S(\mu))$. Since $t_S(\mu)>T_{k,D}$, \eqref{ac-S-epsilon-estimates} proves $\norm{f_S(t_S(\mu))-f_S^\infty}_X<\frac13$. Moreover, \eqref{ac-S-uniform-estimates} and \eqref{ac-tS-lower-bound} imply
\[
 \begin{aligned}
 \bigl\|\lambda_{\mu,S}^{-1/(2(p-2))}\mathcal Ru_{\mu,S}-f_S^\infty\bigr\|_X
 &=\norm{f_S(t_S(\mu))-f_S^\infty}_X\\
 &\leq Ce^{-t_S(\mu)\ell_*/2}\\
 &\leq
 \begin{cases}
  \displaystyle C\exp\!\left[-\frac{\ell_*}{2}\left(\frac{1}{4\mu}\right)^{(p-2)/(p-4)}\right],&p>4,\\[2mm]
  \displaystyle C\exp\!\left[-\frac{\ell_*}{2}\left(\frac{\mu}{2kD^{p/(p-2)}}\right)^{(p-2)/(4-p)}\right],&2<p<4,
 \end{cases}
 \end{aligned}
\]
where $C=C(p,\ell_*,k,D)>0$ is independent of $S$ and $\mu$. Consequently, uniformly over all nonempty finite $S\subset\V_0$ satisfying $\#S\leq k$ and $\max_{\vv\in S}\deg(\vv)\leq D$,
\begin{equation*}
 \begin{aligned}
  \lambda_{\mu,S}&\longrightarrow+\infty,\\
  \bigl\|\lambda_{\mu,S}^{-1/(2(p-2))}\mathcal Ru_{\mu,S}-f_S^\infty\bigr\|_X&\longrightarrow0,
 \end{aligned}
 \qquad\text{as}\quad
 \begin{cases}
  \mu\to0^+,&p>4,\\
  \mu\to+\infty,&2<p<4.
 \end{cases}
\end{equation*}

It remains to prove pairwise distinctness. If $S_1\ne S_2$ are nonempty subsets of $W$, then
\[
 \norm{f_{S_1}^\infty-f_{S_2}^\infty}_X^2
 =\sum_{\vv\in S_1\setminus S_2}\deg(\vv)^{p/(p-2)}
  +\sum_{\vv\in S_2\setminus S_1}\deg(\vv)^{p/(p-2)}\geq1.
\]
Suppose that $u_{\mu,S_1}=u_{\mu,S_2}=:u$. On every edge,
\[
 u''=\lambda_{\mu,S_1}u=\lambda_{\mu,S_2}u,\qquad
 (\lambda_{\mu,S_1}-\lambda_{\mu,S_2})u=0.
\]
Since $u>0$, we have $\lambda_{\mu,S_1}=\lambda_{\mu,S_2}$ and hence $t_{S_1}(\mu)=\sqrt{\lambda_{\mu,S_1}}=\sqrt{\lambda_{\mu,S_2}}=t_{S_2}(\mu)=:t>T_{k,D}$. Equality of the vertex restrictions gives $f_{S_1}(t)=\lambda_{\mu,S_1}^{-1/(2(p-2))}\mathcal Ru_{\mu,S_1}=\lambda_{\mu,S_2}^{-1/(2(p-2))}\mathcal Ru_{\mu,S_2}=f_{S_2}(t)$, and hence
\[
 1\leq\norm{f_{S_1}^\infty-f_{S_2}^\infty}_X
 \leq\norm{f_{S_1}^\infty-f_{S_1}(t)}_X+\norm{f_{S_2}(t)-f_{S_2}^\infty}_X
 <\frac23,
\]
a contradiction. The $2^{\#W}-1$ solutions are therefore pairwise distinct.
\end{proof}
We are now in a position to prove Theorem~\ref{th2} and Theorem~\ref{cor-ac-infinite-orbits}.
\begin{proof}[Proofs of Theorems~\ref{th2} and \ref{cor-ac-infinite-orbits}]
For a given nonempty finite $W\subset\V_0$, take $k:=\#W$ and $D:=\max_{\vv\in W}\deg(\vv)$ in Theorem~\ref{th-ac-uniform-solutions} and put $\mu_W:=\mu_{k,D}$. All conclusions of Theorem~\ref{th2} follow.

To prove Theorem~\ref{cor-ac-infinite-orbits}, let $\V_*:=\{\vv\in\V_0:\deg(\vv)\leq D\}$, which is an infinite set by assumption, and fix $\vv_*\in\V_*$. Since $\G$ is locally finite, $\ell_*>0$ and $\V_*$ is infinite, we can choose vertices $\vv_n\in\V_*\setminus\{\vv_*\}$ such that
\begin{equation}\label{ac-vertex-distances}
 \operatorname{dist}(\vv_*,\vv_{n+1})>\operatorname{dist}(\vv_*,\vv_n),
 \qquad n\geq1.
\end{equation}
Set $S_n:=\{\vv_*,\vv_n\}$. For every $n\geq1$,
\[
 \#S_n=2,\qquad\max_{\vv\in S_n}\deg(\vv)\leq D,\qquad
 2\leq\norm{f_{S_n}^\infty}_X^2\leq2D^{p/(p-2)}.
\]
Apply Theorem~\ref{th-ac-uniform-solutions} with $k=2$ and select the solution indexed by $S=W=S_n$. Its thresholds $\mu_{2,D}>0$ and $T_{2,D}>1$, depending only on $p$, $\ell_*$ and $D$, apply to every $n$. Set $\mu_*:=\mu_{2,D}$. For every $\mu$ with $0<\mu<\mu_*$ if $p>4$ and $\mu>\mu_*$ if $2<p<4$, Theorem~\ref{th-ac-uniform-solutions} provides $t_n:=t_{S_n}(\mu)>T_{2,D}$ with $\mu_{S_n}(t_n)=\mu$, and the resulting solutions satisfy
\[
 u_{\mu,n}:=u_{\mu,S_n}=u_{S_n}(t_n),\qquad
 \lambda_{\mu,n}:=t_n^2,\qquad
 \norm{f_{S_n}(t_n)-f_{S_n}^\infty}_X<\frac13.
\]
The energy signs and both Morse indices follow from Theorem~\ref{th-ac-uniform-solutions}.
Applying \eqref{ac-tS-lower-bound} to $t_n=\sqrt{\lambda_{\mu,n}}$ gives
\[
 \lim_{\mu\to0^+}\inf_{n\geq1}\lambda_{\mu,n}=+\infty
 \quad\text{if }p>4,\qquad
 \lim_{\mu\to+\infty}\inf_{n\geq1}\lambda_{\mu,n}=+\infty
 \quad\text{if }2<p<4 .
\]

We show that no two of these solutions are $\Gamma$-equivalent. Every $\gamma\in\Gamma$ is an isometry of $\G$, preserves vertex degrees, and induces isometries on $H^1(\G)$, $L^2(\G)$ and $X$. Moreover,
\[
 \mathcal R(\gamma u)=\gamma\mathcal Ru,\qquad
 \gamma f_S^\infty=f_{\gamma S}^\infty.
\]
Suppose that $u_{\mu,n}=\gamma u_{\mu,m}$ for some $\gamma\in\Gamma$. Since $\gamma$ preserves edge lengths and $\V_0$, the function $\gamma u_{\mu,m}$ satisfies $(\gamma u_{\mu,m})''=t_m^2\,\gamma u_{\mu,m}$ on every edge $\e \in \E$, whereas $u_{\mu,n}''=t_n^2u_{\mu,n}$ on every edge $\e \in \E$. Hence
\[
 (t_n^2-t_m^2)u_{\mu,n}=0
 \quad\text{on every edge } \e \in \E,
\]
and $u_{\mu,n}>0$ on $\G$ gives $t_n=t_m=:t$. 
Since $t_n=t_m$ and $\mathcal R(\gamma u)=\gamma\mathcal Ru$, taking vertex restrictions yields 
\[f_{S_n}(t)=\lambda_{\mu,S_n}^{-1/(2(p-2))}\mathcal Ru_{\mu,S_n}=\lambda_{\mu,S_m}^{-1/(2(p-2))}\gamma\mathcal Ru_{\mu,S_m}=\gamma f_{S_m}(t).\]
Argue by contradiction, and suppose that $S_n\ne\gamma S_m$. Then
\[
 \norm{f_{S_n}^\infty-f_{\gamma S_m}^\infty}_X^2
 =\sum_{\vv\in S_n\triangle\gamma S_m}\deg(\vv)^{p/(p-2)}\geq1,
\]
where $S_n\triangle\gamma S_m:=(S_n\setminus\gamma S_m)\cup(\gamma S_m\setminus S_n)$ is the symmetric difference of $S_n$ and $\gamma S_m$. Moreover, the uniform profile estimate \eqref{ac-S-epsilon-estimates}, applied to $S_n$ and to $S_m$, gives
\[
 \begin{aligned}
 \norm{f_{S_n}^\infty-f_{\gamma S_m}^\infty}_X
 &\leq\norm{f_{S_n}^\infty-f_{S_n}(t)}_X
       +\norm{\gamma f_{S_m}(t)-\gamma f_{S_m}^\infty}_X\\
 &<\frac13+\frac13=\frac23,
 \end{aligned}
\]
which is a contradiction. Thus $S_n=\gamma S_m$. Since $\gamma$ preserves distances, the two vertices of $S_n$ and the two vertices of $\gamma S_m$ satisfy $\operatorname{dist}(\vv_*,\vv_n)=\operatorname{dist}(\gamma\vv_*,\gamma\vv_m)=\operatorname{dist}(\vv_*,\vv_m)$, and the strict increase \eqref{ac-vertex-distances} forces $n=m$. This proves that the solutions are pairwise $\Gamma$-inequivalent, as asserted in Theorem~\ref{cor-ac-infinite-orbits}.
\end{proof}
Now, we are in a position to prove Theorem~\ref{cor-ac-supercritical-all-mass}.
\begin{proof}[Proof of Theorem~\ref{cor-ac-supercritical-all-mass}]
Since the degrees of infinitely many defect vertices lie in $\{1,\ldots,D\}$, there are $q\in\{1,\ldots,D\}$ and an infinite set $W\subset\V_0$ such that $\deg(\vv)=q$ for every $\vv\in W$.
For $t\geq1$ and a bounded edge $\e\succ\vv$, we have 
\[
 \coth(t\ell_\e)-1+\operatorname{csch}(t\ell_\e)
 =\frac{\cosh(t\ell_\e)+1}{\sinh(t\ell_\e)}-1
 =\coth(t\ell_\e/2)-1
 \leq\coth(t\ell_*/2)-1.
\]
Therefore, by \eqref{ac-DtN-coordinate}, for $g\in X$ and $t \geq 1$,  using $\coth(t\ell_*/2)-1=2e^{-t\ell_*}/(1-e^{-t\ell_*})$, we obtain
\begin{equation}\label{ac-all-mass-DtN-sup}
 \norm{(A_t-I_X)g}_\infty
 \leq\bigl[\coth(t\ell_*/2)-1\bigr]\norm g_\infty
 \leq\frac{2e^{-t\ell_*}}{1-e^{-\ell_*}}\norm g_\infty.
\end{equation}
Fix a nonempty finite $S\subset W$. Then $\norm{f_S^\infty}_\infty=q^{1/(p-2)}$ by \eqref{ac-limit-profile}. Put
\[
 C=(p-1)(p-2)\bigl[q^{1/(p-2)}+1\bigr]^{p-3},
 \qquad
 \delta=\min\bigl\{1/3,1/(4C)\bigr\}.
\]
For $g,h\in X$ with $\norm g_\infty,\norm h_\infty\leq\delta$, the fundamental theorem of calculus and \eqref{ac-N-derivative} give, for every $\vv\in\V$,
\[
 \begin{aligned}
 &\bigl[\mathcal N(f_S^\infty+g)-\mathcal N(f_S^\infty+h)
                  -\mathcal N'(f_S^\infty)(g-h)\bigr](\vv)\\
 &\quad=\int_0^1
   \bigl[(\mathcal N'(f_S^\infty+(1-\theta)h+\theta g)-\mathcal N'(f_S^\infty))(g-h)\bigr](\vv)\,d\theta\\
 &\quad=\frac{(p-1)\mathbf1_{\V_0}(\vv)}{\deg(\vv)}[g(\vv)-h(\vv)]
    \int_0^1
    \Bigl[\abs{f_S^\infty(\vv)+(1-\theta)h(\vv)+\theta g(\vv)}^{p-2}
          -\abs{f_S^\infty(\vv)}^{p-2}\Bigr]\,d\theta .
 \end{aligned}
\]
Since $p>4$,
\begin{equation}\label{ac-all-mass-scalar-derivative}
 \left|\frac{d}{ds}\abs{s}^{p-2}\right|
 =(p-2)\abs{s}^{p-3}
 \leq(p-2)\bigl[q^{1/(p-2)}+1\bigr]^{p-3}
 \quad\text{for } |s|\leq q^{1/(p-2)}+1.
\end{equation}
Therefore, by $\norm{f_S^\infty}_\infty=q^{1/(p-2)}$ and $\abs{(1-\theta)h(\vv)+\theta g(\vv)}\leq\delta\leq1$, the mean value theorem gives
\[
 \begin{aligned}
 \left|\abs{f_S^\infty(\vv)+(1-\theta)h(\vv)+\theta g(\vv)}^{p-2}
           -\abs{f_S^\infty(\vv)}^{p-2}\right|&\leq(p-2)\bigl[q^{1/(p-2)}+1\bigr]^{p-3}
               \abs{(1-\theta)h(\vv)+\theta g(\vv)}\\
 &\leq(p-2)\bigl[q^{1/(p-2)}+1\bigr]^{p-3}\delta .
 \end{aligned}
\]
Since $\mathbf1_{\V_0}(\vv)/\deg(\vv)\leq1$, we obtain
\begin{equation}\label{ac-all-mass-pointwise-lipschitz}
 \abs{\bigl[\mathcal N(f_S^\infty+g)-\mathcal N(f_S^\infty+h)
                  -\mathcal N'(f_S^\infty)(g-h)\bigr](\vv)}
 \leq C\delta\abs{g(\vv)-h(\vv)}.
\end{equation}

Let $C_0,C_1$ be the constants in \eqref{ac-DtN-exponential} and choose $T>\max\{1,2/\ell_*\}$ with
\[
 C_0e^{-T\ell_*/2}\leq\frac14,\qquad
 \frac{2q^{1/(p-2)}e^{-T\ell_*}}{1-e^{-\ell_*}}\leq\frac{\delta}{2}.
\]
These inequalities hold with $T$ replaced by every $t\geq T$, and $T$ depends only on $p$, $\ell_*$ and $q$. On the set $\{g\in X:\norm g_\infty\leq\delta\}$, define
\[
 \Phi_{t,S}(g):=(L_S^\infty)^{-1}
 \bigl[-(A_t-I_X)(f_S^\infty+g)+\mathcal N(f_S^\infty+g)
       -\mathcal N(f_S^\infty)-\mathcal N'(f_S^\infty)g\bigr].
\]
By \eqref{ac-DtN-exponential}, \eqref{ac-all-mass-DtN-sup} and the choice of $T$, for every $t\geq T$ and $g\in X$,
\[
 \norm{A_t-I_X}_{\mathcal L(X)}\leq\frac14,
 \qquad
 \norm{(A_t-I_X)g}_\infty\leq\frac14\norm g_\infty .
\]
By \eqref{ac-limit-inverse-coordinate} and $p>4$,
\[
 \norm{(L_S^\infty)^{-1}}_{\mathcal L(X)}\leq1,
 \qquad
 \norm{(L_S^\infty)^{-1}g}_\infty\leq\norm g_\infty
 \quad\text{for every }g\in X .
\]
Combining these bounds with \eqref{ac-all-mass-pointwise-lipschitz} and $C\delta\leq1/4$, we obtain, for every $t\geq T$ and $g,h\in X$ with $\norm g_\infty,\norm h_\infty\leq\delta$,
\[
 \norm{\Phi_{t,S}(g)-\Phi_{t,S}(h)}_X
 \leq\frac12\norm{g-h}_X,\qquad
 \norm{\Phi_{t,S}(g)-\Phi_{t,S}(h)}_\infty
 \leq\frac12\norm{g-h}_\infty .
\]
By the definition of $\Phi_{t,S}$, $\Phi_{t,S}(0)=-(L_S^\infty)^{-1}(A_t-I_X)f_S^\infty$. Using \eqref{ac-limit-inverse-coordinate}, \eqref{ac-all-mass-DtN-sup} and the choice of $T$, we obtain for every $t\geq T$,
\[
 \begin{aligned}
 \norm{\Phi_{t,S}(0)}_\infty\leq\norm{(A_t-I_X)f_S^\infty}_\infty\leq\frac{2q^{1/(p-2)}e^{-t\ell_*}}{1-e^{-\ell_*}}
 \leq\frac{\delta}{2}.
 \end{aligned}
\]
For every $g\in X$ with $\norm g_\infty\leq\delta$, the contraction estimate gives
\[
 \norm{\Phi_{t,S}(g)}_\infty
 \leq\norm{\Phi_{t,S}(g)-\Phi_{t,S}(0)}_\infty+\norm{\Phi_{t,S}(0)}_\infty
 \leq\frac12\norm g_\infty+\frac{\delta}{2}\leq\delta .
\]
Thus
\[
 \Phi_{t,S}\bigl(\{g\in X:\norm g_\infty\leq\delta\}\bigr)
 \subseteq\{g\in X:\norm g_\infty\leq\delta\}.
\]
The set $\{g\in X:\norm g_\infty\leq\delta\}$ is closed in $X$ by \eqref{ac-norm-ordering}. 
The contraction theorem therefore gives a fixed point $g_{t,S}$ satisfying
\begin{equation}\label{ac-all-mass-fixed-point}
 g_{t,S}=\Phi_{t,S}(g_{t,S})\qquad\text{for every }t\geq T.
\end{equation}
By the definition of $\Phi_{t,S}$, $f_S(t):=f_S^\infty+g_{t,S}$ satisfies \eqref{ac-vertex-equation}. Arguing as in the proof of Lemma~\ref{ac-IFT-branch}, see \eqref{ac-uniform-fixedpoint-error}, and using \eqref{ac-DtN-exponential}, we obtain
\begin{equation}\label{ac-all-mass-fixed-point-bound}
 \norm{g_{t,S}}_X\leq2\norm{\Phi_{t,S}(0)}_X
 \leq2C_0 e^{-t\ell_*/2}\norm{f_S^\infty}_X.
\end{equation}
For $t,s\geq T$, evaluating the two maps at the same argument $g_{s,S}$ gives
\[
 \begin{aligned}
 \Phi_{t,S}(g_{s,S})-\Phi_{s,S}(g_{s,S})
 &=-(L_S^\infty)^{-1}\bigl[(A_t-I_X)-(A_s-I_X)\bigr]
       (f_S^\infty+g_{s,S})\\
 &=-(L_S^\infty)^{-1}(A_t-A_s)f_S(s).
 \end{aligned}
\]
By \eqref{ac-all-mass-fixed-point}, the contraction estimate in $X$ and the bound $\norm{(L_S^\infty)^{-1}}_{\mathcal L(X)}\leq1$, we have
\[
 \begin{aligned}
 \norm{g_{t,S}-g_{s,S}}_X
 &\leq\norm{\Phi_{t,S}(g_{t,S})-\Phi_{t,S}(g_{s,S})}_X
      +\norm{\Phi_{t,S}(g_{s,S})-\Phi_{s,S}(g_{s,S})}_X\\
 &\leq\frac12\norm{g_{t,S}-g_{s,S}}_X
      +\norm{(A_t-A_s)f_S(s)}_X .
 \end{aligned}
\]
Since $f_S(t)-f_S(s)=g_{t,S}-g_{s,S}$, we have $\norm{f_S(t)-f_S(s)}_X\leq2\norm{(A_t-A_s)f_S(s)}_X$. By Lemma~\ref{ac-DtN-estimates}, $t\mapsto A_t$ is continuous in $\mathcal L(X)$, so $f_S$ is continuous in $X$.

Write $f_S^-(t)=\min\{f_S(t),0\}$, so $\norm{f_S^-(t)}_\infty\leq\delta$. Testing \eqref{ac-vertex-equation} with $f_S^-(t)$ as in the proof of Lemma~\ref{ac-IFT-branch} and using \eqref{ac-DtN-coercivity}, we obtain
\[
 \tanh(t\ell_*/2)\norm{f_S^-(t)}_X^2
 \leq\sum_{\vv\in\V_0}\abs{f_S^-(t)(\vv)}^p
 \leq\delta^{p-2}\norm{f_S^-(t)}_X^2.
\]
Since $T>2/\ell_*$, $p>4$ and $\delta\leq1/3$, we have $f_S(t)\geq0$. Moreover, by \eqref{ac-all-mass-fixed-point-bound} and the choice of $T$,
\[
 \norm{f_S(t)}_X
 \geq\norm{f_S^\infty}_X-\norm{g_{t,S}}_X
 \geq\bigl(1-2C_0e^{-t\ell_*/2}\bigr)\norm{f_S^\infty}_X
 \geq\frac12\norm{f_S^\infty}_X>0 .
\]
Arguing as in the positivity part of the proof of Lemma~\ref{ac-IFT-branch}, we conclude that $U_tf_S(t)>0$ on $\G$.

Set $L_S(t):=A_t-\mathcal N'(f_S(t))$. Since $\norm{f_S^\infty}_\infty=q^{1/(p-2)}$ and $\norm{f_S(t)-f_S^\infty}_\infty=\norm{g_{t,S}}_\infty\leq\delta\leq1$, the values $f_S(t)(\vv)$ and $f_S^\infty(\vv)$ lie in the interval in \eqref{ac-all-mass-scalar-derivative}. By \eqref{ac-N-derivative} and the mean value theorem, for every $h\in X$ and $\vv\in\V$,
\[
 \begin{aligned}
 \bigl|\bigl[(\mathcal N'(f_S(t))-\mathcal N'(f_S^\infty))h\bigr](\vv)\bigr|&=(p-1)\frac{\mathbf1_{\V_0}(\vv)}{\deg(\vv)}
    \bigl|\abs{f_S(t)(\vv)}^{p-2}-\abs{f_S^\infty(\vv)}^{p-2}\bigr|\,|h(\vv)|\\
 &\leq C|f_S(t)(\vv)-f_S^\infty(\vv)|\,|h(\vv)|
 \leq C\delta|h(\vv)|.
 \end{aligned}
\]
Taking the $X$ norm and the supremum over $\norm h_X=1$ gives $\norm{\mathcal N'(f_S(t))-\mathcal N'(f_S^\infty)}_{\mathcal L(X)}\leq C\delta$. Hence, by \eqref{ac-limit-linearization}, \eqref{ac-DtN-exponential} and the choice of $T$,
\[
 \norm{L_S(t)-L_S^\infty}_{\mathcal L(X)}
 \leq C_0e^{-t\ell_*/2}+C\delta\leq\frac12
 \qquad\text{for every }t\geq T .
\]
As in the proof of Lemma~\ref{ac-IFT-branch}, write
\[
 L_S(t)=L_S^\infty(I_X+K_S(t)),\qquad
 K_S(t):=(L_S^\infty)^{-1}\bigl(L_S(t)-L_S^\infty\bigr).
\]
By \eqref{ac-limit-inverse-bound} and the bound on $L_S(t)-L_S^\infty$,
\[
 \norm{K_S(t)}_{\mathcal L(X)}
 \leq\norm{(L_S^\infty)^{-1}}_{\mathcal L(X)}
      \norm{L_S(t)-L_S^\infty}_{\mathcal L(X)}\leq\frac12.
\]
Thus $I_X+K_S(t)$ is invertible by \cite[Exercise~6.14]{Ha}, and
\begin{equation}\label{ac-all-mass-inverse-bound}
 \norm{L_S(t)^{-1}}_{\mathcal L(X)}
 \leq\frac{\norm{(L_S^\infty)^{-1}}_{\mathcal L(X)}}
           {1-\norm{K_S(t)}_{\mathcal L(X)}}\leq2 .
\end{equation}
We apply the implicit function theorem to \eqref{ac-vertex-equation} along the continuous branch as in the proof of Lemma~\ref{ac-IFT-branch} and obtain $f_S\in C^1((T,+\infty),X)$. By \eqref{ac-all-mass-fixed-point-bound} and $C_0e^{-T\ell_*/2}\leq1/4$, we have $\norm{f_S(t)}_X\leq\frac32\norm{f_S^\infty}_X$ for $t\geq T$. Differentiating \eqref{ac-vertex-equation}, we obtain
\[
 L_S(t)f_S'(t)=-(\partial_t A_t)f_S(t),
 \]
 hence, \eqref{ac-DtN-exponential} and \eqref{ac-all-mass-inverse-bound} lead to 
 \begin{equation}\label{ac-all-mass-branch-derivative-bound}
 \norm{f_S'(t)}_X
 \leq2C_1e^{-t\ell_*/2}\norm{f_S(t)}_X
 \leq3C_1e^{-t\ell_*/2}\norm{f_S^\infty}_X .
\end{equation}
Define $u_S(t)$ and $\mu_S(t)$ by \eqref{ac-physical-branch}. By Lemma~\ref{ac-restriction-extension} and the regularity of $f_S$, $u_S$ is $C^1$ on $(T,+\infty)$.

Use $\eta_S(t)=\langle(A_t+t\partial_tA_t)f_S(t),f_S(t)\rangle_X-\norm{f_S^\infty}_X^2$ as in the proof of Proposition~\ref{prop-ac-branches}. By \eqref{ac-DtN-exponential}, \eqref{ac-all-mass-fixed-point-bound} and $\norm{f_S(t)}_X\leq\frac32\norm{f_S^\infty}_X$, for $t>T$,
\[
 \begin{aligned}
 \bigl|\langle A_tf_S(t),f_S(t)\rangle_X-\norm{f_S^\infty}_X^2\bigr|
 &\leq C_0e^{-t\ell_*/2}\norm{f_S(t)}_X^2+\bigl(\norm{f_S(t)}_X+\norm{f_S^\infty}_X\bigr)
              \norm{f_S(t)-f_S^\infty}_X\\
 &\leq\frac{29}{4}C_0e^{-t\ell_*/2}\norm{f_S^\infty}_X^2.
 \end{aligned}
\]
Moreover, since $t|\langle(\partial_t A_t)f_S(t),f_S(t)\rangle_X|\leq tC_1e^{-t\ell_*/2}\norm{f_S(t)}_X^2$ by \eqref{ac-DtN-exponential}, we obtain
\[
 |\eta_S(t)|\leq\frac{29C_0+9tC_1}{4}e^{-t\ell_*/2}\norm{f_S^\infty}_X^2 .
\]
For the derivative, the identity in \eqref{ac-mass-quadratic-derivative}, together with \eqref{ac-DtN-exponential}, gives
\[
 \begin{aligned}
 |\eta_S'(t)|
 &\leq(2C_1+tC_2)e^{-t\ell_*/2}\norm{f_S(t)}_X^2+2\bigl[1+(C_0+tC_1)e^{-t\ell_*/2}\bigr]
        \norm{f_S(t)}_X\norm{f_S'(t)}_X\\
 &\leq\left[\frac94(2C_1+tC_2)+9C_1\bigl[1+(C_0+tC_1)e^{-t\ell_*/2}\bigr]\right]
           e^{-t\ell_*/2}\norm{f_S^\infty}_X^2 ,
 \end{aligned}
\]
where the last inequality uses \eqref{ac-all-mass-branch-derivative-bound} and  $\norm{f_S(t)}_X\leq\frac32\norm{f_S^\infty}_X$. Dividing these estimates by $\norm{f_S^\infty}_X^2$ and using $t>T>1$ yields
\[
 \begin{aligned}
 \frac{|\eta_S(t)|+|\eta_S'(t)|}{\norm{f_S^\infty}_X^2}
 &\leq C_{\mathrm{rem}}te^{-t\ell_*/2},\\
 \frac{\bigl|\langle A_tf_S(t),f_S(t)\rangle_X-\norm{f_S^\infty}_X^2\bigr|}
      {\norm{f_S^\infty}_X^2}
 &\leq C_{\mathrm{rem}}e^{-t\ell_*/2}
 \end{aligned}
 \qquad\text{for every }t>T ,
\]
where $C_{\mathrm{rem}}>0$ depends only on $\ell_*$. By \eqref{ac-vertex-equation}, $\sum_{\vv\in\V_0}\abs{f_S(t)(\vv)}^p=\langle A_tf_S(t),f_S(t)\rangle_X$. By \eqref{ac-limit-profile} and $\deg(\vv)=q$ for $\vv\in S$, $\norm{f_S^\infty}_X^2=q^{p/(p-2)}\#S$. Substituting the preceding estimates into the mass, mass derivative and energy calculations in the proof of Proposition~\ref{prop-ac-branches}, we obtain, for every nonempty finite $S\subset W$ and every $t>T$,
\begin{align}
 \mu_S(t)
 &=\frac{q^{p/(p-2)}\#S}{2}\,t^{(4-p)/(p-2)}
       \bigl(1+\varepsilon_{S,1}(t)\bigr),\label{ac-all-mass-mass-asymptotic}\\
 \mu_S'(t)
 &=-\frac{p-4}{2(p-2)}\,q^{p/(p-2)}\#S\,t^{(4-p)/(p-2)-1}
       \bigl(1+\varepsilon_{S,2}(t)\bigr),\label{ac-all-mass-slope-asymptotic}\\
 E(u_S(t),\G)
 &=\frac{p-4}{4p}\,q^{p/(p-2)}\#S\,t^{p/(p-2)}
       \bigl(1+\varepsilon_{S,3}(t)\bigr),\label{ac-all-mass-energy-asymptotic}
\end{align}
where the relative errors satisfy
\[
 \max_{1\leq j\leq3}|\varepsilon_{S,j}(t)|
 \leq C_*t^2e^{-t\ell_*/2}\qquad\text{for every }t>T ,
\]
where $C_*>0$ depends only on $p$ and $\ell_*$. Choose $\widehat T>\max\{T,4/\ell_*\}$ such that $C_*\widehat T^2 e^{-\widehat T\ell_*/2}\leq1/2$. Thus, for every $t\geq\widehat T$,
\[
 \mu_S'(t)<0,\qquad E(u_S(t),\G)>0,\qquad
 m(u_S(t))=\hat m_0(u_S(t))=\#S .
\]
The first two conclusions follow from \eqref{ac-all-mass-slope-asymptotic} and \eqref{ac-all-mass-energy-asymptotic}. For the indices, argue as in the proof of Proposition~\ref{prop-ac-branches}, using \eqref{ac-Hessian-splitting}, \eqref{ac-limit-form} and $\norm{L_S(t)-L_S^\infty}_{\mathcal L(X)}\leq1/2$ to obtain $m(u_S(t))=\#S$, and \eqref{ac-correction-formula} with $\mu_S'(t)<0$ to obtain $\hat m_0(u_S(t))=\#S$.

For each integer $n\geq1$, choose $S_n\subset W$ with $\#S_n=n$. For every $\mu>0$, choose an integer $n_\mu\geq1$ with
\[
 \frac14 q^{p/(p-2)}n_\mu\widehat T^{(4-p)/(p-2)}>\mu .
\]
Then \eqref{ac-all-mass-mass-asymptotic} gives $\mu_{S_n}(\widehat T)>\mu$ for every $n\geq n_\mu$ and $\mu_{S_n}(t)\to0$ as $t\to+\infty$.  Inverting the mass curve as in the proof of Theorem~\ref{th-ac-uniform-solutions}, we obtain a unique $t_n>\widehat T$ with $\mu_{S_n}(t_n)=\mu$. The solution $u_{\mu,n}:=u_{S_n}(t_n)$ has multiplier $\lambda_{\mu,n}:=t_n^2$ and both Morse indices equal to $n$. 

As in the proof of Theorem~\ref{cor-ac-infinite-orbits}, each $\gamma\in\Gamma$ acts isometrically and preserves the multiplier, with $\gamma(\V_0)=\V_0$. Thus \eqref{ac-Q-definition} gives $m(\gamma u_{\mu,n})=m(u_{\mu,n})=n$ for every $\gamma\in\Gamma$, so the solutions are pairwise $\Gamma$-inequivalent. Also, \eqref{ac-all-mass-mass-asymptotic} at $t=t_n$ gives
\[
 \lambda_{\mu,n}\geq
 \left(\frac{q^{p/(p-2)}n}{4\mu}\right)^{2(p-2)/(p-4)} .
\]
Hence $\lambda_{\mu,n}\to+\infty$ as $n\to+\infty$ for each fixed $\mu>0$. The same lower bound with $n\geq1$ gives $\inf_{n\geq n_\mu}\lambda_{\mu,n}\to+\infty$ as $\mu\to0^+$.
\end{proof}

\begin{remark}\label{rem-ac-all-mass-limitation}
The proof of Theorem~\ref{cor-ac-supercritical-all-mass} uses supports $S_n$ of $n$ vertices of a common degree $q$ and a frequency threshold $\widehat T$ independent of $n$. By \eqref{ac-all-mass-mass-asymptotic}, for every $t\geq\widehat T$,
\[
 \frac14q^{p/(p-2)}n\,t^{(4-p)/(p-2)}
 \leq\mu_{S_n}(t)
 \leq\frac34q^{p/(p-2)}n\,t^{(4-p)/(p-2)}.
\]
For $p>4$, \eqref{ac-all-mass-mass-asymptotic} gives
\[
 \mu_{S_n}(\widehat T)\longrightarrow+\infty\quad\text{as }n\to+\infty,
 \qquad
 \mu_{S_n}(t)\longrightarrow0\quad\text{as }t\to+\infty\text{ for each fixed }n.
\]
Thus, for every $\mu>0$ and every $n$ with $\mu_{S_n}(\widehat T)>\mu$, there exists $t_n>\widehat T$ such that $\mu_{S_n}(t_n)=\mu$.

For $2<p<4$, let $S\subset\V_0$ be nonempty and finite with $\#S\leq k$ and $\max_{\vv\in S}\deg(\vv)\leq D$. At the threshold $T_{k,D}$ in the proof of Theorem~\ref{th-ac-uniform-solutions}, \eqref{ac-mass-asymptotic} and \eqref{ac-S-epsilon-estimates} give
\[
 \mu_S(t)\geq\frac14t^{(4-p)/(p-2)}
        \sum_{\vv\in S}\deg(\vv)^{p/(p-2)}
 \geq\frac14\#S\qquad\text{for every }t\geq T_{k,D}>1.
\]
A prescribed mass $\mu$ on these branches therefore requires $\#S\leq4\mu$. Thus, any prescribed mass attained on these branches must satisfy $\mu\geq\#S/4\geq1/4$, even if the same lower bound were available for every $t\geq1$.
\end{remark}

\section[Proofs of Theorems~\ref{th-ac-critical} and \ref{cor-ac-critical-infinite-orbits}]{Proofs of Theorems~\ref{th-ac-critical} and \ref{cor-ac-critical-infinite-orbits}}
\label{sect-critical}

Throughout this section, $p=4$, $\G\in\mathbf G$ and $\varnothing\ne\V_0\subset\V$. We use the branches defined in \eqref{ac-physical-branch}. Proposition~\ref{prop-ac-critical-branches} gives general estimates for each branch and relates the constrained Morse index to the sign of $\mu_S'(t)$. Under local degree and edge-length bounds, Proposition~\ref{prop-ac-critical-local-sign} gives expansions of the mass, its derivative and the energy. Proposition~\ref{prop-ac-critical-sign-conditions} uses these expansions to give sufficient conditions for the direction of approach to the limiting mass. Inverting the increasing mass curves on a common interval gives Theorems~\ref{th-ac-critical} and \ref{cor-ac-critical-infinite-orbits}. Appendix~\ref{app-ac-critical-all-mass-tree} proves Proposition~\ref{prop-ac-critical-all-mass-tree} by constructing continuous branches on a suitable metric tree.

At $p=4$, the factor $t^{(4-p)/(p-2)}$ in the mass identity~\eqref{ac-mass-exact} equals one. Thus the limiting vertex profile determines a finite positive limiting mass.

For every nonempty finite $S\subset\V_0$, the critical mass is
\[
 \mu_S^\infty:=\frac12\norm{f_S^\infty}_X^2=\frac12\sum_{\vv\in S}\deg(\vv)\bigl|f_S^\infty(\vv)\bigr|^2=\frac12\sum_{\vv\in S}\deg(\vv)^{2/(p-2)+1}=\frac12\sum_{\vv\in S}\deg(\vv)^2.
\]
The following proposition gives quantitative mass and energy estimates for each fixed branch $u_S$ and determines its constrained Morse index from the sign of $\mu_S'(t)$.

\begin{proposition}\label{prop-ac-critical-branches}
Assume $p=4$. For every finite nonempty $S\subset\V_0$, there is $C_S>0$, depending only on $\ell_*$, $\#S$ and $\max_{\vv\in S}\deg(\vv)$, such that for every $t>T_S$,
\begin{equation}\label{ac-critical-main-estimates}
 \begin{aligned}
 \abs{\mu_S(t)-\mu_S^\infty}+\abs{\mu_S'(t)}
   &\leq C_St e^{-t\ell_*/2},\\
 \abs{E(u_S(t),\G)}&\leq C_St^3e^{-t\ell_*/2},
 \end{aligned}
\end{equation}
where $T_S$ is the threshold in Lemma~\ref{ac-IFT-branch}. In particular, $\mu_S(t)\to\mu_S^\infty$ and $E(u_S(t),\G)\to0$ as $t\to+\infty$. There is $\widehat T_S>T_S$, depending only on $\ell_*$, $\#S$ and $\max_{\vv\in S}\deg(\vv)$, such that, for every $t\geq\widehat T_S$,
\begin{equation}\label{ac-critical-main-index}
 m(u_S(t))=\#S,\qquad
 \hat m_0(u_S(t))=
 \begin{cases}
  \#S,&\mu_S'(t)<0,\\
  \#S-1,&\mu_S'(t)\geq0.
 \end{cases}
\end{equation}
\end{proposition}
\begin{proof}
Fix a nonempty finite $S\subset\V_0$ and put $p=4$. Lemma~\ref{ac-IFT-branch} and \eqref{ac-physical-branch} give $u_S(t)=t^{1/2}U_tf_S(t)$ for every $t>T_S$. By \eqref{ac-branch-weak-identity}, $(u_S(t), t^2)$ is a positive solution of \eqref{eqpd}. We use the constant $C$ of \eqref{ac-branch-estimates}, which depends only on $\ell_*$, $\#S$ and $\max_{\vv\in S}\deg(\vv)$. For every $t>T_S>1$,
\begin{equation}\label{ac-critical-fS-bounds}
 \norm{f_S(t)-f_S^\infty}_X+\norm{f_S'(t)}_X\leq Ce^{-t\ell_*/2}, \qquad  \norm{f_S(t)}_X\leq\norm{f_S^\infty}_X+C.
\end{equation}
Let $C_0,C_1,C_2$ be the constants of \eqref{ac-DtN-exponential}. We have
\begin{equation}\label{ac-critical-bounds}
 \begin{aligned}
 \norm{A_t-I_X}_{\mathcal L(X)}
 &\leq C_0 e^{-t\ell_*/2},\\
\norm{\partial_tA_t}_{\mathcal L(X)}&\leq C_1e^{-t\ell_*/2},\\
\norm{\partial_t^2A_t}_{\mathcal L(X)}&\leq C_2e^{-t\ell_*/2},\\
 \bigl|\norm{f_S(t)}_X^2-\norm{f_S^\infty}_X^2\bigr|
 &=\bigl|\norm{f_S(t)}_X-\norm{f_S^\infty}_X\bigr|\bigl(\norm{f_S(t)}_X+\norm{f_S^\infty}_X\bigr)\\
 &\leq\norm{f_S(t)-f_S^\infty}_X\bigl(2\norm{f_S^\infty}_X+C\bigr)\\
 &\leq C\bigl(2\norm{f_S^\infty}_X+C\bigr)e^{-t\ell_*/2}.
 \end{aligned}
\end{equation}
Then \eqref{ac-critical-bounds} and \eqref{ac-mass-fixed-restriction} give a fixed constant $\widehat C_S=\widehat C_S(\ell_*,\#S,\max_{\vv\in S}\deg(\vv))>0$ such that
\begin{equation}\label{ac-critical-mass-expansion}
 \begin{aligned}
 \mu_S(t)
 &=\frac12\bigl\langle(A_t+t\partial_tA_t)f_S(t),f_S(t)\bigr\rangle_X,\\
 \mu_S'(t)
 &=\frac12\bigl\langle(2\partial_tA_t+t\partial_t^2A_t)f_S(t),f_S(t)\bigr\rangle_X
   +\bigl\langle(A_t+t\partial_tA_t)f_S(t),f_S'(t)\bigr\rangle_X,\\
 |\mu_S(t)-\mu_S^\infty|
 &\leq\tfrac12\norm{A_t+t\partial_tA_t-I_X}_{\mathcal L(X)}\norm{f_S(t)}_X^2
   +\tfrac12\bigl|\norm{f_S(t)}_X^2-\norm{f_S^\infty}_X^2\bigr|,\\
 |\mu_S'(t)|
 &\leq\tfrac12\norm{2\partial_tA_t+t\partial_t^2A_t}_{\mathcal L(X)}\norm{f_S(t)}_X^2
   +\norm{A_t+t\partial_tA_t}_{\mathcal L(X)}\norm{f_S(t)}_X\norm{f_S'(t)}_X,\\
 |\mu_S(t)-\mu_S^\infty|+|\mu_S'(t)|
 &\leq \widehat C_St e^{-t\ell_*/2}.
 \end{aligned}
\end{equation}
Testing \eqref{ac-vertex-equation} with $f_S(t)$ gives
\[
 \sum_{\vv\in\V_0}\abs{f_S(t)(\vv)}^4
 =\langle A_tf_S(t),f_S(t)\rangle_X.
\]
Since $u_S(t)=t^{1/2}U_tf_S(t)$ and $\mathcal RU_t=I_X$, we have $u_S(t)(\vv)=t^{1/2}f_S(t)(\vv)$ for every $\vv\in\V$. Moreover, testing \eqref{ac-branch-weak-identity} with $u_S(t)$ then yields $\norm{(u_S(t))'}_2^2+t^2\mu_S(t)=t^2\langle A_tf_S(t),f_S(t)\rangle_X$. Consequently, \eqref{ac-critical-fS-bounds}, \eqref{ac-critical-bounds} and the first line of \eqref{ac-critical-mass-expansion} give
\begin{equation}\label{ac-critical-energy-exact}
 \begin{aligned}
 E(u_S(t),\G)
 &=\frac{t^2}{4}\langle A_tf_S(t),f_S(t)\rangle_X-\frac{t^2}{2}\mu_S(t)\\
 &=-\frac{t^3}{4}\langle\partial_tA_tf_S(t),f_S(t)\rangle_X,\\
 \abs{E(u_S(t),\G)}
 &=\frac{t^3}{4}\bigl|\langle\partial_tA_tf_S(t),f_S(t)\rangle_X\bigr|
   \leq\frac{t^3}{4}\norm{\partial_tA_t}_{\mathcal L(X)}\norm{f_S(t)}_X^2\\
 &\leq\frac{C_1}{4}\bigl(\norm{f_S^\infty}_X+C\bigr)^2t^3e^{-t\ell_*/2}.
 \end{aligned}
\end{equation}
In particular, $\mu_S(t)\to\mu_S^\infty$ and $E(u_S(t),\G)\to0$ as $t\to+\infty$.

We next compute the Morse indices. Recall the bounded self-adjoint operator $L_S(t)=A_t-\mathcal N'(f_S(t)):X\to X$ and the limit $L_S^\infty=-2P_S+(I_X-P_S)$. By \eqref{ac-N-derivative}, \eqref{ac-norm-ordering} and \eqref{ac-critical-fS-bounds},
\[
 \begin{aligned}
 \norm{\mathcal N'(f_S(t))-\mathcal N'(f_S^\infty)}_{\mathcal L(X)}
 &\leq3\Bigl(\sum_{\vv\in\V_0}\frac{\bigl(\abs{f_S(t)(\vv)}^2-\abs{f_S^\infty(\vv)}^2\bigr)^2}{\deg(\vv)}\Bigr)^{1/2}\\
 &\leq3\norm{f_S(t)-f_S^\infty}_X\sup_{\vv\in\V_0}\bigl(2\abs{f_S^\infty(\vv)}+C\bigr)\\
 &\leq3\bigl(2\norm{f_S^\infty}_X+C\bigr)\norm{f_S(t)-f_S^\infty}_X,
 \end{aligned}
\]
and hence, by \eqref{ac-critical-bounds},
\[
 \begin{aligned}
 \norm{L_S(t)-L_S^\infty}_{\mathcal L(X)}
 &=\norm{(A_t-I_X)-\bigl(\mathcal N'(f_S(t))-\mathcal N'(f_S^\infty)\bigr)}_{\mathcal L(X)}\\
 &\leq\norm{A_t-I_X}_{\mathcal L(X)}+\norm{\mathcal N'(f_S(t))-\mathcal N'(f_S^\infty)}_{\mathcal L(X)}\\
 &\leq C_0e^{-t\ell_*/2}+3\bigl(2\norm{f_S^\infty}_X+C\bigr)\norm{f_S(t)-f_S^\infty}_X\\
 &\leq\bigl[C_0+3C\bigl(2\norm{f_S^\infty}_X+C\bigr)\bigr]e^{-t\ell_*/2}.
 \end{aligned}
\]
Choose $\widehat T_S>\max\{T_S,1\}$ such that
\[
 \bigl[C_0+3C\bigl(2\norm{f_S^\infty}_X+C\bigr)\bigr]
 e^{-\widehat T_S\ell_*/2}<\frac12.
\]
This choice depends only on $\ell_*$, $\#S$ and $\max_{\vv\in S}\deg(\vv)$. For every $t\geq\widehat T_S$, we have $\norm{L_S(t)-L_S^\infty}_{\mathcal L(X)}<1/2$. Arguing as in the proof of Proposition~\ref{prop-ac-branches}, using \eqref{ac-Hessian-splitting} and $L_S^\infty=-2P_S+(I_X-P_S)$, we obtain
\[
 m(u_S(t))=\#S.
\]
Differentiating \eqref{ac-branch-weak-identity} gives \eqref{ac-Q-branch-derivative} and \eqref{ac-positive-branch-direction} also at $p=4$. At each $t\geq\widehat T_S$ with $\mu_S'(t)\ne0$, arguing as in the proof of Proposition~\ref{prop-ac-branches}, using the tangent-space decomposition \eqref{ac-exact-tangent-correction}--\eqref{ac-correction-formula}, we obtain
\[
 \hat m_0(u_S(t))=
 \begin{cases}
 \#S,&\mu_S'(t)<0,\\
 \#S-1,&\mu_S'(t)>0.
 \end{cases}
\]
It remains to consider $t\geq\widehat T_S$ with $\mu_S'(t)=0$. Then $\mu_S'(t)=2(u_S(t),u_S'(t))_2$, and \eqref{ac-Q-branch-derivative} at $h=u_S(t)$, together with \eqref{ac-positive-branch-direction}, gives
\begin{equation}\label{ac-zero-slope-identities}
 \begin{aligned}
 (u_S(t),u_S'(t))_2&=0,\qquad Q_{u_S(t)}[u_S'(t),u_S'(t)]=0,\\
 Q_{u_S(t)}[u_S'(t),u_S(t)]&=-2t\mu_S(t)\ne0.
 \end{aligned}
\end{equation}
In particular, $u_S'(t)$ and $u_S(t)$ are linearly independent. Put $Z:=\operatorname{span}\{u_S'(t),u_S(t)\}$.
Then, for every $a\in\R$, \eqref{ac-zero-slope-identities} gives
\[
 \begin{aligned}
 Q_{u_S(t)}[u_S(t)+a u_S'(t),u_S(t)+a u_S'(t)]=Q_{u_S(t)}[u_S(t),u_S(t)]-4at\mu_S(t).
 \end{aligned}
\]
Since $t\mu_S(t)>0$, the right-hand side takes both positive and negative values as $a$ varies in $\R$. Thus $Q_{u_S(t)}|_Z$ has one negative and one positive direction, because $\dim Z=2$.

Define
\[
 Z^{\perp}:=\{h\in H^1(\G):Q_{u_S(t)}[h,z]=0\text{ for every }z\in Z\}.
\]
For every $h\in H^1(\G)$, the equations
\[
 \begin{aligned}
 -2t\mu_S(t)b&=Q_{u_S(t)}[h,u_S'(t)],\\
 -2t\mu_S(t)a+bQ_{u_S(t)}[u_S(t),u_S(t)]
 &=Q_{u_S(t)}[h,u_S(t)]
 \end{aligned}
\]
uniquely determine $b \in \R$ and then $a\in \R$, since $t\mu_S(t)>0$. For such $a$ and $b$, $h_0:=h-au_S'(t)-bu_S(t)$ belongs to $Z^\perp$. The uniqueness of $a,b$ therefore gives
\begin{equation}\label{ac-space-splitting}
 H^1(\G)=Z\oplus Z^{\perp},
\end{equation}
where $\oplus$ denotes the algebraic direct sum and the decomposition is orthogonal for $Q_{u_S(t)}$. Also, $Q_{u_S(t)}[u_S'(t),h]=-2t(u_S(t),h)_2$ implies $Z^{\perp}\subset T_{u_S(t)}H^1_{\mu_S(t)}(\G)$. If a tangent vector $h \in T_{u_S(t)}H^1_{\mu_S(t)}(\G)$ is written as $h=au_S'(t)+bu_S(t)+h_0$ with $h_0\in Z^{\perp}$, then $0=(u_S(t),h)_2=b\mu_S(t)$, hence $b=0$. Therefore
\begin{equation}\label{ac-tangent-splitting}
 T_{u_S(t)}H^1_{\mu_S(t)}(\G)=\operatorname{span}\{u_S'(t)\}\oplus Z^{\perp},
\end{equation}
where the decomposition is again orthogonal for $Q_{u_S(t)}$. Then \eqref{ac-space-splitting} gives $m(Q_{u_S(t)}|_{Z^{\perp}})=m(u_S(t))-1=\#S-1$, and then \eqref{ac-tangent-splitting} together with the second identity in \eqref{ac-zero-slope-identities} gives $\hat m_0(u_S(t))=\#S-1$, where $m(Q_{u_S(t)}|_{Z^{\perp}})$ is the maximal dimension of a subspace of $Z^{\perp}$ on which $Q_{u_S(t)}$ is negative definite.  This proves \eqref{ac-critical-main-index}.

Set $C_S:=\max\{\widehat C_S,\frac{C_1}{4}(\norm{f_S^\infty}_X+C)^2\}$. By \eqref{ac-critical-mass-expansion} and \eqref{ac-critical-energy-exact}, this constant gives \eqref{ac-critical-main-estimates} and depends only on $\ell_*$, $\#S$ and $\max_{\vv\in S}\deg(\vv)$.
\end{proof}

\begin{remark}\label{rem-ac-critical-fixed-multiplier}
Let $p=4$ and let $W\subset\V_0$ be nonempty and finite. There is $T_W>\max_{\varnothing\ne S\subset W}\widehat T_S$, depending only on $\ell_*$, $\#W$ and $\max_{\vv\in W}\deg(\vv)$, such that, for every $t>T_W$, the edge and vertex equations in \eqref{eqpd} have $2^{\#W}-1$ pairwise distinct positive solutions $u_S(t)$, indexed by $\varnothing\ne S\subset W$, all with multiplier $t^2$. Indeed, apply Proposition~\ref{prop-ac-critical-branches} to these sets $S$ and choose $T_W$ so large that \eqref{ac-branch-estimates} gives $\norm{f_S(t)-f_S^\infty}_X<1/3$ for every $t>T_W$ and every $\varnothing\ne S\subset W$. For distinct nonempty $S_1,S_2\subset W$, we have $\norm{f_{S_1}^\infty-f_{S_2}^\infty}_X\geq1$, and hence
\[
 \norm{f_{S_1}(t)-f_{S_2}(t)}_X
 \geq\norm{f_{S_1}^\infty-f_{S_2}^\infty}_X
       -\norm{f_{S_1}(t)-f_{S_1}^\infty}_X
       -\norm{f_{S_2}(t)-f_{S_2}^\infty}_X>\frac13.
\]
Since $t^{-1/2}\mathcal Ru_S(t)=f_S(t)$ by \eqref{ac-physical-branch} and $\mathcal RU_t=I_X$, the functions $u_S(t)$ are pairwise distinct.

By \eqref{ac-critical-main-estimates}, these solutions satisfy
\[
 \norm{u_S(t)}_2^2=\mu_S(t)\longrightarrow\mu_S^\infty
 =\frac12\sum_{\vv\in S}\deg(\vv)^2
 \qquad\text{as }t\to+\infty.
\]
For sets $S_1,S_2$ with $\mu_{S_1}^\infty\ne\mu_{S_2}^\infty$, we have $\mu_{S_1}(t)-\mu_{S_2}(t)\to\mu_{S_1}^\infty-\mu_{S_2}^\infty\ne0$. Thus the solutions with the same multiplier need not have the same mass. Moreover, these estimates alone do not determine the sign of $\mu_S'(t)$ or give a common prescribed-mass interval for the branches. Consequently, on a general graph $\G\in\mathbf G$, these estimates do not yield a prescribed-mass multiplicity result analogous to Theorem~\ref{th2}.
\end{remark}
 We use the neighbor set $N(W)$ defined in the Introduction. The space $X=\ell^2(\V,\deg)$ and its norm are those of Section~\ref{sect5}; in particular, $\norm f_\infty\leq\norm f_X$ by \eqref{ac-norm-ordering}.

For every nonempty finite $S\subset\V_0$ and $t>0$, set
\begin{equation}\label{ac-critical-weighted-size}
 \sigma_S(t)^2:=\sum_{\vv\in S}\sum_{\e\in\E_b,\,\e\succ\vv}e^{-2t\ell_\e}
\leq e^{-2t\ell_*}\sum_{\vv\in S}\deg(\vv).
\end{equation}
Only bounded edges enter $\sigma_S(t)$, while $\deg(\vv)$ still counts all edges incident at $\vv$.
The set of internal edges is
\[
 \E(S):=\{\e\in\E_b:\vv_\e,\w_\e\in S\}.
\]
For $\vv\in\V$, denote by
\[
 \E_S(\vv):=\{\e\in\E_b:\e\succ\vv,\ \w_{\e,\vv}\in S\}
\]
the set of edges joining $\vv$ to a vertex of $S$. For every $\vv\in\V$ and every $(\e,\e')\in\E_S(\vv)\times\E_S(\vv)$, define
\[
 \omega_\vv(\e,\e')
 :=\frac{\sqrt{\deg(\w_{\e,\vv})\deg(\w_{\e',\vv})}}{\deg(\vv)}.
\]

The next proposition expands $\mu_S(t)$, $\mu_S'(t)$ and $E(u_S(t),\G)$, with remainder estimates expressed in terms of $\sigma_S(t)$.
\begin{proposition}\label{prop-ac-critical-local-sign}
Let $\G\in\mathbf G$, let $p=4$, and let $\varnothing\ne\V_0\subset\V$. Fix integers $q\geq1$, $k\geq1$ and a number $\ell_+\geq\ell_*$. There are $T_0>\max\{1,3/\ell_*\}$ and $C_f,C_r>0$, depending only on $q$, $k$, $\ell_*$ and $\ell_+$, with the following property. Let $S\subset\V_0$ consist of $k$ vertices, with $\deg(\vv)\leq q$ for every $\vv\in S$, and suppose every bounded edge incident at $S$ has length at most $\ell_+$. Then the branch $u_S$ in \eqref{ac-physical-branch} is defined for $t\geq T_0$, and its vertex profile satisfies $\norm{f_S(t)-f_S^\infty}_X\leq C_f\sigma_S(t)$ for every $t\geq T_0$.
For $t\geq T_0$, the mass and its derivative have the expansions
\begin{equation}\label{ac-critical-general-mass-expansion}
 \begin{aligned}
 \mu_S(t)={}&\mu_S^\infty
       +2\sum_{\e\in\E(S)}\sqrt{\deg(\vv_\e)\deg(\w_\e)} (t\ell_\e-2)e^{-t\ell_\e}\\
 &-2\sum_{\vv\in S}\sum_{\e\in\E_b,\,\e\succ\vv}\deg(\vv) (t\ell_\e-1)e^{-2t\ell_\e}\\
 &-\sum_{\vv\in S}\sum_{\e,\e'\in\E_S(\vv)}\omega_\vv(\e,\e')
       [t(\ell_\e+\ell_{\e'})-2]e^{-t(\ell_\e+\ell_{\e'})}\\
 &+2\sum_{\vv\in N(S)}
       \sum_{\e,\e'\in\E_S(\vv)}\omega_\vv(\e,\e')
       [t(\ell_\e+\ell_{\e'})-2]e^{-t(\ell_\e+\ell_{\e'})}
       +\mathcal E_S(t),
 \end{aligned}
\end{equation}
\begin{equation}\label{ac-critical-general-mass-slope}
 \begin{aligned}
 \mu_S'(t)={}&2\sum_{\e\in\E(S)}\sqrt{\deg(\vv_\e)\deg(\w_\e)}
          \ell_\e(3-t\ell_\e)e^{-t\ell_\e}\\
 &+2\sum_{\vv\in S}\sum_{\e\in\E_b,\,\e\succ\vv}\deg(\vv)
          \ell_\e(2t\ell_\e-3)e^{-2t\ell_\e}\\
 &+\sum_{\vv\in S}\sum_{\e,\e'\in\E_S(\vv)}\omega_\vv(\e,\e')
          (\ell_\e+\ell_{\e'})
          [t(\ell_\e+\ell_{\e'})-3]e^{-t(\ell_\e+\ell_{\e'})}\\
 &-2\sum_{\vv\in N(S)}
          \sum_{\e,\e'\in\E_S(\vv)}\omega_\vv(\e,\e')(\ell_\e+\ell_{\e'})[t(\ell_\e+\ell_{\e'})-3]e^{-t(\ell_\e+\ell_{\e'})}
          +\mathcal E_S'(t),
 \end{aligned}
\end{equation}
where $\mathcal E_S\in C^1([T_0,+\infty),\R)$  defined in \eqref{ac-critical-general-mass-remainder} satisfies
\begin{equation}\label{ac-critical-general-mass-error}
 |\mathcal E_S(t)|+|\mathcal E_S'(t)|
 \leq C_rt e^{-t\ell_*}\sigma_S(t)^2
 \qquad\text{for every }t\geq T_0.
\end{equation}
Moreover, for $t\geq T_0$, the energy satisfies the expansion
\begin{equation}\label{ac-critical-general-energy-expansion}
 \begin{aligned}
 E(u_S(t),\G)
 &=-t^3\sum_{\e\in\E(S)}\sqrt{\deg(\vv_\e)\deg(\w_\e)}\,\ell_\e e^{-t\ell_\e}
   +t^3\sum_{\vv\in S}\sum_{\e\in\E_b,\,\e\succ\vv}\deg(\vv)\ell_\e e^{-2t\ell_\e}\\
 &\quad+\frac{t^3}{2}\sum_{\vv\in S}\sum_{\e,\e'\in\E_S(\vv)}\omega_\vv(\e,\e')
                 (\ell_\e+\ell_{\e'})e^{-t(\ell_\e+\ell_{\e'})}\\
 &\quad-t^3\sum_{\vv\in N(S)}\sum_{\e,\e'\in\E_S(\vv)}\omega_\vv(\e,\e')
                 (\ell_\e+\ell_{\e'})e^{-t(\ell_\e+\ell_{\e'})}
       -\frac{t^3}{4}\sum_{j=0}^2R_{S,j}(t),
 \end{aligned}
\end{equation}
where $R_{S,j}\in C^1([T_0,+\infty),\R)$ defined in \eqref{ac-critical-grid-scalar-remainder-definitions} satisfy $\frac{t^3}{4}\left|\sum_{j=0}^2R_{S,j}(t)\right|\leq C_rt^3e^{-t\ell_*}\sigma_S(t)^2$.
\end{proposition}

\begin{proof}
Fix $S\subset\V_0$ consisting of $k$ vertices, with $\deg(\vv)\leq q$ for every $\vv\in S$, and suppose every bounded edge incident at $S$ has length at most $\ell_+$.  Apply Lemma~\ref{ac-IFT-branch} to each possible value $1,\ldots,q$ of $\max_{\vv\in S}\deg(\vv)$. Choose $T_b=T_b(q,k,\ell_*)>1$ strictly above the largest branch threshold and $C_b=C_b(q,k,\ell_*)>0$ equal to the largest estimate constant. Then $f_S$ is defined for $t\geq T_b$ and
\begin{equation}\label{ac-critical-initial-branch-bound}
 \norm{f_S(t)-f_S^\infty}_X+\norm{f_S'(t)}_X
 \leq C_b e^{-t\ell_*/2}.
\end{equation}

For $t>0$, define bounded self-adjoint operators $B_t,D_t:X\to X$ by
\[
   (B_tf)(\vv):=\frac1{\deg(\vv)}\sum_{\substack{\e\in\E_b\\\e\succ\vv}}e^{-t\ell_\e}f(\w_{\e,\vv}),
 \qquad
 (D_tf)(\vv):=\frac1{\deg(\vv)}
       \left[\sum_{\substack{\e\in\E_b\\\e\succ\vv}}e^{-2t\ell_\e}\right]f(\vv).
\]
Formula \eqref{ac-DtN-coordinate}, including its half-line contributions, gives
\begin{equation}\label{ac-critical-grid-operator}
 A_t=I_X-2B_t+2D_t+R_A(t),
\end{equation}
where $R_A(t):X\to X$ is the bounded self-adjoint operator

\begin{equation}\label{ac-critical-grid-DtN-expansion}
 \begin{aligned}
 (R_A(t)f)(\vv)={}&\frac2{\deg(\vv)}
 \sum_{\substack{\e\in\E_b\\\e\succ\vv}}
 \frac{e^{-4t\ell_\e}}{1-e^{-2t\ell_\e}}f(\vv) -\frac2{\deg(\vv)}\sum_{\substack{\e\in\E_b\\\e\succ\vv}}
  \frac{e^{-3t\ell_\e}}{1-e^{-2t\ell_\e}}f(\w_{\e,\vv}).
 \end{aligned}
\end{equation}

Set $\varepsilon(t):=e^{-t\ell_*}$. There is a single constant $C_A=C_A(q,\ell_*,\ell_+)\geq1$ such that, for $j=0,1,2$ and $t\geq\max\{1,3/\ell_*\}$, the operators satisfy
\begin{equation}\label{ac-critical-grid-operator-bounds}
 \begin{aligned}
 \norm{\partial_t^jB_t}_{\mathcal L(X)}&\leq C_A\varepsilon,&
 \norm{\partial_t^jD_t}_{\mathcal L(X)}&\leq C_A\varepsilon^2,\\
 \norm{\partial_t^j(A_t-I_X)}_{\mathcal L(X)}&\leq C_A\varepsilon,&
 \norm{\partial_t^jR_A(t)}_{\mathcal L(X)}&\leq C_A\varepsilon^3,
 \end{aligned}
\end{equation}
and their action on the limiting profile satisfies
\begin{equation}\label{ac-critical-local-operator-bounds}
 \begin{aligned}
 \norm{\partial_t^jB_tf_S^\infty}_X
 +\norm{\partial_t^j(A_t-I_X)f_S^\infty}_X&\leq C_A\sigma_S,\\
 \norm{\partial_t^jD_tf_S^\infty}_X&\leq C_A\varepsilon\sigma_S,\\
 \norm{\partial_t^jR_A(t)f_S^\infty}_X&\leq C_A\varepsilon^2\sigma_S,\\
 |\langle\partial_t^jR_A(t)f_S^\infty,f_S^\infty\rangle_X|
 &\leq C_A\varepsilon\sigma_S^2.
 \end{aligned}
\end{equation}

We first prove \eqref{ac-critical-grid-operator-bounds}. For any bounded family of real coefficients $(a_\e)_{\e\in\E_b}$, the Cauchy--Schwarz inequality and $\#\{\e\in\E_b:\e\succ\vv\}\leq\deg(\vv)$ give
\begin{equation}\label{ac-critical-edge-coefficient-bound}
 \begin{aligned}
 \sum_{\vv\in\V}\deg(\vv)
 \left|\frac1{\deg(\vv)}\sum_{\substack{\e\in\E_b\\\e\succ\vv}}a_\e f(\w_{\e,\vv})\right|^2
 &\leq \sum_{\vv\in\V}\sum_{\substack{\e\in\E_b\\\e\succ\vv}}
              |a_\e|^2|f(\w_{\e,\vv})|^2\\
 &=\sum_{\w\in\V}\sum_{\substack{\e\in\E_b\\\e\succ\w}}
             |a_\e|^2|f(\w)|^2
 \leq\sup_{\e\in\E_b}|a_\e|^2\norm f_X^2 .
 \end{aligned}
\end{equation}
Taking $a_\e=e^{-t\ell_\e}$ in \eqref{ac-critical-edge-coefficient-bound} gives $\norm{B_tf}_X\leq\varepsilon(t)\norm f_X$. The same estimate with $f(\vv)$ in place of $f(\w_{\e,\vv})$ and $a_\e=e^{-2t\ell_\e}$ gives $\norm{D_tf}_X\leq\varepsilon(t)^2\norm f_X$.

For $j=0,1,2$ and $t\geq3/\ell_*$, the function $\ell\mapsto\ell^je^{-t\ell}$ decreases on $[\ell_*,+\infty)$, and hence
\begin{equation}\label{ac-critical-sup-bound}
 \sup_{\ell\geq\ell_*}\ell^je^{-t\ell}\leq\ell_*^j e^{-t\ell_*}.
\end{equation}
To estimate $\partial_t^jB_t$, take $a_\e=(-\ell_\e)^je^{-t\ell_\e}$ in \eqref{ac-critical-edge-coefficient-bound}. For $\partial_t^jD_t$, take $a_\e=(-2\ell_\e)^je^{-2t\ell_\e}$ and use $f(\vv)$ in place of $f(\w_{\e,\vv})$. By \eqref{ac-critical-sup-bound} and $e^{-t\ell_\e}\leq\varepsilon(t)$, the absolute values of these coefficients satisfy
\[
 \begin{aligned}
 \ell_\e^je^{-t\ell_\e}
 &\leq\ell_*^j\varepsilon(t),\\
 (2\ell_\e)^je^{-2t\ell_\e}
 &=2^j\bigl(\ell_\e^je^{-t\ell_\e}\bigr)e^{-t\ell_\e}
 \leq(2\ell_*)^j\varepsilon(t)^2.
 \end{aligned}
\]
For the two coefficients in \eqref{ac-critical-grid-DtN-expansion}, put
\[
 c_m(s):=\frac{e^{-ms}}{1-e^{-2s}},\qquad m=3,4.
\]
Direct differentiation gives
\[
 \begin{aligned}
 c_m'(s)
 &=-e^{-ms}\left[\frac{m}{1-e^{-2s}}+\frac{2e^{-2s}}{(1-e^{-2s})^2}\right],\\
 c_m''(s)
 &=e^{-ms}\left[\frac{m^2}{1-e^{-2s}}
 +\frac{4(m+1)e^{-2s}}{(1-e^{-2s})^2}
 +\frac{8e^{-4s}}{(1-e^{-2s})^3}\right].
 \end{aligned}
\]
Since $\partial_t^j c_m(t\ell)=\ell^jc_m^{(j)}(t\ell)$ and $e^{-2t\ell}\leq e^{-6}$ for $t\ell\geq3$, we obtain
\begin{equation}\label{ac-critical-remainder-coefficients}
 \begin{aligned}
 \left|\partial_t^j c_m(t\ell)\right|
 &\leq M\ell^je^{-mt\ell},\qquad j=0,1,2,\\
 M&:=\frac{16}{1-e^{-6}}+\frac{20e^{-6}}{(1-e^{-6})^2}+\frac{8e^{-12}}{(1-e^{-6})^3}\geq16.
 \end{aligned}
\end{equation}
Since the sum at each vertex in \eqref{ac-critical-grid-DtN-expansion} is finite, the coordinate formula for $\partial_t^jR_A(t)$ is obtained by differentiating its coefficients.

By \eqref{ac-critical-sup-bound}, $\sup_{\ell\geq\ell_*}\ell^je^{-mt\ell}\leq\ell_*^j\varepsilon(t)^m$ for $m\geq1$. Apply \eqref{ac-critical-edge-coefficient-bound} to the two sums in \eqref{ac-critical-grid-DtN-expansion}, with $a_\e=2\partial_t^jc_4(t\ell_\e)$ and $f(\vv)$ in the first sum and $a_\e=-2\partial_t^jc_3(t\ell_\e)$ and $f(\w_{\e,\vv})$ in the second. Using \eqref{ac-critical-remainder-coefficients} and $(a+b)^2\leq2(a^2+b^2)$ gives
\[
 \norm{\partial_t^j R_A(t)}_{\mathcal L(X)}^2
 \leq8M^2\ell_*^{2j}\bigl[\varepsilon(t)^8+\varepsilon(t)^6\bigr]
 \leq16M^2\ell_*^{2j}\varepsilon(t)^6.
\]
The bound for $\partial_t^j(A_t-I_X)$ follows from \eqref{ac-critical-grid-operator}.

We next prove \eqref{ac-critical-local-operator-bounds}. By \eqref{ac-limit-profile}, $f_S^\infty(\vv)=\sqrt{\deg(\vv)}\mathbf1_S(\vv)$. Using $\norm{f_S^\infty}_\infty\leq\sqrt q$ and $\#\E_S(\w)\leq\deg(\w)$ in the Cauchy--Schwarz step of \eqref{ac-critical-edge-coefficient-bound} gives
\begin{equation}\label{ac-critical-local-B-bound}
 \begin{aligned}
 \norm{\partial_t^jB_tf_S^\infty}_X^2
 &=\sum_{\w\in\V}\frac1{\deg(\w)}
       \left[\sum_{\e\in\E_S(\w)}\ell_\e^je^{-t\ell_\e}f_S^\infty(\w_{\e,\w})\right]^2\\
 &\leq q\ell_+^{2j}\sum_{\w\in\V}\sum_{\e\in\E_S(\w)}e^{-2t\ell_\e}
 =q\ell_+^{2j}\sigma_S(t)^2 .
 \end{aligned}
\end{equation}
For $D_t$, using $|f_S^\infty(\vv)|^2=\deg(\vv)$ on $S$ and applying the same step at each $\vv\in S$ gives
\[
 \begin{aligned}
 \norm{\partial_t^jD_tf_S^\infty}_X^2
 &=\sum_{\vv\in S}\left[\sum_{\e\in\E_b,\,\e\succ\vv}(2\ell_\e)^je^{-2t\ell_\e}\right]^2\\
 &\leq q(2\ell_+)^{2j}\sum_{\vv\in S}\sum_{\e\in\E_b,\,\e\succ\vv}e^{-4t\ell_\e}
 \leq q(2\ell_+)^{2j}\varepsilon(t)^2\sigma_S(t)^2 .
 \end{aligned}
\]
In the coordinate formula for $\partial_t^jR_A(t)f_S^\infty$, a term is nonzero only if $\vv\in S$ or $\w_{\e,\vv}\in S$. Thus $\ell_\e\leq\ell_+$ for each nonzero term. Using \eqref{ac-critical-remainder-coefficients}, $e^{-4t\ell_\e}\leq\varepsilon(t)^2e^{-2t\ell_\e}$ and $e^{-3t\ell_\e}\leq\varepsilon(t)^2e^{-t\ell_\e}$, we obtain
\[
 \begin{aligned}
 \bigl|\bigl(\partial_t^jR_A(t)f_S^\infty\bigr)(\vv)\bigr|
 &\leq\frac{2M\ell_+^j}{\deg(\vv)}
 \sum_{\substack{\e\in\E_b\\\e\succ\vv}}
 \bigl[e^{-4t\ell_\e}f_S^\infty(\vv)
       +e^{-3t\ell_\e}f_S^\infty(\w_{\e,\vv})\bigr]\\
 &\leq 2M\ell_+^j\varepsilon(t)^2
 \bigl[(D_tf_S^\infty)(\vv)+(B_tf_S^\infty)(\vv)\bigr].
 \end{aligned}
\]
Taking the $X$ norm and using \eqref{ac-critical-local-B-bound} and the subsequent estimate for $D_t$, both at $j=0$, gives
\[
 \begin{aligned}
 \norm{\partial_t^jR_A(t)f_S^\infty}_X
 &\leq 2M\ell_+^j\varepsilon(t)^2
       \bigl[\norm{D_tf_S^\infty}_X+\norm{B_tf_S^\infty}_X\bigr]\\
 &\leq 2M\sqrt q\,\ell_+^j\varepsilon(t)^2
        \bigl[1+\varepsilon(t)\bigr]\sigma_S(t)
 \leq 4M\sqrt q\,\ell_+^j\varepsilon(t)^2\sigma_S(t).
 \end{aligned}
\]

Taking the scalar product of the differentiated coordinate formula \eqref{ac-critical-grid-DtN-expansion} with $f_S^\infty$ gives
\[
 \begin{aligned}
 \langle\partial_t^jR_A(t)f_S^\infty,f_S^\infty\rangle_X
 &=2\sum_{\vv\in S}\sum_{\e\in\E_b,\,\e\succ\vv}\deg(\vv)\partial_t^jc_4(t\ell_\e)\\
 &\quad-2\sum_{\vv\in S}\sum_{\e\in\E_S(\vv)}
       f_S^\infty(\vv)f_S^\infty(\w_{\e,\vv})\partial_t^jc_3(t\ell_\e)\\
 &=2\sum_{\vv\in S}\sum_{\e\in\E_b,\,\e\succ\vv}\deg(\vv)\partial_t^jc_4(t\ell_\e)
   -4\sum_{\e\in\E(S)}\sqrt{\deg(\vv_\e)\deg(\w_\e)}\,\partial_t^jc_3(t\ell_\e).
 \end{aligned}
\]
Using \eqref{ac-critical-remainder-coefficients}, $\deg(\vv)\leq q$ on $S$, $\sqrt{\deg(\vv_\e)\deg(\w_\e)}\leq q$ on $\E(S)$, and $\displaystyle 2\sum_{\e\in\E(S)}e^{-2t\ell_\e}=\sum_{\vv\in S}\sum_{\e\in\E(S),\,\e\succ\vv}e^{-2t\ell_\e}\leq\sum_{\vv\in S}\sum_{\e\in\E_b,\,\e\succ\vv}e^{-2t\ell_\e}=\sigma_S(t)^2$, we find
\[
 \begin{aligned}
 \bigl|\langle\partial_t^jR_A(t)f_S^\infty,f_S^\infty\rangle_X\bigr|
 &\leq2qM\ell_+^j\left[
   \varepsilon(t)^2\sum_{\vv\in S}\sum_{\e\in\E_b,\,\e\succ\vv}e^{-2t\ell_\e}
   +2\varepsilon(t)\sum_{\e\in\E(S)}e^{-2t\ell_\e}
   \right]\\
 &\leq2qM\ell_+^j\bigl[\varepsilon(t)^2+\varepsilon(t)\bigr]\sigma_S(t)^2.
 \end{aligned}
\]
Finally, \eqref{ac-critical-grid-operator} gives
\[
 \begin{aligned}
 \norm{\partial_t^jB_tf_S^\infty}_X+
 \norm{\partial_t^j(A_t-I_X)f_S^\infty}_X
 &\leq 3\norm{\partial_t^jB_tf_S^\infty}_X+
        2\norm{\partial_t^jD_tf_S^\infty}_X+
        \norm{\partial_t^jR_A(t)f_S^\infty}_X\\
 &\leq \sqrt q\left[3\ell_+^j+2(2\ell_+)^j\varepsilon(t)
                  +4M\ell_+^j\varepsilon(t)^2\right]\sigma_S(t).
 \end{aligned}
\]
Take $C_A=4M(6\sqrt q+q)\max\{1,\ell_+^2\}$ in \eqref{ac-critical-grid-operator-bounds} and \eqref{ac-critical-local-operator-bounds}. Since $\ell_+^j\leq\max\{1,\ell_+^2\}$, $(2\ell_+)^j\leq4\max\{1,\ell_+^2\}$ and $\varepsilon(t)\leq1$, the coefficient of $\sigma_S(t)$ in the last inequality is bounded by
\[
 \sqrt q(11+4M)\max\{1,\ell_+^2\}
 \leq15M\sqrt q\max\{1,\ell_+^2\}
 \leq C_A,
\]
where the first inequality uses $M\geq16$ from \eqref{ac-critical-remainder-coefficients}.

The vertex equation \eqref{ac-vertex-equation} is
\begin{equation}\label{ac-critical-grid-vertex-equation}
 A_tf=\mathcal N(f),\qquad
 \mathcal N(f)(\vv)=\frac{\mathbf1_{\V_0}(\vv)}{\deg(\vv)}f(\vv)^3,
 \qquad f\in X.
\end{equation}

Put $v_S(t):=f_S(t)-f_S^\infty$ and recall $L_S^\infty=I_X-\mathcal N'(f_S^\infty)$ from \eqref{ac-limit-linearization}.
By \eqref{ac-critical-grid-vertex-equation}, $A_t(f_S^\infty+v_S)=\mathcal N(f_S^\infty+v_S)$. Using $\mathcal N(f_S^\infty)=f_S^\infty$, we obtain
\[
 \begin{aligned}
 v_S&=f_S^\infty+v_S-\mathcal N(f_S^\infty)\\
    &=-(A_t-I_X)(f_S^\infty+v_S)
      +\mathcal N(f_S^\infty+v_S)-\mathcal N(f_S^\infty).
 \end{aligned}
\]
At each vertex $\vv\in\V$, expanding the cubic in \eqref{ac-critical-grid-vertex-equation} gives
\begin{equation}\label{ac-critical-cubic-expansion}
 \begin{aligned}
 &\bigl(\mathcal N(f_S^\infty+v_S)-\mathcal N(f_S^\infty)\bigr)(\vv)\\
 &\qquad=\frac{\mathbf1_{\V_0}(\vv)}{\deg(\vv)}
  \bigl[3f_S^\infty(\vv)^2v_S(\vv)+3f_S^\infty(\vv)v_S(\vv)^2+v_S(\vv)^3\bigr]\\
 &\qquad=3(P_Sv_S)(\vv)+\bigl(\mathcal N'(v_S)f_S^\infty\bigr)(\vv)+\mathcal N(v_S)(\vv).
 \end{aligned}
\end{equation}
Here $f_S^\infty(\vv)^2=\deg(\vv)\mathbf1_S(\vv)$, $P_S$ is the projection in \eqref{ac-coordinate-subspace}, and the quadratic term $3\frac{\mathbf1_{\V_0}(\vv)}{\deg(\vv)}f_S^\infty(\vv)v_S(\vv)^2$ is $\mathcal N'(v_S)f_S^\infty$ by \eqref{ac-N-derivative}. By \eqref{ac-limit-linearization}, $L_S^\infty=I_X-3P_S$ at $p=4$, so moving $3P_Sv_S$ to the left yields
\begin{equation}\label{ac-critical-grid-linearized-identity}
 \begin{aligned}
 L_S^\infty v_S
 &=-(A_t-I_X)(f_S^\infty+v_S)
   +\mathcal N'(v_S)f_S^\infty+\mathcal N(v_S).
 \end{aligned}
\end{equation}
By \eqref{ac-limit-inverse-bound}, $\norm{(L_S^\infty)^{-1}}_{\mathcal L(X)}\leq1$. Since $\deg(\vv)\geq1$ and $\norm{f_S^\infty}_\infty\leq\sqrt q$, the two terms \(\mathcal N'(v_S)f_S^\infty\) and \(\mathcal N(v_S)\) are bounded pointwise by $3\sqrt q\,|v_S(\vv)|^2$ and $|v_S(\vv)|^3$, respectively. By \eqref{ac-norm-ordering},
\begin{equation}\label{ac-critical-grid-nonlinear-remainder-bound}
 \left\|\mathcal N'(v_S)f_S^\infty+\mathcal N(v_S)\right\|_X
 \leq3\sqrt q\,\norm{v_S}_\infty\norm{v_S}_X
       +\norm{v_S}_\infty^2\norm{v_S}_X \leq3\sqrt q\,\norm{v_S}_X^2+\norm{v_S}_X^3.
\end{equation}
Choose $T_0>\max\{T_b,1,3/\ell_*\}$ such that
\[
 C_Ae^{-T_0\ell_*}\leq\frac14,\qquad
 C_be^{-T_0\ell_*/2}
 \leq\min\left\{1,\frac1{4(6\sqrt q+3)}\right\}.
\]
This choice depends only on $q$, $k$, $\ell_*$ and $\ell_+$. For $t\geq T_0$, applying $(L_S^\infty)^{-1}$ to \eqref{ac-critical-grid-linearized-identity} and using its norm bound \eqref{ac-limit-inverse-bound}, the operator estimates \eqref{ac-critical-grid-operator-bounds}--\eqref{ac-critical-local-operator-bounds} at $j=0$, the nonlinear estimate \eqref{ac-critical-grid-nonlinear-remainder-bound}, and \eqref{ac-critical-initial-branch-bound} together with the choice of $T_0$, we obtain
\[
 \norm{v_S}_X
 \leq C_A\sigma_S+C_A\varepsilon\norm{v_S}_X
       +\left(3\sqrt q+1\right)\norm{v_S}_X^2
 \leq C_A\sigma_S+\frac12\norm{v_S}_X.
\]
Consequently, with $C_f:=2C_A$,
\begin{equation}\label{ac-critical-grid-profile-bound}
 \norm{v_S(t)}_X\leq C_f\sigma_S(t)\qquad\text{for }t\geq T_0.
\end{equation}
Set $L_S(t):=A_t-\mathcal N'(f_S(t))$. By \eqref{ac-N-derivative} at $p=4$ and $f_S=f_S^\infty+v_S$, for every $g\in X$ and every $\vv\in\V$,
\[
 \bigl[(\mathcal N'(f_S)-\mathcal N'(f_S^\infty))g\bigr](\vv)
 =\frac{3\mathbf1_{\V_0}(\vv)}{\deg(\vv)}
   \bigl[2f_S^\infty(\vv)v_S(\vv)+v_S(\vv)^2\bigr]g(\vv).
\]
By \eqref{ac-norm-ordering}, \eqref{ac-critical-initial-branch-bound} and the choice of $T_0$, we have $\norm{v_S}_\infty\leq\norm{v_S}_X\leq1$ for $t\geq T_0$. Since $\deg(\vv)\geq1$ and $\norm{f_S^\infty}_\infty\leq\sqrt q$, taking the $X$ norm and the supremum over $\norm g_X=1$ gives
\[
 \norm{\mathcal N'(f_S)-\mathcal N'(f_S^\infty)}_{\mathcal L(X)}
 \leq 6\sqrt q\,\norm{v_S}_\infty+3\norm{v_S}_\infty^2\leq (6\sqrt q+3)\norm{v_S}_X.
\]
By \eqref{ac-limit-linearization}, $L_S(t)-L_S^\infty=(A_t-I_X)-[\mathcal N'(f_S)-\mathcal N'(f_S^\infty)]$. Using \eqref{ac-critical-grid-operator-bounds} at $j=0$, \eqref{ac-critical-grid-profile-bound}, and $\sigma_S(t)\leq\sqrt{kq}\,\varepsilon$ by \eqref{ac-critical-weighted-size}, we obtain
\begin{equation}\label{ac-critical-grid-linearization-bound}
 \begin{aligned}
 \norm{L_S(t)-L_S^\infty}_{\mathcal L(X)}
 &\leq C_A\varepsilon+\left(6\sqrt q+3\right)\norm{v_S}_X \leq C_A\varepsilon+\left(6\sqrt q+3\right)C_f\sigma_S(t)
 \leq C_L\varepsilon,\\
 C_L&:=C_A+\left(6\sqrt q+3\right)C_f\sqrt{kq}.
 \end{aligned}
\end{equation}
The first inequality and the choice of $T_0$ also show that $\norm{L_S(t)-L_S^\infty}_{\mathcal L(X)}\leq1/2$. Applying \cite[Exercise~6.14]{Ha} to $I_X+(L_S^\infty)^{-1}(L_S(t)-L_S^\infty)$ gives
\begin{equation}\label{ac-critical-grid-inverse-bound}
 \norm{L_S(t)^{-1}}_{\mathcal L(X)}\leq2\qquad\text{for }t\geq T_0.
\end{equation}
Using \eqref{ac-branch-derivative}, \eqref{ac-critical-grid-inverse-bound}, \eqref{ac-critical-grid-operator-bounds}--\eqref{ac-critical-local-operator-bounds} and \eqref{ac-critical-grid-profile-bound}, we obtain
\begin{equation}\label{ac-critical-grid-branch-derivative-bound}
 \norm{f_S'}_X\leq2\norm{(\partial_tA_t)f_S}_X
 \leq2C_A(1+C_f)\sigma_S=:C_{\mathrm{der}}\sigma_S.
\end{equation}

Define the first correction and its remainder by
\[
 g_S(t):=2(L_S^\infty)^{-1}B_tf_S^\infty,\qquad
 r_S(t):=v_S(t)-g_S(t).
\]
Since $(L_S^\infty)^{-1}=-P_S/2+(I_X-P_S)$,
\[
 g_S=-P_SB_tf_S^\infty+2(I_X-P_S)B_tf_S^\infty,\qquad
 L_S^\infty g_S=2B_tf_S^\infty,\qquad
 L_S^\infty g_S'=2\partial_tB_tf_S^\infty.
\]
By \eqref{ac-critical-local-B-bound} and $\norm{(L_S^\infty)^{-1}}\leq1$, we have $\norm{g_S}_X\leq2\sqrt q\,\sigma_S$ and $\norm{g_S'}_X\leq2\sqrt q\,\ell_+\sigma_S$. Substitution into \eqref{ac-critical-grid-linearized-identity} gives
\begin{equation}\label{ac-critical-grid-remainder-identity}
 \begin{aligned}
 L_S^\infty r_S
 &=-(A_t-I_X)f_S^\infty-L_S^\infty g_S-(A_t-I_X)v_S
      +\mathcal N'(v_S)f_S^\infty+\mathcal N(v_S)\\
 &=-\left(2D_t+R_A(t)\right)f_S^\infty
      -(A_t-I_X)v_S
      +\mathcal N'(v_S)f_S^\infty+\mathcal N(v_S).
 \end{aligned}
\end{equation}
By \eqref{ac-limit-linearization}, $3P_Sv_S=\mathcal N'(f_S^\infty)v_S$, so \eqref{ac-critical-cubic-expansion} gives $\mathcal N'(v_S)f_S^\infty+\mathcal N(v_S)=\mathcal N(f_S)-\mathcal N(f_S^\infty)-\mathcal N'(f_S^\infty)v_S$. Differentiating the second expression in \eqref{ac-critical-grid-remainder-identity} and using $v_S'=r_S'+g_S'$, we obtain
\begin{equation*}
 \begin{aligned}
 L_S^\infty r_S'
 &=-\left(2\partial_tD_t+\partial_tR_A(t)\right)f_S^\infty-(\partial_tA_t)v_S-(A_t-I_X)v_S'+\bigl[\mathcal N'(f_S)-\mathcal N'(f_S^\infty)\bigr]v_S'\\
 &=-\left(2\partial_tD_t+\partial_tR_A(t)\right)f_S^\infty-(\partial_tA_t)v_S-(L_S(t)-L_S^\infty)(r_S'+g_S'),
 \end{aligned}
\end{equation*}
which leads to
\begin{equation}\label{ac-critical-grid-remainder-derivative-identity}
 L_S(t)r_S'
 =-\left(2\partial_tD_t+\partial_tR_A(t)\right)f_S^\infty
       -(\partial_tA_t)v_S-(L_S(t)-L_S^\infty)g_S'.
\end{equation}
By \eqref{ac-critical-grid-remainder-identity}, the inverse bound  $\norm{(L_S^\infty)^{-1}}\leq1$ from \eqref{ac-limit-inverse-bound}, \eqref{ac-critical-grid-operator-bounds}--\eqref{ac-critical-local-operator-bounds} at $j=0$, and \eqref{ac-critical-grid-nonlinear-remainder-bound}--\eqref{ac-critical-grid-profile-bound}, we have, using $\varepsilon\leq1$ and $\sigma_S\leq\sqrt{kq}\,\varepsilon$ from \eqref{ac-critical-weighted-size},
\[
 \begin{aligned}
 \norm{r_S}_X
 &\leq 2\norm{D_tf_S^\infty}_X+\norm{R_A(t)f_S^\infty}_X
       +\norm{A_t-I_X}_{\mathcal L(X)}\norm{v_S}_X+\norm{\mathcal N'(v_S)f_S^\infty+\mathcal N(v_S)}_X\\
 &\leq 2C_A\varepsilon\sigma_S+C_A\varepsilon^2\sigma_S
       +C_AC_f\varepsilon\sigma_S+(3\sqrt q+1)C_f^2\sigma_S^2\\
 &\leq\left[C_A\left(3+C_f\right)
       +\left(3\sqrt q+1\right)C_f^2\sqrt{kq}\right]\varepsilon\sigma_S.
 \end{aligned}
\]
For $r_S'$, \eqref{ac-critical-grid-remainder-derivative-identity}, the inverse bound \eqref{ac-critical-grid-inverse-bound}, \eqref{ac-critical-grid-operator-bounds}--\eqref{ac-critical-local-operator-bounds} at $j=1$, \eqref{ac-critical-grid-profile-bound}--\eqref{ac-critical-grid-linearization-bound}, and $\norm{g_S'}_X\leq2\sqrt q\,\ell_+\sigma_S$  give
\[
 \begin{aligned}
 \norm{r_S'}_X
 &\leq 2\Bigl[2\norm{(\partial_tD_t)f_S^\infty}_X
      +\norm{(\partial_tR_A(t))f_S^\infty}_X
      +\norm{\partial_tA_t}_{\mathcal L(X)}\norm{v_S}_X+\norm{L_S(t)-L_S^\infty}_{\mathcal L(X)}\norm{g_S'}_X\Bigr]\\
 &\leq 2\bigl[2C_A\varepsilon\sigma_S+C_A\varepsilon^2\sigma_S
       +C_AC_f\varepsilon\sigma_S+2\sqrt q\,\ell_+C_L\varepsilon\sigma_S\bigr]\\
 &\leq2\left[C_A\left(3+C_f\right)+2\sqrt q\,\ell_+C_L\right]\varepsilon\sigma_S.
 \end{aligned}
\]
Thus
\begin{equation}\label{ac-critical-grid-branch-remainder-bound}
 \norm{r_S}_X+\norm{r_S'}_X\leq C_{\mathrm{res}}\varepsilon\sigma_S,
 \quad
 C_{\mathrm{res}}:=3C_A\left(3+C_f\right)
 +\left(3\sqrt q+1\right)C_f^2\sqrt{kq}+4\sqrt q\,\ell_+C_L.
\end{equation}

Define
\begin{equation}\label{ac-critical-grid-HJ-definitions}
 H_S(t):=\langle A_tf_S(t),f_S(t)\rangle_X,\qquad
 J_S(t):=\langle\partial_tA_tf_S(t),f_S(t)\rangle_X.
\end{equation}
The cubic equation \eqref{ac-critical-grid-vertex-equation} and \eqref{ac-branch-derivative} give $L_S(t)f_S=-2A_tf_S$ and $L_S(t)f_S'=-(\partial_tA_t)f_S$. By self-adjointness,
\[
 -2\langle A_tf_S,f_S'\rangle_X
 =\langle L_S(t)f_S,f_S'\rangle_X
 =\langle f_S,L_S(t)f_S'\rangle_X=-J_S(t),
\]
which implies that $H_S'=2J_S$.
Consequently, \eqref{ac-critical-mass-expansion} and \eqref{ac-critical-energy-exact} imply
\begin{equation}\label{ac-critical-grid-mass-energy-identities}
 \mu_S=\frac12(H_S+tJ_S),\qquad
 \mu_S'=\frac12(3J_S+tJ_S'),\qquad
 E(u_S(t),\G)=-\frac{t^3}{4}J_S.
\end{equation}
We compute $J_S$ and $J_S'$ using the following scalar remainders:
\begin{equation}\label{ac-critical-grid-scalar-remainder-definitions}
 \begin{aligned}
 R_{S,0}:={}&-4\langle\partial_tB_tf_S^\infty,r_S\rangle_X
           -2\langle\partial_tB_tv_S,v_S\rangle_X,\\
 R_{S,1}:={}&4\langle\partial_tD_tf_S^\infty,v_S\rangle_X
           +2\langle\partial_tD_tv_S,v_S\rangle_X,\\
 R_{S,2}:={}&\langle\partial_tR_A(t)f_S,f_S\rangle_X\\
 ={}&\langle\partial_tR_A(t)f_S^\infty,f_S^\infty\rangle_X
       +2\langle\partial_tR_A(t)f_S^\infty,v_S\rangle_X
       +\langle\partial_tR_A(t)v_S,v_S\rangle_X.
 \end{aligned}
\end{equation}
Inserting $f_S=f_S^\infty+g_S+r_S$ into the definition of $J_S$ in \eqref{ac-critical-grid-HJ-definitions} gives the exact identity
\begin{equation}\label{ac-critical-general-J}
 J_S=-2\langle\partial_tB_tf_S^\infty,f_S^\infty\rangle_X
       -4\langle\partial_tB_tf_S^\infty,g_S\rangle_X+2\langle\partial_tD_tf_S^\infty,f_S^\infty\rangle_X
       +\sum_{j=0}^2R_{S,j}.
\end{equation}
Differentiating the remainders in \eqref{ac-critical-grid-scalar-remainder-definitions}, with $v_S'=f_S'$, gives
\begin{equation}\label{ac-critical-grid-scalar-remainder-derivatives}
 \begin{aligned}
 R_{S,0}'&=-4\left(\langle\partial_t^2B_tf_S^\infty,r_S\rangle_X
                    +\langle\partial_tB_tf_S^\infty,r_S'\rangle_X\right)-2\left(\langle\partial_t^2B_tv_S,v_S\rangle_X
                      +2\langle\partial_tB_tv_S,f_S'\rangle_X\right),\\
 R_{S,1}'&=4\left(\langle\partial_t^2D_tf_S^\infty,v_S\rangle_X
                    +\langle\partial_tD_tf_S^\infty,f_S'\rangle_X\right)+2\left(\langle\partial_t^2D_tv_S,v_S\rangle_X
                      +2\langle\partial_tD_tv_S,f_S'\rangle_X\right),\\
 R_{S,2}'&=\langle\partial_t^2R_A(t)f_S,f_S\rangle_X
             +2\langle\partial_tR_A(t)f_S,f_S'\rangle_X\\
 &=\langle\partial_t^2R_A(t)f_S^\infty,f_S^\infty\rangle_X
       +2\langle\partial_t^2R_A(t)f_S^\infty,v_S\rangle_X
       +\langle\partial_t^2R_A(t)v_S,v_S\rangle_X\\
 &\quad+2\langle\partial_tR_A(t)f_S^\infty,f_S'\rangle_X
       +2\langle\partial_tR_A(t)v_S,f_S'\rangle_X.
 \end{aligned}
\end{equation}
Therefore, $\varepsilon\leq 1$, \eqref{ac-critical-grid-operator-bounds}, \eqref{ac-critical-local-operator-bounds}, \eqref{ac-critical-grid-profile-bound}, \eqref{ac-critical-grid-branch-derivative-bound} and \eqref{ac-critical-grid-branch-remainder-bound} imply
\[
 \begin{aligned}
 |R_{S,0}|+|R_{S,0}'|
 &\leq 4C_A\sigma_S\bigl[2\norm{r_S}_X+\norm{r_S'}_X\bigr]
       +4C_AC_f^2\varepsilon\sigma_S^2+4C_AC_fC_{\mathrm{der}}\varepsilon\sigma_S^2\\
 &\leq C_A(8C_{\mathrm{res}}+4C_f^2+4C_fC_{\mathrm{der}})\varepsilon\sigma_S^2,\\
 |R_{S,1}|+|R_{S,1}'|
 &\leq 8C_AC_f\varepsilon\sigma_S^2+4C_AC_f^2\varepsilon^2\sigma_S^2
       +4C_AC_{\mathrm{der}}\varepsilon\sigma_S^2
       +4C_AC_fC_{\mathrm{der}}\varepsilon^2\sigma_S^2\\
 &\leq C_A(8C_f+4C_{\mathrm{der}}+4C_f^2+4C_fC_{\mathrm{der}})\varepsilon\sigma_S^2,\\
|R_{S,2}|+|R_{S,2}'|
 &\leq 2C_A\varepsilon\sigma_S^2
       +4C_AC_f\varepsilon^2\sigma_S^2+2C_AC_f^2\varepsilon^3\sigma_S^2+2C_AC_{\mathrm{der}}\varepsilon^2\sigma_S^2
       +2C_AC_fC_{\mathrm{der}}\varepsilon^3\sigma_S^2\\
 &\leq C_A(2+4C_f+2C_f^2+2C_{\mathrm{der}}
       +2C_fC_{\mathrm{der}})\varepsilon\sigma_S^2.
 \end{aligned}
\]
Then there exists a constant $C_J=C_J(q,k,\ell_*,\ell_+)>0$ such that
\begin{equation}\label{ac-critical-grid-scalar-remainders}
 \sum_{j=0}^2\bigl(|R_{S,j}(t)|+|R_{S,j}'(t)|\bigr)
 \leq C_J\varepsilon(t)\sigma_S(t)^2\qquad\text{for }t\geq T_0.
\end{equation}

 To estimate the inner products in \eqref{ac-critical-general-J}, we use the definitions of $B_t$ and $D_t$ and $f_S^\infty(\vv)=\sqrt{\deg(\vv)}\mathbf1_S(\vv)$ from \eqref{ac-limit-profile}. 
In particular,
\begin{equation}\label{ac-critical-B-profile-coordinates}
 \begin{aligned}
 (B_tf_S^\infty)(\vv)
 &=\frac1{\deg(\vv)}\sum_{\e\in\E_S(\vv)}
     \sqrt{\deg(\w_{\e,\vv})}\,e^{-t\ell_\e},\\
 ((\partial_tB_t)f_S^\infty)(\vv)
 &=-\frac1{\deg(\vv)}\sum_{\e\in\E_S(\vv)}
     \sqrt{\deg(\w_{\e,\vv})}\,\ell_\e e^{-t\ell_\e}.
 \end{aligned}
\end{equation}
Moreover, these definitions imply
\[
 \begin{aligned}
 \langle\partial_tB_tf_S^\infty,f_S^\infty\rangle_X
 &=-\sum_{\vv\in S}\sum_{\e\in\E_S(\vv)}
       \sqrt{\deg(\vv)\deg(\w_{\e,\vv})}\,\ell_\e e^{-t\ell_\e}\\
 &=-2\sum_{\e\in\E(S)}\sqrt{\deg(\vv_\e)\deg(\w_\e)}\,\ell_\e e^{-t\ell_\e},\\
 \langle\partial_tD_tf_S^\infty,f_S^\infty\rangle_X
 &=-2\sum_{\vv\in S}\sum_{\e\in\E_b,\,\e\succ\vv}\deg(\vv)\ell_\e e^{-2t\ell_\e}.
 \end{aligned}
\]
Using \eqref{ac-critical-B-profile-coordinates} in the two components of $g_S=-P_SB_tf_S^\infty+2(I_X-P_S)B_tf_S^\infty$, we obtain
\[
 \begin{aligned}
 \langle\partial_tB_tf_S^\infty,P_SB_tf_S^\infty\rangle_X
 &=-\sum_{\vv\in S}\sum_{\e,\e'\in\E_S(\vv)}
        \omega_\vv(\e,\e')\ell_\e e^{-t(\ell_\e+\ell_{\e'})}\\
 &=-\frac12\sum_{\vv\in S}\sum_{\e,\e'\in\E_S(\vv)}
        \omega_\vv(\e,\e')(\ell_\e+\ell_{\e'})e^{-t(\ell_\e+\ell_{\e'})},\\
 \langle\partial_tB_tf_S^\infty,(I_X-P_S)B_tf_S^\infty\rangle_X
 &=-\sum_{\vv\in N(S)}\sum_{\e,\e'\in\E_S(\vv)}
        \omega_\vv(\e,\e')\ell_\e e^{-t(\ell_\e+\ell_{\e'})}\\
 &=-\frac12\sum_{\vv\in N(S)}\sum_{\e,\e'\in\E_S(\vv)}
        \omega_\vv(\e,\e')(\ell_\e+\ell_{\e'})e^{-t(\ell_\e+\ell_{\e'})}.
 \end{aligned}
\]
Substitution into \eqref{ac-critical-general-J} yields
\begin{equation}\label{ac-critical-general-edge-J}
 \begin{aligned}
 J_S(t)={}&4\sum_{\e\in\E(S)}\sqrt{\deg(\vv_\e)\deg(\w_\e)}\,\ell_\e e^{-t\ell_\e}
       -4\sum_{\vv\in S}\sum_{\e\in\E_b,\,\e\succ\vv}\deg(\vv)\ell_\e e^{-2t\ell_\e}\\
 &-2\sum_{\vv\in S}\sum_{\e,\e'\in\E_S(\vv)}
        \omega_\vv(\e,\e')(\ell_\e+\ell_{\e'})e^{-t(\ell_\e+\ell_{\e'})}\\
 &+4\sum_{\vv\in N(S)}\sum_{\e,\e'\in\E_S(\vv)}
        \omega_\vv(\e,\e')(\ell_\e+\ell_{\e'})e^{-t(\ell_\e+\ell_{\e'})}
       +\sum_{j=0}^2R_{S,j}(t).
 \end{aligned}
\end{equation}
By \eqref{ac-critical-grid-operator-bounds} and \eqref{ac-critical-grid-profile-bound}, $H_S(t)\to\norm{f_S^\infty}_X^2=\sum_{\vv\in S}\deg(\vv)^2=2\mu_S^\infty$.
Using $H_S'=2J_S$ and $\mu_S=(H_S+tJ_S)/2$ from \eqref{ac-critical-grid-mass-energy-identities}, we obtain
\[
 \mu_S(t)=\mu_S^\infty-\int_t^\infty J_S(s)\,ds+\frac t2J_S(t).
\]
For every $a>0$, a term $a e^{-ta}$ in \eqref{ac-critical-general-edge-J} contributes
\[
 -\int_t^\infty a e^{-sa}\,ds+\frac t2 a e^{-ta}
 =-e^{-ta}+\frac{ta}{2}e^{-ta}
 =\frac{ta-2}{2}e^{-ta}
\]
to $\mu_S(t)-\mu_S^\infty$. Applying this identity with $a=\ell_\e$, $2\ell_\e$ and $\ell_\e+\ell_{\e'}$ in \eqref{ac-critical-general-edge-J} proves the mass expansion \eqref{ac-critical-general-mass-expansion}. The remainder and its derivative are
\begin{equation}\label{ac-critical-general-mass-remainder}
 \begin{aligned}
 \mathcal E_S(t)&:=\frac t2\sum_{j=0}^2R_{S,j}(t)
       -\int_t^\infty\sum_{j=0}^2R_{S,j}(s)\,ds,\\
 \mathcal E_S'(t)&=\frac12\sum_{j=0}^2
                     [3R_{S,j}(t)+tR_{S,j}'(t)].
 \end{aligned}
\end{equation}

To prove \eqref{ac-critical-general-mass-slope}, differentiate \eqref{ac-critical-general-edge-J} and use $\mu_S'=(3J_S+tJ_S')/2$ from \eqref{ac-critical-grid-mass-energy-identities}. For each exponential term,
\[
 \frac12\left[3\ell e^{-t\ell}
       +t\partial_t(\ell e^{-t\ell})\right]
 =\frac{\ell}{2}(3-t\ell)e^{-t\ell}.
\]
Applying this identity with $\ell=\ell_\e$, $2\ell_\e$ and $\ell_\e+\ell_{\e'}$ gives the displayed coefficients in \eqref{ac-critical-general-mass-slope}, and the contribution of $\sum_{j=0}^2R_{S,j}$ is precisely $\mathcal E_S'(t)$ in \eqref{ac-critical-general-mass-remainder}.

Since $\ell_*+2\ell_\e\geq3\ell_*$, \eqref{ac-critical-grid-scalar-remainders} gives
\[
 \begin{aligned}
 \int_t^\infty\sum_{j=0}^2|R_{S,j}(s)|\,ds
 &\leq C_J\int_t^\infty e^{-s\ell_*}\sigma_S(s)^2\,ds\\
 &=C_J\sum_{\vv\in S}\sum_{\e\in\E_b,\,\e\succ\vv}
          \int_t^\infty e^{-s(\ell_*+2\ell_\e)}\,ds\\
 &=C_J\sum_{\vv\in S}\sum_{\e\in\E_b,\,\e\succ\vv}
          \frac{e^{-t(\ell_*+2\ell_\e)}}{\ell_*+2\ell_\e}
 \leq\frac{C_J}{3\ell_*}e^{-t\ell_*}\sigma_S(t)^2.
  \end{aligned}
\]
Combining this integral estimate with
\eqref{ac-critical-grid-scalar-remainders} and
\eqref{ac-critical-general-mass-remainder}, we obtain\[
 |\mathcal E_S(t)|+|\mathcal E_S'(t)|
 \leq C_J\left[\frac t2+\frac32+\frac1{3\ell_*}\right]
                  e^{-t\ell_*}\sigma_S(t)^2
 \leq C_rt e^{-t\ell_*}\sigma_S(t)^2.
\]
This proves \eqref{ac-critical-general-mass-error} for $t\geq T_0>1$ and $C_r:=C_J(3+1/\ell_*)$, where the constant $C_r$ depends only on $q$, $k$, $\ell_*$ and $\ell_+$. By \eqref{ac-critical-grid-mass-energy-identities}, $E(u_S(t),\G)=-t^3J_S(t)/4$. Thus \eqref{ac-critical-general-edge-J} gives \eqref{ac-critical-general-energy-expansion}. By \eqref{ac-critical-grid-scalar-remainders},
\[
 \frac{t^3}{4}\left|\sum_{j=0}^2R_{S,j}(t)\right|
 \leq\frac{C_J}{4}t^3e^{-t\ell_*}\sigma_S(t)^2
 \leq C_rt^3e^{-t\ell_*}\sigma_S(t)^2.
\]
\end{proof}
  We now study the signs of $\mu_S'(t)$.
If $S$ has no bounded incident edge, then $\sigma_S(t)=0$ and the mass expansion gives $\mu_S(t)=\mu_S^\infty$ for every $t\geq T_0$.

For a set $S$ consisting of $k$ vertices of degree $q$, define
\[
\Xi_S:=\sum_{\vv\in\V}\#\E_S(\vv)-2\sum_{\vv\in N(S)}\frac{(\#\E_S(\vv))^2}{\deg(\vv)}.
\]

The next proposition uses the expansions in Proposition~\ref{prop-ac-critical-local-sign} to determine the signs of the mass gap, the mass derivative and the energy under three sets of geometric conditions.

\begin{proposition}\label{prop-ac-critical-sign-conditions}
Assume the hypotheses of Proposition~\ref{prop-ac-critical-local-sign}, and suppose that $\deg(\vv)=q$ for all $\vv\in S$. Assume also that at least one bounded edge is incident at $S$, so $\sigma_S(t)>0$ by \eqref{ac-critical-weighted-size}. There is $T_{\rm loc}>T_0$, depending only on $q$, $k$, $\ell_*$ and $\ell_+$, with the following properties.
\begin{enumerate}
\renewcommand{\labelenumi}{\textup{(\roman{enumi})}}
\item If $\G$ is simple, $\operatorname{dist}_{\rm c}(\vv,\w)\geq3$ for all distinct $\vv,\w\in S$ and $\deg(\vv)\geq3$ for every $\vv\in N(S)$, then for every $t\geq T_{\rm loc}$,
\begin{equation}\label{ac-critical-grid-two-sided-gap}
 \frac q3\ell_*t\sigma_S(t)^2
 \leq\frac{kq^2}{2}-\mu_S(t)
 \leq3q\ell_+t\sigma_S(t)^2,
\end{equation}
\begin{equation}\label{ac-critical-grid-positive-slope-energy}
 \mu_S'(t)\geq\frac{5q}{9}\ell_*^2t\sigma_S(t)^2>0,
 \qquad
 E(u_S(t),\G)\geq\frac{2q}{9}\ell_*t^3\sigma_S(t)^2>0.
\end{equation}
\item Suppose that $\E(S)\ne\varnothing$ and that there is $L>0$ such that $\ell_\e=L$ for every $\vv\in S$ and every $\e\in\E_b$ with $\e\succ\vv$. Then for every $t\geq T_{\rm loc}$,
\[
 \frac q2 L\,\#\E(S)t e^{-tL}
 \leq\mu_S(t)-\frac{kq^2}{2}
 \leq\frac{5q}{2}L\,\#\E(S)t e^{-tL},
\]
\[
 -\mu_S'(t)\geq\frac q2L^2\,\#\E(S)t e^{-tL}>0,
 \qquad
 -E(u_S(t),\G)\geq\frac q2L\,\#\E(S)t^3e^{-tL}>0.
\]
\item Suppose that $\E(S)=\varnothing$, $\Xi_S\ne0$ and that there is $L>0$ such that $\ell_\e=L$ for every $\vv\in S$ and every $\e\in\E_b$ with $\e\succ\vv$. There is $T_{\Xi,S}\geq T_{\rm loc}$, depending only on $q$, $k$, $\ell_*$, $\ell_+$ and $|\Xi_S|^{-1}$, such that for every $t\geq T_{\Xi,S}$,
\[
 qL|\Xi_S|t e^{-2tL}
 \leq(\operatorname{sgn}\Xi_S)\left[\frac{kq^2}{2}-\mu_S(t)\right]
 \leq\frac{7q}{3}L|\Xi_S|t e^{-2tL},
\]
\[
 \begin{aligned}
 (\operatorname{sgn}\Xi_S)\mu_S'(t)&\geq\frac{5q}{3}L^2|\Xi_S|t e^{-2tL}>0,\\
 (\operatorname{sgn}\Xi_S)E(u_S(t),\G)&\geq\frac{2q}{3}L|\Xi_S|t^3e^{-2tL}>0.
 \end{aligned}
\]
\end{enumerate}
\end{proposition}

\begin{proof}
  Assume the hypotheses of Proposition~\ref{prop-ac-critical-local-sign}, and suppose that $\deg(\vv)=q$ for all $\vv\in S$.
By the equal-degree assumption, $\mu_S^\infty=kq^2/2$, $\sqrt{\deg(\vv_\e)\deg(\w_\e)}=q$ on $\E(S)$, and $\omega_\vv(\e,\e')=q/\deg(\vv)$ for every $\e,\e'\in\E_S(\vv)$. Let $C_r$ be the constant in Proposition~\ref{prop-ac-critical-local-sign}, and put
\[
 C_{\rm eq}:=8kq^2\left[1+4\ell_++2\ell_+^2\right]+kqC_r.
\]
Choose
\[
 T_{\rm loc}>\max\left\{T_0,\frac6{\ell_*},
 \frac1{\ell_*}\log\frac{9C_r}{q\min\{\ell_*,\ell_*^2\}},
 \frac1{\ell_*}\log\frac{2C_{\rm eq}}{q\min\{\ell_*,\ell_*^2\}}\right\}.
\]
In part~\textup{(i)}, separation gives $\E(S)=\varnothing$ and $\E_S(\vv)=\varnothing$ for every $\vv\in S$, and each $\vv\in N(S)$ has exactly one neighbor in $S$. Since $\G$ is simple, $\#\E_S(\vv)=1$ for every $\vv\in N(S)$. Thus, for every $\phi:(0,+\infty)\to\R$,
\begin{equation}\label{ac-critical-separated-sum-reduction}
 \sum_{\vv\in N(S)}\frac1{\deg(\vv)}\sum_{\e,\e'\in\E_S(\vv)}\phi(\ell_\e+\ell_{\e'})
 =\sum_{\vv\in S}\sum_{\e\in\E_b,\,\e\succ\vv}\frac{\phi(2\ell_\e)}{\deg(\w_{\e,\vv})}.
\end{equation}
For $\vv\in S$ and $\e\succ\vv$ with $\e\in\E_b$, define
\begin{equation}\label{ac-critical-local-coefficient}
 \kappa_{\vv,\e}:=\deg(\vv)\left[1-\frac2{\deg(\w_{\e,\vv})}\right].
\end{equation}
Since $\E(S)=\varnothing$ and $\E_S(\vv)=\varnothing$ for every $\vv\in S$, the corresponding sums vanish in \eqref{ac-critical-general-mass-expansion} and \eqref{ac-critical-general-mass-slope}. Reducing the sums over $N(S)$ by \eqref{ac-critical-separated-sum-reduction} and then combining coefficients by \eqref{ac-critical-local-coefficient}, we obtain
\begin{equation}\label{ac-critical-grid-mass-expansion}
 \begin{aligned}
 \mu_S(t)
 &=\frac{kq^2}{2}-2q\sum_{\vv\in S}\sum_{\e\in\E_b,\,\e\succ\vv}
      (t\ell_\e-1)e^{-2t\ell_\e}+2q\sum_{\vv\in S}\sum_{\e\in\E_b,\,\e\succ\vv}
      \frac{(2t\ell_\e-2)e^{-2t\ell_\e}}{\deg(\w_{\e,\vv})}
      +\mathcal E_S(t)\\
 &=\frac{kq^2}{2}-2\sum_{\vv\in S}\sum_{\e\in\E_b,\,\e\succ\vv}
      \kappa_{\vv,\e}(t\ell_\e-1)e^{-2t\ell_\e}+\mathcal E_S(t),\\
 \mu_S'(t)
 &=2q\sum_{\vv\in S}\sum_{\e\in\E_b,\,\e\succ\vv}
      \ell_\e(2t\ell_\e-3)e^{-2t\ell_\e}-2q\sum_{\vv\in S}\sum_{\e\in\E_b,\,\e\succ\vv}
      \frac{2\ell_\e(2t\ell_\e-3)e^{-2t\ell_\e}}{\deg(\w_{\e,\vv})}
      +\mathcal E_S'(t)\\
 &=2\sum_{\vv\in S}\sum_{\e\in\E_b,\,\e\succ\vv}
      \kappa_{\vv,\e}\ell_\e(2t\ell_\e-3)e^{-2t\ell_\e}+\mathcal E_S'(t).
 \end{aligned}
\end{equation}
The same reduction in \eqref{ac-critical-general-energy-expansion} gives
\begin{equation}\label{ac-critical-grid-energy-expansion}
 \left|E(u_S(t),\G)-t^3\sum_{\vv\in S}\sum_{\e\in\E_b,\,\e\succ\vv}\kappa_{\vv,\e}
                  \ell_\e e^{-2t\ell_\e}\right|
 \leq C_r t^3e^{-t\ell_*}\sigma_S(t)^2.
\end{equation}
In part (i), since $\deg(\w)\geq3$ for every $\w\in N(S)$, we have $q/3\leq\kappa_{\vv,\e}\leq q$. For $t\geq T_{\rm loc}$, we have $t\ell_\e\geq t\ell_*\geq3$ for every bounded edge incident at $S$. Hence $t\ell_\e-1\geq2t\ell_\e/3$ and $2t\ell_\e-3\geq t\ell_\e$. The choice of $T_{\rm loc}$ also gives $C_re^{-t\ell_*}\leq(q/9)\min\{\ell_*,\ell_*^2\}$. By \eqref{ac-critical-grid-mass-expansion}, \eqref{ac-critical-general-mass-error} and the definition of $\sigma_S$ in \eqref{ac-critical-weighted-size},
\[
 \begin{aligned}
 \frac{4q}{9}\ell_*t\sigma_S(t)^2
 &\leq2\sum_{\vv\in S}\sum_{\e\in\E_b,\,\e\succ\vv}
       \kappa_{\vv,\e}(t\ell_\e-1)e^{-2t\ell_\e}
 \leq2q\ell_+t\sigma_S(t)^2,\\
 |\mathcal E_S(t)|
 &\leq C_rt e^{-t\ell_*}\sigma_S(t)^2
 \leq\frac q9\ell_*t\sigma_S(t)^2.
 \end{aligned}
\]
These bounds give \eqref{ac-critical-grid-two-sided-gap}. The derivative in \eqref{ac-critical-grid-mass-expansion} and the energy estimate \eqref{ac-critical-grid-energy-expansion} give
\[
 \begin{aligned}
 \mu_S'(t)
 &\geq2\sum_{\vv\in S}\sum_{\e\in\E_b,\,\e\succ\vv}
       \kappa_{\vv,\e}\ell_\e(2t\ell_\e-3)e^{-2t\ell_\e}
       -C_rt e^{-t\ell_*}\sigma_S(t)^2\\
 &\geq\frac{2q}{3}\ell_*^2t\sigma_S(t)^2-\frac q9\ell_*^2t\sigma_S(t)^2
 =\frac{5q}{9}\ell_*^2t\sigma_S(t)^2,\\
 E(u_S(t),\G)
 &\geq t^3\sum_{\vv\in S}\sum_{\e\in\E_b,\,\e\succ\vv}
       \kappa_{\vv,\e}\ell_\e e^{-2t\ell_\e}
       -C_rt^3e^{-t\ell_*}\sigma_S(t)^2\\
 &\geq\frac q3\ell_*t^3\sigma_S(t)^2-\frac q9\ell_*t^3\sigma_S(t)^2
 =\frac{2q}{9}\ell_*t^3\sigma_S(t)^2.
 \end{aligned}
\]
Since $\sigma_S(t)>0$, this proves \eqref{ac-critical-grid-positive-slope-energy}.

For parts (ii) and (iii), every bounded edge incident at $S$ has length $L$. In $\sum_{\vv\in\V}\#\E_S(\vv)$, each such edge is counted once for each endpoint in $S$. Hence $\sum_{\vv\in S}\sum_{\e\in\E_b,\,\e\succ\vv}1=\sum_{\vv\in\V}\#\E_S(\vv)$. Since $\deg(\w)=q$ for every $\w\in S$, we have $\omega_\vv(\e,\e')=q/\deg(\vv)$ for every $\e,\e'\in\E_S(\vv)$. The three terms of order $e^{-2tL}$ in \eqref{ac-critical-general-mass-expansion} therefore sum to
\[
 \begin{aligned}
 &-2q\sum_{\vv\in\V}\#\E_S(\vv)(tL-1)e^{-2tL}
   -\sum_{\vv\in S}(\#\E_S(\vv))^2(2tL-2)e^{-2tL}+2q\sum_{\vv\in N(S)}\frac{(\#\E_S(\vv))^2}{\deg(\vv)}(2tL-2)e^{-2tL}\\
 &=-2\left[q\left(\sum_{\vv\in\V}\#\E_S(\vv)
           -2\sum_{\vv\in N(S)}\frac{(\#\E_S(\vv))^2}{\deg(\vv)}\right)
           +\sum_{\vv\in S}(\#\E_S(\vv))^2\right](tL-1)e^{-2tL}\\
 &=-2\left[q\Xi_S+\sum_{\vv\in S}(\#\E_S(\vv))^2\right](tL-1)e^{-2tL}.
 \end{aligned}
\]
With $\ell_\e=L$ on every bounded edge incident at $S$, \eqref{ac-critical-general-mass-expansion}, \eqref{ac-critical-general-mass-slope} and \eqref{ac-critical-general-energy-expansion} give
\begin{equation}\label{ac-critical-equal-length-expansions}
 \begin{aligned}
 \mu_S(t)
 &=\frac{kq^2}{2}+2q\,\#\E(S)(tL-2)e^{-tL}\\
 &\quad-2\left[q\Xi_S+\sum_{\vv\in S}(\#\E_S(\vv))^2\right](tL-1)e^{-2tL}
      +\mathcal E_S(t),\\
 \mu_S'(t)
 &=2qL\,\#\E(S)(3-tL)e^{-tL}\\
 &\quad+2L\left[q\Xi_S+\sum_{\vv\in S}(\#\E_S(\vv))^2\right](2tL-3)e^{-2tL}
      +\mathcal E_S'(t),\\
 E(u_S(t),\G)
 &=-qL\,\#\E(S)t^3e^{-tL}\\
 &\quad+L\left[q\Xi_S+\sum_{\vv\in S}(\#\E_S(\vv))^2\right]t^3e^{-2tL}
      -\frac{t^3}{4}\sum_{j=0}^2R_{S,j}(t).
 \end{aligned}
\end{equation}
Since $\sum_{\vv\in\V}\#\E_S(\vv)=\sum_{\vv\in S}\sum_{\e\in\E_b,\,\e\succ\vv}1\leq kq$ and $\#\E_S(\vv)\leq\deg(\vv)$, we have 
\[
 \left|q\Xi_S+\sum_{\vv\in S}(\#\E_S(\vv))^2\right|
 \leq kq^2+2q\sum_{\w\in N(S)}\#\E_S(\w)+q\sum_{\vv\in S}\#\E_S(\vv)
 \leq4kq^2 .
\]
For part (ii), this bound, $L\leq\ell_+$ and \eqref{ac-critical-general-mass-error} give
\begin{equation}\label{ac-critical-adjacent-mass}
 \begin{aligned}
 &\left|\mu_S(t)-\frac{kq^2}{2}-2q\,\#\E(S)(tL-2)e^{-tL}\right|\\
 &\quad+\left|\mu_S'(t)-2qL\,\#\E(S)(3-tL)e^{-tL}\right|
 \leq C_{\rm eq}t e^{-2tL}.
 \end{aligned}
\end{equation}
Indeed, the sum of the absolute values of the coefficients of $e^{-2tL}$ in the expansions \eqref{ac-critical-equal-length-expansions} of $\mu_S(t)$ and $\mu_S'(t)$ is at most $8kq^2[t\ell_++1+\ell_+(2t\ell_++3)]$.
By \eqref{ac-critical-general-mass-error} and $\sigma_S(t)^2\leq kq e^{-2tL}$,
\[
 \abs{\mathcal E_S(t)}+\abs{\mathcal E_S'(t)}
 \leq C_rt e^{-t\ell_*}\sigma_S(t)^2
 \leq kqC_rt e^{-t(2L+\ell_*)}
 \leq kqC_rt e^{-2tL}.
\]
For $t\geq T_{\rm loc}$, we have $tL\geq6$ and hence $tL\leq2(tL-2)\leq2tL$ and $tL\leq2(tL-3)$. Since $\#\E(S)\geq1$, the choice of $T_{\rm loc}$ gives
\[
 C_{\rm eq}te^{-2tL}\leq\frac q2\min\{\ell_*,\ell_*^2\}te^{-tL}
 \leq\frac q2\,\#\E(S)\min\{L,L^2\}\,te^{-tL}.
\]
Using this error bound, with $L$ for the mass and $L^2$ for the derivative, in \eqref{ac-critical-adjacent-mass} proves the mass and derivative estimates in part (ii).
The bound on the coefficient of $e^{-2tL}$ in the expansion of $E(u_S(t),\G)$ and the remainder estimate after \eqref{ac-critical-general-energy-expansion} give
\[
 \begin{aligned}
 \left|E(u_S(t),\G)+qL\,\#\E(S)t^3e^{-tL}\right|
 &\leq(4kq^2\ell_++kqC_r)t^3e^{-2tL}\\
 &\leq C_{\rm eq}t^3e^{-2tL}
 \leq\frac q2\ell_*t^3e^{-tL}
 \qquad\text{for }t\geq T_{\rm loc}.
 \end{aligned}
\]
Since $qL\,\#\E(S)\geq q\ell_*$, this gives the energy bound in part (ii).

For part (iii), $\E(S)=\varnothing$ gives $\#\E_S(\vv)=0$ on $S$. Hence the coefficients of $e^{-2tL}$ in the expansions \eqref{ac-critical-equal-length-expansions} of $\mu_S(t)$ and $\mu_S'(t)$, together with \eqref{ac-critical-general-mass-error}, give
\begin{equation}\label{ac-critical-independent-mass}
 \begin{aligned}
 &\left|\frac{kq^2}{2}-\mu_S(t)-2q\Xi_S(tL-1)e^{-2tL}\right|\\
 &\quad+\left|\mu_S'(t)-2qL\Xi_S(2tL-3)e^{-2tL}\right|
 \leq kq C_r t e^{-t(2L+\ell_*)}.
 \end{aligned}
\end{equation}
If $\Xi_S\ne0$, choose
\[
 T_{\Xi,S}>\max\left\{T_{\rm loc},
 \frac1{\ell_*}\log\frac{3kC_r}{|\Xi_S|\min\{\ell_*,\ell_*^2\}}\right\}.
\]
For $t\geq T_{\Xi,S}$, we have $tL\geq3$ by the choice of $T_{\rm loc}$. Hence the signed leading terms in \eqref{ac-critical-independent-mass} satisfy
\[
 \begin{aligned}
 \frac{4q}{3}|\Xi_S|tLe^{-2tL}
 &\leq 2q|\Xi_S|(tL-1)e^{-2tL}
 \leq 2q|\Xi_S|tLe^{-2tL},\\
 2q|\Xi_S|L(2tL-3)e^{-2tL}
 &\geq 2q|\Xi_S|tL^2e^{-2tL}.
 \end{aligned}
\]
By \eqref{ac-critical-general-mass-error}, $\sigma_S(t)^2\leq kq e^{-2tL}$ and the choice of $T_{\Xi,S}$, $\abs{\mathcal E_S(t)}\leq(q/3)L|\Xi_S|te^{-2tL}$ and $\abs{\mathcal E_S'(t)}\leq(q/3)L^2|\Xi_S|te^{-2tL}$. This gives the mass and derivative estimates in part (iii). The energy expansion for $\E(S)=\varnothing$ in \eqref{ac-critical-equal-length-expansions} gives
\[
 \left|E(u_S(t),\G)-qL\Xi_St^3e^{-2tL}\right|
 \leq kqC_rt^3e^{-t(2L+\ell_*)}
 \leq\frac q3|\Xi_S|\ell_*t^3e^{-2tL}
 \qquad\text{for }t\geq T_{\Xi,S}.
\]
This proves the energy bound in part (iii). 
\end{proof}

The conclusions of Proposition~\ref{prop-ac-critical-sign-conditions} apply to a broad class of metric graphs. In particular, we give two examples illustrating how to determine the sign of the mass derivative $\mu_S'(t)$.

 \begin{remark}\label{rem-ac-critical-star}
 Let $q\geq3$ and let $\G=(\V,\E)$ be a star graph with center $\vv_0$. For some integer $1\leq m\leq q$, its vertex and edge sets are
\[
 \V=\{\vv_0,\vv_1,\ldots,\vv_m\},\qquad
 \E=\{\e_1,\ldots,\e_q\}.
\]
For $1\leq j\leq m$, the edge $\e_j$ is a bounded interval joining $\vv_0$ to $\vv_j$. For $m<j\leq q$, the edge $\e_j$ is a half-line starting at $\vv_0$. Assume that $\vv_0\in\V_0$ and set $S=\{\vv_0\}$. The reduction \eqref{ac-critical-separated-sum-reduction} in the proof of part~\textup{(i)} of Proposition~\ref{prop-ac-critical-sign-conditions} applies. Each bounded edge $\e$ has a degree-one endpoint other than $\vv_0$, so $\kappa_{\vv_0,\e}=-q$ by \eqref{ac-critical-local-coefficient}. With $k=1$ and $\ell_+$ the maximum length of the bounded edges, \eqref{ac-critical-grid-mass-expansion} and the error bound \eqref{ac-critical-general-mass-error} give
 \[
  \mu_S(t)>\frac{q^2}{2},\qquad \mu_S'(t)<0
  \qquad\text{for every }t\geq T_{\rm loc},
\]
with $T_{\rm loc}$ chosen as in the proof of Proposition~\ref{prop-ac-critical-sign-conditions}. Thus $\mu_S(t)$ decreases to $q^2/2$ from above. 
\end{remark}

\begin{remark}\label{rem-ac-critical-adjacent-supports}
On the metric grid $\mathcal G^d$, all vertices have degree $q=2d$ and all edges have length one.

Assume $d\geq2$ and let $S=\{\vv,\w\}\subset\V_0$ be a pair at combinatorial distance at least three. Every vertex of $N(S)$ has degree $q=2d\geq4$, so part~\textup{(i)} of Proposition~\ref{prop-ac-critical-sign-conditions} gives
\[
 \mu_S(t)<q^2,\qquad \mu_S'(t)>0
 \qquad\text{for }t\geq T_{\rm loc}.
\]
Thus the mass $\mu_S(t)$ approaches $q^2$ from below as $t \to +\infty$.

For every $d\geq1$, let $S=\{\vv,\w\}\subset\V_0$ be an adjacent pair. Then $k=2$ and $\#\E(S)=1$. Hence part~\textup{(ii)} of Proposition~\ref{prop-ac-critical-sign-conditions} gives
\[
 \mu_S(t)>q^2,\qquad \mu_S'(t)<0
 \qquad\text{for }t\geq T_{\rm loc}.
\]
Thus the mass $\mu_S(t)$ approaches $q^2$ from above as $t \to +\infty$.

On $\mathcal G^1$, let $S\subset\V_0$ consist of $k$ vertices with $\E(S)=\varnothing$. Since every vertex has degree two and $\#\E_S(\vv)\in\{1,2\}$ for $\vv\in N(S)$,
\[
 \begin{aligned}
 \Xi_S&=\sum_{\vv\in N(S)}\bigl[\#\E_S(\vv)-(\#\E_S(\vv))^2\bigr]\\
 &=-2\#\{\vv\in N(S):\#\E_S(\vv)=2\}.
 \end{aligned}
\]
Suppose that two vertices of $S$ have combinatorial distance two. Then $\Xi_S<0$, and part~\textup{(iii)} of Proposition~\ref{prop-ac-critical-sign-conditions}, with $q=2$ and $L=1$, gives
\[
 \mu_S(t)>2k,\qquad \mu_S'(t)<0
 \qquad\text{for }t\geq T_{\Xi,S}.
\]
Thus $\mu_S(t)$ approaches $2k$ from above as $t\to+\infty$. With $\eta_S:=\mu_S(T_{\Xi,S})-2k>0$, inverting this mass curve gives a positive solution $(u_{\mu,S},\lambda_{\mu,S})$ of \eqref{eqpd} for every $\mu\in(2k,2k+\eta_S)$.
If all pairwise combinatorial distances in $S$ are at least three, then $\Xi_S=0$ and the coefficient of $e^{-2t}$ in \eqref{ac-critical-equal-length-expansions} vanishes.
\end{remark}
We are now in a position to prove Theorem~\ref{th-ac-critical} and Theorem~\ref{cor-ac-critical-infinite-orbits}.
\begin{proof}[Proofs of Theorems~\ref{th-ac-critical} and \ref{cor-ac-critical-infinite-orbits}]
We first prove Theorem~\ref{th-ac-critical}. As in the proof of Theorem~\ref{th-ac-uniform-solutions}, we invert the mass curves on a common interval. Every vertex of $\G$ is incident to a bounded edge. Indeed, if a vertex were incident only to half-lines, connectedness would make $\G$ a star graph with those half-lines as its edges. Fix $W$ as in Theorem~\ref{th-ac-critical} and let $S\subset W$ with $\#S=k$. Then $S$ satisfies the hypotheses of part~\textup{(i)} of Proposition~\ref{prop-ac-critical-sign-conditions}, so \eqref{ac-critical-grid-two-sided-gap} and \eqref{ac-critical-grid-positive-slope-energy} hold for $t\geq T_{\rm loc}$. Proposition~\ref{prop-ac-critical-local-sign} gives the profile bound \eqref{ac-critical-grid-profile-bound}. The thresholds and constants depend only on $q$, $k$, $\ell_*$ and $\ell_+$, uniformly over all such $W$ and $S$. By Proposition~\ref{prop-ac-critical-branches}, the Morse index conclusions hold for $t\geq T_1$, where $T_1=T_1(q,k,\ell_*)>1$. Choose
\[
 T_{q,k}>\max\left\{T_{\rm loc},T_1,1,\frac{3(q-2)}{\ell_*},
                   \frac1{\ell_*}\log(3C_f\sqrt{kq})\right\}.
\]
By \eqref{ac-critical-weighted-size}, $\sigma_S(t)^2\leq kq e^{-2t\ell_*}$. The choice of $T_{q,k}$ gives $C_f\sigma_S(t)<1/3$ for $t\geq T_{q,k}$. To obtain a common mass interval, define
\begin{equation}\label{ac-critical-grid-common-mass-threshold}
 \eta_{q,k}:=\frac{kq}{3}\ell_*T_{q,k}e^{-2\ell_+T_{q,k}}.
\end{equation}
Since $q \geq 3$, $\ell_+\geq\ell_*$ and $x\mapsto xe^{-2x}$ decreases for $x\geq3$,
\[
 0<\ell_*T_{q,k}e^{-2\ell_+T_{q,k}}<3e^{-6}<1,
 \qquad
 \frac{\eta_{q,k}}{kq^2/2}<\frac{2e^{-6}}q<1.
\]
Also $\eta_{q,k}<kq/3<kq(q-2)$. Since each vertex of $S$ is incident to a bounded edge and every bounded edge incident at $S$ has length at most $\ell_+$, \eqref{ac-critical-weighted-size} gives
\[
 \sigma_S(t)^2
 \geq e^{-2t\ell_+}\sum_{\vv\in S}\#\{\e\in\E_b:\e\succ\vv\}
 \geq k e^{-2t\ell_+}
 \qquad\text{for every }t\geq T_{q,k}.
\]
By \eqref{ac-critical-grid-two-sided-gap} at $t=T_{q,k}$, $\mu_S(T_{q,k})\leq kq^2/2-\eta_{q,k}$. Proposition~\ref{prop-ac-critical-branches} gives $\lim_{t\to+\infty}\mu_S(t)=kq^2/2$. Since $\mu_S'(t)>0$ by \eqref{ac-critical-grid-positive-slope-energy},
\[
 \mu_S\bigl((T_{q,k},+\infty)\bigr)
 =\bigl(\mu_S(T_{q,k}),kq^2/2\bigr)
 \supset\bigl(kq^2/2-\eta_{q,k},kq^2/2\bigr).
\]
For each prescribed mass in this common interval $\bigl(kq^2/2-\eta_{q,k},kq^2/2\bigr)$, there is therefore a unique $t_S(\mu)>T_{q,k}$ with $\mu_S(t_S(\mu))=\mu$. Define $u_{\mu,S}:=u_S(t_S(\mu))$ and $\lambda_{\mu,S}:=t_S(\mu)^2$. By \eqref{ac-physical-branch}, \eqref{ac-branch-weak-identity} and Lemma~\ref{ac-IFT-branch}, $(u_{\mu,S},\lambda_{\mu,S})\in H^1_\mu(\G)\times(0,+\infty)$ is a positive solution of \eqref{eqpd}. By \eqref{ac-critical-grid-positive-slope-energy}, $E(u_{\mu,S},\G)>0$ and $\mu_S'(t_S(\mu))>0$. Proposition~\ref{prop-ac-critical-branches} then gives
\[
 m(u_{\mu,S})=k,\qquad\hat m_0(u_{\mu,S})=k-1.
\]
The profile bound \eqref{ac-critical-grid-profile-bound} and the choice of $T_{q,k}$ give
\[
 \bigl\|\lambda_{\mu,S}^{-1/4}\mathcal Ru_{\mu,S}-f_S^\infty\bigr\|_X
 =\norm{v_S(t_S(\mu))}_X<\frac13.
\]
This and the separation $\norm{f_{S_1}^\infty-f_{S_2}^\infty}_X\geq1$ for $S_1\ne S_2$ give pairwise distinctness by the argument in the proof of Theorem~\ref{th-ac-uniform-solutions}.

Finally, $\sigma_S(t)^2\geq k e^{-2t\ell_+}$, \eqref{ac-critical-grid-two-sided-gap} and $t_S(\mu) > T_{q,k}>3(q-2)/\ell_*$ give
\[
 \begin{aligned}
 \frac{kq^2}{2}-\mu
 &\geq\frac q3\ell_*t_S(\mu)\sigma_S(t_S(\mu))^2
 \geq q(q-2)\sigma_S(t_S(\mu))^2\\
 &\geq kq(q-2)e^{-2\ell_+t_S(\mu)}.
 \end{aligned}
\]
Since $kq^2/2-\mu<\eta_{q,k}<kq(q-2)$, taking logarithms gives $2\ell_+t_S(\mu)\geq\log[kq(q-2)/(kq^2/2-\mu)]>0$. Hence
\begin{equation}\label{ac-critical-grid-mass-parameter-estimates}
 \lambda_{\mu,S}\geq\frac1{4\ell_+^2}
 \left[\log\frac{kq(q-2)}{kq^2/2-\mu}\right]^2.
\end{equation}
This proves the multiplier limit in Theorem~\ref{th-ac-critical}.

For Theorem~\ref{cor-ac-critical-infinite-orbits}, fix $k\geq2$ and let $W_\infty$ be as in its hypotheses. By local finiteness, every ball of finite combinatorial radius contains finitely many vertices. Choose $F:=\{\vv_1,\ldots,\vv_{k-1}\}\subset W_\infty$ satisfying $\operatorname{dist}_{\rm c}(\vv_i,\vv_j)\geq3$ for all $1\leq i<j\leq k-1$. For every nonempty finite $A\subset\V$, define $\operatorname{diam}_{\rm c}A:=\max_{\vv,\w\in A}\operatorname{dist}_{\rm c}(\vv,\w)$. Choose $\w_n\in W_\infty$ and set $S_n:=F\cup\{\w_n\}$ such that
\[
 \min_{\vv\in F}\operatorname{dist}_{\rm c}(\vv,\w_1)\geq3,\qquad
 \min_{\vv\in F}\operatorname{dist}_{\rm c}(\vv,\w_{n+1})
 >\operatorname{diam}_{\rm c}S_n
 \quad\text{for every }n\geq1.
\]
Since $W_\infty$ is infinite, $\G$ has infinitely many vertices and is therefore not a star graph whose edges are all half-lines. Each $S_n$ has $k$ vertices and satisfies the assumptions of Theorem~\ref{th-ac-critical}. The choice of $\w_{n+1}$ also gives $\operatorname{diam}_{\rm c}S_{n+1}>\operatorname{diam}_{\rm c}S_n$. Fix $\mu\in(kq^2/2-\eta_{q,k},kq^2/2)$ and use the preceding construction with $S=W=S_n$. The same $\eta_{q,k}$ and $T_{q,k}$ apply to every $n$. We obtain $u_{\mu,n}:=u_{\mu,S_n}$ and $t_n:=t_{S_n}(\mu)>T_{q,k}$ with
\[
 \norm{f_{S_n}(t_n)-f_{S_n}^\infty}_X<\frac13.
\]
Arguing as in the proof of Theorem~\ref{cor-ac-infinite-orbits}, we obtain $S_n=\gamma S_m$ whenever $u_{\mu,n}=\gamma u_{\mu,m}$ for some $\gamma\in\Gamma$. Since $\gamma$ preserves combinatorial diameters and $\operatorname{diam}_{\rm c}S_{n+1}>\operatorname{diam}_{\rm c}S_n$, we have $n=m$. This proves the required $\Gamma$-inequivalence. The energy signs and Morse indices are those given by the construction in the proof of Theorem~\ref{th-ac-critical}. Finally, the lower bound for $\lambda_{\mu,S}$ in \eqref{ac-critical-grid-mass-parameter-estimates} applies with $S=S_n$ and is independent of $n$, so $\inf_{n\geq1}\lambda_{\mu,n}\to+\infty$ as $\mu\to(kq^2/2)^-$.
\end{proof}

\begin{remark}\label{rem-ac-critical-equal-length-normalized}
Under part~\textup{(ii)} of Proposition~\ref{prop-ac-critical-sign-conditions}, the mass curve is strictly decreasing to $kq^2/2$ for $t\geq T_{\rm loc}$. Set $\eta_S:=\mu_S(T_{\rm loc})-kq^2/2>0$. Inverting the mass curve gives a positive solution $(u_{\mu,S},\lambda_{\mu,S})$ of \eqref{eqpd} with mass $\mu$ for every $\mu\in(kq^2/2,kq^2/2+\eta_S)$.

Let $W\subset\V_0$ be finite with $\#W\geq k$. Suppose every $k$-element subset $S\subset W$ satisfies part~\textup{(iii)} of Proposition~\ref{prop-ac-critical-sign-conditions}, and the $\Xi_S$ have the same sign with $|\Xi_S|\geq\xi>0$. At a common threshold $T$, that part gives $|\mu_S(T)-kq^2/2|\geq q\ell_*\xi T e^{-2\ell_+T}$. The mass inversion and distinctness argument in the proof of Theorem~\ref{th-ac-critical} then gives $\binom{\#W}{k}$ distinct solutions on a common mass interval immediately below $kq^2/2$ if $\Xi_S>0$, and immediately above $kq^2/2$ if $\Xi_S<0$. In the latter case $E(u_{\mu,S},\G)<0$ and $\hat m_0(u_{\mu,S})=k$.
\end{remark}

\appendix
\section[Existence for every prescribed mass]{Existence for every prescribed mass}
\label{app-th1}

In this appendix, we prove the following existence result. The proof uses the minimax principle of Borthwick, Chang, Jeanjean and Soave \cite[Theorem~1.10]{Bort0}, which combines the monotonicity trick with approximate Morse index information on Palais--Smale sequences. In \cite{Bort}, this principle gives solutions of parameterized problems for almost every parameter, and a blow-up analysis is used to pass to the limit. For point defects, we instead bound the Lagrange multipliers directly, by testing equation \eqref{eqpdn}; see Lemma~\ref{lembu}.

\begin{theorem}\label{th1}
Let $\G=(\V,\E)\in\mathbf G$ contain at least one half-line, let $p>4$, and let $\varnothing\ne\V_0\subset\V$ be finite. Then, for every prescribed mass $\mu>0$, problem \eqref{eqpd} has a positive solution $(u,\lambda)\in H^1_\mu(\G)\times(0,+\infty)$ with constrained Morse index $\hat m_0(u)\leq1$, as defined in \eqref{ac-constrained-index-definition}.
\end{theorem}

\subsection{The parameterized problems}\label{subsec-app-parameterized}
Throughout Subsections~\ref{subsec-app-parameterized}--\ref{sect4}, $\G\in\mathbf G$ contains at least one half-line, $p>4$, and $\varnothing\ne\V_0\subset\V$ is finite, as in Theorem~\ref{th1}. The edge set may be infinite. Finiteness of $\V_0$ is used for the convergence of the nonlinear terms and for the multiplier estimate in Lemma~\ref{lembu}.

As in \cite{Bort}, for $\rho\in[1/2,1]$, we consider
\begin{equation}\label{eqpdrho}
       \begin{cases}
u'' = \lambda u & \text{on every }\e \in \E, \\
\int_\G\abs{u}^2\, dx = \mu & \\
\displaystyle\sum_{\e \succ \vv} u'_\e(\vv) = -\rho\abs{u(\vv)}^{p-2}u(\vv) & \text{at every } \vv \in \V_0,\\
\displaystyle\sum_{\e \succ \vv} u'_\e(\vv) = 0 & \text{at every } \vv \in \V \setminus \V_0.
\end{cases}
\end{equation}
A pair $(u,\lambda)\in H^1_\mu(\G)\times\R$ solves \eqref{eqpdrho} if and only if $u$ is a critical point of the energy functional
\begin{equation}\label{funcrho}
E_\rho(u,\G) := \frac{1}{2}\norm{u'}^{2}_2-\frac{\rho}{p}\sum_{\vv \in \V_{0}}\abs{u(\vv)}^p,\qquad u\in H^1(\G),
\end{equation}
constrained on $H^1_\mu(\G)$, with Lagrange multiplier $\lambda$. For $\rho=1$, $E_\rho(\cdot,\G)$ is the energy $E(\cdot,\G)$ in \eqref{func}. Since $p>4$, the maps $E_\rho'(\cdot,\G)$ and $E_\rho''(\cdot,\G)$ are Lipschitz continuous on bounded subsets of $H^1(\G)$, uniformly for $\rho\in[1/2,1]$. By \cite[Eq.~(2.7)]{Do0}, for every noncompact $\G\in\mathbf G$ the vertex values satisfy the Gagliardo--Nirenberg inequality
\begin{equation}\label{eqgn}
 \sum_{\vv\in\V_0}\abs{u(\vv)}^p
 \leq C_{p,\V_0}\norm{u'}_2^{p/2}\norm u_2^{p/2}
 \qquad\text{for }u\in H^1(\G),
 \qquad C_{p,\V_0}:=2^{p/2}\#\V_0.
\end{equation}

For $u\in H^1_\mu(\G)$, let $T_uH^1_\mu(\G)=\{v\in H^1(\G): (u,v)_2=0\}$ be the tangent space of $H^1_\mu(\G)$ at $u$. Following \cite[Definitions~1.3 and~1.4]{Bort0}, for $\theta\geq0$, the approximate Morse index of $u$ with respect to $E_\rho(\cdot,\G)$ is
\[
 \begin{aligned}
 \hat m_\theta(u):=\sup\bigl\{\dim L: \ &L \text{ is a subspace of } T_uH^1_\mu(\G) \text{ such that}\\
 &D^2E_\rho(u)[\varphi,\varphi]<-\theta\norm{\varphi}^2 \text{ for every } \varphi\in L\setminus\{0\}\bigr\},
 \end{aligned}
\]
where $D^2E_\rho(u):=E_\rho''(u,\G)-\mu^{-1}E_\rho'(u,\G)[u]\,(\cdot,\cdot)_2$. If $(u,\lambda)$ solves \eqref{eqpdrho}, then $E_\rho'(u,\G)[u]=-\lambda\mu$, so $D^2E_\rho(u)=E_\rho''(u,\G)+\lambda(\cdot,\cdot)_2$. In particular, for $\rho=1$, $\hat m_0(u)$ is the constrained Morse index in \eqref{ac-constrained-index-definition}.

\subsection{Boundedness of the sequence of solutions to the parameterized equations}\label{subsec-app-multipliers}
In this subsection, we establish a key compactness lemma for passing from the parameterized problems \eqref{eqpdrho} to the original problem \eqref{eqpd}. More precisely, if $\rho_n\to1^-$ and $\{(u_{\rho_n},\lambda_{\rho_n})\}\subset H^1(\G)\times(0,+\infty)$ is a sequence of normalized solutions of \eqref{eqpdrho} with bounded energy, then $\{u_{\rho_n}\}$ is a bounded Palais--Smale sequence for $E(\cdot,\G)$ constrained on $H^1_\mu(\G)$.
\begin{lemma}\label{lembu}
Let $\G\in\mathbf G$ contain a half-line, let $\varnothing\ne\V_0\subset\V$ be finite, let $p>4$, and let $\{\rho_n\}\subset[1/2,1]$ satisfy $\rho_n\to1^-$. Suppose that
\[\{(u_{\rho_n},\lambda_{\rho_n})\} \subset H^1(\G)\times (0,+\infty)\]
is a sequence of solutions  to 
\begin{equation}\label{eqpdn}
    \begin{cases} 
u_{\rho_n}'' = \lambda_{\rho_n} u_{\rho_n} & \text{on every }\e \in \E, \\ 
\int_\G\abs{u_{\rho_n}}^2\, dx = \mu & \\
\displaystyle\sum_{\e \succ \vv} u'_{\rho_n,\e}(\vv) = -\rho_n\abs{u_{\rho_n}(\vv)}^{p-2}u_{\rho_n}(\vv) & \text{at every } \vv \in \V_0,\\
 \displaystyle\sum_{\e \succ \vv} u'_{\rho_n,\e}(\vv) = 0 & \text{at every } \vv \in \V \setminus \V_0,
\end{cases}
\end{equation}
where $u_{\rho_n,\e}'(\vv)$ stands for $u_{\rho_n}|_\e'(0)$ or $-u_{\rho_n}|_\e'(\ell_\e)$, according to whether the vertex $\vv$ is identified with $0$ or $\ell_\e$. Assume additionally  that the sequence of energy levels $\{c_{\rho_n} := E_{{\rho_n}}(u_{\rho_n}, \G)\}$ is bounded. Then, $\{\lambda_{\rho_n}\} \subset (0,+\infty)$ is bounded from above. Moreover, up to a subsequence, $\{u_{\rho_n}\}$ satisfies the following properties.
	\begin{enumerate}[label=\rm(\roman*)]
		\item $E(u_{\rho_n},\G) \to c$, where  $\displaystyle c:= \lim_{n \to \infty}c_{\rho_n}$;\label{bps1}
		\item $E'(u_{\rho_n},\G)+\lambda_{\rho_n}(u_{\rho_n},\cdot)_2 \to 0 $ \text{in the dual of $H^1(\mathcal{G})$};
		\item $\{u_{\rho_n}\}$ is bounded in $ H^1(\G) $.\label{bps3}
	\end{enumerate}
\end{lemma}
\begin{proof}
Let $\{(u_{\rho_n},\lambda_{\rho_n})\} \subset H^1(\G)\times (0,+\infty)$ be a sequence of solutions to \eqref{eqpdn}, where $ (\rho_n)$ is a sequence in $[\frac{1}{2},1]$ satisfying $\rho_n \to 1^{-}$ as $n \to \infty$. We first show that  $\{\lambda_{\rho_n}\}$ cannot diverge to $+\infty$. Assume by contradiction that $\lambda_{\rho_n} \to +\infty$ as $n \to \infty$. By a straightforward computation, we have
\[
\begin{aligned}
    2c_{\rho_n}&=2E_{\rho_n}(u_{\rho_n},\G)-E_{\rho_n}'(u_{\rho_n} ,\G)[u_{\rho_n}]-\lambda_{\rho_n}(u_{\rho_n},u_{\rho_n})_2\\
    &=\frac{\rho_n(p-2)}{p}\sum_{\vv \in \V_0}\abs{u_{\rho_n}(\vv)}^p-\lambda_{\rho_n}\mu\\
    &\geq \frac{p-2}{2p}\sum_{\vv \in \V_0}\abs{u_{\rho_n}(\vv)}^p-\lambda_{\rho_n}\mu,
\end{aligned}
\]
and thus, writing $C_c:=\sup_n|c_{\rho_n}|<+\infty$ and $C_1:=[2p(2C_c+\mu)/(p-2)]^{1/p}$, we have, for every $n\in\mathbb{N}^+$ with $\lambda_{\rho_n}\geq1$,
\begin{equation}\label{eqc1lambdan0}
    \abs{u_{\rho_n}(\vv)} \leq \left(\sum_{\vv \in \V_0}\abs{u_{\rho_n}(\vv)}^p\right)^\frac{1}{p} \leq C_1\lambda_{\rho_n}^\frac{1}{p}, \quad \forall \vv \in \V_0.
\end{equation}
Set $t_n:=\sqrt{\lambda_{\rho_n}}$. Passing to a subsequence, we may assume that $t_n>\max\{1,3/\ell_*\}$ for every $n\in\mathbb{N}^+$. For a bounded edge $\e$, let $\vv_\e$ and $\w_\e$ be the vertices identified with $0$ and $\ell_\e$, and define
\[
 w_{n,\e}(x):=\mathbf1_{\V_0}(\vv_\e)u_{\rho_n}(\vv_\e)\bigl(1-t_nx\bigr)_+
 +\mathbf1_{\V_0}(\w_\e)u_{\rho_n}(\w_\e)\bigl(1-t_n(\ell_\e-x)\bigr)_+.
\]
For a half-line $\e$, with endpoint $\vv_\e$,
\[
 w_{n,\e}(x):=\mathbf1_{\V_0}(\vv_\e)u_{\rho_n}(\vv_\e)\bigl(1-t_nx\bigr)_+,
\]
where $(s)_+:=\max\{s,0\}$. Set $w_n:=w_{n,\e}$ on $I_\e$ and $w_n:=0$ on the rest of $\G$. Since $\sum_{\vv\in\V_0}\deg(\vv)<+\infty$ and $2/t_n<\ell_*\leq\ell_\e$, we know that $w_n\in H^1(\G)$ and $w_n(\vv)=u_{\rho_n}(\vv)$ for every $\vv\in\V_0$.

The weak form of \eqref{eqpdn}, with $\lambda_{\rho_n}=t_n^2$, is
\[
 \int_\G\bigl(u_{\rho_n}'\eta'+t_n^2u_{\rho_n}\eta\bigr)\,dx
 =\rho_n\sum_{\vv\in\V_0}\abs{u_{\rho_n}(\vv)}^{p-2}u_{\rho_n}(\vv)\eta(\vv)
 \qquad \text{for every }\eta\in H^1(\G).
\]
Taking $\eta=u_{\rho_n}-w_n$, we have $\eta(\vv)=0$ for every $\vv\in\V_0$, and hence
\[
 \int_\G\bigl(|u_{\rho_n}'|^2+t_n^2|u_{\rho_n}|^2\bigr)\,dx
 =\int_\G\bigl(u_{\rho_n}'w_n'+t_n^2u_{\rho_n}w_n\bigr)\,dx.
\]
Consequently, by the Cauchy--Schwarz inequality for $\int_\G\bigl(u_{\rho_n}'w_n'+t_n^2u_{\rho_n}w_n\bigr)\,dx$,
\begin{equation}\label{eq-localized-test-energy}
 \begin{aligned}
 t_n^2\mu
 &\leq\int_\G\bigl(|u_{\rho_n}'|^2+t_n^2|u_{\rho_n}|^2\bigr)\,dx\\
 &\leq\int_\G\bigl(|w_n'|^2+t_n^2|w_n|^2\bigr)\,dx\\
 &=\sum_{\vv\in\V_0}\deg(\vv)|u_{\rho_n}(\vv)|^2
     \int_0^{1/t_n}\bigl(t_n^2+t_n^2(1-t_nx)^2\bigr)\,dx\\
 &=\frac43t_n\sum_{\vv\in\V_0}\deg(\vv)|u_{\rho_n}(\vv)|^2.
 \end{aligned}
\end{equation}
Using \eqref{eqc1lambdan0} and $p>4$, we therefore obtain
\[
 0<\mu\leq\frac43 C_1^2
       \left(\sum_{\vv\in\V_0}\deg(\vv)\right)
       \lambda_{\rho_n}^{2/p-1/2}\longrightarrow0.
\] This contradiction proves that the multipliers are bounded above.

Next, we compute that 
\[
pc_{\rho_n}=pE_{\rho_n}(u_{\rho_n},\G)-E_{\rho_n}'(u_{\rho_n} ,\G)[u_{\rho_n}]-\lambda_{\rho_n}(u_{\rho_n},u_{\rho_n})_2=(\frac{p}{2}-1)\norm{u_{\rho_n}'}_2^2-\lambda_{\rho_n}\mu,
\]
which shows that $\{u_{\rho_n}\}$ is bounded in $H^1(\G)$. Thus,
\[
E(u_{\rho_n},\G)-E_{\rho_n}(u_{\rho_n},\G)\to 0  \quad \text{and}\quad E'(u_{\rho_n},\G)-E'_{\rho_n}(u_{\rho_n},\G)\to 0 \quad \text{as }n \to \infty.
\]
Therefore, up to a subsequence, $\{u_{\rho_n}\}$ is a bounded Palais--Smale sequence for $E(\cdot, \G)$ constrained on $H_\mu^1(\G)$, satisfying properties \ref{bps1}--\ref{bps3}.
\end{proof}

\subsection{Compactness of bounded Palais--Smale sequences}\label{sec3}
In this subsection, we assume that $\rho \in [\frac{1}{2},1]$. Let $\{u_n\} \subset H_\mu^1(\mathcal{G})$ be a bounded Palais--Smale sequence for $E_\rho(\cdot, \G)$ constrained on $H_\mu^1(\G)$ at some level $c>0$, that is, $\{u_n\}$ is bounded in $H^1(\mathcal{G})$ and satisfies
\begin{equation}\label{PS1}
E_\rho(u_n, \mathcal{G}) \to c 
\end{equation}
and 
\begin{equation}\label{PS2}
E^{ \prime}_{\rho}(u_n, \mathcal{G}) + \lambda_n (u_n, \cdot)_2 \to 0 \quad \text{in the dual of $H^1(\mathcal{G})$},
\end{equation}
where we choose  the sequence of Lagrange multipliers $\{\lambda_n\} \subset \R$  as
\[
 \lambda_n:=-\frac{E_\rho'(u_n,\G)[u_n]}{\mu}.
\]
Indeed, if $\{\widetilde\lambda_n\}$ is any multiplier sequence satisfying \eqref{PS2}, testing \eqref{PS2} with $u_n$ and using the boundedness of $\{u_n\}$ in $H^1(\G)$ gives
\[
 \widetilde\lambda_n-\lambda_n=\frac{E_\rho'(u_n,\G)[u_n]+\widetilde\lambda_n\mu}{\mu}\longrightarrow0.
\]
Since $\{u_n\}$ is bounded, \eqref{PS2} also holds with this choice of $\lambda_n$.

We now prove the compactness of the bounded Palais--Smale sequences $\{u_n\} \subset H^1(\mathcal{G})$.
\begin{lemma}\label{lemurhosolve}
There exists \( u_\rho \in H^1(\mathcal{G}) \) such that, up to a subsequence,
\[
u_n \rightharpoonup  u_\rho \quad \text{in } H^1(\mathcal{G})
\]
and
\[
u_n(x) \to u_\rho(x), \quad \forall x \in \G. 
\]
Moreover, $u_\rho$ is a solution of equation
\begin{equation}\label{eqpdnomu}
\begin{cases} 
u_{\rho}''=\lambda_\rho u_\rho &\quad \text{on every } \e\in \E,\\
\sum_{\e\succ\vv}u_{\rho,\e}'(\vv)=-\rho\abs{u_{\rho}(\vv)}^{p-2}u_{\rho}(\vv)  &\quad \text{at every } \vv\in\V_0,\\
\sum_{\e\succ\vv}u_{\rho,\e}'(\vv)=0  &\quad \text{at every } \vv\in \V \setminus \V_0,
\end{cases}
\end{equation}
where $\displaystyle \lambda_\rho=\lim_{n \to +\infty} \lambda_n $ and $u_{\rho,\e}'(\vv)$ stands for $u_\rho|_\e'(0)$ or $-u_\rho|_\e'(\ell_\e)$, according to whether the vertex $\vv$ is identified with $0$ or $\ell_\e$.
\end{lemma}
\begin{proof}
Since \(\{u_n\} \subset H^1(\mathcal{G})\) is bounded, by \cite[Theorem~8.8]{Ha},  we may assume that, up to a subsequence, there exists \( u_\rho \in H^1(\mathcal{G}) \) such that
\[
u_n \rightharpoonup  u_\rho \quad \text{in } H^1(\mathcal{G})
\]
and
\[
u_n(x) \to u_\rho(x), \quad \forall x \in \G.
\]
Observe also that, from \eqref{PS2} and the boundedness of \(\{u_n\} \subset H^1(\mathcal{G})\), we obtain
\[
   o(1)= E_{\rho}'(u_n ,\G)[u_n]+\lambda_n(u_n,u_n)_2 =E_{\rho}'(u_n ,\G)[u_n]+\lambda_n \mu,
\]
which implies that $\{\lambda_n\}$ is bounded. Thus, up to a subsequence, there exists \(\lambda_\rho \in \mathbb{R}\) such that 
\[
\lambda_n \rightarrow \lambda_\rho,\,\,\,\text{as}\,\,\,n \to +\infty.
\]
 Next, we show that $u_{\rho}\in H^{1}(\G)$ is a weak solution of \eqref{eqpdnomu}. Indeed, using \eqref{PS2} and  $\displaystyle \lim_{n \to +\infty} \lambda_n = \lambda_\rho$, we have, for every \(\eta\in H^{1}(\mathcal{G})\),
\begin{equation}\label{solution u}
\begin{split}
0&=\lim_{n\to\infty}\left(E^{ \prime}_{\rho}(u_{n},\mathcal{G})+\lambda_{n}(u_{n},\cdot)_2\right)[\eta]\\
&=\lim_{n\to\infty}\left(\int_{\mathcal{G}}u_{n}^{\prime}\eta^{ \prime}\,dx+\lambda_{n}\int_{\mathcal{G}}u_{n}\eta\,dx-\rho\sum_{\vv\in \V_0} \abs{u_{n}(\vv)}^{p-2 }u_{n}(\vv)\eta(\vv)\right)\\
&=\int_{\mathcal{G}}(u_{\rho})^{\prime}\eta^{\prime}\,dx+\lambda_{\rho}\int_{\mathcal{G}}u_{\rho}\eta\,dx-\rho\sum_{\vv\in \V_0}\abs{u_{\rho}(\vv)}^{p-2}u_{\rho}(\vv)\eta(\vv).
\end{split}
\end{equation}
This completes the proof of the lemma.
\end{proof}
Next, we focus on proving that $u_{n} \to u_\rho$ in $H^1(\G)$, as $n \to \infty$, to ensure that the limit \(u_\rho\) satisfies the mass constraint, that is, \(u_\rho\in H^1_\mu(\mathcal{G})\).

\begin{lemma}\label{lemurhoconverge}
 Assume that $\lambda_\rho \geq 0$ in  Lemma~\ref{lemurhosolve}. Then,
\begin{equation}\label{strong conv}
\int_{\mathcal{G}}|(u_{n}-u_{\rho})^{\prime}|^{2}\,dx+\lambda_{\rho}\int_{\mathcal{G}}|u_{n}-u_{\rho}|^{2}\,dx\to 0 \quad \text{as }n \to \infty.
\end{equation}
In particular, if \(\lambda_{\rho}>0\), the sequence \(\{u_{n}\}\) converges strongly to $u_{\rho}$ in \(H^{1}(\mathcal{G})\).
\end{lemma}

\begin{proof}
Using \eqref{PS2} and the fact that $\displaystyle\lim_{n \to +\infty} \lambda_n = \lambda_\rho$, for every \(\eta\in H^{1}(\mathcal{G})\), we have
\begin{align*}
o(1)\|\eta\|&=\left(E^{ \prime}_{\rho}(u_{n},\mathcal{G})+\lambda_{n}(u_{n},\cdot)_2\right)[\eta]\\
&= \int_{\mathcal{G}}u_{n}^{\prime}\eta^{\prime}\,dx-\rho\sum_{\vv\in \V_0}\abs{u_n(\vv)}^{p-2}u_n(\vv)\eta(\vv)+\lambda_{n}\int_{\mathcal{G}}u_{n}\eta\,dx\\
&=\int_{\mathcal{G}}u_{n}^{\prime}\eta^{\prime}\,dx-\rho\sum_{\vv\in \V_0}\abs{u_{\rho}(\vv)}^{p-2}u_{\rho}(\vv)\eta(\vv)+\lambda_{\rho}\int_{\mathcal{G}}u_{n}\eta\,dx+(\lambda_{n}-\lambda_{\rho})\int_{\mathcal{G}}u_{n}\eta\,dx+o(1)\|\eta\|,
\end{align*}
that is,
\begin{equation}\label{PS21}
\int_{\mathcal{G}}u_{n}^{\prime}\eta^{\prime}\,dx-\rho\sum_{\vv\in \V_0}\abs{u_{\rho}(\vv)}^{p-2}u_{\rho}(\vv)\eta(\vv)+\lambda_{\rho}\int_{\mathcal{G}}u_{n}\eta\,dx=o(1)\|\eta\|.
\end{equation}

Now, taking the difference between \eqref{PS21} and \eqref{solution u}, choosing \(\eta=\eta_{n}:=u_{n}-u_{\rho}\), we obtain
\begin{align*}
o(1)=o(1)\|u_{n}-u_{\rho}\|
&=\int_{\mathcal{G}}(u_{n}^{\prime}-u_{\rho}^{\prime})\eta_{n}^{\prime}\,dx+\lambda_{\rho}\int_{\mathcal{G}}(u_{n}-u_{\rho})\eta_{n}\,dx\\
&=\int_{\mathcal{G}}\big{|}(u_{n}-u_{\rho})^{\prime}\big{|}^{2 }\,dx+\lambda_{\rho}\int_{\mathcal{G}}|u_{n}-u_{\rho}|^{2}\,dx,
\end{align*}
which yields \eqref{strong conv}.
\end{proof}

\begin{lemma}\label{lem 5.2}
For any \(\lambda<0\) and \(d\in\mathbb{N}^+\), there exists a subspace \(Y\subset H^{1}(\mathcal{G})\) with \(\dim Y =d\) such that
\[
E''_\rho(u_n, \mathcal{G})[w, w] + \lambda \|w\|_2^2 \leq\frac{\lambda}{2}\|w\|^{2}, \qquad\forall\,w\in Y.
\]
\end{lemma}
\begin{proof}
Choose $\varphi\in C_c^\infty(\mathbb{R})$ such that $\operatorname{supp}\varphi\subset(0,1)$ and $\int_0^{+\infty}|\varphi|^2\,dx=1$. For \(t\in\mathbb{R}^{+}\), define the function \(\varphi^{t}\) by
\begin{align}\label{scal}
	\varphi^{t}(x):=t^{1/2}\varphi(tx).
\end{align}
We view $\varphi$ as a function in $H^1(\G)$ whose support is contained in a half-line identified with $[0,\infty)$ and define
\[
\varphi_{1}:=\varphi^{\tau},
\]
where $\tau>0$ is chosen small enough so that
\begin{equation}\label{5.9}
\tau^{2}\|\varphi^{\prime}\|_{L^{2}(\mathbb{R})}^{2}+\lambda\leq\frac{\lambda}{2}(\tau^{2}\|\varphi^{\prime}\|_{L^{2}(\mathbb{R})}^{2}+1).
\end{equation}
We have
\[
\|\varphi_{1}\|_{L^{2}({\G})}^{2}=1,\qquad\|\varphi^{\prime}_{1}\|_{L^{2}({\G})}^{2}=\tau^{2}\|\varphi^{\prime}\|_{L^{2}(\mathbb{R})}^{2}.
\]
For $i\geq 2$, define
\[
\varphi_{i}(x):=\varphi_{1}\left(x-\frac{i-1}{\tau}\right).
\]
Since $\operatorname{supp}(\varphi_{i})\subset\left(\frac{i-1}{\tau},\frac{i}{\tau}\right)$, the functions $\varphi_{i}$ have mutually disjoint supports. Let $Y=\text{span}\{\varphi_{1},\ldots,\varphi_{d}\}$. Every $w\in Y$ can be written as
\[
w=\sum_{i=1}^{d}\theta_{i}\varphi_{i},\,\,\,\theta_{i}\in \mathbb{R}.
\]
By a direct calculation, we have
\begin{align*}
\int_{\G}|w^{\prime}|^{2}\,dx+\lambda\int_{\G}|w|^{2}\,dx 
&= \tau^{2}\Bigg{(}\sum_{i=1}^{d}\theta_{i}^{2}\|\varphi^{\prime}\|_{L^{2}(\mathbb{R})}^{2}\Bigg{)}+\lambda\left(\sum_{i=1}^{d}\theta_{i}^{2}\right)\\
&=(\tau^{2}\|\varphi^{\prime}\|_{L^{2}(\mathbb{R})}^{2}+\lambda)\sum_{i=1}^{d}\theta_{i}^{2},
\end{align*}
and similarly, 
\[\|w\|^{2}=(\tau^{2}\|\varphi^{\prime}\|_{L^{2}(\mathbb{R})}^{2}+1)\sum_{i=1}^{d}\theta_{i}^{2}.\]
Therefore, \eqref{5.9} and $\operatorname{supp}(\varphi_{i})\subset\left(\frac{i-1}{\tau},\frac{i}{\tau}\right)$ imply that
\begin{align*}
\int_\mathcal{G} |w'|^2 \,dx+ \lambda\int_\mathcal{G} w^2 \,dx-(p-1)\rho\sum _{\vv\in \V_0}\abs{u_n(\vv)}^{p-2}w^{2}(\vv)
&= (\tau^{2}\|\varphi^{\prime}\|_{L^{2}(\mathbb{R})}^{2}+\lambda)\sum_{i=1}^{d}\theta_{i}^{2}\\
&\leq \frac{\lambda}{2}(\tau^{2}\|\varphi^{\prime}\|_{L^{2}(\mathbb{R})}^{2}+1)\sum_{i=1}^{d}\theta_{i}^{2}\\
&= \frac{\lambda}{2}\|w\|^{2}.
\end{align*}
This completes the proof of Lemma~\ref{lem 5.2}.
\end{proof}

\begin{lemma}\label{lamuna+}
Assume that there are an integer $N\geq0$, independent of $n$, and a sequence $\{\zeta_n\} \subset \mathbb{R}^+$ with $\zeta_n \to 0^+$ such that  $ \hat{\mathrm{m}}_{\zeta_n}(u_n) \leq N $ for every $n$, that is, if a subspace $W_n\subset T_{u_n}H^1_\mu(\mathcal{G})$ satisfies
\begin{equation}\label{5.4}
\int_\mathcal{G} |\varphi'|^2 \,dx+ \lambda_n\int_\mathcal{G} \varphi^2  \,dx-(p-1)\rho\sum _{\vv\in \V_0}\abs{u_n(\vv)}^{p-2}\varphi^{2}(\vv)= E''_\rho(u_n, \mathcal{G})[\varphi, \varphi] + \lambda_n \|\varphi\|_2^2 < -\zeta_n \|\varphi\|^2 
\end{equation}
 for every \(\varphi\in W_n\setminus\{0\}\), then \(\dim W_n\leq N\). Then $\lambda_\rho\geq0$.
\end{lemma}
\begin{proof}
    Observe that the codimension of $T_{u_{n}}H^{1}_{\mu}({\G})$ in $H^{1}({\G})$ is one. Therefore, if inequality \eqref{5.4} holds for all $\varphi\in W_{n}\setminus\{0\}$ in a subspace $W_{n}\subset H^{1}({\G})$, then  \(\dim W_n \leq  N+1 \). Let $\lambda<0$ and $Y$ be the subspace of dimension $d=N+2$ provided by Lemma~\ref{lem 5.2}. Applying \cite[Lemma~2.7]{Ca}, we conclude that
\begin{equation}\label{5.10}
\lambda_{\rho}\geq 0.
\end{equation}
\end{proof}

\subsection{Proof of Theorem~\ref{th1}}\label{sect4}
In this subsection, we prove Theorem~\ref{th1} by applying the abstract minimax principle of \cite[Theorem~1.10]{Bort0}. To ensure that the family $\{E_{\rho}(\cdot,\G)\}$ fits into the framework of \cite[Theorem~1.10]{Bort0}, we first show that it has a mountain pass geometry on $H^{1}_{\mu}(\G)$ uniformly with respect to $\rho\in [\frac{1}{2},1]$.
\begin{lemma}\label{lem3.1}
For every \(\mu>0\), there exist \(w_{1},w_{2}\in H_{\mu}^{1}(\G)\) with $w_1,w_2\geq 0$ independent of \(\rho\in\left[\frac{1}{2},1\right]\) such that  
\[
c_{\rho}:=\inf_{\gamma\in\mathcal P_\mu}\max_{t\in[0,1]}E_{\rho}(\gamma(t), \mathcal{G})>\max\{E_{\rho}(w_{1},\mathcal{G}),E_{\rho}(w_{2}, \mathcal{G})\},\qquad\forall\rho\in\left[\frac{1}{2},1\right],
\]
where \(\mathcal P_\mu:=\left\{\gamma\in C([0,1],H^{1}_{\mu}(\mathcal{G})): \gamma(0)=w_{1},\gamma(1)=w_{2}\right\}\).
\end{lemma}
\begin{proof}
For any \(\mu,k>0\), denote
\begin{equation*}
A_{\mu,k}:=\{u\in H^{1}_{\mu}(\mathcal{G}): \int_{\mathcal{G}}|u^{\prime}|^{2}\,dx<k\}, \,\,\,  \partial A_{\mu,k}:=\{u\in H^{1}_{\mu}(\mathcal{G}): \int_{\mathcal{G}}|u^{\prime}|^{2}\,dx=k\}.
\end{equation*}
By assumption, $\G$ contains a half-line. Choose $v\in C_c^\infty((0,\infty))$ with $v\geq0$ and $\|v\|_2^2=\mu$, regard it as a function supported on that half-line, and define $v^t(x)=t^{1/2}v(tx)$ for $t>0$. Then $\|v^t\|_2^2=\mu$ and $\|(v^t)'\|_2^2=t^2\|v'\|_2^2$. Hence $v^t\in A_{\mu,k}$ for all sufficiently small $t>0$, while the choice $t=\sqrt{k}/\|v'\|_2$ gives $v^t\in\partial A_{\mu,k}$. Therefore, both $A_{\mu,k}$ and $\partial A_{\mu,k}$ are nonempty.

Now, by the Gagliardo--Nirenberg inequality \eqref{eqgn}, for any $u\in H^{1}_{\mu}(\G)$, we obtain
\[
E_{\rho}(u,\G)\geq\frac{1}{2}\|u^{\prime}\|_2^2-\rho\frac{C_{p,\V_0} }{p}\mu^{\frac{p}{4}}\norm{u'}_2^{\frac{p}{2}}.
\]
Then, for any \(\mu>0\) and \(u\in\partial A_{\mu,k_{0}}\) with \(k_{0}:=\frac{1}{2}(\frac{p}{2C_{p,\V_0}})^{\frac{4}{p-4}}\mu^{-\frac{p}{p-4}}\), we have
\[
\inf_{u\in\partial A_{\mu,k_{0}}}E_{\rho}(u,\mathcal{G})\geq k_{0}\left(\frac{1}{2}-\frac{C_{p,\V_0}}{p}\mu^{\frac{p}{4}}k_{0}^{\frac{p-4}{4}}\right)=\colon\alpha>0,
\]
for every \(\rho\in[\frac{1}{2},1]\). Next, observe that for any \(u\in A_{\mu,k}\)
\[
E_{\rho}(u,\mathcal{G})\leq\frac{1}{2}\int_{\mathcal{G}}|u^{\prime}|^{2}\,dx\leq\frac{1}{2}k.
\]
Thus, recalling that \(A_{\mu,k}\neq\varnothing\) for all \(\mu,k>0\), it is possible to choose a function \(w_{1}\in A_{\mu,k}\)  with $w_1\geq 0$ for $k>0$ sufficiently small, such that
\begin{equation}\label{eq3.3}
\|w_{1}^{\prime}\|_2^2<k_{0}\quad\text{and}\quad E_{\rho}(w_{1},\mathcal{G})<\frac{\alpha}{2}\quad\forall\rho\in\left[\frac{1}{2},1\right].
\end{equation}

Next, take a vertex \(\vv \in \V_0\), which we identify in the following with $0$, and denote by \(\e_j\), \(j = 1, \ldots, K\), the edges incident at \(\vv\). Set
\begin{equation*}
\ell :=\begin{cases}
    \min_{1 \leq j \leq K} |\e_j|&\text{if }\e_j \text{ is bounded for some }j = 1, ..., K,\\
    1\quad &\text{ otherwise}.
\end{cases} 
\end{equation*}
Define \(u_{s_{0}}\in H^1(\mathcal{G})\) by
\[
	u_{s_0}(x) = 
	\begin{cases} 
		s_0(\frac{\ell}{K} - x) & \text{on }  \e_j \cap [0,\frac{\ell}{K}], \quad \forall \, j = 1, \ldots, K,\\
		0 & \text{elsewhere,}
	\end{cases}
\]
with $s_0>0$ chosen so that \(\|u_{s_{0}}\|_2^{2}=\mu\), and denote  $w:=u_{s_{0}}$. For \(t>1\), set \(w^t(x):=t^{1/2}w(tx)\). Then \(w^t\in H^1_\mu(\mathcal{G})\), and
\[
E_{\rho}(w^{t},\mathcal{G})= \frac{t^2}{2}\norm{w'}^2_2-\frac{\rho}{p} t^{p/2} \left ( \frac{s_0\ell}{K}  \right )^{p},
\]
for every \(\rho\in[\frac{1}{2},1]\). Since \(p>4\), $E_{\rho}(w^{t},\mathcal{G})\rightarrow -\infty$ as \(t\to+\infty\). Hence,  there exists a sufficiently large \(t_2>0\) such that,  for all \(t>t_{2}\) 
\[t^{2}\|w^{\prime}\|_2^2>2k_{0} \,\,\, \text{and}\,\,\,E_{\rho}(w^{t},\mathcal{G})<0,\,\,\,\rho\in[\frac{1}{2},1].\]
 We set \(w_{2}:=w^{2t_{2}}\geq 0 \), so that
\begin{equation}\label{eq3.4}
\|w_{2}^{\prime}\|_2^2>2k_{0}\quad\text{and}\quad E_{\rho}(w_{2},\mathcal{G})<0,\quad\forall\rho\in\left[\frac{1}{2},1\right].
\end{equation}
Finally, let \(\mathcal P_\mu\) and \(c_{\rho}\) be defined as in the statement of the lemma for our choice of \(w_{1}\) and \(w_{2}\); clearly, \(\mathcal P_\mu\neq \varnothing\), since the path
\[
\gamma_{0}(t):=\frac{\mu^{1/2}}{\|(1-t)w_{1}+tw_{2}\|_2}\left((1-t)w_{1}+tw_{2}\right),\quad  t\in[0,1]
\]
belongs to \(\mathcal P_\mu\). By continuity and \eqref{eq3.3}--\eqref{eq3.4}, for every \(\gamma \in \mathcal P_\mu\), there exists a constant \(t_\gamma \in [0,1]\) such that \(\gamma(t_\gamma) \in \partial A_{\mu,k_0}\). Therefore, for any \(\rho \in [\frac{1}{2},1]\), we have, for every \(\gamma \in \mathcal P_\mu\),
\[
\max_{t \in [0,1]} E_\rho(\gamma(t),\mathcal{G}) \geq E_\rho(\gamma(t_\gamma),\mathcal{G}) \geq \inf_{u \in \partial A_{\mu,k_0}} E_\rho(u,\mathcal{G}) \geq \alpha.
\]
It follows that \(c_\rho \geq \alpha\),  while
\[
\max\{E_\rho(w_1,\mathcal{G}),E_\rho(w_2,\mathcal{G})\} < \frac{\alpha}{2}
\]
by \eqref{eq3.3} and \eqref{eq3.4}, which completes the proof of the lemma.
\end{proof}

\begin{proof}[Proof of Theorem~\ref{th1}]
We use \cite[Theorem~1.10]{Bort0} with $E=H^1(\G)$, $H=L^2(\G)$, $I=[1/2,1]$ and $\Phi_\rho:=E_\rho(\cdot,\G)=A-\rho B$, where $A(u):=\frac12\norm{u'}_2^2$ and $B(u):=\frac1p\sum_{\vv\in\V_0}\abs{u(\vv)}^p\geq0$. Since $\norm u_2^2=\mu$ on $H^1_\mu(\G)$, the condition $\norm u\to+\infty$ implies $A(u)\to+\infty$. By Subsection~\ref{subsec-app-parameterized}, $\Phi_\rho'$ and $\Phi_\rho''$ are Lipschitz continuous on bounded sets, that is, $\alpha$-H\"older continuous on bounded sets with $\alpha=1$ in the sense of \cite[Definition~2.1]{Bort}. By Lemma~\ref{lem3.1} and \cite[Theorem~1.10 and Remarks~1.6 and~1.7]{Bort0}, we deduce that for almost every $\rho\in[\frac{1}{2},1]$, there exists a bounded sequence \(\{u_{\rho,n}\}_{n=1}^{\infty} \subset H^1_\mu(\mathcal{G})\) such that
\begin{equation*}
E_\rho(u_{\rho,n}, \mathcal{G}) \to c_\rho,
\end{equation*}
and
\begin{equation*}
E'_\rho(u_{\rho,n}, \mathcal{G}) + \lambda_{\rho,n} (u_{\rho,n}, \cdot)_2 \to 0 \quad \text{in the dual of } H^1(\mathcal{G}),
\end{equation*}
where
\begin{equation*}
\lambda_{\rho,n} := -\frac{1}{\mu} E'_\rho(u_{\rho,n}, \mathcal{G})[u_{\rho,n}]. 
\end{equation*}
Moreover, there exists a sequence \(\{\zeta_{\rho,n}\} \subset \mathbb{R}^+\) with \(\zeta_{\rho,n}\to 0^+\) such that, if 
\begin{equation*}
 E''_\rho(u_{\rho,n}, \mathcal{G})[\varphi, \varphi] + \lambda_{\rho,n} \|\varphi\|_2^2 < -\zeta_{\rho,n} \|\varphi\|^2 
\end{equation*}
holds for any \( \varphi \in W_n \setminus \{0\} \) in a subspace \( W_n \) of \( T_{u_{\rho,n}} H^1_\mu(\mathcal{G}) \), then   \(\dim W_n \leq  1 \). In addition, as noted in \cite[Remark~1.8]{Bort0}, since $u\in H^{1}_{\mu}(\mathcal{G})\implies|u|\in H^{1}_{\mu}(\mathcal{G})$, $w_{1},w_{2}\geq 0$ (see Lemma~\ref{lem3.1}), the map $u\mapsto|u|$ is continuous, and $E_{\rho}(u)=E_{\rho}(|u|)$, it is possible to choose $\{u_{\rho,n}\}$ with the property that $u_{\rho,n}\geq 0$ on $\mathcal{G}$. Moreover, the monotonicity of $c_\rho$ as a function of $\rho\in [\frac{1}{2},1]$ implies that $\{c_\rho\}$ is bounded, with $c_1\leq c_\rho \leq c_\frac{1}{2}$.

By Lemma~\ref{lamuna+}, up to a subsequence, $\displaystyle \lambda_{\rho}:= \lim_{n \to \infty}\lambda_{\rho,n}\geq 0$. Moreover, by  Lemmas~\ref{lemurhosolve} and \ref{lemurhoconverge}, there exists a solution $u_\rho\in H^1(\G)$ to \eqref{eqpdnomu} such that $u_{\rho,n}(x) \to u_\rho(x)$ for all $x \in \G$,  and
\begin{equation*}
\int_{\mathcal{G}}|(u_{\rho,n}-u_{\rho})^{\prime}|^{2}\,dx+\lambda_{\rho}\int_{\mathcal{G}}|u_{\rho,n}-u_{\rho}|^{2}\,dx \to 0, \quad \text{as }n \to \infty.
\end{equation*} 
Since $u_{\rho,n}\geq 0$ on $\mathcal{G}$, we have $u_\rho \geq 0$ on $\mathcal{G}$.

Suppose, by contradiction, that $u_{\rho}\equiv 0$. Then
\begin{equation}\label{3.17}
\int_{\mathcal{G}}|u_{\rho,n}^{\prime}|^{2}\,dx+\lambda_{\rho}\int_{\mathcal{G}}|u_{\rho,n}|^{2}\,dx\to 0,\quad \text{as }n \to \infty.
\end{equation}
If $\lambda_{\rho}>0$, this is impossible since $\|u_{\rho,n}\|^{2}_2=\mu>0$. If $\lambda_{\rho}=0$, then, by the Gagliardo--Nirenberg inequality \eqref{eqgn}, \eqref{3.17} contradicts the fact that  $c_{\rho}>0$. Hence, $u_{\rho}\not\equiv 0$.
By the vertex conditions and the uniqueness theorem for ODEs, we conclude that $u_{\rho}>0$ in $\mathcal{G}$. Indeed, suppose that there exists $x_0\in\G$ such that $u_\rho(x_0)=0$. If $x_0$ lies in the interior of an edge, then the strong
maximum principle implies that $u_\rho$ vanishes identically on that edge.
In particular, $u_\rho$ vanishes at an endpoint of that edge. Therefore, it
is enough to consider the case $x_0\in\V$. Since $u_\rho\geq0$ on $\G$ and $u_\rho(x_0)=0$, one has $u_{\rho,\e}'(x_0)\geq0$ for each $\e\succ x_0$,
where $(u_{\rho,\e})'(x_0)$ stands for $u_{\rho}|_\e'(0)$ or $-u_{\rho}|_\e'(\ell_\e)$, according to whether the vertex $x_0$ is identified with $0$ or $\ell_\e$. Moreover, by \eqref{eqpdnomu}, we obtain
 \[\sum_{\e\succ x_0}u_{\rho,\e}'(x_0)=0,
 \]
 and thus, $u_{\rho,\e}'(x_0)=0$ for all $\e \succ x_0$. Hence, by the uniqueness of solutions to second-order ODEs (see, e.g., {\cite[Theorem~10 in Section~3.1]{AD}}), $u_\rho \equiv 0$ on any edge containing $x_0$, and repeating this argument over every edge $\e \in \E$, we deduce that $u_{\rho}\equiv 0$ on $\mathcal{G}$, which leads to a contradiction.

Next, we claim that $\lambda_\rho > 0$. If $\lambda_\rho = 0$, consider an unbounded edge $\e_1 \in \E$ identified with $[0, +\infty)$. Then,  $u_\rho$ is $C^2$ on $[0, +\infty)$ and satisfies $u_\rho''(x) = 0$ for all $x \in (0, +\infty)$. It follows from \cite[Remark~4 in Section~3.3]{AD} that $u_\rho > 0$ on $[0, +\infty)$ and
\[
u_\rho(x) = u_\rho(0) + u_\rho'(0)x,
\]
which leads to a contradiction with $u_\rho \in L^2(\mathcal{G})$. Hence the claim holds. 

By Lemma~\ref{lemurhoconverge}, we conclude that $(u_\rho,\lambda_\rho) \in H^1(\G)\times (0,\infty)$ is a solution to equation \eqref{eqpdrho} and satisfies $E_{\rho}(u_{\rho},\G)=c_{\rho}$.

Finally, there exists a sequence $\{\rho_n\} \subset [\frac{1}{2},1]$ with $\rho_n \to 1^-$ such that each $(u_{\rho_n},\lambda_{\rho_n}) \in H^1(\G) \times (0, +\infty)$ solves 
\begin{align*}
\begin{cases} 
u_{\rho_n}'' = \lambda_{\rho_n} u_{\rho_n} & \text{on every }\e \in \E, \\ 
\int_\G\abs{u_{\rho_n} }^2\, dx = \mu & \\
\displaystyle\sum_{\e \succ \vv} u_{\rho_n,\e}' (\vv)= -\rho_n\abs{u_{\rho_n}(\vv)}^{p-2}u_{\rho_n}(\vv)  & \text{at every } \vv \in \V_0,\\
\displaystyle\sum_{\e \succ \vv} u_{\rho_n,\e}' (\vv)= 0  & \text{at every } \vv \in \V \setminus \V_0
\end{cases}
\end{align*}
with $E_{\rho_n}(u_{\rho_n},\G)=c_{\rho_n}$, where $u_{\rho_n,\e}'(\vv)$ stands for $u_{\rho_n}|_\e'(0)$ or $-u_{\rho_n}|_\e'(\ell_\e)$, according to whether the vertex $\vv$ is identified with $0$ or $\ell_\e$. In view of $c_1\leq c_{\rho_n}\leq c_{1/2}$ and Lemma~\ref{lembu}, $\{u_{\rho_{n}}\} \subset H^{1}(\G)$ is a bounded sequence satisfying the assumptions \eqref{PS1} and \eqref{PS2}. 
By repeating the previous arguments, we conclude that, up to a subsequence, $\displaystyle\lambda =\lim_{n \to \infty}\lambda_{\rho_{n}}>0$, $\{u_{\rho_{n}}\} $ converges to $u$, and $(u,\lambda)\in H_\mu^1(\G)\times(0,+\infty)$ is a positive solution of \eqref{eqpd}.

Arguing as in \cite[Remark~2.6]{Ca},  the solution $(u, \lambda)$ of Theorem~\ref{th1} can be chosen with $\hat m_0(u)\leq1$. Proposition~\ref{prop-ac-mp-single-vertex} identifies such solutions with single-vertex branches for all sufficiently small masses.
\end{proof}

\section[Proof of Proposition~\ref{prop-ac-mp-single-vertex} and examples on metric grids]{Proof of Proposition~\ref{prop-ac-mp-single-vertex} and examples on metric grids}\label{app-exact-count}
In this appendix, we prove Proposition~\ref{prop-ac-mp-single-vertex}. We also show, by examples on metric grids, that its growth assumption for $2<p<4$ cannot be dropped, that the restriction to large masses in Theorem~\ref{th2} cannot be removed for $2<p<4$, and that the half-line in Theorem~\ref{th1} cannot be replaced by the condition $\lambda(\G)=0$.
\subsection[Small-mass and large-mass solutions]{Small-mass and large-mass solutions}\label{app-small-mass-solutions}

The following lemma identifies the branch of a high-frequency sequence once its rescaled vertex restrictions converge in $X$.

\begin{lemma}\label{lem-ac-high-frequency-classification}
Let $\G\in\mathbf G$, $p>2$ and $\varnothing\ne\V_0\subset\V$. Let $\{(u_n,t_n^2)\}$ be a sequence of positive solutions of \eqref{eqpd} in $H^1(\G)\times(0,+\infty)$, with $t_n\to+\infty$ and
\[
 t_n^{-1/(p-2)}\mathcal Ru_n\longrightarrow f
 \qquad\text{in }X.
\]
Then $f=f_S^\infty$ for some nonempty finite $S\subset\V_0$, and there is $n_0 \in \mathbb{N}^+$ such that $u_n=u_S(t_n)$ for every $n\geq n_0$.
\end{lemma}

\begin{proof}
Set $f_n:=t_n^{-1/(p-2)}\mathcal Ru_n$. By \eqref{ac-vertex-equation}, $A_{t_n}f_n=\mathcal N(f_n)$. The estimate \eqref{ac-DtN-exponential} gives $A_{t_n}\to I_X$ in $\mathcal L(X)$, and \eqref{ac-N-remainder} implies that $\mathcal N:X\to X$ is continuous. Hence
\[
 f=\lim_{n\to\infty}A_{t_n}f_n
  =\lim_{n\to\infty}\mathcal N(f_n)=\mathcal N(f).
\]
Moreover, \eqref{ac-norm-ordering} and $u_n>0$ on $\G$ give, for every $\vv\in\V$,
\[
 \begin{aligned}
 f(\vv)&=\lim_{n\to\infty}t_n^{-1/(p-2)}u_n(\vv)\geq0.
 \end{aligned}
\]
Testing \eqref{ac-vertex-equation} with $f_n$ and using \eqref{ac-DtN-coercivity}, we obtain
\[
 \begin{aligned}
 \tanh(t_n\ell_*/2)\norm{f_n}_X^2
 &\leq\langle A_{t_n}f_n,f_n\rangle_X
 =\sum_{\vv\in\V_0}\abs{f_n(\vv)}^p\\
 &\leq\norm{f_n}_\infty^{p-2}\norm{f_n}_X^2.
 \end{aligned}
\]
Therefore $\norm{f_n}_\infty^{p-2}\geq\tanh(t_n\ell_*/2)\to1$. By \eqref{ac-norm-ordering},
\[
 \norm f_\infty=\lim_{n\to\infty}\norm{f_n}_\infty\geq1.
\]

Since $f\geq0$, the equation $f=\mathcal N(f)$ gives
\[
 f(\vv)=0\quad \text{for }\vv\notin\V_0,\qquad
 f(\vv)\in\{0,\deg(\vv)^{1/(p-2)}\}\quad \text{for }\vv\in\V_0.
\]
Thus $f=f_S^\infty$ for $S:=\{\vv\in\V_0:f(\vv)>0\}\ne\varnothing$. This set $S$ is finite, because
\[
 \#S\leq\sum_{\vv\in S}\deg(\vv)^{p/(p-2)}=\norm f_X^2<+\infty.
\]
Since $t_n\to+\infty$ and $f_n\to f_S^\infty$ in $X$, there is $n_0$ such that $t_n>T_S$ and $\norm{f_n-f_S^\infty}_X<\delta_S$ for every $n\geq n_0$. Lemma~\ref{ac-IFT-branch} gives $f_n=f_S(t_n)$ for every $n\geq n_0$. By Lemma~\ref{ac-restriction-extension} and \eqref{ac-physical-branch}, $u_n=t_n^{1/(p-2)}U_{t_n}f_n=u_S(t_n)$ for every $n\geq n_0$.
\end{proof}

\begin{lemma}\label{lem-ac-bounded-index-classification}
Let $\G\in\mathbf G$ be noncompact with $\lambda(\G)=0$, and let $\varnothing\ne\V_0\subset\V$ be finite.
\begin{enumerate}[label=\rm(\roman*)]
\item Let $p>2$. Every positive solution $(u,\lambda)$ of \eqref{eqpd} satisfies $\lambda\geq0$, and $\lambda>0$ if $\G$ contains a half-line. If $\lambda>0$, then $\hat m_0(u)\leq m(u)\leq\#\V_0$.
\item Let $p>4$, and let $M\geq1$ be an integer. There is $\mu_0=\mu_0(p,\G,\V_0,M)>0$ such that every positive solution $(u,\lambda)\in H^1_\mu(\G)\times\R$ of \eqref{eqpd} with $0<\mu<\mu_0$ and $\hat m_0(u)\leq M$ satisfies
\[
 (u,\lambda)=(u_{\mu,S},\lambda_{\mu,S})
 \quad\text{for some nonempty } S\subset\V_0,
 \qquad 1\leq\#S=\hat m_0(u)\leq M.
\]
\item Let $2<p<4$, let $M\geq1$ be an integer, and assume that there are $o\in\G$ and $A>0$ such that $\abs{\overline B_r(o)}\leq A(1+r^{4/p})$ for every $r\geq0$. There is $\mu_0=\mu_0(p,\G,\V_0,M)>0$ such that every positive solution $(u,\lambda)\in H^1_\mu(\G)\times\R$ of \eqref{eqpd} with $\mu>\mu_0$ and $\hat m_0(u)\leq M$ satisfies
\[
 (u,\lambda)=(u_{\mu,S},\lambda_{\mu,S})
 \quad\text{for some nonempty } S\subset\V_0,
 \qquad 1\leq\#S=\hat m_0(u)+1\leq M+1.
\]
\end{enumerate}
Here $(u_{\mu,S},\lambda_{\mu,S})$ is the solution indexed by $S$ in Theorem~\ref{th2}.
\end{lemma}%
\begin{proof}
In this proof, $\G\in\mathbf G$ is noncompact with $\lambda(\G)=0$. Recall that the inequality \eqref{eqgn} holds for every noncompact $\G\in\mathbf G$, and that the cutoff functions in the proof of Lemma~\ref{lembu} are defined on every $\G\in\mathbf G$.

We first prove \textup{(i)}. We adapt the argument in the proof of \cite[Lemma~2.5(a)]{Do0}, where the same bound is proved when the nonlinearity acts on the whole graph. Let $(u,\lambda)$ be a positive solution of \eqref{eqpd}, and put $c_\vv:=u(\vv)^{p-2}$ for $\vv\in\V_0$. Let $\varphi\in H^1(\G)$ have compact support, and put $g:=\varphi/u$. Since $u$ is continuous and positive, it is bounded below by a positive constant on $\operatorname{supp}\varphi$, so $g\in H^1(\G)$ and $\operatorname{supp}g=\operatorname{supp}\varphi$. We have $\varphi'=u'g+ug'$ and $(u'g)^2+2uu'gg'=u'(ug^2)'$ on every edge, so
\[
 \norm{\varphi'}_2^2=\int_\G u'(ug^2)'\,dx+\int_\G u^2\abs{g'}^2\,dx .
\]
 On every edge $\e$, $u''=\lambda u$ and $u^2g^2=\varphi^2$, so integrating $u'(ug^2)'$ by parts on $\e$ gives
\[
 \int_\e u'(ug^2)'\,dx=-\lambda\int_\e\varphi^2\,dx-\sum_{\vv\in\V,\,\e\succ\vv}u(\vv)g(\vv)^2u'_\e(\vv).
\]
Summing these identities over all edges, we obtain
\[
 \int_\G u'(ug^2)'\,dx=-\lambda\norm{\varphi}_2^2-\sum_{\vv\in\V}u(\vv)g(\vv)^2\sum_{\e\succ\vv}u'_\e(\vv)=-\lambda\norm{\varphi}_2^2+\sum_{\vv\in\V_0}c_\vv\abs{\varphi(\vv)}^2,
\]
where the last equality follows from the vertex conditions in \eqref{eqpd}. Hence
\begin{equation}\label{ac-ground-state-identity}
 \norm{\varphi'}_2^2+\lambda\norm{\varphi}_2^2-\sum_{\vv\in\V_0}c_\vv\abs{\varphi(\vv)}^2=\int_\G u^2\abs{g'}^2\,dx\geq0 .
\end{equation}
In particular, $\lambda\norm{\varphi}_2^2\geq-\norm{\varphi'}_2^2$. Since functions with compact support are dense in $H^1(\G)$ (see, e.g., \cite[Lemma~2.4]{DAHJ}), this inequality holds for every $\varphi\in H^1(\G)$, and hence $\lambda\geq-\lambda(\G)=0$. If $\G$ contains a half-line $\e$, then $u$ solves $u''=\lambda u$ on $\e$, and the solution formulas in \cite[Remark~4 in Section~3.3]{AD} show that, for $\lambda\leq0$, the only solution in $L^2(\e)$ is $u|_\e\equiv0$. Thus $\lambda>0$ in this case.

Now let $\lambda>0$, and write $\lambda=t^2$ with $t>0$. Let $L\subset H^1(\G)$ be a subspace such that $Q_u[\varphi,\varphi]<0$ for every $\varphi\in L\setminus\{0\}$, where $Q_u$ is defined in \eqref{ac-Q-definition}. If $\varphi\in L$ satisfies that $\varphi(\vv_0)=0$ for every $\vv_0\in\V_0$, then
\[
 Q_u[\varphi,\varphi]=\norm{\varphi'}_2^2+t^2\norm{\varphi}_2^2\geq0,
\]
so $\varphi=0$. Hence the linear map $L\to\R^{\V_0}$, $\varphi\mapsto(\varphi(\vv))_{\vv\in\V_0}$, is injective, and $\dim L\leq\#\V_0$. Taking the supremum over all such $L$ gives $m(u)\leq\#\V_0$. Finally, $\hat m_0(u)\leq m(u)$ by the definitions \eqref{ac-free-index-definition} and \eqref{ac-constrained-index-definition}.

We next prove \textup{(ii)}.
Fix an integer $M\geq1$, and set $D:=\max_{\vv\in\V_0}\deg(\vv)$. Let $\mu_{M,D}$ and $T_{M,D}$ be as in Theorem~\ref{th-ac-uniform-solutions}. Suppose that the conclusion fails for arbitrarily small masses: for each $n\geq1$, choose a solution $(u_n,\lambda_n)\in H^1_{\mu_n}(\G)\times\R$ of \eqref{eqpd} satisfying
\[
 0<\mu_n<\min\{1/n,\mu_{M,D}\},\qquad
 u_n>0\text{ on }\G,\qquad \hat m_0(u_n)\leq M,
\]
and, for every nonempty $S\subset\V_0$ with $\#S\leq M$,
\[
 (u_n,\lambda_n)\ne(u_{\mu_n,S},\lambda_{\mu_n,S})
 \qquad\text{or}\qquad
 \hat m_0(u_n)\ne\#S.
\]

By \textup{(i)}, $\lambda_n\geq0$. We first prove $\lambda_n\to+\infty$. Set $a_n:=\norm{u_n'}_2^2>0$. Testing the equation \eqref{eqpd} with $u_n$ and using \eqref{eqgn} gives
\[
 0<a_n\leq a_n+\lambda_n\mu_n
 =\sum_{\vv\in\V_0}\abs{u_n(\vv)}^p
 \leq C_{p,\V_0}\mu_n^{p/4}a_n^{p/4}.
 \]
Since $p>4$ and $\mu_n\to0^+$, it follows that
\begin{equation}\label{ac-small-mass-kinetic-lower-bound}
 a_n\geq C_{p,\V_0}^{-4/(p-4)}\mu_n^{-p/(p-4)}\longrightarrow+\infty\quad\text{as }n\to+\infty.
\end{equation}
If a subsequence of $\{u_n\}$ satisfies $\lambda_n\leq L<+\infty$, put $r:=\ell_*/3$ and define $w_n$ as in the proof of Lemma~\ref{lembu}, replacing $u_{\rho_n}$ by $u_n$ and the cutoff width $1/t_n$ by $r$. Since $2r<\ell_*$, the endpoint segments are disjoint and $u_n-w_n$ vanishes on $\V_0$. The energy comparison in \eqref{eq-localized-test-energy}, applied to this cutoff, gives
\[
 \begin{aligned}
 a_n
 &\leq a_n+\lambda_n\mu_n
 \leq\norm{w_n'}_2^2+\lambda_n\norm{w_n}_2^2\\
 &=(r^{-1}+\lambda_nr/3)\sum_{\vv\in\V_0}\deg(\vv)|u_n(\vv)|^2\\
 &\leq2(r^{-1}+Lr/3)\sqrt{\mu_na_n}\sum_{\vv\in\V_0}\deg(\vv),
 \end{aligned}
\]
Here the last line uses \cite[Eq.~(2.7)]{Do0}, namely
\[
|u_n(\vv)|^2\leq\norm{u_n}_\infty^2\leq2\norm{u_n}_2\norm{u_n'}_2=2\sqrt{\mu_na_n}\qquad\text{for every }\vv\in\V.
\]
Dividing by $\sqrt{a_n}$ gives $\sqrt{a_n}\to0$ because $\mu_n\to0$ and the degree sum over $\V_0$ is finite. This contradicts \eqref{ac-small-mass-kinetic-lower-bound}. Thus, $\lambda_n\to+\infty$.

Put $t_n:=\sqrt{\lambda_n}$ and $f_n:=t_n^{-1/(p-2)}\mathcal Ru_n$. By \eqref{ac-vertex-equation}, $A_{t_n}f_n=\mathcal N(f_n)$. For $\vv\in\V_0$, \eqref{ac-DtN-coordinate} and $f_n>0$ give
\begin{equation}\label{ac-small-mass-vertex-pointwise-bound}
 \begin{aligned}
 f_n(\vv)^{p-1}
 &=\deg(\vv)(A_{t_n}f_n)(\vv)\\
 &=\sum_{\e\in\E_b,\,\e\succ\vv}
   \bigl[\coth(t_n\ell_\e)f_n(\vv)
    -\operatorname{csch}(t_n\ell_\e)f_n(\w_{\e,\vv})\bigr]
    +\sum_{\e\in\E_\infty,\,\e\succ\vv}f_n(\vv)\\
 &\leq\left[\sum_{\e\in\E_b,\,\e\succ\vv}\coth(t_n\ell_\e)
           +\sum_{\e\in\E_\infty,\,\e\succ\vv}1\right]f_n(\vv)\\
 &\leq\deg(\vv)\coth(t_n\ell_*)f_n(\vv).
 \end{aligned}
\end{equation}
Dividing by $f_n(\vv)>0$ gives $f_n(\vv)^{p-2}\leq\deg(\vv)\coth(t_n\ell_*)$. With $P_{\V_0}$ as in \eqref{ac-coordinate-subspace}, the support of $\mathcal N(f_n)$ gives
\[
 (I_X-P_{\V_0})A_{t_n}f_n=(I_X-P_{\V_0})\mathcal N(f_n)=0,
 \]
 which implies
 \[
 (I_X-P_{\V_0})f_n
 =(I_X-P_{\V_0})(I_X-A_{t_n})f_n.
 \]
 For every $n$ with $t_n\geq1$ and $C_0e^{-t_n\ell_*/2}\leq1/2$, \eqref{ac-DtN-exponential} and $\norm{I_X-P_{\V_0}}_{\mathcal L(X)}\leq1$ give
\begin{equation}\label{ac-small-mass-complement-bound}
 \norm{(I_X-P_{\V_0})f_n}_X
 \leq C_0e^{-t_n\ell_*/2}\norm{f_n}_X
 \leq\frac12\norm{f_n}_X,
\end{equation}
which leads to
\[
 \norm{f_n}_X
 \leq\norm{P_{\V_0}f_n}_X+\norm{(I_X-P_{\V_0})f_n}_X
 \leq\norm{P_{\V_0}f_n}_X+\frac12\norm{f_n}_X.
\]
Hence $\norm{f_n}_X\leq2\norm{P_{\V_0}f_n}_X$.
 Moreover, by \eqref{ac-small-mass-vertex-pointwise-bound}, we have
 \[
 \norm{P_{\V_0}f_n}_X^2
 \leq\coth(t_n\ell_*)^{2/(p-2)}
               \sum_{\vv\in\V_0}\deg(\vv)^{p/(p-2)}.
\]
Hence $(f_n)$ is bounded. Since $t_n\to+\infty$, \eqref{ac-small-mass-complement-bound} gives $(I_X-P_{\V_0})f_n\to0$ in $X$. Since $\V_0$ is finite, a subsequence of $P_{\V_0}f_n$ converges in $X$, and so does $f_n$.

Lemma~\ref{lem-ac-high-frequency-classification} now gives a nonempty finite $S\subset\V_0$ such that $f_n\to f_S^\infty$ and $u_n=u_S(t_n)$ for every sufficiently large $n$.
Since $p>4$, Proposition~\ref{prop-ac-branches} gives
\[
 1\leq\#S=\hat m_0(u_n)\leq M
\]
for all sufficiently large $n$. For such $n$ with $t_n>T_{M,D}$, we have $\mu_S(t_n)=\mu_n<\mu_{M,D}$. The uniqueness of the mass parameter in Theorem~\ref{th-ac-uniform-solutions} gives $t_n=t_S(\mu_n)$ and hence $(u_n,\lambda_n)=(u_{\mu_n,S},\lambda_{\mu_n,S})$.

This contradicts the choice of $(u_n,\lambda_n)$ and proves \textup{(ii)}.

We next prove \textup{(iii)}. Let $2<p<4$.
Choose $o\in\G$ and $A>0$ as in the assumption of \textup{(iii)}. Since $\V_0$ is nonempty and finite, set
\[
 R_0:=\max_{\vv\in\V_0}\operatorname{dist}(o,\vv)<+\infty.
\]
For every $r\geq0$ and $x\in\G$ with $\operatorname{dist}(x,\V_0)\leq r$, there is $\vv\in\V_0$ such that $\operatorname{dist}(x,\vv)=\operatorname{dist}(x,\V_0)$, and the triangle inequality gives
\[
 \operatorname{dist}(x,o)
 \leq\operatorname{dist}(x,\vv)+\operatorname{dist}(\vv,o)
 \leq r+R_0.
\]
Since $4/p>1$, for every $r\geq0$ we obtain
\[
 \begin{aligned}
 \abs{\{x\in\G:\operatorname{dist}(x,\V_0)\leq r\}}
 &\leq\abs{\overline B_{r+R_0}(o)}
 \leq A\bigl[1+(r+R_0)^{4/p}\bigr]\\
 &\leq A\bigl[1+2^{4/p-1}(r^{4/p}+R_0^{4/p})\bigr]
 \leq A_0(1+r^{4/p}),
 \end{aligned}
\]
where $A_0:=A\bigl[1+2^{4/p-1}(1+R_0^{4/p})\bigr]$ depends only on $p$, $A$ and $R_0$. 

We first show that every positive solution $(u,\lambda)$ of \eqref{eqpd} satisfies $\lambda>0$ and \eqref{ac-large-mass-multiplier-bound} below. Put $d(x):=\operatorname{dist}(x,\V_0)$, so that $\abs{\{d\leq r\}}\leq A_0(1+r^{4/p})$ for every $r\geq0$. Let $(u,\lambda)\in H^1_\mu(\G)\times\R$ be a positive solution of \eqref{eqpd}. For $R>0$, the function $\varphi_R:=\min\{1,\max\{0,2-d/R\}\}$ belongs to $H^1(\G)$, equals $1$ on $\V_0$, vanishes outside $\{d\leq2R\}$ and satisfies $\abs{\varphi_R'}\leq R^{-1}$. Testing \eqref{eqpd} with $\varphi_R$ gives
\[
 \sum_{\vv\in\V_0}u(\vv)^{p-1}=\int_\G u'\varphi_R'\,dx+\lambda\int_\G u\varphi_R\,dx ,
\]
and $\abs{\int_\G u'\varphi_R'\,dx}\leq R^{-1}(A_0(1+(2R)^{4/p}))^{1/2}\norm{u'}_2\to0$ as $R\to+\infty$ because $4/p<2$. Since the left-hand side is positive, $\lambda\ne0$, so $\lambda>0$ by \textup{(i)}, and monotone convergence gives 
\begin{equation}\label{ac-mass-integral-identity}
 \int_\G u\,dx=\lambda^{-1}\sum_{\vv\in\V_0}u(\vv)^{p-1}.
\end{equation}
Write $\lambda=t^2$ with $t>0$. For $R>0$, put $w_R:=e^{t\min\{d,R\}/2}$. Then $uw_R^2\in H^1(\G)$, $w_R=1$ on $\V_0$ and $\abs{w_R'}\leq tw_R/2$. Testing \eqref{eqpd} with $uw_R^2$ and using $u'(uw_R^2)'=((uw_R)')^2-u^2(w_R')^2$, we obtain $\frac34t^2\int_\G u^2w_R^2\,dx\leq\sum_{\vv\in\V_0}u(\vv)^p$. Letting $R\to+\infty$ gives
\[
 \frac34t^2\int_\G u^2e^{td}\,dx\leq\sum_{\vv\in\V_0}u(\vv)^p .
\]
Since $u\in H^1(\G)$ and $\ell_*>0$, for every $\varepsilon>0$ the set $\{u\geq\varepsilon\}$ meets only finitely many edges and is bounded, so the maximum $m:=\max_\G u$ is attained. Since $u''=t^2u>0$ on every edge, $m$ is attained at a vertex $\vv_*$. If $\vv_*\notin\V_0$, then $u'_\e(\vv_*)\leq0$ for every $\e\succ\vv_*$, and these derivatives add up to zero by the Kirchhoff condition in \eqref{eqpd}, so $u_\e'(\vv_*)=0$ for every $\e \succ \vv_*$, however, since $u''>0$, this gives $u>m$ near $\vv_*$, which is impossible. Thus $\vv_*\in\V_0$, and therefore $\sum_{\vv\in\V_0}u(\vv)^p\leq m\sum_{\vv\in\V_0}u(\vv)^{p-1}=mt^2\int_\G u\,dx$ by \eqref{ac-mass-integral-identity}, which leads to $ \int_\G u^2e^{td}\,dx \leq \frac43 m\int_\G u\,dx$. By the Cauchy--Schwarz inequality and
\[
 \int_\G e^{-td}\,dx=\int_\G\Bigl(\int_{d(x)}^{+\infty}te^{-tr}\,dr\Bigr)dx=t\int_0^{+\infty}e^{-tr}\abs{\{d\leq r\}}\,dr\leq A_0\bigl(1+2t^{-4/p}\bigr),
\]
where the last inequality uses $t\int_0^{+\infty}e^{-tr}\,dr=1$ and, by H\"older's inequality, $t\int_0^{+\infty}e^{-tr}r^{4/p}\,dr=t^{-4/p}\int_0^{+\infty}e^{-s}s^{4/p}\,ds\leq t^{-4/p}\bigl(\int_0^{+\infty}e^{-s}s^2\,ds\bigr)^{2/p}=2^{2/p}t^{-4/p}<2t^{-4/p}$, we have
\[
 \Bigl(\int_\G u\,dx\Bigr)^2\leq\int_\G u^2e^{td}\,dx\int_\G e^{-td}\,dx\leq\frac43 A_0\bigl(1+2t^{-4/p}\bigr) m\int_\G u\,dx .
\]
Hence $\int_\G u\,dx\leq\frac43 A_0(1+2t^{-4/p})m$, and by \eqref{ac-mass-integral-identity}, 
\[
 \begin{aligned}
 m^{p-1}&\leq\sum_{\vv\in\V_0}u(\vv)^{p-1}=t^2\int_\G u\,dx\leq\frac43 A_0(t^2+2t^{2-4/p})m,\\
 \mu&\leq\int_\G u^2e^{td}\,dx\leq\frac43m\int_\G u\,dx\leq\frac{16}9 A_0\bigl(1+2t^{-4/p}\bigr)m^2 .
 \end{aligned}
\]
Put $C_{A_0,p}:=\frac{16}3A_0(4A_0)^{2/(p-2)}$.
Eliminating $m$ from these bounds and using $1+2t^{-4/p}=t^{-4/p}(t^{4/p}+2)$ and $t^2+2t^{2-4/p}=t^{2-4/p}(t^{4/p}+2)$, we obtain
\[
 \mu\leq\frac{16}9A_0\Bigl(\frac43A_0\Bigr)^{2/(p-2)}\bigl(t^{4/p}+2\bigr)^{p/(p-2)}.
\]
If $t\leq1$, then $t^{4/p}+2\leq3$, so $\mu\leq C_{A_0,p}$. If $t>1$, then $t^{4/p}+2<3t^{4/p}$, so $\mu<C_{A_0,p}t^{4/(p-2)}$.%
Consequently,
\begin{equation}\label{ac-large-mass-multiplier-bound}
 \lambda> C_{A_0,p}^{-(p-2)/2}\mu^{(p-2)/2}\qquad\text{for }\mu>C_{A_0,p} .
\end{equation}

Let $D:=\max_{\vv\in\V_0}\deg(\vv)$, and let $T_{M+1,D}$ and $\mu_{M+1,D}$ be as in Theorem~\ref{th-ac-uniform-solutions}. Suppose that the conclusion of \textup{(iii)} fails for arbitrarily large masses: for each $n\geq1$, there is a positive solution $(u_n,\lambda_n)\in H^1_{\mu_n}(\G)\times\R$ of \eqref{eqpd} with $\mu_n>\max\{n,\mu_{M+1,D}\}$ and $\hat m_0(u_n)\leq M$ such that, for every nonempty $S\subset\V_0$ with $\#S\leq M+1$, either $(u_n,\lambda_n)\ne(u_{\mu_n,S},\lambda_{\mu_n,S})$ or $\hat m_0(u_n)\ne\#S-1$. By \eqref{ac-large-mass-multiplier-bound}, $\lambda_n\to+\infty$. 
 After passing to a subsequence, repeating the argument in the proof of \textup{(ii)} gives a nonempty $S\subset\V_0$ and an integer $n_0$ such that $u_n=u_S(t_n)$ for every $n\geq n_0$. Since $2<p<4$, Proposition~\ref{prop-ac-branches} gives $1\leq\#S=\hat m_0(u_n)+1\leq M+1$ for all sufficiently large $n \in \mathbb{N}^+$. For such $n$ with $t_n>T_{M+1,D}$, we have $\mu_S(t_n)=\mu_n>\mu_{M+1,D}$, and the uniqueness of the mass parameter in Theorem~\ref{th-ac-uniform-solutions} gives $t_n=t_S(\mu_n)$. Hence $(u_n,\lambda_n)=(u_{\mu_n,S},\lambda_{\mu_n,S})$ and $\hat m_0(u_n)=\#S-1$, a contradiction.
\end{proof}
Now, we are in a position to prove Proposition~\ref{prop-ac-mp-single-vertex}.
\begin{proof}[Proof of Proposition~\ref{prop-ac-mp-single-vertex}]
Put $q:=\#\V_0$. We begin with the case $p>4$ and show that there is $\mu_1=\mu_1(p,\G,\V_0)>0$ such that every positive solution $(u,\lambda)\in H^1_\mu(\G)\times\R$ of \eqref{eqpd} with $0<\mu<\mu_1$ satisfies $\lambda>0$. Indeed, if there are such solutions $(u_n,\lambda_n)$ with $\mu_n\to0^+$ and $\lambda_n\leq0$, then, the first step of the proof of Lemma~\ref{lem-ac-bounded-index-classification}\textup{(ii)} gives $\lambda_n\to+\infty$, which is a contradiction. Let $\mu_0(p,\G,\V_0,q)$ be given by Lemma~\ref{lem-ac-bounded-index-classification}\textup{(ii)} with $M=q$, and let $\mu_{\V_0}$ be given by Theorem~\ref{th2} with $W=\V_0$. Set $\bar\mu:=\min\{\mu_1,\mu_0(p,\G,\V_0,q),\mu_{\V_0}\}$, and let $0<\mu<\bar\mu$. By Theorem~\ref{th2}, the solutions $(u_{\mu,S},\lambda_{\mu,S})$, $\varnothing\ne S\subset\V_0$, are $2^q-1$ pairwise distinct positive solutions of \eqref{eqpd} with $\hat m_0(u_{\mu,S})=\#S$. Conversely, let $(u,\lambda)\in H^1_\mu(\G)\times\R$ be a positive solution of \eqref{eqpd}. Then $\lambda>0$ because $\mu<\mu_1$, so Lemma~\ref{lem-ac-bounded-index-classification}\textup{(i)} gives $\hat m_0(u)\leq q$, and Lemma~\ref{lem-ac-bounded-index-classification}\textup{(ii)} with $M=q$ gives $(u,\lambda)=(u_{\mu,S},\lambda_{\mu,S})$ for some nonempty $S\subset\V_0$. Hence these $2^q-1$ solutions are all the positive solutions at the mass $\mu$. The solutions with constrained Morse index $k$ are those with $\#S=k$, and there are $\binom{q}{k}$ of them. This proves the proposition for $p>4$.

Now let $2<p<4$, and let $o\in\G$ and $A>0$ satisfy the growth assumption of the proposition.
By the proof of Lemma~\ref{lem-ac-bounded-index-classification}\textup{(iii)}, every positive solution $(u,\lambda)$ of \eqref{eqpd} satisfies $\lambda>0$, so Lemma~\ref{lem-ac-bounded-index-classification}\textup{(i)} gives $\hat m_0(u)\leq q$. Let $\mu_0(p,\G,\V_0,q)$ be given by Lemma~\ref{lem-ac-bounded-index-classification}\textup{(iii)} with $M=q$, and let $\mu_{\V_0}$ be given by Theorem~\ref{th2} with $W=\V_0$. Set $\bar\mu:=\max\{\mu_0(p,\G,\V_0,q),\mu_{\V_0}\}$, and let $\mu>\bar\mu$. By Theorem~\ref{th2}, the solutions $(u_{\mu,S},\lambda_{\mu,S})$, $\varnothing\ne S\subset\V_0$, are $2^q-1$ pairwise distinct positive solutions of \eqref{eqpd} with $\hat m_0(u_{\mu,S})=\#S-1$. Conversely, Lemma~\ref{lem-ac-bounded-index-classification}\textup{(iii)} with $M=q$ shows that every positive solution $(u,\lambda)\in H^1_\mu(\G)\times\R$ of \eqref{eqpd} equals $(u_{\mu,S},\lambda_{\mu,S})$ for some nonempty $S\subset\V_0$. The solutions with constrained Morse index $k-1$ are those with $\#S=k$, and there are $\binom{q}{k}$ of them. This completes the proof of Proposition~\ref{prop-ac-mp-single-vertex}.
\end{proof}
\subsection{Examples on metric grids}\label{app-grid-examples}
In this subsection, $\mathcal G^d$, $d \geq 2$, is the $d$-dimensional metric grid defined in the Introduction, and $0$ denotes its vertex at the origin.  Then, $\mathcal G^d$ belongs to $\mathbf G$, is noncompact and satisfies $\lambda(\mathcal G^d)=0$, but it contains no half-line. 
For $\mathcal G^1$, which is excluded here, \eqref{ac-star-single-defect-branch} with $q=2$ gives a branch of positive solutions for $\V_0=\{0\}$, since the degree-two Kirchhoff vertices outside $\V_0$ simply subdivide the two half-lines of the star graph $\G_2$; if $p\ne4$, their masses $2^{2/(p-2)}t^{(4-p)/(p-2)}$ take every value in $(0,+\infty)$, so every $\mu>0$ carries a positive solution, and the small-mass nonexistence result in the following proposition fails for $d=1$.
For simplicity, apart from the nonexistence of solutions at small masses for $d\in\{2,3,4\}$ and $2<p<4$, we only consider the case $\V_0=\{0\}$. For a general nonempty finite set $\V_0$, Theorem~\ref{th2} with $W=\V_0$ still gives $2^{\#\V_0}-1$ positive solutions on $\mathcal G^d$ at every sufficiently large mass if $2<p<4$ and at every sufficiently small mass if $p>4$.
The following proposition shows that, if $d\in\{2,3,4\}$ and $2<p<4$, for every nonempty finite set $\V_0$ of vertices, small masses carry no positive solution, and large masses carry at least two for $\V_0=\{0\}$. Hence the restriction to large masses in Theorem~\ref{th2} cannot be removed for $2<p<4$, and the growth assumption of Proposition~\ref{prop-ac-mp-single-vertex} for $2<p<4$ cannot be dropped; for $d=2$, this complements \cite[Theorem~1.3]{BDS} (see Remark~\ref{rem-grid-bds}). For $\V_0=\{0\}$, if $d\in\{2,3,4\}$ and $p>4$, every prescribed mass carries such a solution, although Theorem~\ref{th1} does not apply, and if $d\geq5$ and $p>4$, large masses carry no such solution, so the half-line in Theorem~\ref{th1} cannot be replaced by the condition $\lambda(\G)=0$, even though Theorem~\ref{th2} gives such a solution for every sufficiently small $\mu>0$.
\begin{proposition}\label{prop-grid-examples}
Let $d\geq2$ be an integer, let $\G=\mathcal G^d$, and let $\varnothing\ne\V_0\subset\V$ be finite.
\begin{enumerate}[label=\rm(\roman*)]
\item If $d\leq4$ and $2<p<4$, there is $\mu_\flat>0$ such that, for $0<\mu<\mu_\flat$, problem \eqref{eqpd} has no positive solution $(u,\lambda)\in H^1_\mu(\G)\times\R$. If, in addition, $\V_0=\{0\}$, then problem \eqref{eqpd} has at least one positive solution $(u,\lambda)\in H^1_{\mu_\flat}(\G)\times(0,+\infty)$ for $\mu=\mu_\flat$, and at least two positive solutions $(u,\lambda)\in H^1_\mu(\G)\times(0,+\infty)$ for every $\mu>\mu_\flat$.
\item If $d\leq4$, $p>4$ and $\V_0=\{0\}$, then for every $\mu>0$, problem \eqref{eqpd} has a positive solution $(u,\lambda)\in H^1_\mu(\G)\times(0,+\infty)$.
\item If $d\geq5$, $p>4$ and $\V_0=\{0\}$, there is $\mu^\sharp>0$ such that, for every $\mu>\mu^\sharp$, problem \eqref{eqpd} has no positive solution $(u,\lambda)\in H^1_\mu(\G)\times\R$.
\end{enumerate}
\end{proposition}
\begin{proof}
We first let $\V_0=\{0\}$ and describe, for a fixed multiplier, the positive solutions $(u,\lambda)$ of \eqref{eqpd}. Identify every edge of $\G$ with $[0,1]$, write $\w\sim\vv$ if the vertices $\w$ and $\vv$ are joined by an edge, and let $(\mathcal Af)(\vv):=\sum_{\w\sim\vv}f(\w)$ be the adjacency operator on $\ell^2(\mathbb Z^d)$. For $f\in\ell^2(\mathbb Z^d)$ and every vertex $\vv \in \V$, the Cauchy--Schwarz inequality gives $\abs{(\mathcal Af)(\vv)}^2\leq2d\sum_{\w\sim\vv}\abs{f(\w)}^2$. Summing over $\vv$ and using that every vertex has $2d$ neighbors, we obtain $\norm{\mathcal A}\leq2d$ on $\ell^2(\mathbb Z^d)$. Fix $t>0$. Since $\norm{(2d\cosh t)^{-1}\mathcal A}\leq1/\cosh t<1$, the operator $2d\cosh t-\mathcal A=2d\cosh t\,\bigl(I-(2d\cosh t)^{-1}\mathcal A\bigr)$ is invertible on $\ell^2(\mathbb Z^d)$ by \cite[Exercise~6.14]{Ha}, and 
\begin{equation}\label{ac-grid-neumann}
 g_t:=(2d\cosh t-\mathcal A)^{-1}\mathbf 1_{\{0\}}=\sum_{n\geq0}(2d\cosh t)^{-n-1}\mathcal A^n\mathbf 1_{\{0\}}
\end{equation}
 defines a map $t\mapsto g_t$ in $C((0,+\infty),\ell^2(\mathbb Z^d))$: the partial sums of the series are continuous in $t$ and, by \cite[Exercise~6.14]{Ha}, converge in $\ell^2(\mathbb Z^d)$ uniformly for $t\geq a$, for every $a>0$, since $\norm{(2d\cosh t)^{-1}\mathcal A}\leq1/\cosh a<1$ for $t\geq a$. Moreover, since $\mathcal A^n\mathbf 1_{\{0\}}\geq0$ for every $n$ and $(\mathcal A^n\mathbf 1_{\{0\}})(\vv)>0$ for $n=\operatorname{dist}_{\rm c}(0,\vv)$, we have $g_t>0$. Since the dual group of $\mathbb Z^d$ is the $d$-dimensional torus \cite[Section~1.2.7 and Theorem~2.2.2]{Ru}, the Plancherel theorem \cite[Theorem~1.6.1 and Section~1.6.2]{Ru} gives a Fourier transform from $\ell^2(\mathbb Z^d)$ onto $L^2([-\pi,\pi]^d)$. For $f\in\ell^2(\mathbb Z^d)$, we denote this transform by $\hat f$ and write $\hat f(k):=\sum_{\vv\in\mathbb Z^d}f(\vv)e^{-\mathrm ik\cdot\vv}$, where the series is understood in $L^2([-\pi,\pi]^d)$. Moreover, the following Parseval's identity holds
\begin{equation}\label{ac-grid-parseval}
 \sum_{\vv\in\mathbb Z^d}\abs{f(\vv)}^2=(2\pi)^{-d}\int_{[-\pi,\pi]^d}\abs{\hat f(k)}^2\,dk.
\end{equation}
Applying the polarized form of this identity \cite[Section~1.6.2]{Ru} to $f$ and $\mathbf 1_{\{\vv\}}$, whose transform is $e^{-\mathrm ik\cdot\vv}$, we obtain the inversion formula
\begin{equation}\label{ac-grid-inversion}
 f(\vv)=\frac1{(2\pi)^d}\int_{[-\pi,\pi]^d}\hat f(k)e^{\mathrm ik\cdot\vv}\,dk\qquad\text{for every }\vv\in\mathbb Z^d.
\end{equation}
Let $\mathbf e_1,\ldots,\mathbf e_d$ be the unit vectors of $\mathbb Z^d$. Then $(\mathcal Af)(\vv)=\sum_{j=1}^d\bigl(f(\vv+\mathbf e_j)+f(\vv-\mathbf e_j)\bigr)$, and, for finitely supported $f$, shifting the summation index gives
\begin{equation}\label{ac-grid-multiplier}
 \begin{aligned}
 \widehat{\mathcal Af}(k)&=\sum_{j=1}^d\sum_{\vv\in\mathbb Z^d}f(\vv)\bigl(e^{-\mathrm ik\cdot(\vv-\mathbf e_j)}+e^{-\mathrm ik\cdot(\vv+\mathbf e_j)}\bigr)=\sum_{j=1}^d\bigl(e^{\mathrm ik_j}+e^{-\mathrm ik_j}\bigr)\hat f(k)=2\sum_{j=1}^d\cos k_j\,\hat f(k).
 \end{aligned}
\end{equation}
By Parseval's identity and $\norm{\mathcal A}\leq2d$, for every $f\in\ell^2(\mathbb Z^d)$,
\[
 \begin{aligned}
 \norm{\widehat{\mathcal Af}}_{L^2([-\pi,\pi]^d)}&=(2\pi)^{d/2}\norm{\mathcal Af}_{\ell^2(\mathbb Z^d)}\leq2d(2\pi)^{d/2}\norm{f}_{\ell^2(\mathbb Z^d)},\\
 \Bigl\|2\sum_{j=1}^d\cos k_j\,\hat f\Bigr\|_{L^2([-\pi,\pi]^d)}&\leq2d\norm{\hat f}_{L^2([-\pi,\pi]^d)}=2d(2\pi)^{d/2}\norm{f}_{\ell^2(\mathbb Z^d)}.
 \end{aligned}
\]
Since finitely supported functions are dense in $\ell^2(\mathbb Z^d)$, the identity $\widehat{\mathcal Af}=2\sum_{j=1}^d\cos k_j\,\hat f$ holds for every $f\in\ell^2(\mathbb Z^d)$. Applied to $(2d\cosh t-\mathcal A)g_t=\mathbf 1_{\{0\}}$, it gives $\bigl(2d\cosh t-2\sum_{j=1}^d\cos k_j\bigr)\hat g_t(k)=1$, so, by the inversion formula \eqref{ac-grid-inversion}, the function $g_t$ has the Fourier representation
\begin{equation}\label{ac-grid-fourier}
 g_t(\vv)=\frac1{(2\pi)^d}\int_{[-\pi,\pi]^d}\frac{e^{\mathrm ik\cdot\vv}}{\delta_t(k)}\,dk,
 \qquad
 \delta_t(k):=2d\cosh t-2\sum_{j=1}^d\cos k_j .
\end{equation}
Let $h_t:=U_t\bigl(g_t/g_t(0)\bigr)$ with $U_t$ as in \eqref{ac-edge-extension}, that is, $h_t=g_t/g_t(0)$ at the vertices and, for every edge $\e$ with endpoints $\vv,\w$, identified with $[0,1]$ in such a way that $\vv$ is identified with $0$ and $\w$ with $1$,
\[
 (h_t)_\e(x)=h_t(\vv)\frac{\sinh\bigl(t(1-x)\bigr)}{\sinh t}+h_t(\w)\frac{\sinh(tx)}{\sinh t},\qquad x\in[0,1].
\]
Then $h_t\in H^1(\G)$, $h_t''=t^2h_t$ on every edge and, for every vertex $\vv$,
\[
 \sum_{\e\succ\vv}(h_t)'_\e(\vv)=\frac{t}{\sinh t}\bigl((\mathcal A-2d\cosh t)h_t\bigr)(\vv)=-D(t)\mathbf 1_{\{0\}}(\vv),
 \qquad
 D(t):=\frac{t}{g_t(0)\sinh t}.
\]
Hence $u_t:=D(t)^{1/(p-2)}h_t>0$ satisfies $\sum_{\e\succ0}(u_t)'_\e(0)=-u_t(0)^{p-1}$ and  $\sum_{\e\succ \vv}(u_t)'_\e(\vv)=0$ at every vertex $\vv\ne0$. Thus $(u_t,t^2)$ solves \eqref{eqpd} with mass $\mu(t):=\norm{u_t}_2^2=D(t)^{2/(p-2)}\norm{h_t}_2^2$, which depends continuously on $t>0$. 

Conversely, let $(u,t^2)$ be a positive solution of \eqref{eqpd}. Since $u\in H^1(\G)$ and $\deg(\vv)=2d$ for every $\vv\in\V$, Lemma~\ref{ac-restriction-extension} shows that the vertex values $f:=u|_{\mathbb Z^d}$ belong to $\ell^2(\mathbb Z^d)$. Since $u''=t^2u$ on every edge, $u=U_tf$, so the vertex conditions in \eqref{eqpd} give $(2d\cosh t-\mathcal A)f=(\sinh t/t)\,u(0)^{p-1}\mathbf 1_{\{0\}}$. Since $2d\cosh t-\mathcal A$ is invertible, $f=(\sinh t/t)\,u(0)^{p-1}g_t$. This leads to $u(0)^{p-2}=D(t)$, hence $f=D(t)^{1/(p-2)}g_t/g_t(0)$ and $u=U_tf=D(t)^{1/(p-2)}h_t=u_t$.

Next, let $t\geq1$. Since $\norm{(2d\cosh t)^{-1}\mathcal A}\leq1/\cosh t$, \cite[Exercise~6.14]{Ha} gives $\norm{(2d\cosh t-\mathcal A)^{-1}}\leq(2d\cosh t)^{-1}(1-1/\cosh t)^{-1}=(2d(\cosh t-1))^{-1}$. Applied to $\mathbf 1_{\{0\}}$, this gives $g_t(0)^2\leq\sum_{\vv\in\mathbb Z^d}g_t(\vv)^2\leq(2d(\cosh t-1))^{-2}$, while the first term of the series in \eqref{ac-grid-neumann}, whose terms are nonnegative, gives $g_t(0)\geq(2d\cosh t)^{-1}$. Hence
\begin{equation}\label{ac-grid-D-bounds}
2dt\tanh(t/2)=2dt\frac{\cosh t-1}{\sinh t}\leq D(t)\leq2dt\coth t.
\end{equation}
For every edge $\e \in \E$ with endpoints $\vv,\w$, identified with $[0,1]$ as above, the inequality $(\alpha+\beta)^2\leq2\alpha^2+2\beta^2$ and $\int_0^1\sinh^2(ts)\,ds=\sinh t\cosh t/(2t)-1/2\leq\sinh t\cosh t/(2t)$ show that 
\[\int_0^1(h_t)_\e(x)^2\,dx\leq(h_t(\vv)^2+h_t(\w)^2)\coth t/t.\] 
For each of the $2d$ edges $\e\succ0$, identified with $[0,1]$ in such a way that $0$ is identified with $0$, by $h_t(0)=1$ and $h_t\geq0$, we have $(h_t)_\e(x)\geq\sinh(t(1-x))/\sinh t$. Moreover, since $t/\sinh^2t$ is decreasing, we obtain $\int_0^1(h_t)_\e(x)^2\,dx\geq\coth t/(2t)-1/(2\sinh^2t)\geq(1-\operatorname{csch}^2 1)/(2t)$. Since every vertex lies on $2d$ edges, we obtain
\[
 \frac{d(1-\operatorname{csch}^2 1)}{t}\leq\norm{h_t}_2^2\leq\frac{2d\coth t}{t}\sum_{\vv\in\mathbb Z^d}\frac{g_t(\vv)^2}{g_t(0)^2}\leq\frac{2d\coth(t)\cosh^2t}{t(\cosh t-1)^2}.
\]
Together with the bounds \eqref{ac-grid-D-bounds} for $D(t)$, since $\tanh(t/2)$ is increasing while $\coth t$ and $\cosh t/(\cosh t-1)$ are decreasing, this shows 
\begin{equation}\label{ac-grid-mu-large}
 d(1-\operatorname{csch}^2 1)\bigl(2d\tanh(1/2)\bigr)^{2/(p-2)}\leq t^{-(4-p)/(p-2)}\mu(t)\leq\frac{2d\coth1\,\cosh^21}{(\cosh1-1)^2}\,(2d\coth1)^{2/(p-2)}.
\end{equation}
In particular, as $t\to+\infty$, $\mu(t)\to0$ if $p>4$, and $\mu(t)\to+\infty$ if $2<p<4$.

Now let $d\leq4$ and $0<t\leq1$. Since, by Taylor's formula and by the fact that $s\mapsto(1-\cos s)/s^2$ is decreasing on $(0,\pi]$, $dt^2\leq2d(\cosh t-1)\leq2dt^2$ and $4s^2/\pi^2\leq2-2\cos s\leq s^2$ for $\abs{s}\leq\pi$, we have $dt^2+4\abs{k}^2/\pi^2\leq\delta_t(k)\leq2dt^2+\abs{k}^2$ on $[-\pi,\pi]^d$. Hence the Fourier representation \eqref{ac-grid-fourier} and Parseval's identity \eqref{ac-grid-parseval}, together with polar coordinates, give
\[
 \begin{aligned}
 g_t(0)&\leq\gamma_d(t):=\frac1{(2\pi)^d}\int_{\abs{k}\leq\pi\sqrt d}\frac{dk}{dt^2+4\abs{k}^2/\pi^2}=\frac{\omega_d}{(2\pi)^d}\int_0^{\pi\sqrt d}\frac{r^{d-1}\,dr}{dt^2+4r^2/\pi^2},\\
 \sum_{\vv\in\mathbb Z^d}g_t(\vv)^2&\geq\sigma_d(t):=\frac1{(2\pi)^d}\int_{\abs{k}\leq\pi}\frac{dk}{(2dt^2+\abs{k}^2)^2}=\frac{\omega_d}{(2\pi)^d}\int_0^{\pi}\frac{r^{d-1}\,dr}{(2dt^2+r^2)^2},
 \end{aligned}
\]
where $\omega_d$ is the area of the unit sphere in $\mathbb R^d$, so that $\omega_2=2\pi$, $\omega_3=4\pi$ and $\omega_4=2\pi^2$. For $d=2$, the last integral in the first line equals $\frac{\pi^2}8\log(1+4/t^2)$, and for $d\in\{3,4\}$ we use $dt^2 +4r^2/\pi^2\geq 4r^2/\pi^2$ in it, hence
\[
 \gamma_2(t)\leq\frac\pi{16}\log\Bigl(1+\frac4{t^2}\Bigr),\qquad
 \gamma_3(t)\leq\frac{\sqrt3\,\pi}8,\qquad
 \gamma_4(t)\leq\frac{\pi^2}{16}.
\]
Since $2dt^2+r^2\leq2r^2$ for $r\geq\sqrt{2d}\,t$, and $\sqrt{2d}\,t\leq\sqrt8<\pi$, the last integral in the second line satisfies
\[
 \begin{aligned}
 \int_0^\pi\frac{r\,dr}{(4t^2+r^2)^2}&=\frac{\pi^2}{8t^2(4t^2+\pi^2)}&&\text{if }d=2,\\
 \int_0^\pi\frac{r^2\,dr}{(6t^2+r^2)^2}&\geq\int_{\sqrt6\,t}^\pi\frac{dr}{4r^2}=\frac14\Bigl(\frac1{\sqrt6\,t}-\frac1\pi\Bigr)&&\text{if }d=3,\\
 \int_0^\pi\frac{r^3\,dr}{(8t^2+r^2)^2}&\geq\int_{\sqrt8\,t}^\pi\frac{dr}{4r}=\frac14\log\frac{\pi}{2\sqrt2\,t}&&\text{if }d=4,
 \end{aligned}
\]
and therefore, since $t\leq1$,
\[
 \sigma_2(t)\geq\frac{\pi}{16(4+\pi^2)t^2},\qquad
 \sigma_3(t)\geq\frac{1-\sqrt6/\pi}{8\sqrt6\,\pi^2t},\qquad
 \sigma_4(t)\geq\frac1{32\pi^2}\log\frac{\pi}{2\sqrt2\,t}.
\]
Moreover, $\sinh(ts)\geq ts$ and $\sinh t=\int_0^t\cosh s\,ds\leq t\cosh1$ for $0<t\leq1$ give $\sinh(ts)/\sinh t\geq s/\cosh1$ for $s\in[0,1]$. Hence, for every edge $\e$ with endpoints $\vv,\w$, identified with $[0,1]$ as above,
\[
 (h_t)_\e(x)\geq\frac{h_t(\vv)(1-x)+h_t(\w)x}{\cosh1},\qquad x\in[0,1],
\]
and
\[
 \int_0^1(h_t)_\e(x)^2\,dx\geq\frac{h_t(\vv)^2+h_t(\w)^2}{3\cosh^21}.
\]
Since every vertex lies on $2d$ edges and $D(t)\geq1/(g_t(0)\cosh1)$, we obtain
\begin{equation}\label{ac-grid-mu-small}
 \mu(t)\geq\frac{2d}{3\cosh^21}\cdot\frac{\sigma_d(t)}{(\gamma_d(t)\cosh1)^{2/(p-2)}\gamma_d(t)^2}\longrightarrow+\infty
 \qquad\text{as }t\to0^+
\end{equation}
for $d\in\{2,3,4\}$.

We now prove \textup{(i)}. Let $d\leq4$ and $2<p<4$, put $q:=\#\V_0$, and let $(u,\lambda)\in H^1_\mu(\G)\times\R$ be a positive solution of \eqref{eqpd}. By Lemma~\ref{lem-ac-bounded-index-classification}\textup{(i)}, $\lambda\geq0$. Suppose that $\lambda=0$. Then $u''=\lambda u=0$ on every edge by \eqref{eqpd}, so $u$ is affine on every edge, and, since $\deg(\vv)=2d$ for every $\vv\in\V$, Lemma~\ref{ac-restriction-extension} shows that its vertex values $f:=u|_{\mathbb Z^d}$ belong to $\ell^2(\mathbb Z^d)$. By the vertex conditions in \eqref{eqpd}, $(2d-\mathcal A)f=\sum_{\vv\in\V_0}u(\vv)^{p-1}\mathbf 1_{\{\vv\}}$. In Fourier variables, this reads $\delta_0\hat f=P$ on $[-\pi,\pi]^d$, where $P(k):=\sum_{\vv\in\V_0}u(\vv)^{p-1}e^{-\mathrm ik\cdot\vv}$ and $\delta_0(k)=2\sum_{j=1}^d(1-\cos k_j)\leq\abs{k}^2$. Since $P$ is continuous and $P(0)>0$, we have $\abs{\hat f(k)}\geq P(0)\abs{k}^{-2}/2$ near $k=0$, contradicting $f\in\ell^2(\mathbb Z^d)$ since $d\le4$. 
Thus $\lambda=t^2$ with $t>0$. 

As in the first part of the proof, $u$ is given on every edge by its vertex values $f\in\ell^2(\mathbb Z^d)$, and the vertex conditions give $(2d\cosh t-\mathcal A)f=(\sinh t/t)\sum_{\vv\in\V_0}u(\vv)^{p-1}\mathbf 1_{\{\vv\}}$, so, by \eqref{ac-grid-neumann}, we have
\[
 f=\frac{\sinh t}{t}\sum_{\vv\in\V_0}u(\vv)^{p-1}g_t(\cdot-\vv).
\]
Choose $\vv_*\in\V_0$ with $u(\vv_*)=\max_{\V_0}u$. Since $g_t\in\ell^2(\mathbb Z^d)$ is positive, it attains its maximum, and only at $0$: at a maximum point $\vv\ne0$ we would have $2d\cosh t\,g_t(\vv)=(\mathcal Ag_t)(\vv)\le2d\,g_t(\vv)$, which is impossible for $t>0$. Hence $0<g_t(\vv)\leq g_t(0)$ for every $\vv \in \V$, so, since $f(\vv_*)=u(\vv_*)$, the formula for $f$ gives $u(\vv_*)\leq q(\sinh t/t)g_t(0)u(\vv_*)^{p-1}$, that is, $u(\vv_*)^{p-2}\geq D(t)/q$. Moreover, for every $\vv \in \V$, $f(\vv)\geq(\sinh t/t)\,u(\vv_*)^{p-1}g_t(\vv-\vv_*)$, and the right-hand side equals $u(\vv_*)^{p-1}D(t)^{-1}h_t(\vv-\vv_*)$. Since $u=U_tf$ and $h_t(\cdot-\vv_*)$ are given on every edge by \eqref{ac-edge-extension}, whose coefficients are nonnegative, we have $u \geq u(\vv_*)^{p-1}D(t)^{-1}h_t(\cdot-\vv_*)$, and therefore
\[
 \norm{u}_2^2\geq u(\vv_*)^{2(p-1)}D(t)^{-2}\norm{h_t}_2^2\geq q^{-2(p-1)/(p-2)}D(t)^{2/(p-2)}\norm{h_t}_2^2=q^{-2(p-1)/(p-2)}\mu(t).
\]
The function $\mu(\cdot)$ is continuous and positive on $(0,+\infty)$ and tends to $+\infty$ as $t\to0^+$ by \eqref{ac-grid-mu-small} and as $t\to+\infty$ by \eqref{ac-grid-mu-large}, since $2<p<4$. Hence it attains its minimum at some $t_0>0$, and $\mu_\flat:=q^{-2(p-1)/(p-2)}\mu(t_0)>0$ has the required property. 

Now let $\V_0=\{0\}$, so that $q=1$, $\mu_\flat=\mu(t_0)$ and $(u_t,t^2)\in H^1_{\mu(t)}(\G)\times(0,+\infty)$ is a positive solution of \eqref{eqpd} for every $t>0$. At $t=t_0$ this gives a positive solution $(u_{t_0},t_0^2)$ of mass $\mu_\flat$. If $\mu>\mu_\flat$, the intermediate value theorem gives $t_1\in(0,t_0)$ and $t_2\in(t_0,+\infty)$ with $\mu(t_1)=\mu(t_2)=\mu$, and $(u_{t_1},t_1^2)$ and $(u_{t_2},t_2^2)$ are two distinct positive solutions of \eqref{eqpd}.

We next prove \textup{(ii)}. Let $d\leq4$, $p>4$ and $\V_0=\{0\}$. The function $\mu(\cdot)$ is continuous on $(0,+\infty)$, tends to $+\infty$ as $t\to0^+$ by \eqref{ac-grid-mu-small} and tends to $0$ as $t\to+\infty$ by \eqref{ac-grid-mu-large}, since $p>4$, so $\mu((0,+\infty))=(0,+\infty)$. For $\mu>0$, the pair $(u_t,t^2)$ with $\mu(t)=\mu$ is the required solution.

We finally prove \textup{(iii)}. Let $d\geq5$, $p>4$ and $\V_0=\{0\}$. By Taylor's formula for the upper bound and by the monotonicity of $s\mapsto(1-\cos s)/s^2$ on $(0,\pi]$ for the lower one, $\frac{4}{\pi^2}\abs{k}^2\leq\delta_0(k)=2\sum_{j=1}^d(1-\cos k_j)\leq\abs{k}^2$ on $[-\pi,\pi]^d$. Since $\int_{[-\pi,\pi]^d}\abs{k}^{-4}\,dk<+\infty$, we have $1/\delta_0\in L^2([-\pi,\pi]^d)$, and we define
\[
 G_0(\vv):=\frac1{(2\pi)^d}\int_{[-\pi,\pi]^d}\frac{e^{\mathrm ik\cdot\vv}}{\delta_0(k)}\,dk,\qquad\vv\in\mathbb Z^d.
\]
By the Plancherel theorem \cite[Theorem~1.6.1]{Ru} and the inversion formula \eqref{ac-grid-inversion}, $G_0\in\ell^2(\mathbb Z^d)$ and $\widehat{G_0}=1/\delta_0$. Hence the identity $\widehat{\mathcal Af}=2\sum_{j=1}^d\cos k_j\,\hat f$, which holds for every $f\in\ell^2(\mathbb Z^d)$ by \eqref{ac-grid-multiplier} and the density argument following it, gives $\widehat{(2d-\mathcal A)G_0}=\delta_0\widehat{G_0}=1=\widehat{\mathbf 1_{\{0\}}}$, and, since $f\mapsto\hat f$ is injective on $\ell^2(\mathbb Z^d)$ by Parseval's identity \eqref{ac-grid-parseval}, $(2d-\mathcal A)G_0=\mathbf 1_{\{0\}}$.

By the Neumann series \eqref{ac-grid-neumann}, $t\mapsto g_t(\vv)$ is decreasing on $(0,+\infty)$ for every $\vv\in\mathbb Z^d$. Moreover, $\widehat{g_t}=1/\delta_t$ by \eqref{ac-grid-fourier}, and $1/\delta_t(k)\nearrow1/\delta_0(k)$ as $t\searrow0$ for every $k\in[-\pi,\pi]^d\setminus\{0\}$, where $1/\delta_0\in L^2([-\pi,\pi]^d)$. Since $(1/\delta_0-1/\delta_t)^2\leq1/\delta_0^2-1/\delta_t^2$, monotone convergence and Parseval's identity \eqref{ac-grid-parseval} give $g_t\to G_0$ in $\ell^2(\mathbb Z^d)$. Hence $g_t(\vv)\nearrow G_0(\vv)$ as $t\searrow0$ for every $\vv\in\mathbb Z^d$; in particular, $G_0>0$. Let $h_0$ be the function on $\G$ such that $h_0(\vv)=G_0(\vv)/G_0(0)$ for every $\vv \in \mathbb Z^d$ and, for every edge $\e$ with endpoints $\vv,\w$, identified with $[0,1]$ as above,
\[
 (h_0)_\e(x)=h_0(\vv)(1-x)+h_0(\w)x,\qquad x\in[0,1],
\]
and put $D(0):=1/G_0(0)$ and $u_0:=D(0)^{1/(p-2)}h_0$. Since $G_0\in\ell^2(\mathbb Z^d)$, every vertex lies on $2d$ edges and, for every edge $\e$ with endpoints $\vv,\w$, $\int_0^1(h_0)_\e(x)^2\,dx\leq(h_0(\vv)^2+h_0(\w)^2)/2$ and $\int_0^1(h_0)'_\e(x)^2\,dx=(h_0(\w)-h_0(\vv))^2\leq2(h_0(\vv)^2+h_0(\w)^2)$, we have $h_0\in H^1(\G)$. Then $\sum_{\e\succ\vv}(h_0)'_\e(\vv)=((\mathcal A-2d)h_0)(\vv)=-D(0)\mathbf 1_{\{0\}}(\vv)$ for every vertex $\vv$, so $(u_0,0)$ is a positive solution of \eqref{eqpd}. Moreover, as $t\to0^+$, we have $D(t)\to D(0)$ and $f_t:=g_t/g_t(0)\to f_0:=G_0/G_0(0)$ in $\ell^2(\mathbb Z^d)$. Integrating over each edge, by \eqref{ac-edge-extension} for $h_t$ and since $h_0$ is affine on every edge, and summing over the edges, we obtain
\[
 \begin{aligned}
 \norm{h_t}_2^2&=2d\,\alpha_t\norm{f_t}_{\ell^2(\mathbb Z^d)}^2+\beta_t\sum_{\vv\in\mathbb Z^d}f_t(\vv)(\mathcal Af_t)(\vv),\\
 \norm{h_0}_2^2&=\frac{2d}3\norm{f_0}_{\ell^2(\mathbb Z^d)}^2+\frac16\sum_{\vv\in\mathbb Z^d}f_0(\vv)(\mathcal Af_0)(\vv),
 \end{aligned}
\]
where
\[
 \alpha_t:=\frac1{\sinh^2t}\int_0^1\sinh^2(tx)\,dx,\qquad\beta_t:=\frac1{\sinh^2t}\int_0^1\sinh\bigl(t(1-x)\bigr)\sinh(tx)\,dx.
\]
Since $\sinh(tx)/\sinh t\to x$ uniformly on $[0,1]$, we have $\alpha_t\to1/3$ and $\beta_t\to1/6$, and hence $\norm{h_t}_2\to\norm{h_0}_2$ as $t\to0^+$. Thus,  $\mu(\cdot)$, extended by $\mu(0):=\norm{u_0}_2^2>0$, is continuous on $[0,+\infty)$ and tends to $0$ as $t\to+\infty$ by \eqref{ac-grid-mu-large} and $p>4$, so $\mu^\sharp:=\max_{t\geq0}\mu(t)$ is finite.

Now let $(u,\lambda)\in H^1_\mu(\G)\times\R$ be a positive solution of \eqref{eqpd}. By Lemma~\ref{lem-ac-bounded-index-classification}\textup{(i)}, $\lambda\geq0$. If $\lambda=t^2>0$, then $u=u_t$ by the first part of the proof. If $\lambda=0$, then $u''=0$ on every edge by \eqref{eqpd}, so $u$ is affine on every edge and its vertex values $f:=u|_{\mathbb Z^d}$, which belong to $\ell^2(\mathbb Z^d)$ by Lemma~\ref{ac-restriction-extension} as in the proof of \textup{(i)}, satisfy $(2d-\mathcal A)f=u(0)^{p-1}\mathbf 1_{\{0\}}$. If $\psi\in\ell^2(\mathbb Z^d)$ and $(2d-\mathcal A)\psi=0$, then $\delta_0\hat\psi=0$ in Fourier variables, so $\hat\psi=0$ because $\delta_0>0$ on $[-\pi,\pi]^d\setminus\{0\}$, and $\psi=0$. Thus $2d-\mathcal A$ is injective on $\ell^2(\mathbb Z^d)$, and, since $(2d-\mathcal A)G_0=\mathbf 1_{\{0\}}$, $f=u(0)^{p-1}G_0$, and $f(0)=u(0)$ then gives $u=u_0$. In both cases $\mu\leq\mu^\sharp$.
\end{proof}
\begin{remark}\label{rem-grid-bds}
Let $d=2$ and $2<p<4$. The graph $\mathcal G^2$ is the two-dimensional square grid of \cite{BDS}. For every nonempty finite set $\V_0$ of vertices, \cite[Theorem~1.3]{BDS} gives a critical mass $\mu_p^*>0$ such that ground states of $E(\cdot,\G)$ in $H^1_\mu(\G)$ exist for $\mu>\mu_p^*$ and do not exist for $\mu<\mu_p^*$, whereas Proposition~\ref{prop-grid-examples}\textup{(i)} excludes, for $0<\mu<\mu_\flat$, all positive solutions $(u,\lambda)$ of \eqref{eqpd}, not only ground states. Now let $\V_0=\{0\}$. By the proof of \cite[Proposition~1.4]{BDS}, for every ground state $u$, up to a choice of sign, $(u,\lambda)$ is a positive solution of \eqref{eqpd} for some $\lambda>0$. By the first part of the proof of Proposition~\ref{prop-grid-examples}, $(u,\lambda)=(u_t,t^2)$ for some $t>0$, so $\mu_\flat\leq\mu_p^*$. Thus Proposition~\ref{prop-grid-examples}\textup{(i)} gives at least two positive solutions $(u,\lambda)$ for every $\mu>\mu_\flat$, hence for every $\mu>\mu_p^*$, whereas \cite[Theorem~1.3]{BDS} gives a ground state. Proposition~\ref{prop-grid-examples}\textup{(i)} also covers the grids $\mathcal G^3$ and $\mathcal G^4$, which are not considered in \cite{BDS}.
\end{remark}

\section[All-mass multiplicity on a tree]{All-mass multiplicity on a tree}
\label{app-ac-critical-all-mass-tree}

In this appendix, $p=4$. We construct a metric tree on which problem \eqref{eqpd} has infinitely many solutions at every prescribed mass. We construct the branches on a frequency interval extending down to a point where their mass is smaller than the prescribed value.

\begin{proposition}\label{prop-ac-critical-all-mass-tree}
There exist a noncompact metric tree $\G\in\mathbf G$, with every edge of length one and every vertex of degree at most three, and an infinite defect set $\V_0$ consisting of leaves, such that the following holds. For $p=4$ and every $\mu>0$, problem \eqref{eqpd} has infinitely many positive solutions $(u,\lambda)\in H^1_\mu(\G)\times(0,+\infty)$, pairwise inequivalent under the metric graph automorphisms preserving $\V_0$.
\end{proposition}

\begin{proof}
\noindent\textit{Step~1: construction of the tree.}
We first fix two sequences of integers. 
Choose $\{L_j\}_{j\in\mathbb Z}$ such that
\[
 L_j\in\mathbb Z,\qquad L_j\geq2, \qquad
 L_i\ne L_j\ \text{whenever }i\ne j,\qquad
 \lim_{|j|\to\infty}L_j=+\infty .
\]
For example, $L_j:=2j+2$ for $j\geq0$ and $L_j:=1-2j$ for $j<0$. Let $d_j$, $j\in\mathbb Z$, be positive integers such that every positive integer occurs infinitely often in the sequence $\{d_j\}$; for example, 
\[
 d_j:=1+\max\bigl\{m\in\mathbb Z_{\geq0}:2^m\mid (|j|+1)\bigr\},
\]
where, for $a,b\in\mathbb Z$, we write
\[
 a\mid b\quad\Longleftrightarrow\quad
 \text{there exists }k\in\mathbb Z\text{ such that }b=ka.
\]

\emph{The chain $\G^{(0)}$.} Choose the integers $b_j$, $j\in\mathbb Z$, by $b_0=0$ and $b_{j+1}-b_j=L_j$. Let $\G^{(0)}=(\V^{(0)},\E^{(0)})$ be the bi-infinite chain in $\R$ with vertex set and edge set
\[
 \V^{(0)}:=\{b_j:j\in\mathbb Z\},\qquad
 \E^{(0)}:=\{\e_j:j\in\mathbb Z\},\qquad
 \e_j=(b_j,b_{j+1}).
\]
Each chain edge $\e_j$ has length $L_j$ and is identified with $I_{\e_j}=[0,L_j]$.

\emph{The tree $\G^{(1)}$.}
For an integer $d\geq1$, set
\[
 \mathcal W_d:=\bigcup_{k=0}^{d}\{0,1\}^k,
 \qquad \{0,1\}^0:=\{\varnothing\},
\]
where $\varnothing$ is the empty word, i.e., the unique element of $\mathcal W_d$ of length zero. For $\omega\in\{0,1\}^k$, we write $|\omega|:=k$, so in particular $|\varnothing|=0$. For $a\in\{0,1\}$, set
\[
 \varnothing a:=(a).
\]
For $\omega=(\omega_1,\ldots,\omega_k)\in\{0,1\}^k$ with $k\geq1$ and $a\in\{0,1\}$, set
\[
 \omega a:=(\omega_1,\ldots,\omega_k,a).
\]

Define $\G^{(1)}=(\V^{(1)},\E^{(1)})$ by
\[
 \V^{(1)}
 =\V^{(0)}\cup\{r_{j,\omega}:j\in\mathbb Z,\ \omega\in\mathcal W_{d_j}\},\qquad
 \E^{(1)}
 =\E^{(0)}\cup\{\e_{j,\omega}:j\in\mathbb Z,\ \omega\in\mathcal W_{d_j}\}.
\]
For each $j\in\mathbb Z$, the added edges are given by
\[
 \begin{aligned}
 \e_{j,\varnothing}&=(b_j,r_{j,\varnothing}),\\
 \e_{j,\omega a}&=(r_{j,\omega},r_{j,\omega a})
 \qquad\text{for }\omega\in\mathcal W_{d_j},\ |\omega|<d_j,\ a\in\{0,1\}.
 \end{aligned}
\]
Every edge in $\E^{(1)}\setminus\E^{(0)}$ is identified with $[0,1]$.
Thus, for each $j$, the vertices $r_{j,\omega}$ and the edges between them form a complete binary tree of depth $d_j$ with root $r_{j,\varnothing}$, and the edge \(\mathrm e_{j,\varnothing}\) joins the root \(r_{j,\varnothing}\) to the chain vertex \(b_j\). For $k=0,\ldots,d_j$, the vertices $r_{j,\omega}$ with $|\omega|=k$ form the $k$-th level of this binary tree. The leaves are the vertices of the last level, and we set
\[
 S_j:=\{r_{j,\omega}:\omega\in\mathcal W_{d_j},\ |\omega|=d_j\},\qquad \V_0:=\bigcup_{j\in\mathbb Z}S_j .
\]
So $\#S_j=2^{d_j}$, and $\V_0$ is the set of all leaves of $\G^{(1)}$. In $\G^{(1)}$, every $b_j$ has degree three (two chain edges  $\e_{j-1},\e_j$ and one stem edge  $\e_{j,\varnothing}$), every $r_{j,\omega}$ with $|\omega|<d_j$ has degree three, and every leaf has degree one. All edges are bounded, and their lengths are at least one. Hence $\G^{(1)}\in\mathbf G$ with $\ell_*=1$. Moreover, $\G^{(1)}$ is connected, has no cycles and has infinitely many edges, so it is a noncompact tree. A part of $\G^{(1)}$ is shown in Figure~\ref{fig-ac-all-mass-tree-chain}\textup{(a)}.

\emph{The tree $\G$.} Subdivide each chain edge $\e_j$ at the vertices $c_{j,i}:=b_j+i$ for $1\leq i<L_j$, and set $c_{j,0}:=b_j$ and $c_{j,L_j}:=b_{j+1}$. The resulting unit-edge tree $\G=(\V,\E)$ is given by
\[
 \begin{aligned}
 \V&=\V^{(1)}\cup\{c_{j,i}:j\in\mathbb Z,\ 1\leq i<L_j\},\\
 \E&=\bigl(\E^{(1)}\setminus\E^{(0)}\bigr)
 \cup\{\e_j^i:j\in\mathbb Z,\ 0\leq i<L_j\}.
 \end{aligned}
\]
The unit edges $\e_j^i$ have endpoints
\[
\e_j^i=(c_{j,i},c_{j,i+1})
 \qquad\text{for }j\in\mathbb Z,\ 0\leq i<L_j.
\]
The defect set is still $\V_0$. The added vertices $c_{j,i}$ with $1\leq i<L_j$ have degree two and lie outside $\V_0$ (see Figure~\ref{fig-ac-all-mass-tree-chain}\textup{(b)}), so every vertex of $\G$ has degree at most three.

Restriction to the subdivided edges and gluing at the subdivision vertices $c_{j,i}$ identify \(H^1(\G^{(1)})\) and \(H^1(\G)\) isometrically and preserve the \(L^2\) norm. At each subdivision vertex \(c_{j,i}\), the Kirchhoff condition is equivalent to the matching of the one-sided derivatives along the chain. Therefore a function $u\in H^1(\G)$ satisfies the edge equations and the vertex conditions of \eqref{eqpd} on $\G$, with multiplier $t^2$, if and only if it satisfies them on $\G^{(1)}$. We construct the solutions on $\G^{(1)}$.
\begin{figure}[H]
\centering
\begin{tikzpicture}[x=.93cm,y=.87cm,line cap=round,line join=round,every node/.style={font=\small,inner sep=2pt}]
 \draw[densely dotted] (-4.65,0)--(-3.25,0) (3.25,0)--(4.65,0);
 \draw (-3.25,0)--(0,0)--(3.25,0);
 \node[below] at (-3.25,0) {$b_{j-1}$};
 \node[below] at (0,0) {$b_j$};
 \node[below] at (3.25,0) {$b_{j+1}$};
 \node[above] at (-1.625,0) {$L_{j-1}$};
 \node[above] at (1.625,0) {$L_j$};
 \foreach \xx in {-3.25,3.25}{
   \begin{scope}[shift={(\xx,0)}]
    \draw (0,0)--(0,.75)--(-.55,1.45) (0,.75)--(.55,1.45);
    \fill (0,0) circle[radius=1.4pt] (0,.75) circle[radius=1.4pt];
    \draw[fill=white] (-.55,1.45) circle[radius=2pt] (.55,1.45) circle[radius=2pt];
   \end{scope}
 }
 \draw (0,0)--(0,.65)--(-.63,1.30)--(-.99,2)
       (-.63,1.30)--(-.27,2) (0,.65)--(.63,1.30)--(.27,2)
       (.63,1.30)--(.99,2);
 \fill (0,0) circle[radius=1.4pt] (0,.65) circle[radius=1.4pt]
       (-.63,1.30) circle[radius=1.4pt] (.63,1.30) circle[radius=1.4pt];
 \foreach \xx in {-.99,-.27,.27,.99}
   \draw[fill=white] (\xx,2) circle[radius=2pt];
 \node[right] at (-3.25,.375) {$\e_{j-1,\varnothing}$};
 \node[right] at (0,.325) {$\e_{j,\varnothing}$};
 \node[right] at (3.25,.375) {$\e_{j+1,\varnothing}$};
 \node at (0,-.93) {\textup{(a)}};
 \draw (5.55,1.6)--(9.55,1.6);
 \fill (5.55,1.6) circle[radius=1.4pt] (9.55,1.6) circle[radius=1.4pt];
 \node[above] at (7.55,1.6) {$L_i=4$};
 \draw[->] (7.55,1.28)--(7.55,.58);
 \draw (5.55,0)--(9.55,0);
 \fill (5.55,0) circle[radius=1.4pt] (9.55,0) circle[radius=1.4pt];
 \foreach \xx in {6.55,7.55,8.55}
   \fill (\xx,0) circle[radius=1.1pt];
 \foreach \xx in {6.05,7.05,8.05,9.05}
   \node[above] at (\xx,0) {$1$};
 \node[below] at (5.55,0) {$b_i$};
 \node[below] at (9.55,0) {$b_{i+1}$};
 \node at (7.55,-.93) {\textup{(b)}};
\end{tikzpicture}
\caption{\textup{(a)} A part of the tree $\G^{(1)}$: the chain vertices $b_{j-1}$, $b_j$, $b_{j+1}$ and the binary trees attached to them by their stem edges $ \e_{j-1,\varnothing}$, $ \e_{j,\varnothing}$, $ \e_{j+1,\varnothing}$. Open circles mark the defect vertices. \textup{(b)} In $\G$, a chain edge is divided into edges of length one; the drawing shows $L_i=4$. The three new vertices $c_{i,1}$, $c_{i,2}$, $c_{i,3}$ have degree two and do not belong to $\V_0$.}
\label{fig-ac-all-mass-tree-chain}
\end{figure}

\noindent\textit{Step~2: the branch on a binary tree with a stem edge and two half-lines.} Fix $d\geq1$ and put $N:=2^d$. 
Let $K_N=(\V_{K_N},\E_{K_N})$ be the finite metric tree given by
\[
 \V_{K_N}:=\{o\}\cup\{r_\omega: \omega\in\mathcal W_d\},\qquad
 \E_{K_N}:=\{\e_\omega: \omega\in\mathcal W_d\}.
\]
All edges have length one, and we set
\[
 \begin{aligned}
 \e_\varnothing&=(o,r_\varnothing),\\
 \e_{\omega a}&=(r_\omega,r_{\omega a})
 \qquad\text{for }\omega\in\mathcal W_d,\ |\omega|<d,\ a\in\{0,1\}.
 \end{aligned}
\]
Thus $K_N$ is a copy of the binary tree attached at a vertex $b_j$ with $d_j=d$, together with its stem edge $\e_\varnothing$, and $o$ corresponds to $b_j$. 
Let $\e_-$ and $\e_+$ be two distinct half-lines attached at origin $o$. Define $K_N^+=(\V_{K_N^+},\E_{K_N^+})$ by
\[
 \V_{K_N^+}:=\V_{K_N},\qquad
 \E_{K_N^+}:=\E_{K_N}\cup\{\e_-,\e_+\}.
\]
This is the middle component in Figure~\ref{fig-ac-all-mass-tree-cut}.
Then $K_N^+\in\mathbf G$ with $\ell_*=1$.
We take the $N$ leaves of $K_N^+$ as its defect set and write
\[
 \V_{N,0}:=S:=\{r_\omega: \omega\in\mathcal W_d,\ |\omega|=d\}.
\]
The vertex $o$ and the $N-1$ vertices $r_\omega$ with $|\omega|<d$ have degree three, and the leaves have degree one. Thus $K_N^+$ has $2N$ vertices, and the sum of their degrees is $4N$. We write $A^N_t$ for the normalized Dirichlet-to-Neumann operator of $K_N^+$ given by \eqref{ac-DtN-identity} on the weighted space $X(K_N^+):=\ell^2(\V_{K_N^+},\deg)$, as in Section~\ref{sect5}.

Define the admissible affine set by
\[
 H^1_N:=\{v\in H^1(K_N^+):v=1\text{ on }S\}.
\]
For $t>0$, consider the minimization of $\mathfrak b_t(v,v)=\int_{K_N^+}(\abs{v'}^2+t^2\abs v^2)\,dx$ over $H^1_N$.
The functional $v\mapsto\mathfrak b_t(v,v)$ is coercive by \eqref{ac-form-coercivity} and weakly lower semicontinuous as a sum of squared $L^2$ norms. If $v_n\rightharpoonup v$ in $H^1(K_N^+)$, then $v_n|_{K_N}\to v|_{K_N}$ in $C(K_N)$ by the compact embedding $H^1(K_N)\hookrightarrow C(K_N)$. Hence $v_n=1$ on $S$ implies $v=1$ on $S$, so the admissible affine set is closed under weak limits of sequences. The standard minimization method gives a minimizer $v_t\in H^1_N$. The set $H^1_N$ is convex, so strict convexity of the functional $v\mapsto\mathfrak b_t(v,v)$ gives uniqueness.
For every $\phi\in H^1(K_N^+)$ with $\phi=0$ on $S$, we have $v_t+s\phi\in H^1_N$ for every $s\in\R$. Therefore, since $v_t$ is the minimizer, we have
\[
 0=\frac12\left.\frac{d}{ds}\mathfrak b_t(v_t+s\phi,v_t+s\phi)\right|_{s=0}
 =\int_{K_N^+}\bigl(v_t'\phi'+t^2v_t\phi\bigr)\,dx .
\]
Testing first in the interior of each edge gives the edge equation, and integrating by parts with tests supported near a vertex outside $S$ gives the Kirchhoff condition. Together with $v_t\in H^1_N$, this yields
\[
 \begin{cases}
 -v_t''+t^2v_t=0 & \text{on every }\e\in\E_{K_N^+},\\
 v_t(\vv)=1 & \text{at every }\vv\in S,\\
 \displaystyle\sum_{\e\succ\vv}(v_t)'_\e(\vv)=0
 & \text{at every }\vv\in\V_{K_N^+}\setminus S.
 \end{cases}
\]
On each half-line, $v_t(x)=v_t(o)e^{-tx}$, which contributes $tv_t(o)^2$ to $\mathfrak b_t(v_t,v_t)$ by \eqref{ac-halfline-energy}. Hence the restriction of $v_t$ to $K_N$ minimizes
\[
 \int_{K_N}\bigl(\abs{v'}^2+t^2\abs v^2\bigr)\,dx+2t\abs{v(o)}^2
\]
among all $v\in H^1(K_N)$ with $v=1$ on $S$.
Put $\widetilde v_t:=\min\{\max\{v_t,0\},1\}$ on $K_N^+$. Then, $\widetilde v_t=1$ on $S$, and the chain rule on each edge gives
\[
 \abs{\widetilde v_t}\leq\abs{v_t},\qquad
 \abs{\widetilde v_t'}\leq\abs{v_t'}\qquad\text{a.e. on }K_N^+ .
\]
Thus $\widetilde v_t\in H^1_N$. Moreover, we have
\[
 \mathfrak b_t(v_t,v_t)\leq\mathfrak b_t(\widetilde v_t,\widetilde v_t)
 =\int_{K_N^+}\bigl(\abs{\widetilde v_t'}^2+t^2\abs{\widetilde v_t}^2\bigr)\,dx
 \leq\mathfrak b_t(v_t,v_t).
\]
Uniqueness gives $v_t=\widetilde v_t$, and hence $0\leq v_t\leq1$ on $K_N^+$.

Every automorphism of $K_N$ fixing $o$ extends to an automorphism of $K_N^+$ by acting as the identity on the stem edge $ \e_\varnothing$ and on the two half-lines $ \e_+, \e_-$. Each such automorphism $\gamma$ preserves the leaf set $S$, so $(v_t\circ\gamma)(\vv)=v_t(\gamma(\vv))=1$ for every $\vv\in S$. Thus $v_t\circ\gamma\in H^1_N$. Changing variables on the edges gives
\[
 \mathfrak b_t(v_t\circ\gamma,v_t\circ\gamma)
 =\int_{K_N^+}\bigl(\abs{v_t'}^2+t^2\abs{v_t}^2\bigr)\circ\gamma\,dx=\mathfrak b_t(v_t,v_t).
\]
By uniqueness, we obtain $v_t\circ\gamma=v_t$. Therefore, we obtain
\begin{equation}\label{app-tree-level-symmetry}
 v_t(r_\omega)=v_t(r_{\omega'}),\qquad
(v_t)'_{\e_\omega}(r_\omega)=(v_t)'_{\e_{\omega'}}(r_{\omega'})\qquad 
 \text{for }\omega,\omega'\in\mathcal W_d\text{ with }|\omega|=|\omega'|.
\end{equation}
Let $g_t:=\mathcal Rv_t\in X(K_N^+)$.  Since $v_t''=t^2v_t$ on every edge, Lemma~\ref{ac-restriction-extension} gives $v_t=U_tg_t$. The Kirchhoff conditions outside $S$ and \eqref{app-tree-level-symmetry} on $S$ give, by \eqref{ac-DtN-flux},
\begin{equation}\label{app-tree-CN-energy}
 A^N_tg_t=C_N(t)\mathbf1_S,\qquad
 NC_N(t)=\langle A^N_tg_t,g_t\rangle_X=\frac1t\mathfrak b_t(v_t,v_t),
\end{equation}
where $C_N(t)$ is the common value of $A^N_tg_t$ at the leaves. Since $g_t=1$ on $S$ and $0\leq g_t\leq1$, we have $N\leq\norm{g_t}_X^2\leq4N$. Applying \eqref{ac-DtN-coercivity} and \eqref{ac-DtN-bound} with $\ell_*=1$ gives
\begin{equation}\label{app-tree-CN-bounds}
 \tanh(t/2)\leq C_N(t)\leq4\coth(t/2)\qquad\text{for every }t>0.
\end{equation}
For $p=4$, define the nonlinear map on $K_N^+$ by
\[
 \mathcal N^N:X(K_N^+)\to X(K_N^+),\qquad
 \mathcal N^N(f):=\mathbf1_S f^3.
\]
This is the nonlinearity in \eqref{ac-vertex-equation} for the defect set $S$, whose vertices have degree one. Put $f^0_t:=\sqrt{C_N(t)}\,g_t$.
For $\vv\in S$, we have $\deg(\vv)=1$ and $g_t(\vv)=1$, so
\[
 (A^N_tf^0_t)(\vv)=\sqrt{C_N(t)}(A^N_tg_t)(\vv)
 =C_N(t)^{3/2}=f^0_t(\vv)^3=(\mathcal N^N(f^0_t))(\vv).
\]
For $\vv\notin S$, we have $(A^N_tg_t)(\vv)=0$ by the Kirchhoff condition and $\mathbf1_S(\vv)=0$ in the definition of $\mathcal N^N$. Hence
\[
 (A^N_tf^0_t)(\vv)=\sqrt{C_N(t)}(A^N_tg_t)(\vv)
 =0=(\mathcal N^N(f^0_t))(\vv).
\]
Thus $A^N_tf^0_t=\mathcal N^N(f^0_t)$, and the function
\[
 u_{N,t}:=t^{1/2}U_tf^0_t=\sqrt{tC_N(t)}\,v_t
\]
satisfies the edge equations and the vertex conditions of \eqref{eqpd} on $K_N^+$, with defects at $S$ and multiplier $t^2$. By \eqref{ac-mass-fixed-restriction}, its mass is
\[
 \mu_N(t):=\norm{u_{N,t}}_2^2=\frac12\bigl\langle(A^N_t+t\partial_tA^N_t)f^0_t,f^0_t\bigr\rangle_X .
\]
As in \eqref{ac-coordinate-subspace}, we define the subspace $Y$ and its orthogonal projection by
\[
 Y:=\{h\in X(K_N^+):h=0\text{ on }S\},\qquad
 P_Y:X(K_N^+)\longrightarrow Y,\qquad
 P_Yh:=\mathbf1_{\V_{K_N^+}\setminus S}h .
\]
For $t>0$, define
\[
 B_t^N:Y\longrightarrow Y,\qquad B_t^Nh:=P_YA^N_th
 \quad\text{for }h\in Y.
\]
By \eqref{ac-DtN-coercivity}, for every $h\in Y$,
\[
 \langle B_t^Nh,h\rangle_X=\langle A^N_th,h\rangle_X
 \geq\tanh(t/2)\norm h_X^2.
\]
Thus, $B_t^N$ is injective and, since $Y$ is finite-dimensional, invertible. The Cauchy--Schwarz inequality gives $\|h\|_X\leq \coth(t/2)\|B_t^Nh\|_X$, hence, $\norm{(B_t^N)^{-1}}_{\mathcal L(Y)}\leq\coth(t/2)$.
Moreover, since $g_t=1$ on $S$ and $A^N_tg_t=C_N(t)\mathbf1_S$, we have $g_t-\mathbf1_S\in Y$ and $P_YA^N_tg_t=0$. Thus, the values of $g_t$ outside $S$ are uniquely determined by the finite linear system
\[
 B_t^N(g_t-\mathbf1_S)=-P_YA^N_t\mathbf1_S
 \qquad\text{in }Y.
\]

Fix $s>0$. Since $g_t-g_s\in Y$ and $P_YA^N_sg_s=0$, we obtain
\[
 B_t^N(g_t-g_s)=-P_Y(A^N_t-A^N_s)g_s \qquad \text{and}\qquad
 \norm{g_t-g_s}_X\leq\coth(t/2)
       \norm{A^N_t-A^N_s}_{\mathcal L(X)}\norm{g_s}_X .
\]
Lemma~\ref{ac-DtN-estimates} gives the continuity of $t\mapsto A^N_t$ in operator norm. The last bound therefore gives $g_t\to g_s$ in $X(K_N^+)$. Moreover, since $\norm{\mathbf1_S}_X=\sqrt N$ and $A^N_tg_t=C_N(t)\mathbf1_S$, we have
\[
 \begin{aligned}
 \sqrt N\,|C_N(t)-C_N(s)|
 &=\norm{A^N_tg_t-A^N_sg_s}_X\\
 &\leq\norm{A^N_t}_{\mathcal L(X)}\norm{g_t-g_s}_X
       +\norm{A^N_t-A^N_s}_{\mathcal L(X)}\norm{g_s}_X
 \longrightarrow0\qquad\text{as }t\to s .
 \end{aligned}
\]
Thus $t\mapsto C_N(t)$ is continuous on $(0,+\infty)$. Since $t\mapsto g_t$ is continuous in $X(K_N^+)$, the map $t\mapsto f_t^0=\sqrt{C_N(t)}\,g_t$ is continuous from $(0,+\infty)$ to $X(K_N^+)$.

To compute the limit as $t\searrow0$, define the test function $\psi_t\in H^1_N$ by
\[
 \psi_t|_{K_N}\equiv1,\qquad
 \psi_t|_{\e_\pm}(x)=e^{-tx}\quad\text{for }x\geq0 .
\]
Since $|K_N|=2N-1$ and $v_t(x)=v_t(o)e^{-tx}$ on each half-line, \eqref{ac-halfline-energy} gives
\[
 \begin{aligned}
 \mathfrak b_t(\psi_t,\psi_t)
 &=t^2|K_N|+4t^2\int_0^\infty e^{-2tx}\,dx=t^2|K_N|+2t,\\
 \mathfrak b_t(v_t,v_t)
 &=\int_{K_N}\bigl(\abs{v_t'}^2+t^2\abs{v_t}^2\bigr)\,dx+2tv_t(o)^2 .
 \end{aligned}
\]
Since $v_t$ minimizes $\mathfrak b_t(v,v)$ over $H^1_N$, the preceding energy identities and \eqref{app-tree-CN-energy} give
\[
 \norm{v_t'}_{L^2(K_N)}^2\leq\mathfrak b_t(v_t,v_t)\leq 2t+t^2|K_N|\qquad
\text{and}\qquad 
 2v_t(o)^2\leq NC_N(t)=\frac1t\mathfrak b_t(v_t,v_t)\leq 2+t|K_N|.
\]
Every $x\in K_N$ can be joined to a leaf by a path of length at most $d+1$. Since $v_t(\vv)=1$ at that leaf $\vv \in S$, integration along the path and the Cauchy--Schwarz inequality give
\[
 |v_t(x)-1|\leq (d+1)^{1/2}\norm{v_t'}_{L^2(K_N)}
 \leq(d+1)^{1/2}\bigl(2t+t^2|K_N|\bigr)^{1/2}.
\]
In particular, $v_t(o)\to1$ and the bounds on $NC_N(t)$ give $C_N(t)\to2/N$. Using $0\leq v_t\leq1$ on $K_N$ and integrating $v_t(o)^2e^{-2tx}$ on the two half-lines, we obtain
\[
 \begin{aligned}
 \mu_N(t)&=tC_N(t)\int_{K_N}\abs{v_t}^2\,dx+
          2tC_N(t)v_t(o)^2\int_0^\infty e^{-2tx}\,dx\\
 &=tC_N(t)\int_{K_N}\abs{v_t}^2\,dx+C_N(t)v_t(o)^2
 \longrightarrow\frac2N
 \qquad\text{as }t\searrow0 .
 \end{aligned}
\]

For the limit as $t\to+\infty$, let $P_S=I-P_Y$. Since $P_Sg_t=\mathbf1_S$ and $A^N_tg_t=C_N(t)\mathbf1_S$, we have
\[
 g_t-\mathbf1_S=(I_X-P_S)(I_X-A^N_t)g_t,\qquad
 (C_N(t)-1)\mathbf1_S=P_S(A^N_t-I_X)g_t.
\]
By \eqref{ac-DtN-exponential}, for every $t\geq1$,
\[
 \max\bigl\{\norm{g_t-\mathbf1_S}_X,\norm{(C_N(t)-1)\mathbf1_S}_X\bigr\}
 \leq C_0 e^{-t/2}\norm{g_t}_X
 \leq2C_0\sqrt N\,e^{-t/2}.
\]
Hence $g_t\to\mathbf1_S$ in $X(K_N^+)$ and $C_N(t)\to1$, so $f^0_t=\sqrt{C_N(t)}\,g_t\to\mathbf1_S$ in $X(K_N^+)$.

\noindent\textit{Step~3: reconnecting the isolated tree.} With $\G$ and $\V_0$ fixed, let $\mu>0$ be the prescribed mass. Choose $d\geq1$ such that $N=2^d$ satisfies $2/N<\mu<N/2$, and then choose $t_->0$ with $\mu_N(t_-)<\mu$, using the limit from Step~2. Put $J:=[t_-,+\infty)$. For this $d$, infinitely many indices $j$ satisfy $d_j=d$. We construct the branches on $J$ at those indices for which $L_{j-1}$ and $L_j$ are sufficiently large.

Let $X:=\ell^2(\V^{(1)},\deg)$ and let $A_t:X\to X$ be the normalized Dirichlet-to-Neumann operator of $\G^{(1)}$ given by \eqref{ac-DtN-identity}. Since the vertices of $\V_0$ have degree one, the nonlinear map in \eqref{ac-vertex-equation} for $p=4$ is
\[
 \mathcal N:X\to X,\qquad \mathcal N(f):=\mathbf1_{\V_0}f^3.
\]

For each $j$ with $d_j=d$ and each $i\in\{j-1,j\}$, let $\e_i^-$ and $\e_i^+$ be two distinct half-lines with origins $b_i$ and $b_{i+1}$, respectively. Define $\widetilde\G_j=(\widetilde\V_j,\widetilde\E_j)$ by
\[
 \widetilde\V_j:=\V^{(1)},\qquad
 \widetilde\E_j:=\bigl(\E^{(1)}\setminus\{\e_{j-1},\e_j\}\bigr)
       \cup\{\e_i^-,\e_i^+:i\in\{j-1,j\}\}.
\]
The graph $\widetilde\G_j$ has the three components shown in Figure~\ref{fig-ac-all-mass-tree-cut}. Each component belongs to $\mathbf G$ with $\ell_*=1$.
Let $\widetilde\G_j^-$ and $\widetilde\G_j^+$ denote the components containing $b_{j-1}$ and $b_{j+1}$, respectively. We identify the middle component with $K_N^+$ by identifying $o$ with $b_j$ and $r_\omega$ with $r_{j,\omega}$. Let $A^\pm_{j,t}$ be the normalized Dirichlet-to-Neumann operator of $\widetilde\G_j^\pm$ on $X(\widetilde\G_j^\pm):=\ell^2(\V_{\widetilde\G_j^\pm},\deg_{\widetilde\G_j^\pm})$, defined by \eqref{ac-DtN-identity}. The corresponding operator on $K_N^+$ is $A_t^N$. Because $\widetilde\G_j$ has the same vertex set and vertex degrees as $\G^{(1)}$, restriction to the three components identifies
\[
 X=X(\widetilde\G_j^-)\oplus X(K_N^+)\oplus X(\widetilde\G_j^+)
\]
as an orthogonal direct sum. On this space we define
\[
 A^0_{j,t}:=A^-_{j,t}\oplus A_t^N\oplus A^+_{j,t}.
\]
Since 
the defect set is still $\V_0$, $\mathcal N$ also gives the nonlinearity on $\widetilde\G_j$. Under the identification of the middle summand of $X$ with $X(K_N^+)$, the restriction of $\mathcal N$ to this $X(K_N^+)$ agrees with $\mathcal N^N$.

Applying \eqref{ac-DtN-coercivity}, \eqref{ac-DtN-bound} and \eqref{ac-DtN-derivative-bound} on each component shows that these estimates hold for $A^0_{j,t}$ with $\ell_*=1$. Let $f^0_j(t)\in X$ be equal to $f^0_t$ on the isolated component  $K_N^+$  and to zero  on  $\widetilde\G_j^-$ and $\widetilde\G_j^+$. Then $f^0_j(t)\geq0$, $A^0_{j,t}f^0_j(t)=\mathcal N(f^0_j(t))$. The limit at the end of Step~2 gives
\begin{equation}\label{app-tree-f0-limit}
 f^0_j(t)\longrightarrow\mathbf 1_{S_j}
 \quad\text{in }X, \qquad \text{as }t\to+\infty .
\end{equation}
Moreover, \eqref{app-tree-CN-bounds} and $\norm{g_t}_X\leq2\sqrt N$ give
\begin{equation}\label{app-tree-f0-bound}
 \norm{f^0_j(t)}_X\leq M_N:=4\bigl(N\coth(t_-/2)\bigr)^{1/2}\qquad\text{for every }t\in J .
\end{equation}

\begin{figure}[H]
\centering
\begin{tikzpicture}[x=.94cm,y=.87cm,line cap=round,line join=round,every node/.style={font=\small,inner sep=2pt}]
 \draw[densely dotted] (-6.70,0)--(-4.85,0) (4.85,0)--(6.70,0);
 \draw[->,thick] (-4.85,0)--node[midway,above] {$\e_{j-1}^-$} (-2.60,0);
 \draw[->,thick] (0,0)--node[midway,above] {$\e_{j-1}^+$} (-1.85,0);
 \draw[->,thick] (0,0)--node[midway,above] {$\e_j^-$} (1.85,0);
 \draw[->,thick] (4.85,0)--node[midway,above] {$\e_j^+$} (2.60,0);
 \foreach \xx in {-2.60,-1.85,1.85,2.60}
   \node[above] at (\xx,0) {$\infty$};
 \node[below] at (-4.85,0) {$b_{j-1}$};
 \node[below] at (0,0) {$o=b_j$};
 \node[below] at (4.85,0) {$b_{j+1}$};
 \foreach \xx in {-4.85,4.85}{
   \begin{scope}[shift={(\xx,0)}]
    \draw (0,0)--(0,.75)--(-.55,1.45) (0,.75)--(.55,1.45);
    \fill (0,0) circle[radius=1.4pt] (0,.75) circle[radius=1.4pt];
    \draw[fill=white] (-.55,1.45) circle[radius=2pt] (.55,1.45) circle[radius=2pt];
   \end{scope}
 }
 \draw (0,0)--(0,.65)--(-.63,1.30)--(-.99,2)
       (-.63,1.30)--(-.27,2) (0,.65)--(.63,1.30)--(.27,2)
       (.63,1.30)--(.99,2);
 \fill (0,0) circle[radius=1.4pt] (0,.65) circle[radius=1.4pt]
       (-.63,1.30) circle[radius=1.4pt] (.63,1.30) circle[radius=1.4pt];
 \foreach \xx in {-.99,-.27,.27,.99}
   \draw[fill=white] (\xx,2) circle[radius=2pt];
 \node at (-4.85,-.85) {Left component $\widetilde\G_j^-$};
 \node at (0,-.85) {Isolated component $K_N^+$};
 \node at (4.85,-.85) {Right component $\widetilde\G_j^+$};
\end{tikzpicture}
\caption{The graph $\widetilde\G_j$: the two chain edges at $b_j$ are removed, and the arrows mark the four new half-lines. The isolated component consists of $K_N$ and two half-lines with origin $o=b_j$. The drawing shows $d=2$ and $N=4$.}
\label{fig-ac-all-mass-tree-cut}
\end{figure}

By \eqref{ac-DtN-form} and \eqref{ac-DtN-identity}, the operator $\Delta_{j,t}:=A_t-A^0_{j,t}$  satisfies that, for $f,g\in X$,
\begin{equation}\label{app-tree-dtn-difference}
 \begin{aligned}
 \langle\Delta_{j,t}f,g\rangle_X=\sum_{i\in\{j-1,j\}}\Bigl[&\bigl(\coth(tL_i)-1\bigr)\bigl(f(b_i)g(b_i)+f(b_{i+1})g(b_{i+1})\bigr)\\
 &-\operatorname{csch}(tL_i)\bigl(f(b_i)g(b_{i+1})+f(b_{i+1})g(b_i)\bigr)\Bigr].
 \end{aligned}
\end{equation}
For $t\geq t_-$ we have $tL_i\geq2t_-$ since $L_i\geq2$. Hence
$ 0\leq\coth(tL_i)-1
 \leq\operatorname{csch}(tL_i)
 =\frac{2e^{-tL_i}}{1-e^{-2tL_i}}
 \leq\frac{2e^{-tL_i}}{1-e^{-4t_-}}
$
and $\coth(tL_i)+\operatorname{csch}(tL_i)=\coth(tL_i/2)\leq\coth t_-$. Moreover, we have $(\abs{f(b_i)}+\abs{f(b_{i+1})})(\abs{g(b_i)}+\abs{g(b_{i+1})})\leq2\norm f_X\norm g_X$. Differentiating in $t$ by means of \eqref{ac-DtN-derivatives}, we obtain
\begin{equation}\label{app-tree-epsilon}
 \begin{aligned}
 \norm{\Delta_{j,t}}_{\mathcal L(X)}
 &\leq\frac{4}{1-e^{-4t_-}}\sum_{i\in\{j-1,j\}}e^{-tL_i},
 \qquad t\in J,\\
 \norm{\partial_t\Delta_{j,t}}_{\mathcal L(X)}
 &\leq\frac{4\coth(t_-)}{1-e^{-4t_-}}\sum_{i\in\{j-1,j\}}L_i e^{-tL_i},
 \qquad t\in J,\\
 \varepsilon_j
 &:=\sup_{t\in J}\bigl[\norm{\Delta_{j,t}}_{\mathcal L(X)}
                    +\norm{\partial_t\Delta_{j,t}}_{\mathcal L(X)}\bigr]\leq C_\Delta\sum_{i\in\{j-1,j\}}(1+L_i)e^{-t_-L_i},
 \end{aligned}
\end{equation}
where $C_\Delta:=8\coth(t_-)/(1-e^{-4t_-})$ depends only on $t_-$. Since $L_i\to+\infty$ as $|i|\to\infty$, we have $\varepsilon_j\to0$ as $|j|\to\infty$. Moreover, for each fixed $j\in \mathbb{Z}$, $\norm{\Delta_{j,t}}_{\mathcal L(X)}\to 0$ as $t\to+\infty$.

Define
\[
 X_j:=\bigl\{h\in X:h(r_{j,\omega})=h(r_{j,\omega'})
 \text{ for all }\omega,\omega'\in\mathcal W_{d_j}\text{ with }|\omega|=|\omega'|\bigr\}.
\]
By \eqref{ac-norm-ordering}, $X_j$ is closed in $X$. By \eqref{ac-DtN-coordinate}, $X_j$ is an invariant subspace for $A^0_{j,t}$, since vertices at the same level have identical edge lengths and the same numbers of neighbors at each level. Thus $A^0_{j,t}|_{X_j}$ is self-adjoint by \cite[Proposition~II.3.7]{Conway}. By \eqref{ac-DtN-coercivity} with $\ell_*=1$ and the Cauchy--Schwarz inequality,
\[
 \tanh(t/2)\norm h_X^2
 \leq\langle A^0_{j,t}h,h\rangle_X
 \leq\norm{A^0_{j,t}h}_X\norm h_X
 \qquad\text{for every }h\in X_j.
\]
The lower bound on $\langle A^0_{j,t}h,h\rangle_X$ and \cite[Proposition~6.9]{Ha} give the invertibility of $A^0_{j,t}|_{X_j}$. Moreover, these inequalities give $\norm{(A^0_{j,t})^{-1}}_{\mathcal L(X_j)}\leq\coth(t/2)$.

Let $g_{j,t}\in X$ be equal to $g_t$ on the isolated component $K_N^+$ and to zero on $\widetilde\G_j^-$ and $\widetilde\G_j^+$. By \eqref{app-tree-level-symmetry}, $g_{j,t}\in X_j$.
As in \eqref{ac-coordinate-subspace}, define the orthogonal projection
\[
 P_{S_j}:X_j\to\operatorname{span}\{\mathbf1_{S_j}\},\qquad
 P_{S_j}h:=\mathbf1_{S_j}h.
\]
Since $h$ is constant on $S_j$ for every $h\in X_j$ and the $N$ vertices of $S_j$ have degree one, we have
\[
 P_{S_j}h=\frac{\langle h,\mathbf1_{S_j}\rangle_X}{N}\mathbf1_{S_j}.
\]
We have $A^0_{j,t}g_{j,t}=C_N(t)\mathbf1_{S_j}$, $P_{S_j}g_{j,t}=\mathbf1_{S_j}$ and $f^0_j(t)=\sqrt{C_N(t)}\,g_{j,t}$. Thus \eqref{ac-N-derivative} with $p=4$ gives
\[
 L_j^0(t):=A^0_{j,t}-\mathcal N'(f^0_j(t))=A^0_{j,t}-3C_N(t)P_{S_j}
 \quad\text{on }X_j.
\]
In particular, $L_j^0(t)g_{j,t}=-2C_N(t)\mathbf1_{S_j}$. For $h\in X_j$, self-adjointness of $(A^0_{j,t})^{-1}$ gives
\[
 \bigl\langle(A^0_{j,t})^{-1}h,\mathbf1_{S_j}\bigr\rangle_X
 =\bigl\langle h,(A^0_{j,t})^{-1}\mathbf1_{S_j}\bigr\rangle_X
 =\frac{\langle h,g_{j,t}\rangle_X}{C_N(t)} .
\]
Substitution into $L_j^0(t)(A^0_{j,t})^{-1}h$ now gives
\[
 L_j^0(t)(A^0_{j,t})^{-1}h=h-\frac3N\langle h,g_{j,t}\rangle_X\,\mathbf1_{S_j}
\]
Define
\[
 R_{j,t}:X_j\to X_j,\qquad
 R_{j,t}h:=(A^0_{j,t})^{-1}h-\frac{3\langle h,g_{j,t}\rangle_X}{2NC_N(t)}\,g_{j,t}.
\]
The preceding identity and $L_j^0(t)g_{j,t}=-2C_N(t)\mathbf1_{S_j}$ give $L_j^0(t)R_{j,t}=I_{X_j}$. The operators $(A^0_{j,t})^{-1}$ and $h\mapsto\langle h,g_{j,t}\rangle_Xg_{j,t}$ are self-adjoint on $X_j$, so $R_{j,t}$ is self-adjoint. Since $L_j^0(t)$ is also self-adjoint, taking adjoints gives
\[
 R_{j,t}L_j^0(t)=\bigl(L_j^0(t)R_{j,t}\bigr)^*=I_{X_j}.
\]
Thus $R_{j,t}=L_j^0(t)^{-1}$. Moreover, $\langle A^0_{j,t}g_{j,t},g_{j,t}\rangle_X=NC_N(t)$, so \eqref{ac-DtN-coercivity} gives $\norm{g_{j,t}}_X^2\leq NC_N(t)\coth(t/2)$. Consequently,
\begin{equation}\label{app-tree-inverse}
 \begin{aligned}
 \norm{L_j^0(t)^{-1}}_{\mathcal L(X_j)}
 &\leq\coth(t/2)+\frac{3\norm{g_{j,t}}_X^2}{2NC_N(t)}\\
 &\leq\frac52\coth(t/2)\leq\Lambda_N:=\frac52\coth(t_-/2)
 \qquad\text{for every }t\in J .
 \end{aligned}
\end{equation}
The bound $\Lambda_N$ depends only on $t_-$.

 To solve $A_tf=\mathcal N(f)$ in $X_j$, we write $f=f^0_j(t)+h$ for $h \in X_j$. Since $A_t=A^0_{j,t}+\Delta_{j,t}$ and $A^0_{j,t}f^0_j(t)=\mathcal N(f^0_j(t))$, this equation is equivalent to $L_j^0(t)h=Q_{j,t}(h)-\Delta_{j,t}(f^0_j(t)+h)$, where
\[
 Q_{j,t}(h):=\mathcal N\bigl(f^0_j(t)+h\bigr)-\mathcal N\bigl(f^0_j(t)\bigr)-\mathcal N'\bigl(f^0_j(t)\bigr)h ,
\]
that is, to the fixed-point equation
\[
 h=\mathcal F_{j,t}(h):=L_j^0(t)^{-1}\Bigl[Q_{j,t}(h)-\Delta_{j,t}\bigl(f^0_j(t)+h\bigr)\Bigr],\qquad h\in X_j .
\]
By \eqref{app-tree-dtn-difference}, $\operatorname{supp}(\Delta_{j,t}f)\subset\{b_{j-1},b_j,b_{j+1}\}$ for every $f\in X$, so $\Delta_{j,t}(X_j)\subset X_j$. Since $\mathcal N(f)=\mathbf1_{\V_0}f^3$ and $\V_0$ contains precisely the leaves of the attached trees, we have $\mathcal N(X_j)\subset X_j$. Moreover, $f_j^0(t)\in X_j$ and $\mathcal N'(f_j^0(t))h=3C_N(t)P_{S_j}h$ for $h\in X_j$. Thus $Q_{j,t}(h)\in X_j$ and $\mathcal F_{j,t}(h)\in X_j$ for every $h\in X_j$.

Let $0<r\leq1$, and let $h_1,h_2\in X_j$ satisfy $\norm{h_1}_X,\norm{h_2}_X\leq r$. The fundamental theorem of calculus gives
\[
 Q_{j,t}(h_1)-Q_{j,t}(h_2)=\int_0^1\Bigl[\mathcal N'\bigl(f^0_j(t)+h_2+\sigma(h_1-h_2)\bigr)-\mathcal N'\bigl(f^0_j(t)\bigr)\Bigr](h_1-h_2)\,d\sigma .
\]
By \eqref{app-tree-f0-bound}, $\norm{f_j^0(t)}_X\leq M_N$, and for every $\sigma\in[0,1]$,
\[
 \begin{aligned}
 \norm{f_j^0(t)+(1-\sigma)h_2+\sigma h_1}_X
 &\leq M_N+(1-\sigma)\norm{h_2}_X+\sigma\norm{h_1}_X\\
 &\leq M_N+r\leq M_N+1=:R.
 \end{aligned}
\]
Hence \eqref{ac-N-lipschitz} with $p=4$, whose constant is $C_R=3\max\{1,2R\}=6(M_N+1)$, gives $\norm{Q_{j,t}(h_1)-Q_{j,t}(h_2)}_X\leq6(M_N+1)\,r\norm{h_1-h_2}_X$. Since $\mathcal F_{j,t}(0)=-L_j^0(t)^{-1}\Delta_{j,t}f^0_j(t)$, \eqref{app-tree-epsilon} and \eqref{app-tree-inverse} give
\begin{equation}\label{app-tree-contraction}
 \begin{gathered}
 \norm{\mathcal F_{j,t}(0)}_X\leq\Lambda_NM_N\norm{\Delta_{j,t}}_{\mathcal L(X)},\\
 \norm{\mathcal F_{j,t}(h_1)-\mathcal F_{j,t}(h_2)}_X\leq\Lambda_N\bigl(6(M_N+1)r+\varepsilon_j\bigr)\norm{h_1-h_2}_X .
 \end{gathered}
\end{equation}
Put $r_j:=2\Lambda_NM_N\varepsilon_j$ and $b_-:=\coth(t_-/2)+t_-\tilde C(t_-,1)$, where $\tilde C(t_-,1)$ is the constant in \eqref{ac-DtN-derivative-bound} for $I=[t_-,t_-+1]$ and $\ell_*=1$; thus $b_-$ depends only on $t_-$. Since $\varepsilon_j\to0$, there is $j_0\geq1$ such that every $j\in\mathbb Z$ with $|j|\geq j_0$ and $d_j=d$ satisfies
\begin{equation}\label{app-tree-smallness}
 \begin{gathered}
 r_j\leq1,\qquad
 \Lambda_N\bigl(6(M_N+1)r_j+\varepsilon_j\bigr)\leq\tfrac12,\qquad
 r_j<\tfrac12\sqrt{\tanh(t_-/2)},\\
 \tfrac12\bigl[\max\{1,t_-\}(M_N+1)^2\varepsilon_j+b_-(2M_N+1)r_j\bigr]<\mu-\mu_N(t_-).
 \end{gathered}
\end{equation}
Fix such a $j\in\mathbb Z$ and define the closed ball
\[
 \overline B_{X_j}(0,r_j):=\{h\in X_j:\norm h_X\leq r_j\}.
\]
For every $t\in J$, \eqref{app-tree-contraction} with $r=r_j$ and \eqref{app-tree-smallness} show that $\mathcal F_{j,t}$ is a $\frac12$-contraction on $\overline B_{X_j}(0,r_j)$ and that
\[
 \norm{\mathcal F_{j,t}(h)}_X
 \leq\norm{\mathcal F_{j,t}(0)}_X+\tfrac12\norm h_X
 \leq\tfrac12r_j+\tfrac12r_j=r_j
 \qquad\text{for every }h\in\overline B_{X_j}(0,r_j).
\]
Thus $\mathcal F_{j,t}:\overline B_{X_j}(0,r_j)\to\overline B_{X_j}(0,r_j)$ has a unique fixed point $h_j(t)$ in the ball $\overline B_{X_j}(0,r_j)$ by the contraction mapping theorem. Put $f_j(t):=f^0_j(t)+h_j(t)$. Since $\norm{h_j(t)}_X\leq\frac12\norm{h_j(t)}_X+\norm{\mathcal F_{j,t}(0)}_X$, we obtain, for every $t\in J$,
\begin{equation}\label{app-tree-close}
 A_tf_j(t)=\mathcal N\bigl(f_j(t)\bigr),\qquad
 \norm{f_j(t)-f^0_j(t)}_X\leq2\Lambda_NM_N\norm{\Delta_{j,t}}_{\mathcal L(X)}\leq r_j .
\end{equation}
By Lemma~\ref{ac-DtN-estimates}, applied to $\G^{(1)}$ and to each component of $\widetilde\G_j$, the maps $t\mapsto A_t$ and $t\mapsto A^0_{j,t}$ are continuous in $\mathcal L(X)$, and so is $t\mapsto\Delta_{j,t}$. The identity $L_j^0(t)=A^0_{j,t}-3C_N(t)P_{S_j}$ and the continuity of $t\mapsto C_N(t)$ give continuity of $t\mapsto L_j^0(t)$ in $\mathcal L(X_j)$. For $t,\tau\in J$, we have
\[
 L_j^0(t)^{-1}-L_j^0(\tau)^{-1}
 =L_j^0(t)^{-1}\bigl[L_j^0(\tau)-L_j^0(t)\bigr]L_j^0(\tau)^{-1}.
\]
Taking norms and applying \eqref{app-tree-inverse} to both inverses gives
\[
 \norm{L_j^0(t)^{-1}-L_j^0(\tau)^{-1}}_{\mathcal L(X_j)}
 \leq\Lambda_N^2\norm{L_j^0(t)-L_j^0(\tau)}_{\mathcal L(X_j)}.
\]
Thus $t\mapsto L_j^0(t)^{-1}$ is continuous in $\mathcal L(X_j)$.

Moreover, for each fixed $h\in X_j$, by the continuity of $t \mapsto f_j^0(t)$ and $\mathcal N \in C(X,X)$ and by the identity $\mathcal N'(f_j^0(t))h=3C_N(t)P_{S_j}h$, the map $t\mapsto Q_{j,t}(h)$ is continuous. The definition of $\mathcal F_{j,t}$ and the continuity of $L_j^0(t)^{-1}$ and $\Delta_{j,t}$ then give continuity of $t\mapsto\mathcal F_{j,t}(h)$.

For $t,\tau\in J$, the fixed-point identities and the contraction estimate give
\[
 \begin{aligned}
 \norm{h_j(t)-h_j(\tau)}_X
 &=\norm{\mathcal F_{j,t}(h_j(t))-\mathcal F_{j,\tau}(h_j(\tau))}_X\\
 &\leq\norm{\mathcal F_{j,t}(h_j(t))-\mathcal F_{j,t}(h_j(\tau))}_X+\norm{\mathcal F_{j,t}(h_j(\tau))-\mathcal F_{j,\tau}(h_j(\tau))}_X\\
 &\leq\tfrac12\norm{h_j(t)-h_j(\tau)}_X
       +\norm{\mathcal F_{j,t}(h_j(\tau))-\mathcal F_{j,\tau}(h_j(\tau))}_X.
 \end{aligned}
\]
Hence
\[
 \norm{h_j(t)-h_j(\tau)}_X
 \leq2\norm{\mathcal F_{j,t}(h_j(\tau))-\mathcal F_{j,\tau}(h_j(\tau))}_X.
\]
Letting $t\to\tau$ gives continuity of $h_j$ on $J$. Hence $f_j=f_j^0+h_j$ is continuous on $J$.

Let $t\in J$ and $f:=f_j(t)$. We show that $f>0$ at every vertex. Put $f^+:=\max\{f,0\}$ and $f^-:=\min\{f,0\}$. All edges of $\G^{(1)}$ are bounded and $f^+f^-=0$, so \eqref{ac-DtN-coordinate} gives
\[
 \langle A_tf^+,f^-\rangle_X=-\sum_{\vv\in\V}f^-(\vv)\sum_{\e\succ\vv}\operatorname{csch}(t\ell_\e)f^+(\w_{\e,\vv})\geq0 .
\]
Since $f=f^++f^-$, \eqref{ac-DtN-coercivity} and the equation in \eqref{app-tree-close} yield
\[
 \begin{aligned}
 \tanh(t_-/2)\norm{f^-}_X^2
 &\leq\langle A_tf^-,f^-\rangle_X
 =\langle A_tf,f^-\rangle_X-\langle A_tf^+,f^-\rangle_X\\
 &\leq\langle A_tf,f^-\rangle_X
 =\langle\mathcal N(f),f^-\rangle_X
 =\sum_{\vv\in\V_0}\abs{f^-(\vv)}^4
 \leq\norm{f^-}_\infty^2\norm{f^-}_X^2 .
 \end{aligned}
\]
Since $f^0_j(t)\geq0$, for every $\vv\in\V$ we have
\[
 \abs{f^-(\vv)}
 =\max\{-f(\vv),0\}
 \leq\max\{f_j^0(t)(\vv)-f(\vv),0\}
 \leq\abs{f_j^0(t)(\vv)-f(\vv)}.
\]
Then \eqref{ac-norm-ordering}, \eqref{app-tree-close} and \eqref{app-tree-smallness} give $\norm{f^-}_\infty^2\leq r_j^2<\tanh(t_-/2)$. Therefore $f^-=0$. At every $\vv\in S_j$, \eqref{ac-norm-ordering}, \eqref{app-tree-CN-bounds}, \eqref{app-tree-close} and \eqref{app-tree-smallness} give
\begin{equation}\label{app-tree-leaf-bound}
 f_j(t)(\vv)\geq\sqrt{C_N(t)}-r_j\geq\sqrt{\tanh(t_-/2)}-r_j>\tfrac12\sqrt{\tanh(t_-/2)} .
\end{equation}
If $f(\vv)=0$ at some vertex $\vv$, then $(\mathcal N(f))(\vv)=0$, \eqref{app-tree-close} and \eqref{ac-DtN-coordinate} give
\[
 0=\deg(\vv)(A_tf)(\vv)=-\sum_{\e\succ\vv}\operatorname{csch}(t\ell_\e)f(\w_{\e,\vv});
\]
hence $f$ vanishes at every neighbor of $\vv$. As $\G^{(1)}$ is connected, $f$ would then vanish at every vertex, which contradicts \eqref{app-tree-leaf-bound}. Thus $f_j(t)>0$ at every vertex. By \eqref{ac-vertex-equation}, the function
\[
 u_j(t):=t^{1/2}U_tf_j(t)
\]
satisfies the edge equations and the vertex conditions of \eqref{eqpd} on $\G^{(1)}$ with multiplier $t^2$. Moreover, $u_j(t)>0$ on $\G^{(1)}$ by \eqref{ac-edge-extension}, since the two coefficients there are nonnegative and not both zero.

By the continuity of $t \mapsto f_j(t)$ and \eqref{ac-extension-estimates}--\eqref{ac-extension-lipschitz}, we have $t \mapsto u_j(t)\in C(J,H^1(\G^{(1)}))$.

\noindent\textit{Step~4: prescribed mass and inequivalence.}
By \eqref{ac-mass-fixed-restriction}, the mass
\[
 \norm{u_j(t)}_2^2=t\norm{U_tf_j(t)}_2^2=\frac12\bigl\langle(A_t+t\partial_tA_t)f_j(t),f_j(t)\bigr\rangle_X
\]
depends continuously on $t\in J$, because $t\mapsto A_t$ is of class $C^2$ by Lemma~\ref{ac-DtN-estimates} and $f_j$ is continuous. We compare it at $t=t_-$ with $\mu_N(t_-)$. Put $B:=A_{t_-}+t_-\partial_tA_{t_-}$, $B^0:=A^0_{j,t_-}+t_-\partial_tA^0_{j,t_-}$. Since $f^0_j(t_-)$ is supported on the isolated component, $\mu_N(t_-)=\frac12\langle B^0f^0_j(t_-),f^0_j(t_-)\rangle_X$. Moreover, $\norm{B-B^0}_{\mathcal L(X)}\leq\max\{1,t_-\}\varepsilon_j$ by \eqref{app-tree-epsilon}, and $\norm{B^0}_{\mathcal L(X)}\leq b_-$ by \eqref{ac-DtN-bound} and \eqref{ac-DtN-derivative-bound} on each component of $\widetilde\G_j$. We write
\[
 \begin{aligned}
 \langle Bf_j(t_-),f_j(t_-)\rangle_X
  -\langle B^0f^0_j(t_-),f^0_j(t_-)\rangle_X
 &=\langle(B-B^0)f_j(t_-),f_j(t_-)\rangle_X\\
 &\quad+\langle B^0[f_j(t_-)-f^0_j(t_-)],f_j(t_-)\rangle_X\\
 &\quad+\langle B^0f^0_j(t_-),f_j(t_-)-f^0_j(t_-)\rangle_X .
 \end{aligned}
\]
Using $\norm{f^0_j(t_-)}_X\leq M_N$ by \eqref{app-tree-f0-bound}, $\norm{f_j(t_-)-f^0_j(t_-)}_X\leq r_j\leq1$ and \eqref{app-tree-smallness}, we get
\[
 \bigl|\norm{u_j(t_-)}_2^2-\mu_N(t_-)\bigr|\leq\tfrac12\bigl[\max\{1,t_-\}(M_N+1)^2\varepsilon_j+b_-(2M_N+1)r_j\bigr]<\mu-\mu_N(t_-),
\]
so $\norm{u_j(t_-)}_2^2<\mu$. As $t\to+\infty$, \eqref{app-tree-close} and the decay of $\norm{\Delta_{j,t}}_{\mathcal L(X)}$ in \eqref{app-tree-epsilon} give $\norm{f_j(t)-f^0_j(t)}_X\to0$. Combining this with \eqref{app-tree-f0-limit} gives $f_j(t)\to\mathbf 1_{S_j}$ in $X$. By \eqref{ac-DtN-exponential} on $\G^{(1)}$, $\norm{A_t+t\partial_tA_t-I_X}_{\mathcal L(X)}\leq(C_0+tC_1)e^{-t/2}$ for $t\geq1$, which tends to zero as $t\to+\infty$. Since the leaves have degree one,
\[
 \lim_{t\to+\infty}\norm{u_j(t)}_2^2=\frac12\norm{\mathbf 1_{S_j}}_X^2=\frac N2>\mu .
\]
By the intermediate value theorem, there is $t_j>t_-$ with $\norm{u_j(t_j)}_2^2=\mu$. Then $(u_j(t_j),t_j^2)\in H^1_\mu(\G^{(1)})\times(0,+\infty)$ is a positive solution of \eqref{eqpd} on $\G^{(1)}$ and hence on $\G$.

Since the trees attached to $\G^{(0)}$ are finite, $\G^{(0)}$ is the unique bi-infinite simple path in $\G$. By construction, we have
\[
 \begin{aligned}
 \V^{(0)}&=\G^{(0)}\cap\overline{\G\setminus\G^{(0)}},\\
 \gamma(\G^{(0)})&=\G^{(0)},\qquad \gamma(\V^{(0)})=\V^{(0)}
 \qquad\text{for every }\gamma\in\Gamma,
 \end{aligned}
\]
where the closure is taken in $\G$. The connected components of $\G^{(0)}\setminus\V^{(0)}$ are the intervals $(b_i,b_{i+1})$, $i\in\mathbb Z$. For each $i\in\mathbb Z$, there is $j\in\mathbb Z$ such that $\gamma((b_i,b_{i+1}))=(b_j,b_{j+1})$. Since $\gamma$ is a homeomorphism that preserves lengths, it also maps the closed intervals onto one another, and
\[
 \gamma([b_i,b_{i+1}])=[b_j,b_{j+1}],
 \qquad L_i=L_j .
\]
The pairwise distinctness of  $L_i$ gives $j=i$.
It follows that, for every $i\in\mathbb Z$,
\[
 \{\gamma(b_i)\}
 =\gamma([b_{i-1},b_i])\cap\gamma([b_i,b_{i+1}])
 =[b_{i-1},b_i]\cap[b_i,b_{i+1}]
 =\{b_i\}.
\]
Let $T_i$ be the closure in $\G$ of the unique component of $\G\setminus\G^{(0)}$ attached at $b_i$. Then, for every $i\in\mathbb Z$,
\[
 \gamma(T_i)=T_i,\qquad
 \gamma(S_i)
=\gamma(\V_0\cap T_i)
=\gamma(\V_0)\cap \gamma(T_i)=\V_0\cap T_i=S_i .
\]

For $t\in J$, \eqref{app-tree-leaf-bound} gives $f_j(t)(\vv)>\tfrac12\sqrt{\tanh(t_-/2)}$ on $S_j$. For $\vv\in\V_0\setminus S_j$, we have $f_j^0(t)(\vv)=0$.
Hence \eqref{ac-norm-ordering}, \eqref{app-tree-close} and
\eqref{app-tree-smallness} give, for every $\vv\in\V_0\setminus S_j$,
\[
 f_j(t)(\vv)\leq\norm{f_j(t)-f^0_j(t)}_\infty
 \leq\norm{f_j(t)-f^0_j(t)}_X\leq r_j<\tfrac12\sqrt{\tanh(t_-/2)}.
\]
Thus, the rescaled vertex values recover $S_j$:
\begin{equation}\label{app-tree-support-characterization}
 S_j=\left\{\vv\in\V_0:t^{-1/2}u_j(t)(\vv)>\tfrac12\sqrt{\tanh(t_-/2)}\right\}.
\end{equation}
If $u_j(t_j)=\gamma u_k(t_k)$ for some $\gamma\in\Gamma$, the edge equations give $(t_j^2-t_k^2)u_j(t_j)=0$. Since $\norm{u_j(t_j)}_2^2=\mu>0$ and $t_j,t_k>0$, we have $t_j=t_k$. Applying \eqref{app-tree-support-characterization} gives $S_j=\gamma(S_k)=S_k$, hence $j=k$. Since $d_j=d$ for infinitely many $j$ with $\abs j\geq j_0$, this completes the proof of Proposition~\ref{prop-ac-critical-all-mass-tree}.
\end{proof}

\subsection*{Conflict of interest}

The authors declare no conflict of interest.

\subsection*{Ethics approval}
 Not applicable.

\subsection*{Data Availability Statement}
Data sharing not applicable to this article as no datasets were generated or analysed during the current study.

\subsection*{Declaration of Generative AI Use}
During the preparation of this manuscript, the authors used ChatGPT to improve the English and the clarity of the exposition. The authors take full responsibility for the content of the paper. 

\subsection*{Acknowledgements}
 C. Ji is partially supported by the National Natural Science Foundation of China (No.12571117).

  \end{document}